\documentclass[onefignum,onetabnum]{siamart190516}

\usepackage{lipsum}
\usepackage{amsfonts}
\usepackage{graphicx}
\usepackage{epstopdf}
\usepackage{algorithmic}
\ifpdf
  \DeclareGraphicsExtensions{.eps,.pdf,.png,.jpg}
\else
  \DeclareGraphicsExtensions{.eps}
\fi

\newsiamremark{remark}{Remark}
\newsiamremark{hypothesis}{Hypothesis}
\crefname{hypothesis}{Hypothesis}{Hypotheses}
\newsiamthm{claim}{Claim}

\headers{DDSM for EIT}{Ruchi Guo and Jiahua Jiang}

\title{
Construct Deep Neural Networks based on Direct Sampling Methods for Solving Electrical Impedance Tomography
\thanks{Submitted to the editors \today.
\funding{The first author was funded by NSF DMS-2012465. The second author was funded by NSF DMS 1654175.}}}

\author{Ruchi Guo\thanks{Department of Applied Mathematics, University of California, Irvine
  (\email{ruchig@uci.edu}).}
\and Jiahua Jiang\thanks{Department of Applied Mathematics, Virginia Tech 
  (\email{jiahua@vt.edu}).}
}

\usepackage{amsopn}

\ifpdf
\hypersetup{
  pdftitle={Construct Deep Neural Networks based on Direct Sampling Methods for Solving Electrical Impedance Tomography},
  pdfauthor={Ruchi Guo and Jiahua Jiang}
}
\fi

\usepackage{geometry}
\allowdisplaybreaks[4]
\usepackage{subcaption}

\usepackage{bbm,mathabx}
\usepackage{array,multirow,booktabs} 
\usepackage{psfrag, graphicx,enumitem,indentfirst,float,geometry,color,graphicx,subcaption}

\makeatletter
\def\BState{\State\hskip-\ALG@thistlm}
\makeatother

\newcommand{\bs}{\boldsymbol}

\newcommand{\jjh}{\color{black}}
\newcommand{\jjhn}{\color{black}}
\newcommand{\grc}{\color{black}}

\begin{document}


\date{}
\maketitle

\begin{abstract}
This work investigates the electrical impedance tomography (EIT) problem when only limited boundary measurements are available, which is known to be challenging due to the extreme ill-posedness. Based on the direct sampling method (DSM) introduced in \cite{chow2014direct}, we propose deep direct sampling methods (DDSMs) to locate inhomogeneous inclusions in which two types of deep neural networks (DNNs) are constructed to approximate the index {\jjh function(functional)}: fully connected neural network(FNN) and convolutional neural network (CNN). The proposed DDSMs {\jjh are} easy to {\jjh be implemented}, capable of {\jjh incorporating} multiple Cauchy data pairs {\jjh to achieve} high-quality reconstruction and highly robust with respect to large noise. Additionally, {\jjh the} implementation {\jjh of DDSMs} adopts offline-online decomposition{\jjh, which} helps {\jjh to reduce a lot of computational costs} and makes {\jjh DDSMs} {\jjh as efficient as} the conventional DSM {\jjh \cite{chow2014direct}}. The numerical experiments are presented to demonstrate the efficacy and show the potential benefits of combining DNN with DSM.
\end{abstract}
\begin{keywords}
Deep Learning, inverse problems, direct sampling methods, electrical impedance tomography, reconstruction algorithm, limited boundary data
\end{keywords}

\section{Introduction}

Electrical impedance tomography (EIT) is a very promising technique for non-invasive, radiation-free type of medical imaging. In short, by alternating electrical current at a set of electrodes and measuring the corresponding voltages on the boundary of an unknown medium(e.g. on the skin), it is possible to reconstruct the internal electrical conductivity distribution image of the medium \cite{holder2004electrical, calderon1980inverse}. EIT has wide applications, such as biomedicine\cite{zou2003review, holder2004electrical} including monitoring of ventilation distribution \cite{victorino2004imbalances}, geophysics \cite{zhdanov1994geoelectrical}, industrial detection \cite{jordana2001electrical}. The particular application considered {\jjh in this work} is to detect {\jjh the} inclusions buried in a known homogeneous background, for example, measuring the resistivity of tissue associated with malignancy, ischemia, and lung water \cite{brown2000relation,frerichs1998electrical}.

In order to describe the isotropic mathematical model for EIT, we consider a bounded domain $\Omega\subseteq\mathbb{R}^n$, $n=2,3$ with $C^1$ boundary $\partial \Omega$ occupied by some conducting materials with the electrical conductivity described by a positive function $\sigma(x)\in L^{\infty}(\Omega)$. {\jjh Let} the background homogeneous material have the conductivity $\sigma_0$, and then the support of $\sigma-\sigma_0$ indicates the inhomogeneous inclusions denoted by $D$. Suppose that $N$ different electrical currents are injected to the boundary $\partial\Omega$, then the resulted electrical potential should satisfy the {\jjh following} $N$ governing equations with the same coefficient but different boundary conditions:
\begin{subequations}
\label{govern_eq_0}
\begin{align}
\label{govern_eq_1}
\nabla \cdot (\sigma\nabla u_\omega) &= 0 \quad~~ \text{in} ~~~~ \Omega, \\
\label{govern_eq_2}
\sigma\frac{\partial u}{\partial n} &= g_\omega \quad \text{on} ~~~ \partial \Omega, ~~~ \omega=1,2,...,N
\end{align}
\end{subequations}
where $g_{\omega}(x) \in H^{-1/2}(\partial \Omega)$ are current density such that 
\begin{equation}
\label{govern_eq_3}
\int_{\partial \Omega} g_\omega(s) ds = 0.
\end{equation}
For the simplicity of notation, we write the surface potential over $\partial \Omega$ as $f_\omega := u_\omega|_{\partial\Omega}$. Mathematically, our inverse problem associated with EIT is to recover $\sigma(x)$ in $\Omega$, specifically the shape {\jjh and the position} of {\jjh the} inclusions, based on the Neumann boundary value (current data) and the observed Dirichlet boundary value (voltage data). These data typically refer to a collection of the  {\jjh Cauchy} data pairs $\{(g_\omega, f_\omega)\}_{\omega=1}^{N}$. 

Let's temporarily write $u=u_{\omega}$ and $g=g_{\omega}$ for general discussion. Then the {\jjh Neumann-to-Dirichlet} (NtD) map associated to \eqref{govern_eq_0} is defined by 
\begin{equation}
\label{eq:NtD}
\Lambda_{\sigma}: H^{-1/2}(\partial \Omega) \rightarrow H^{1/2}(\partial \Omega), \quad \Lambda_{\sigma}g = u|_{\partial \Omega}.
\end{equation}
Theoretically, the conductivity distribution $\sigma$ as a positive $L^{\infty}$ function can be recovered from a full knowledge of the NtD map $\Lambda_{\sigma}$ \cite{astala2006calderon}. However, in practice, the full knowledge of the NtD map needs a matrix approximation {\jjh requiring} a large number of Cauchy data pairs for example $N=64$ in \eqref{govern_eq_0} \cite{chow2014direct} but may {\jjh be} unfeasible in real-world problems. {\jjh Nevertheless} the stability and accuracy of the approximation also need to be resolved. In addition, in many practical situations, the conductivity coefficient $\sigma$ is not as rough as $L^{\infty}$. In fact the piecewise constant conductivity widely appears, and some researches \cite{ammari2007polarization,1990IsakovPowell, ammari2004reconstruction,cedio1998identification} have shown how one or few boundary measurement(s), i.e., the Cauchy data pairs, can reconstruct {\jjh the} inclusions. {\grc So from both the practical and theoretical point of view, the full knowledge of the NtD map not be very necessary, which justifies the limited Cauchy data assumption aforementioned.  }

It is known that a high-quality reconstruction for the EIT problem is challenging due to its severe ill-posedness. Various numerical algorithms have been developed which can be categorized into two family. The first one is based on optimization algorithm which typically constructs a sequences of {\grc coefficient distributions converging to the true one}. Methods in this family include the finite-element-based methods \cite{martin1997fem,vauhkonen1999three,jin2017convergent} and finite-difference-based methods \cite{rondi2001enhanced}, shape optimization methods \cite{GuoLinLin2017,2004ChanTai,chung2005electrical,2001ItoKunischLi} and so on. For these methods, in general a good initial condition and a significant number of iterations are needed to ensure convergence which also rely on suitable regularization techniques such as Tikhonov regularization \cite{vauhkonen1999three,jin2017convergent}, total variation regularization \cite{chung2005electrical} and $L^1$ regularization \cite{jin2012reconstruction}. {\grc Alternatively}, the second family of methods, i.e., the direct methods, have been developed including the sampling and factorization methods \cite{kirsch2008factorization, chow2014direct}, the multiple signal classification algorithms \cite{ammari2007polarization, ammari2004reconstruction}, and the topological sensitivity approaches \cite{ammari2004reconstruction,bonnet2009higher}. These methods are {\jjh non-iterative in nature} and thus tend to be much more efficient than optimization-type methods. Besides the classical approaches, recently, deep neural networks (DNNs) also have shown great potential for solving the EIT problem \cite{adler1994neural, fan2020solving, hamilton2018deep}, for example, combining fast direct reconstruction procedures with CNN \cite{li2019novel, hamilton2018deep, tan2018image}, the radial basis function neural network \cite{michalikova2014image,wang2009implementation}, the multi-layer neural network \cite{klosowski2017using}, the bayesian neural network \cite{lampinen1999using} and improving the neural network training with particle swarm optimization \cite{martin2015electrical, martin2015nonlinear}.

Recently, the authors in \cite{chow2014direct,chow2015direct,chow2018time,chow2020direct,ito2012direct,ito2013direct} develop a {\jjh so-called} direct sampling method (DSM) to solve not only EIT but also a large group of inclusion identification problems which are shown to be easy to implement, computationally economical, highly robust and very effective. The principal ingredient of DSMs is to construct an index function that attains extreme values for the sampling points belonging to the inclusions and thus provides the estimation of the shape. In particular, the index function for EIT is derived exquisitely in \cite{chow2014direct} showing reasonably nice {\jjh indication} for the inclusions. However, we also note the DSM in \cite{chow2014direct} has the potential to be further improved in several aspects. First, it is mainly focused on the case of \textit{single} Cauchy data pair which may hinder its application to the more complicated case due to the limited accuracy. In engineering practice \cite{holder2004electrical}, the number of experimental measurement data is indeed limited but mostly more than one. Second, {\jjh the} index function in \cite{chow2014direct} highly relies on the domain geometry and is only obtained for some basic shapes such as circle, square and open ball. Indeed, the derivation of suitable index functions may become extremely complicated or even unobtainable for multiple Cauchy data pairs and more complex domain shapes. Furthermore, even with singe Cauchy data pair, the optimality may not be guaranteed for the index functions derived in \cite{chow2014direct}. {\grc We believe it may not be easy to overcome these obstacles by conventional mathematical approaches.}

Alternatively, {\jjh to address these hurdles numerically}, in this work we propose using DNNs to replace the classical intricate derivation to construct index functions{\jjh, and call} the resulted methods \textit{deep direct sampling methods} ({\jjh DDSMs}) which have a few remarkable features. First, the {\jjh DDSMs} are capable of {\grc naturally both} incorporating multiple Cauchy data pairs and fitting for any shaped $\Omega$, and thus {\jjh the proposed methods} can break the accuracy limits and result in high-quality reconstruction for both location and shape. Second, {\jjh the DDSMs} also {\jjh inherit} the robustness feature of the conventional DSM with respect to large noise. {\jjh More specifically, these algorithms} can smooth the noise appearing in the Cauchy data and {\jjh handle} noise as large as like 20\%. Furthermore, although the training process costs more computational budget than the conventional DSM since it requires solving an optimization problem, it only needs to be operated once in offline phase. The online computation involved in reconstruction/prediction is only the fast evaluation of the DNN-based index function which {\jjh is as efficient} as the DSM. The offline-online decomposition structure makes the proposed DDSM has more optimal index function without degrading the efficiency, where the optimality is benefited from DNN and the efficiency is kept from DSM. {\grc Therefore, we think the DDSMs combine} the advantages of both optimization methods and direct methods.

In particular, we develop two types of DDSMs: fully connected neural network based DDSM (FNN-DDSM) and convolutional neural network based DDSM (CNN-DDSM). {\jjh Suggested by the conventional DSM, both the networks take the mathematical format of input data} which are the solutions of an elliptic equation with the known background coefficient and the boundary condition give by a certain difference of the Cauchy data pairs. It is noted that the choice of input and output for DNNs are in general crucial for their performance. For FNN-DDSM, we consider the index function as a pointwise indicator that classifies the points in the domain into two categories, inside the inclusions and outside the inclusions, which can be naturally treated as a classification problem. {\grc A softmax layer is then used to normalize the output of the network to a probability distribution over predicted output classes.} We highlight that the probabilistic meaning of the output of FNN-DDSM computationally hinted the connection within the EIT problem, DSM and the probability problem. For CNN-DDSM, we recast the index function into a functional mapping data images to EIT images, so it can be viewed as a semantic image segmentation problem that aims to extract the contour(detect the edge of inclusions) from the image. It can be also considered as a further generalization of the structure of the DSM in \cite{chow2014direct} since it uses more neighbor information of each point to predict their location.

The remainder of this paper is organized as follows. The mathematical background of the DSM \cite{chow2014direct} is prepared in Section \ref{sec:DSM}. Section \ref{sec:DDSMs} is devoted to the development of our novel DDSMs consisting of FNN-DDSM and CNN-DDSM. Our numerical results are presented in Section \ref{sec:num}. Some concluding remarks are given in Section \ref{sec:con}.


\section{Review of Direct Sampling Methods (DSM)}
\label{sec:DSM}
This section summarizes some necessary theoretical background of the DSM \cite{chow2014direct} for solving the aforementioned EIT problem which also serves as the mathematical foundation of the proposed neural networks. The key idea is to derive a certain index function indicating the inclusion shape and locations, which should ideally satisfy
\begin{equation}
\label{idea_index}
\mathcal{I}(x)=
\begin{cases}
      & 1 ~~~~ x\in D \\
      & 0 ~~~~ x\in \Omega\backslash \overline{D}.
\end{cases}
\end{equation}

Since the work in \cite{chow2014direct} mainly focuses on the case of only \textit{single} Cauchy data available, we here just let $(g, f)$ be the pair of Cauchy data measured over the surface $\partial \Omega$ with $f = \Lambda_{\sigma} g$ and the NtD map $\Lambda_{\sigma}$ defined in (\ref{eq:NtD}). {\grc Without loss of generality, we assume the background medium has $\sigma_0=1$.} A fundamental ingredient in the derivation of DSM is the duality product $\langle\cdot, \cdot\rangle_{\gamma,\partial\Omega}$ defined as
\begin{equation}
\label{eq:dual}
\langle\mathcal{X}, \phi\rangle_{\gamma,\partial\Omega} = \int_{\Gamma} (- \Delta_{\partial\Omega})^{\gamma}\mathcal{X}\phi ds \nonumber = \langle (- \Delta_{\partial\Omega})^{\gamma}\mathcal{X}, \phi\rangle_{L^2(\partial\Omega)} \quad \forall \mathcal{X} \in H^{2\gamma}(\partial \Omega), \forall \phi \in L^2(\partial\Omega)
\end{equation}
where $\gamma\ge 0$ and $\Delta_{\partial\Omega}$ is the surface Laplacian operator. Let $| \cdot |_{Y}$ denote certain semi-norm in $H^{2\gamma}(\partial\Omega)$. An essential component to derive the index function is a family of {\jjh probing functions} $\{\eta_{x,d}\}_{x\in \Omega, d\in \mathbb{R}^n} \subset H^{2\gamma}(\partial\Omega)$ satisfying some conditions of which the critical one for us is {\jjh represented} below
\begin{itemize}
\item[(\textbf{C})] The probing functions are almost orthogonal with each other. That is, $\forall x, y \in \Omega, d_x, d_y \in \mathbb{R}^n$, the function
\begin{equation}
\label{assum1}
K_{d_x,d_y}(x,y) := \frac{\langle \eta_{x,d_x}, \eta_{y,d_y}\rangle_{\gamma, \partial\Omega} }{| \eta_{x,d_x}|_{Y} }
\end{equation}
achieves maximum when $x=y$ behaving like a sharply-peaked Gaussian-like distribution, i.e., it is close to a kernel function $e^{\frac{-|x-y|^2}{a^2}}$ with small $a$.
\end{itemize}
We emphasize that the condition (\textbf{C}) is particularly helpful for us to choose suitable activation functions in the design of networks. The construction of suitable probing functions can be challenging in general. {\jjh The probing functions introduced} in \cite{chow2014direct} {\jjh are} based on the dipole potential \cite{el2000inverse, griffiths2007introduction}, {\jjh and they satisfy} 
\begin{equation}
\label{eq:defphi}
\eta_{x,d}(\xi):= w_{x,d}(\xi), \quad \xi \in \partial\Omega,
\end{equation}
where $w_{x,d}$ is the solution of the following problem
\begin{align}
\label{eq:wproblem}
-\Delta w_{x,d}  = -d \cdot  \nabla\delta_x \quad \text{in} \quad \Omega,  \quad \frac{\partial w_{x,d}}{\partial n} = 0 \quad \text{on} \quad \partial \Omega, \quad \int_{\partial \Omega} w_{x,d} ds = 0
\end{align}
with given $x \in \Omega, d \in \mathbb{R}^n$ and $\delta_{x}(\xi)$ being the delta function {\grc for each $x$}.
Empirically the authors have chosen $|\cdot |_Y = | \cdot |_{H^{3/2}(\partial \Omega)}$  {\grc in \eqref{assum1}}.  
With the probing functions, the index function can be then defined as
\begin{equation}
\label{eq:newindex}
\mathcal{I}(x, d_x) := \frac{\langle \eta_{x, d_x}, f - \Lambda_{\sigma_0}g \rangle_{\gamma, \partial\Omega}}{|| f - \Lambda_{\sigma_0}g  ||_{L^2(\partial \Omega)} | \eta_{x, d_x} |_{H^{3/2}(\partial\Omega)}}, \quad x\in\Omega, d_x \in \mathbb{R}^n,
\end{equation}
With the definition of duality product $\langle \cdot, \cdot \rangle_{\gamma, \partial\Omega}$ in \eqref{eq:dual}, the surface Laplacian operator only performs on the probing function $\eta_{x,d}$ that is itself infinitely smooth over the measurement surface. So the index function in \eqref{eq:newindex} is well defined even for very irregular data $f-\Lambda_{\sigma}g$ containing very rough noise, for example those data only in $L^2(\partial \Omega)$. As mentioned in \cite{chow2014direct}, due to this feature the noise appearing in data can be directly smoothed by the duality product over the measurement surface, and it is the reason of the high robustness of the DSM against noise. We also believe it is also one of the reasons that our DDSMs are highly stable with respect to noise as observed in Section \ref{sec:num}.

Note that \eqref{eq:newindex} is not computable yet since the probing direction $d_x  \in \mathbb{R}^n$ remains unknown. In order to find an appropriate direction $d_x$, an alternative characterization of the index function is considered in \cite{chow2014direct}. 
 Define $\phi$ as the solution to the standard elliptic equation with the boundary condition $(- \Delta_{\partial \Omega})^{\gamma} (f - \Lambda_{\sigma_0}g)$:
\begin{align}
\label{eq:phiproblem}
-\Delta \phi  = 0 \quad \text{in} \quad \Omega, \quad  \frac{\partial \phi}{\partial n} = (- \Delta_{\partial \Omega})^{\gamma}(f - \Lambda_{\sigma_0}g)  \quad \text{on} \quad \partial \Omega, \quad \int_{\partial \Omega} \phi ds = 0.
\end{align}
We highlight that the function $\phi$ in \eqref{eq:phiproblem} is particularly important for designing our DNNs and processing given boundary data to the input(images) of the DNNs. 
Through \eqref{eq:phiproblem}, the index function in \eqref{eq:newindex} can be equivalently rewritten as
\begin{equation}
\label{eq:defindexfun}
\mathcal{I}(x, d_x) := \frac{d_x \cdot \nabla\phi(x)}{|| f - \Lambda_{\sigma_0}g ||_{L^2(\partial \Omega)} | \eta_{x, d_x} |_{H^{3/2}(\partial\Omega)}}, \quad x\in\Omega, d_x \in \mathbb{R}^n.
\end{equation}
It is easy to see that the above index function reaches its maximum value at $d_x = \frac{\nabla\phi(x)}{|\nabla\phi(x)|}$ for each $x \in \Omega$ which is the appropriate choice of the probing direction $d_x$. As observed in the experiments in \cite{chow2014direct}, such choice of $d_x$ is crucial since for example, the existing inclusions might be accidentally removed if $d_x$ is wrongly chosen orthogonal to $\nabla \phi(x)$ at a point $x \in \Omega$. Thus $\mathcal{I}(x,d_x) = \mathcal{I}(x,\nabla\phi(x))$ which is the guideline for us to construct DNNs. 

{\jjh Notice} that computing the index function \eqref{eq:defindexfun} with $d_x=\nabla\phi(x)$ is only up to solving the probing functions from \eqref{eq:defphi} and \eqref{eq:newindex}, and this {\jjh is} challenging since {\jjh the} probing functions {\jjh vary} with respect to $x\in\Omega$. Instead of directly solving \eqref{eq:defphi} and \eqref{eq:newindex} that may slow down the DSM, the authors in \cite{chow2014direct}  derive the explicit forms of the probing  functions with very delicate mathematical skills for some special domain:
\begin{itemize}
\item $\Omega$ is a circular domain: 
\begin{equation}
\label{eq:probcir}
\eta_{x,d}(\xi) = \frac{1}{\pi}\frac{(\xi - x) \cdot d}{|x - \xi|^2}, \quad \xi \in \partial\Omega.
\end{equation}
\item $\Omega$ is an open ball:
\begin{equation}
\label{eq:probball}
\eta_{x,d}(\xi) = \frac{d \cdot \xi - \frac{(x -\xi) \cdot d}{|x - \xi|}}{\sqrt{4\pi}(|x - \xi| - x \cdot \xi +1)}, \quad \xi \in \partial\Omega.
\end{equation}
\end{itemize}
The simple and explicit probing functions $\eta_{x,d}(\xi)$ in \eqref{eq:probcir}-\eqref{eq:probball} are then put into \eqref{eq:defindexfun}, and thus the index function $\mathcal{I}(x, \nabla\phi(x))$ can be computed efficiently without solving any PDEs. So it is essentially different from optimization-based iterative methods which require repeatedly solving forward problems and adjoint problems for many times. 




\section{Deep Direct Sampling Methods (DDSMs)}
\label{sec:DDSMs}


The whole design procedure of the index function of the DSM is delicate and intellectual. Despite these features and the successful application, {\jjh the} accuracy and further application {\jjh of the DSM} are limited mainly by the following two aspects. First, the derivation above is only focused on \textit{single} Cauchy data pair.  But when multiple pairs of Cauchy data are applied on $\Omega$ with any geometry, it is very difficult (at least currently unclear to mathematicians) {\jjhn to obtain} such index function by canonical mathematical derivation. Second, the design of $\mathcal{I}(x, d_x)$ is not necessarily optimal. For example, the tuning parameter $\gamma$ and norm $|\cdot|_Y$ need to be chosen empirically, and the format itself may not {\jjh be} the best approximation to the true index function. {\grc So it motivates us to use deep neural {\jjh network} (DNN) models to learn the index functions through a large number of data since DNNs are able to mimic the human's learning process based on physical data.}


Therefore, to make DSMs applicable to more general situations of EIT problems and ameliorate the quality of reconstruction, in this section {\grc we propose {\em Deep Direct Sampling Methods (DDSMs)}.}
The {\jjh major difficulty} is on suitable design of DNN architectures. We build up the architectures based on two different perspectives to understand the mathematical form of the index function \eqref{eq:defindexfun}, which are the topics of the next two subsections. The first algorithm, fully connect neural network(FNN) based DDSM (FNN-DDSM), is from the perspective of function approximation, while the second one, convolutional neural networks(CNN) based DDSM (CNN-DDSM), is from the angle of functional approximation.
FNN-DDSM and CNN-DDSM are intrinsically different both in their network architectures and their mathematical foundations,
yet our numerical experiments show that they both work very well, enhancing the accuracy and stability without sacrificing the efficiency. 



For simplicity we mainly {\grc discuss our methods for} the two-dimensional $\Omega$ and they can be naturally extended to the three-dimensional case. To avoid repetition, we first prepare some notations which will be frequently used in the following discussion.

\begin{itemize}
\item $\mathcal{N}_h$: {\grc the discretization of domain $\Omega$ that consists of a group of discrete points, i.e., $\mathcal{N}_h = \{x^{k}=(x^{k}_1, x^{k}_2) \}^{K}_{k=1}$ where $K$ is the total number of points.}
\item $\phi^{\omega}$ with $1 \le \omega \le N$: the solutions of \eqref{eq:phiproblem} with the boundary value from the $\omega$-th pair of Cauchy data, i.e., $f_{\omega}-\Lambda_{\sigma_0} g_{\omega}$. Since these {\jjh functions} are so critical in both the {\jjh DSM} and the proposed DDSMs, we shall call them \textit{Cauchy difference functions}.
\end{itemize}
If the training data samples are involved in the context, we further define the following notations
\begin{itemize}
\item $S$: the number of inclusion samples (coefficient distribution samples) in the training sets. 
\item $\mathcal{I}^s$ with $1\le s \le S$: the true index function associated with the $s$-th inclusion sample.
\item $\phi^{(s,\omega)}$ with $1\le s \le S$ and $1 \le \omega \le N$: the \textit{Cauchy difference functions} associated with the $\omega-$th pair of Cauchy data generated by the $s$-th inclusion sample.
\end{itemize}
Here $\gamma$ is chosen to be $0$ to generate $\phi^{\omega}$ or $\phi^{(s,\omega)}$ when solving the problems \eqref{eq:phiproblem}. Indeed $\gamma$ can be non-zero, but numerically the boundary values with $\gamma \neq 0$ can be all mapped from those with $\gamma = 0$, and this mapping can be handled by DNNs automatically. So we simply pick $\gamma = 0$ which also makes solving \eqref{eq:phiproblem} easier.

\subsection{FNN-DDSM}
\label{sec:fnn}
We first discuss the construction of FNN-DDSM. The derivation of $\mathcal{I}(x, d_x)$ in (\ref{eq:defindexfun}) suggests the existence of a non-linear mapping from $x$, $d_x$ to the location of $x$, i.e., whether it inside or outside of the inclusions. Furthermore, inspired by the choice of $d_x$, we {\jjh assume} that the direction $d_x$ {\jjh is} a function of $\nabla \phi^{\omega}(x)$ for all the Cauchy pairs. Therefore, we assume that the new index function takes the form
\begin{equation}
\label{eq:defnewindex}
\mathcal{I}(x) = \mathcal{F}_{\text{FNN}}(x, \nabla\phi^{1}(x),\dots,\nabla\phi^{N}(x)), \quad x \in \Omega,
\end{equation} 
where $\mathcal{F}_{\text{FNN}}$ is trained by a FNN described in details below. Mathematically speaking, $\mathcal{F}_{\text{FNN}}$ is a nonlinear high-dimensional function mapping a data point in $\mathbb{R}^{2N+2}$ to the index of the associated point in $\Omega$. Notice that the design of \eqref{eq:defnewindex} is directly motivated by the form of \eqref{eq:defindexfun} with multiple Cauchy data pairs. But at this moment we do not assume $\mathcal{F}_{\text{FNN}}$ takes any a-priori format such as \eqref{eq:defindexfun} and only make the assumptions on its input. So we believe the proposed new index function is more general than the index function \eqref{eq:defindexfun} even with $N=1$. Since there are no special assumptions and prior knowledge about the nonlinear mapping in \eqref{eq:defnewindex}, we believe it is appropriate to train the index function by {\jjh a FNN} {\grc due to its structural agnostic which are universal approximators capable of learning many functions.}
Notice that the input purely depends on $x\in\Omega$ and values of $\nabla \phi^{\omega}(1\le \omega \le N)$ at $x$. So using \eqref{eq:defnewindex} to approximate \eqref{idea_index} can be also treated as a binary classification problem for high-dimensional data points which are associated with each point in $\Omega$.

\begin{figure}[h!]
\centering
\begin{tabular}{c}
\includegraphics[width = 0.9\textwidth]{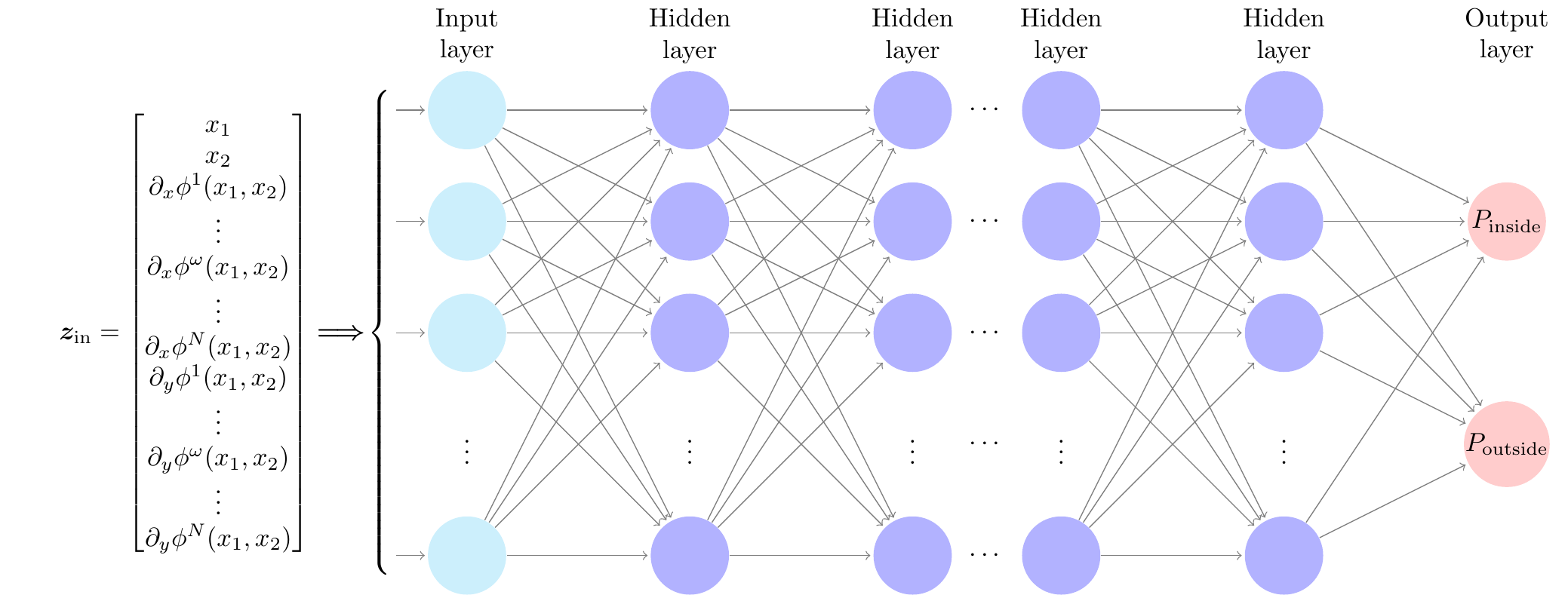}\\
(a)\\
\includegraphics[width = 0.9\textwidth]{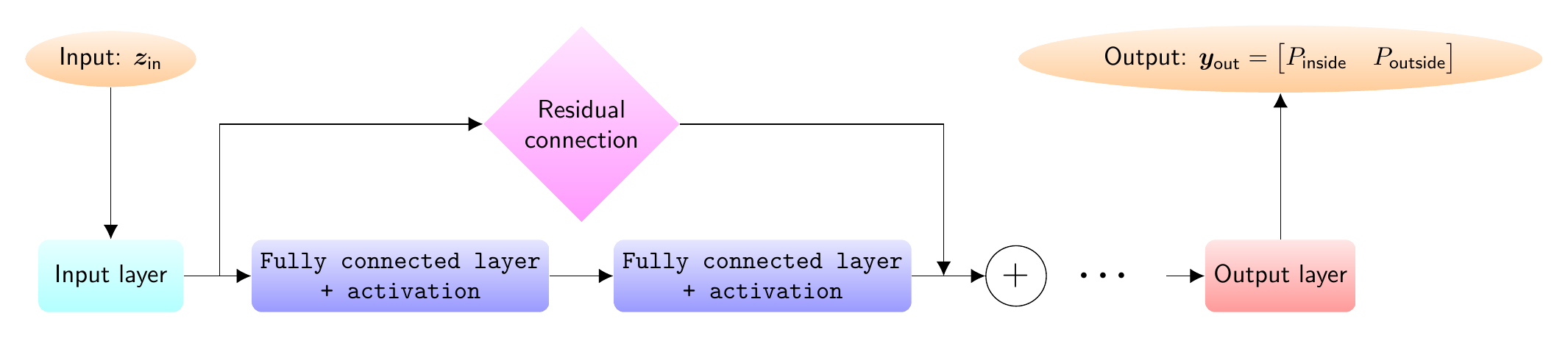}\\
(b)
\end{tabular}
\caption{The structure of FNN. $P_{\text{inside}}$ indicates the probability of a point $x=(x_1,x_2)$ being inside the inclusions. $P_{\text{outside}}$ indicates the probability of a point $x=(x_1,x_2)$ being outside the inclusions.}
\label{fig:fnnstruct}
\end{figure}

Now we proceed to describe the structure of the proposed FNN. Given an inclusion sample and the associated Cauchy difference functions $\phi^{\omega}$, we let the input of {\jjh the} FNN, denoted by $\bs{z}_{\text{in}} \in \mathbb{R}^{2(N+1)}$ (depicted in Figure \ref{fig:fnnstruct}(a)), contain the coordinates of $x$ in the first two components, and the value of $\nabla \phi^{\omega}(x)$ in the next $2N$ components, i.e.,
\begin{subequations}
\label{zinput}
\begin{align}
[\bs{z}_{\text{in}}]_1 &= x_1, \\
[\bs{z}_{\text{in}}]_2 &= x_2, \\
[\bs{z}_{\text{in}}]_{3,\dots,N+2} &= \begin{bmatrix}\partial_x\phi^1(x_1, x_2) & \dots & \partial_x\phi^N(x_1, x_2) \end{bmatrix}^{\top}, \\
[\bs{z}_{\text{in}}]_{N+3,\dots,2N+2} &= \begin{bmatrix}\partial_y\phi^1(x_1, x_2) & \dots & \partial_y\phi^N(x_1, x_2) \end{bmatrix}^{\top}.
\end{align}
\end{subequations}
If the length of the input $\bs{z}_{\text{in}}$ is smaller than the width of the input layer, to handle this discrepancy we can pad $\bs{z}_{\text{in}}$ by a zero vector. The corresponding input linear layer is defined as
\begin{equation}
\psi_{\text{in}}(\bs{z}_{\text{in}}) = W_{\text{in}}\bs{z}_{\text{in}} + \bs{b}_{\text{in}},
\end{equation}
where $W_{\text{in}}$ is the weight matrix associated with the synapses connecting the input layer neurons to the neurons in the current layer and $\bs{b}_{\text{in}}$  is the bias. In the architecture that we use, except input layer and output layer, each layer is associated a block, and each block is constructed by stacking several layers including two linear transformations, two activation functions and a residual connection. In particular, suppose $\bs{z}_i$ is the input of the $i$-th block, then the $i$-th block can be expressed as
\begin{equation}
\label{eq:block}
\tau_i(\bs{z}_i ) = \varphi(W_{i,2} ~ \varphi(W_{i,1}\bs{z}_k  + \bs{b}_{i,1}) + \bs{b}_{i,2}) + \bs{z}_i ,
\end{equation}
where $W_{i,1}, W_{i,2}, \bs{b}_{i,1}, \bs{b}_{i,2}$ are weight matrices and bias vectors. The output linear layer is 
\begin{equation}
\psi_{\text{out}}(\bs{z}_{\text{out}}) =  W_{\text{out}}\bs{z}_{\text{out}} + \bs{b}_{\text{out}},
\end{equation} 
where $W_{\text{out}}, \bs{b}_{\text{out}}$ also denote weight matrices and bias vectors for this layer. Let $\Theta$ denote the collection of the unknown parameters in these matrices and vectors to be {\jjhn learned} in training. The complete FNN with $M$ blocks can now be represented as
\begin{equation}
\label{eq:wholenn}
\bs{y}_{\text{out}}(\bs{z}_{\text{in}})= \kappa\circ\psi_{\text{out}}\circ\tau_{M}\circ\dots\circ\tau_{1}\circ\psi_{\text{in}}(\bs{z}_{\text{in}}),
\end{equation}
where $\kappa$ is the softmax layer. The entire structure including blocks and residual connections is shown in Figure \ref{fig:fnnstruct}. 

Moreover, according to the assumption (\textbf{C}), the inner product $<\cdot,\cdot>_{\gamma,\partial\Omega}$ in \eqref{assum1} of probing functions is very close to a sharply-peaked Gaussian function (small variance), which gives more sensitivity of the index function with respect to the data. Thus, inspired by this mathematical intuition, we employ a Gaussian activation function $\varphi(z) = e^{-\frac{||z||^2}{2a^2}}$($a^2$ is variance chosen to be small in computation) for the first several blocks. For the rest blocks, we instead use the clipped rectified-linear (ReLu) activation function $\varphi(z) = \min\{\max\{0,z\}, 0.1\}$ to mimic the distribution of coefficients. Indeed, our experience has suggested that such a Gaussian activation function makes the network much easier to train and also improves its performance.


In addition, because the softmax layer normalizes the output of a network to a probability distribution over predicted output classes, the output $\bs{y}_{\text{out}}(\bs{z}_{\text{in}}) \in \mathbb{R}^{1\times2}$ consists of two components within $[0,1]$ with the unit sum, i.e.,
\begin{align}
[\bs{y}_{\text{out}}(\bs{z}_{\text{in}})]_1,~ [\bs{y}_{\text{out}}(\bs{z}_{\text{in}})]_2 \in [0, 1], ~~~~~~~~ [\bs{y}_{\text{out}}(\bs{z}_{\text{in}})]_1 + [\bs{y}_{\text{out}}(\bs{z}_{\text{in}})]_2 = 1,
\end{align}
where $[\bs{y}_{\text{out}}(\bs{z}_{\text{in}})]_1$ indicates the probability of the point $x=(x_1, x_2)$ locating inside the inclusion. Namely, the larger $[\bs{y}_{\text{out}}(\bs{z}_{\text{in}})]_1$ (more close to $1$), more possible the point is inside the inclusion, and vice versa. Otherwise, it means that $(x_1, x_2)$ is more likely to be outside of the inclusions. 

It is remarked that the value of the index function \eqref{eq:defindexfun} of the conventional DSM is not exactly $0$ or $1$ either.  It naturally motivates some researchers to think that the DSM is more like evaluating the chance of each mesh point situated within the inclusions, i.e., the possibility of the input point inside the {\jjhn inclusions}. This feature is essentially different from some conventional optimization methods such as the shape optimization from \cite{chung2005electrical,2004ChanTai} that gives the sharp interface of inclusions. But notice that there is no rigorous theoretical analysis on the relationship between the DSM and probability. From this point of view, we think the proposed FNN in certain sense serves as a bridge between the DSM for EIT problems and probability since its output $\bs{y}_{\text{out}}(\bs{z}_{\text{in}})$ has very clear probabilistic meaning for every single input point in $\Omega$. In addition, due to the offline optimization stage of the DDSMs, the estimate to the probabilistic distribution is also much more accurate. We believe this feature makes the proposed algorithm particularly suitable for solving the EIT problems. 



It is well-known that the cross-entropy is a suitable loss function for classification problems, since it minimizes the distance between two probability distributions -- predicted and actual. In order to employ it to evaluate the proposed FNN, based on the notations introduced at the beginning of this section, we first consider the following loss function
\begin{equation}
\label{eq:loss}
\mathcal{L}_{\text{loss}}(\Theta) = \frac{1}{S}\sum_{s = 1}^{S}\big( \frac{1}{K}\sum_{k=1}^{K}\mathcal{I}^{s}(x^{k})\log([\bs{y}_{\text{out}}(\bs{z}^{(s,k)}_{\text{in}})]_1 )   +  (1-\mathcal{I}^{s}(x^{k}))\log([\bs{y}_{\text{out}}(\bs{z}^{(s,k)}_{\text{in}})]_2) \big),
\end{equation} 
where $\mathcal{I}^{s}(x^{k})$ denotes the $s$-th true index function evaluated at the $k$-th point $x^{k}$, $\bs{z}^{(s,k)}_{\text{in}}$ denotes the input formed by the $k$-th point $x^{k}$ and the $s$-th Cauchy difference function, i.e., the $\phi^{\omega}$ in \eqref{zinput} being $\phi^{(s,\omega)}$. However we note that there are in general a large amount of inclusion samples and the discrete points to capture the geometrical details of inclusions. To reduce the computational burden and achieve faster iterations, we only use a random subset of inclusion samples and discrete points at each iteration. Namely, we employ Stochastic Gradient Descent(SGD) \cite{da2014method} to update the weight matrices and bias vectors according to the following formulae
\begin{align}
\Theta^{(j+1)} = \Theta^{(j)} - \alpha \nabla_{\Theta}\mathcal{G}_{\text{loss}}(\Theta^{(j)})
\end{align}
where $\alpha$ indicates the learning rate and $j$ is the number of iteration. $\mathcal{G}_{\text{loss}}(\Theta)$ is defined as
\begin{equation}
\mathcal{G}_{\text{loss}}(\Theta) = \frac{1}{|S_r|}\sum_{s\in S_r}\big( \frac{1}{|K_r|}\sum_{k \in K_r}\mathcal{I}^{s}(x^{k})\log([\bs{y}_{\text{out}}(\bs{z}^{(s,k)}_{\text{in}})]_1 )   +  (1-\mathcal{I}^{s}(x^{k}))\log([\bs{y}_{\text{out}}(\bs{z}^{(s,k)}_{\text{in}})]_2) \big),
\end{equation}
where $S_r$, $K_r$ are random subsets of the integer sets $\{1,2,\cdots,S\}$ and $\{1,2,\cdots,K\}$ respectively. More precisely, for each iteration, we randomly choose $|S_r|$ samples in all $S$ inclusion samples, and $|K_r|$ points from all $K$ discrete points. The gradient of $\mathcal{G}_{\text{loss}}$ is computed by backpropagation \cite{miller1992neural}.

\vspace{0.1in}

\textbf{Implementation Issues}. We emphasize the proposed DDSMs {\jjh are} very easy to implement {\grc based on the current sophisticated DNN packages}. One only needs to prepare the data including the true index functions and Cauchy difference functions for the inclusion samples, and these data are only needed at some discrete points $\mathcal{N}_h$ of $\Omega$. In principle these points do not need to be determined or fixed a-priori, and instead they can be completely randomly chosen during each iteration. But this approach gives extra computation burden to evaluate $\nabla\phi^{\omega}(x)$ since the functions $\phi^{\omega}$ are in general numerically solved without an analytical formula. Therefore, for simplicity's sake, we propose directly using the discretization points involved in the computation of $\phi$ as the discretization points $\mathcal{N}_h$ for {\jjh the FNN} computation since the values of $\nabla\phi^{\omega}(x)$ at these points can be pre-determined and thus do not need be computed during training. For example, if the popular finite element method is used to compute $\phi^{\omega}$ then the underling mesh points can be used as the discretization points for {\jjh the FNN} computation. Furthermore, we mention that the number of these points is related to the resolution of the inclusion geometry. In general, the more complicated geometry requires more points to resolve, but it also makes the training more difficult. 


\subsection{CNN-DDSM}
\label{sec:cnn}

In this subsection, we present the CNN-DDSM which is essentially different from the FNN-DDSM in both the computation and mathematical foundation. Our motivation consists of multiple levels.  First, while being structure agnostic makes FNN broadly applicable, such networks do tend to require a larger number of parameters exacerbated by the fully connected structure. Second, the gradient operator used in \eqref{eq:defnewindex} can be considered as a special convolutional operator utilizing the neighboring pixels for the operation in the context of image process. More specifically, if we define $\bs{A}$ as the source image, the finite difference approximation for horizontal derivative and the vertical derivative are given by
\begin{equation}
\label{eq:grad}
\bs{G}_{x_1} = \frac{1}{h}\begin{bmatrix}0 & 1 \\ -1 & 0\end{bmatrix}*\bs{A} \qquad \text{and} \qquad \bs{G}_{x_2} = \frac{1}{h}\begin{bmatrix}1 & 0 \\ 0 & -1\end{bmatrix}*\bs{A},
\end{equation} 
where $*$ denotes the convolution operation. If we consider $(\phi^1(x), \dots, \phi^N(x))$ as an image with $N$ features(channels) and treat $x$ as pixels, then \eqref{eq:defnewindex} can be understood as the convolution in \eqref{eq:grad} operating on these images. Moreover, from the discretization perspective, gradient or \eqref{eq:grad} only involves the information of the direct neighbor points of a point $x$ to predict the index value at this point. We note that the true situation may be much more complicated than this, i.e., it may involve more neighbor terms and the operation may be more complicated than gradients (weights are not necessary those in \eqref{eq:grad}). Thus, it is reasonable to increase the matrix size in \eqref{eq:grad} and treat their special weights as unknown parameters to be learned. All these considerations motivate us to hybridize CNN and DSM, and our CNN-DDSM naturally arises. 

Mathematically speaking, different from FNN-DDSM, in CNN-DDSM the key is to understand the conventional DSM \eqref{eq:defindexfun} from a novel perspective {\grc that the desired mapping can be viewed as a functional from Cauchy difference functions $\{\phi^{\omega}\}_{\omega=1}^{N}$ to the inclusion distribution instead of a function from spatial variables to the indexes $0$ or $1$, i.e.,}
\begin{equation}
\label{CNN_index}
\mathcal{I} = \mathcal{F}_{\text{CNN}}(x,\phi^1,\dots,\phi^N)~: ~ \left[H^1(\Omega) \right]^{2N+2}  \rightarrow L^2(\Omega).
\end{equation}
By \eqref{CNN_index}, we only assume the index functionals rely on the entire set of Cauchy difference functions rather than their values at a specific point, which is certainly more relaxed than the assumption of FNN-DDSM. Due to the more abstract structure, it is even more difficult to obtain such a functional explicitly. But based on the previous discussion, we can use CNN to approximate \eqref{CNN_index}. Additionally, from image processing perspectives, if we consider $\mathcal{I}$ as the dense prediction of a $(N+2)$-channel image $(x,\phi^1, \dots, \phi^N)$, then the nonlinear functional $\mathcal{F}_{\text{CNN}}$ can be also treated as semantic image segmentation process \cite{chen2017deeplab, hong2015decoupled} that is partitioning a digital image into multiple segments (set of pixels) based on two characteristics: inside or outside the inclusions. {\grc It suggests a relationship between DSM for EIT problems and semantic image segmentation problems which is illustrated in table \ref{tab:compare} for readers from different background.}
\begin{table}[ht]
\centering
\begin{tabular}{|c|c|c|}
\hline
       & DSM for EIT & Image Segmentation \\
\hline
$x$ & mesh point in $\Omega$  & pixel \\
\hline
$\phi^{\omega}(x)$ & Cauchy difference functions  & image features (channels)\\
\hline
$\mathcal{I}$ & index of inclusion distribution & dense prediction \\
\hline
\end{tabular}
\caption{Relationship between DSM for EIT problems and image segmentation problems}
\label{tab:compare}
\end{table}

\begin{figure}[h!]
\hspace{-0.6in}
\includegraphics[width = 1.08\textwidth]{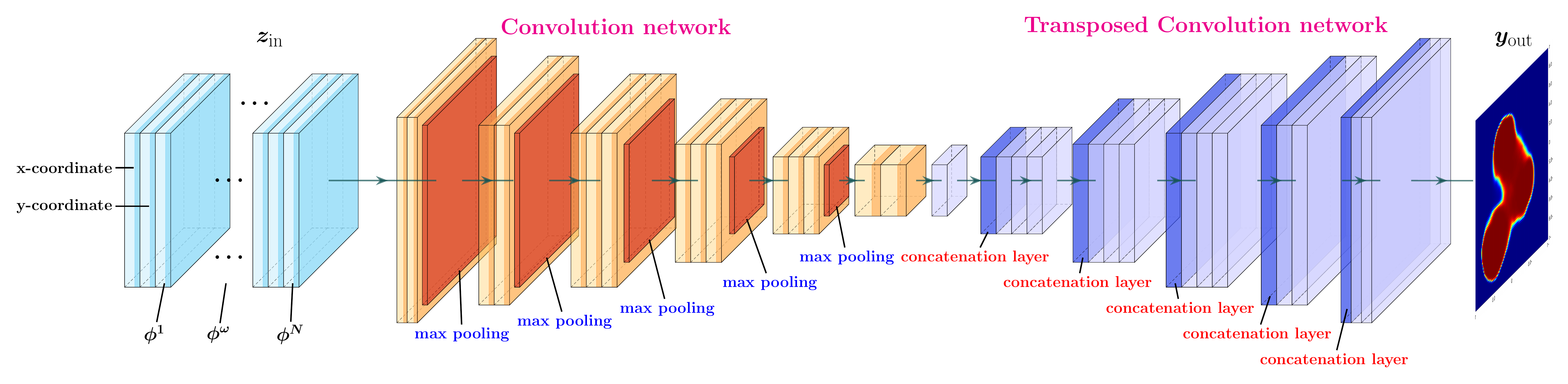}
\caption{The structure of CNN. The input is a 3D matrix in $\mathbb{R}^{n_1\times n_2 \times (N+2)}$ and the output is a matrix in $\mathbb{R}^{n_1\times n_2}$.}
\label{fig:cnnstruct}
\end{figure}

In order to describe the structure of the CNN-DDSM, for {\jjh simplicity we} first assume $\Omega$ has a rectangular shape and leave the general situation to the later discussion about the implementation details. Then we suppose $\Omega$ is discretized by a $n_1\times n_2$ Cartesian grid where $n_1$ and $n_2$ are for the $x_1$ and $x_2$ direction respectively. Based on the previous explanation, the input to the CNN is not a vector but a 3D matrix. In particular, let's focus an inclusion sample with $N$ Cauchy difference functions $\{\phi^\omega\}_{\omega=1}^N$ solved with $\{(g_\omega, f_\omega)\}_{\omega=1}^N$. Then the input denoted by $\bs{z}_{\text{in}} \in \mathbb{R}^{n_x \times n_y \times (N+2)}$ is a stack of $N+2$ matrices in $\mathbb{R}^{n_1\times n_2}${\jjh, where} the first two slices are formed by spatial coordinates $x_1$ and $x_2$ respectively, and the rest $N$ slices correspond to the numerical solutions $\phi^1(x)$, $\phi^2(x)$, ..., $\phi^N(x)$ evaluated at the Cartesian grid points {\jjh and its} pictorial elucidation is provided in Figure \ref{fig:cnnstruct}.

The proposed CNN architecture is composed of convolution networks and transposed convolution networks, where the detailed configuration is illustrated in Figure \ref{fig:cnnstruct}. The convolution part consists of several blocks and each block includes convolution layers, activation layers and max-pooling layers. In particular, the max-pooling layers mainly help in extracting the sharpest features of the input image and reducing the input image size by computing the maximum over each non-overlapping rectangular region. The $i$-th block in the convolution part is expressed as
\begin{equation}
\label{eq:conv_express}
\tau^c_i(\bs{z}_{\text{conv}}) = \mathcal{M}(\zeta(\varrho(W_{\text{conv}}\ast\bs{z}_{i} + \bs{b}_{\text{conv}}))),
\end{equation}
where $*$ denotes the convolution operation, $W_{\text{conv}}$ refers to the convolution filter for the 2D convolutional layer,  $\bs{b}_{\text{conv}}$ denotes the bias, $\bs{z}_{i}$ is the input image, $\mathcal{M}$ is the max-pooling layer, $\zeta$ denotes the activation and $\varrho$ is batch normalization \cite{ioffe2015batch} layer that aims to accelerate the training and reduce the sensitivity of the network initialization. The transposed convolution part also contains several blocks each of which includes a transposed convolution to extrapolate the output of the convolution part to an image with large size (higher resolution). A typical example of the $i$-th block is
\begin{equation}
\label{eq:trconv_express}
\tau^t_i(\bs{z}_{\text{trans}}) = \mathcal{C}(\zeta( {\jjhn \mathcal{T}(\bs{z}_{i}, W_{\text{trans}}, \bs{b}_{\text{trans}})}) ),
\end{equation}
where ${\jjhn \mathcal{T}}$ refers to the transposed convolution {\jjhn operator}, $W_{\text{trans}}$ and $\bs{b}_{\text{trans}}$ are the corresponding {\jjhn transposed convolutional filter and the bias}, $\mathcal{C}$ is a concatenation layer, and other notations are as the same as \eqref{eq:conv_express}. Again $\Theta$ denotes the set of all the unknown parameters {\jjhn including} the convolutional and transposed convolutional filters and {\jjhn the biases} to be learned in training. Notice that convolution and batch normalization layers are also inserted into \eqref{eq:trconv_express} at some blocks. Our experience has suggested choosing sigmoid function as the activation function, 
\begin{equation}
\zeta(z) = \frac{1}{1+e^{-z}}.
\end{equation}
We note that this choice is different from the one used in FNN but it is also motivated by the condition (\textbf{C}) that the non-linearity should involve exponential functions. Then, the full CNN model can be represented as 
\begin{equation}
\label{eq:full-cnn}
\bs{y}_{\text{out}}(\bs{z}_{\text{in}}) = \tau^t_{M_t}\circ\dots\circ\tau^t_{1}\circ\tau^c_{M_c}\circ\dots\circ\tau^c_{1}(\bs{z}_{\text{in}}),
\end{equation}
where the output $\bs{y}_{\text{out}}(\bs{z}_{\text{in}})$ is a $n_1\times n_2$ matrix which is supposed to approximate an inclusion distribution.


To measure the accuracy of the CNN model \eqref{eq:full-cnn}, we employ the mean squared error (MSE) as the loss function
\begin{equation}
\label{eq:loss-mse}
\mathcal{L}_{\text{loss}}(\Theta) = \frac{1}{SK}\sum_{s=1}^{S}\big((\bs{y}_{\text{out}}(\bs{z}^{s}_{\text{in}}) - \mathcal{I}^{s})^\top (\bs{y}_{\text{out}}(\bs{z}^{s}_{\text{in}}) - \mathcal{I}^{s})\big),
\end{equation}
where $\mathcal{I}^{s}$ is the true distribution (index function) corresponding to the $s$-th inclusion sample, $\bs{z}^{s}_{\text{in}}$ is the input image (3D matrix) also corresponding to the $s$-th inclusion sample, namely the Cauchy difference functions $\phi^{\omega}=\phi^{(s,\omega)}$. Similarly, to overcome the infeasibility of gradient descent algorithm when training data size is huge, we apply SGD to find the minimization of the loss function \eqref{eq:loss-mse}
\begin{align}
\Theta^{(j+1)} = \Theta^{(j)} - \alpha \nabla_{\Theta}\mathcal{G}_{\text{loss}}(\Theta^{(j)}) 
\end{align}
where $\alpha$ indicates the learning rate and $j$ is the number of iteration. $\mathcal{G}_{\text{loss}}(\Theta)$ is defined as
\begin{equation}
\mathcal{G}_{\text{loss}}(\Theta) = \frac{1}{|S_r|}\sum_{s\in S_r}\big((\bs{y}_{\text{out}}(\bs{z}^{s}_{\text{in}}) - \mathcal{I}^{s})^\top (\bs{y}_{\text{out}}(\bs{z}^{s}_{\text{in}}) - \mathcal{I}^{s})\big),
\end{equation}
where $S_r$ is a subset of $\{1,2,...,S\}$ randomly chosen for each iteration. 

\vspace{0.1in}

\textbf{Implementation Issues}. Different from FNN-DDSM, the discretization points of $\Omega$ have to be chosen and fixed a-priori. In order to perform convolution, they have to be Cartesian grid points which are natural for rectangular domain. For domain with general shape, we only need to immerse it into a rectangle such that the Cartesian grid can be generated on the whole rectangle. {\grc If Cauchy difference functions $\phi^{\omega}$ are computed by finite element methods on a general triangulation of $\Omega$, then the values $\phi^{\omega}$ at the mesh points can not be directly used for CNN computation though they are already available. Instead we need to recompute their values at the newly generated Cartesian grid points and fill the outside region by zeros.}

\section{Numerical Expeiments}
\label{sec:num}

In this section, we present numerical experiments to demonstrate that our  newly proposed DDSMs are effective and robust for the reconstruction of inhomogeneous inclusions in the EIT problem.

\subsection{Problem Setting and Data Generation}

Let the modeling domain be $\Omega=(-1,1)\times(-1,1)$ which contains two medium with the different conductivity $10$(inclusion) and $1$(background). {\grc Let $\mathcal{U}(a,b)$ be uniform distribution in $[a,b]$.} To verify the efficacy of DDSMs for general inclusion distribution, we explore the following three typical scenarios:


\begin{itemize}
 \setlength{\itemindent}{0.5in}
\item[\textbf{Scenario 1}:] the inclusions are generated by three random circles with the radius sampled from $\mathcal{U}(0.2,0.4)$;

\item[\textbf{Scenario 2}:] the inclusions are generated by five random circles with the radius sampled from $\mathcal{U}(0.2,0.3)$;

\item[\textbf{Scenario 3}:] the inclusions are generated by four random ellipses with the length of the semi-major axis and semi-minor axis sampled from $\mathcal{U}(0.1,0.2)$ and $\mathcal{U}(0.2,0.4)$, respectively, and the rotation angle sampled from $\mathcal{U}(0,2\pi)$. 
\end{itemize}
It is reasonable to require that the inclusions do not touch the boundary \cite{fan2020solving}, so we assume the circles/ellipses have at least distance $0.1$ to the boundary. More precisely, they are uniformly sampled in the square $(-0.9,0.9)\times(-0.9,0.9)$. This sampling strategy has been widely used in solving EIT with deep learning \cite{fan2020solving,he2016deep}. However, as the major difference from these literatures, to make the shape of the inclusions more general and various, we do not require that these circles and ellipses are disjoint from each other, that is, they are free to touch each other. This mechanism will generate much more complicated shapes than basic geometric components, which makes the reconstruction more arduous.  
Let the circles or ellipses be represented by $c_i(x_1,x_2)=0$, $i=1,2,...,N_c$ where $c_i(x_1,x_2)$ are the related level-set functions and $N_c$ denotes the number of circles/ellipses in each configuration. Then the level-set function of their union can be defined as
\begin{equation}
\label{union_set}
c(x_1,x_2) = \min_{i=1,..,N_c}\{c_i(x_1,x_2) \}.
\end{equation}
The homogenous background with the conductivity $1$ fills the subdomain $\Omega\cap\{(x_1,x_2): c(x_1,x_2)>0\}$ and the inhomogeneous inclusions with the conductivity $10$ fills the subdomain $\Omega\cap\{(x_1,x_2): c(x_1,x_2)<0\}$. 
For training set,  \textbf{Scenario 1} and \textbf{2} both have $11200$ samples, and \textbf{Scenario 3} has 14400 samples. For testing set, we use 2000 samples for all the three scenarios. 
For the inserted current data $g_{\omega}$ on the boundary, we follow the idea in \cite{chow2014direct}, 
but here we use  multiple terms,
\begin{equation}
\label{g_cos}
g_{\omega}(x) = \cos(\omega\theta(x)), ~~ x\in\partial\Omega, ~~ \omega=1,2,...,N,
\end{equation}
where $\theta(x)$ is the polar angle of $x$, and $N$ is chosen to be $1,10$ or $20$ to get three different training sets for each scienario. In order to generate the synthetic data, i.e., $u_{\omega}|_{\partial\Omega}$, we need to repeatedly solve the modeling equations \eqref{govern_eq_0}-\eqref{govern_eq_3} with different distribution of discontinuous conductivity coefficients. Here we employ the immersed finite element method (IFEM) \cite{2016GuoLin,2018GuoLinZhuang} which does not require the mesh to resolve the conductivity discontinuity, to efficiently solve all the equations on the same $200\times200$ Cartesian mesh. {\jjhn Standard finite element methods can also be used on the mesh generated by \eqref{union_set}.}
Then we collect all the values of $u_\omega$ at the mesh points on boundary, and use them to generate $\phi^{\omega}$, $\omega=1,2,...,N$, on mesh points by solving \eqref{eq:phiproblem} for each inclusion sample.

According to the practical experience \cite{2019LiZhouWangWangLu}, only very limited number of real data samples can be obtained in the experimental environments. \cite{2019LiZhouWangWangLu} suggests training the DNNs with simulation data rather than experimental data since the simulation model is not subject to the objective factors.
However, as far as we know there are no general rules or empirical suggestions about whether or what types of noise should be added into the training data. On one hand,
as mentioned in \cite{2017MartinChoi}, the noise used in training must be similar to the physiological noise presented in human data varying across patients and hardware systems and usually unknown in practice. On the other hand, the conventional DSM \cite{chow2014direct} is developed in a manner {\grc that the noise on the boundary is smoothed out in the duality product \eqref{eq:dual}}, and thus it can handle relatively large noise in reconstruction. Therefore, due to the insufficient knowledge about the noise and inspired by the robustness feature of the DSM, we here do not include noise in the training set and instead add very large noise in the test set which is similar to the strategy in \cite{2019LiZhouWangWangLu}. This is intentionally to test the robustness of the proposed DNNs with respect to noise if it is not presented in the training set. Then we shall see below that our DNNs can actually handle very large noise in the test set even without any denoising procedure. The noising Cauchy data in the test set are generated by
\begin{equation}
\label{noise_data}
u^{\delta}_{\omega} = ( 1+ \delta G_\omega )u_\omega, ~~~ \omega=1,2,...,N
\end{equation}
for each inclusion sample, where $G_{\omega}$ are Gaussian random variables with zero mean and unity variation, and $\delta=0, 10\%, 20\%$ controls the signal-to-noise ratio. These noising data will be used to compute the noising image features $\phi^{\delta}_\omega$ used for both FNN-DDSM and CNN-DDSM.

\subsection{Numerical Results}

\begin{figure}[htbp]
\centering
\begin{subfigure}{.32\textwidth}
     \includegraphics[width=2.2in]{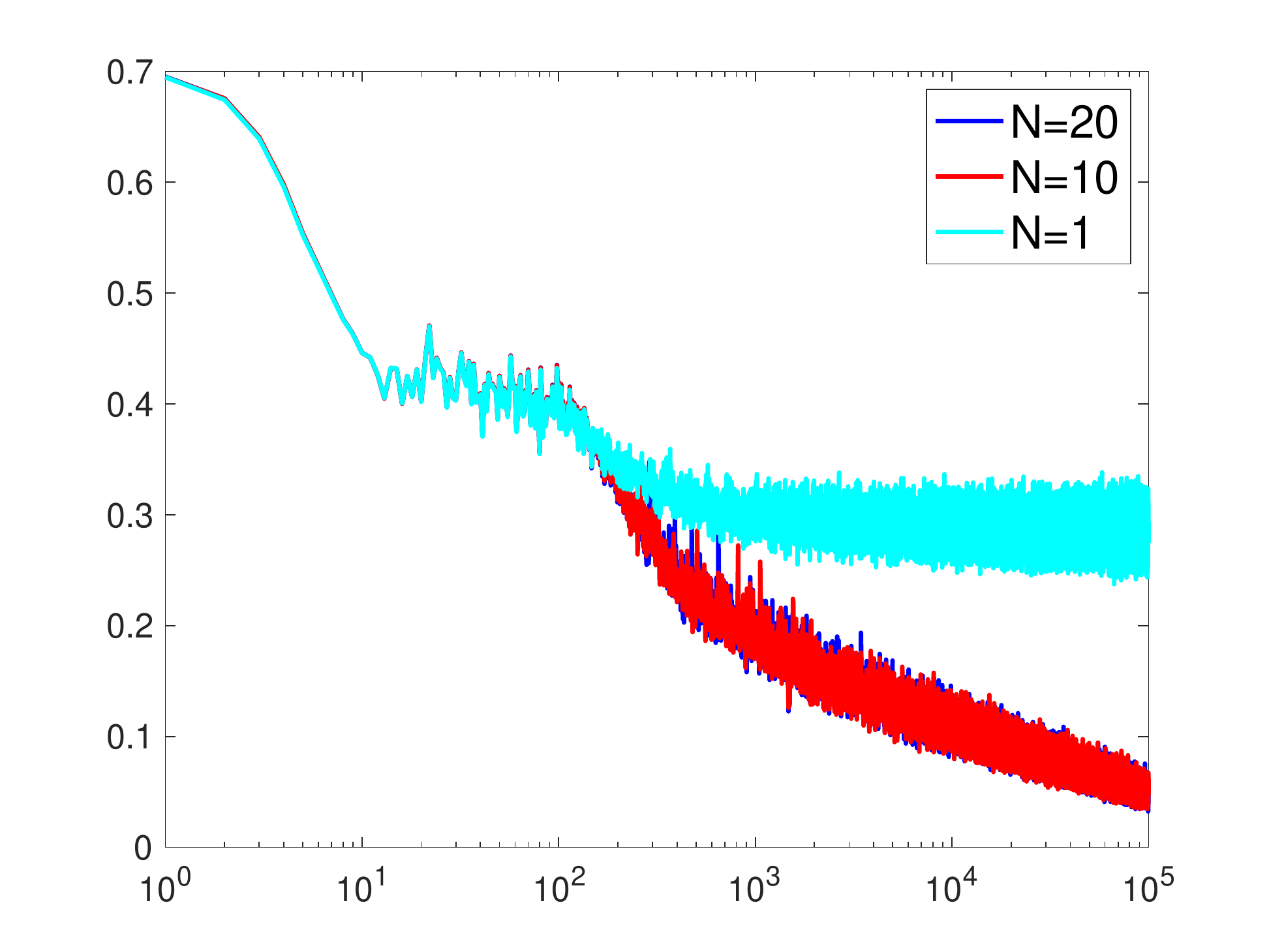}
\end{subfigure}
\begin{subfigure}{.32\textwidth}
     \includegraphics[width=2.2in]{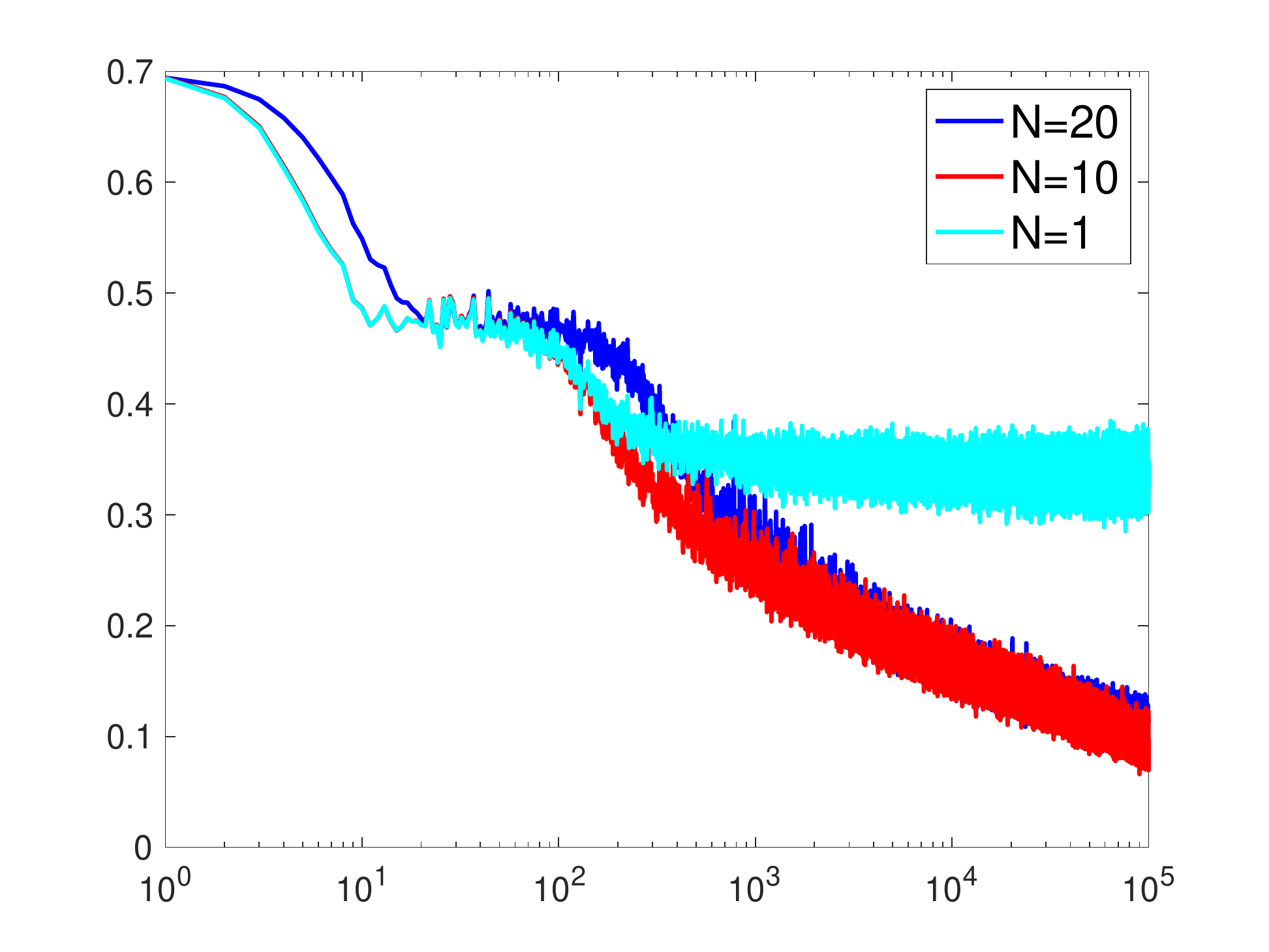}
\end{subfigure}
\begin{subfigure}{.32\textwidth}
     \includegraphics[width=2.2in]{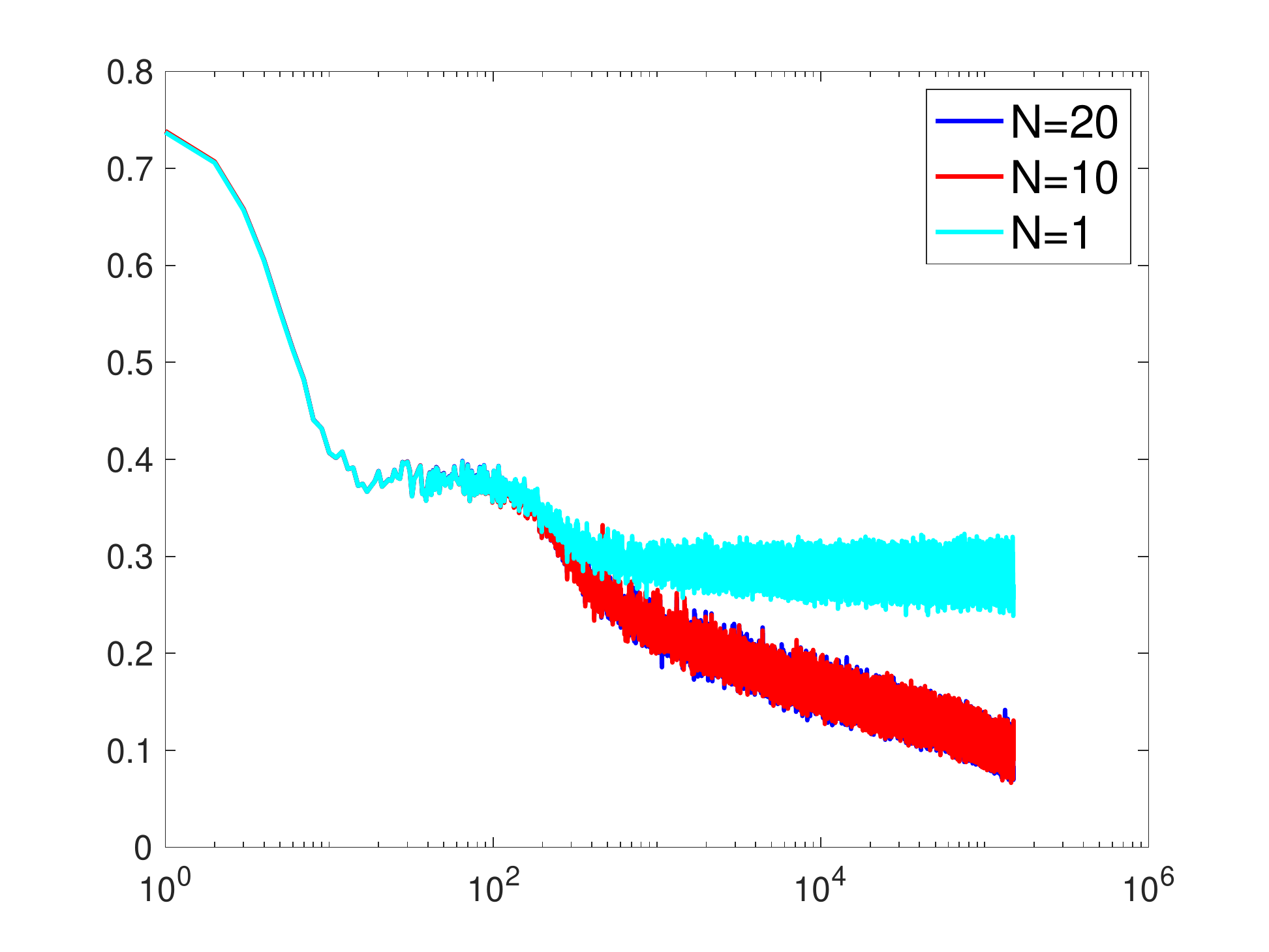}
\end{subfigure}
\caption{Evolution of cross entropy loss values versus training iterations of the FNN-DDSM (from left to right) for the \textbf{Scenario 1}(3 circles), \textbf{Scenario 2}(5 circles) and \textbf{Scenario 3}(4 ellipses)}
  \label{fig:FN_train_process} 
\end{figure}

\begin{figure}[htbp]
\centering
\begin{subfigure}{.32\textwidth}
     \includegraphics[width=2.2in]{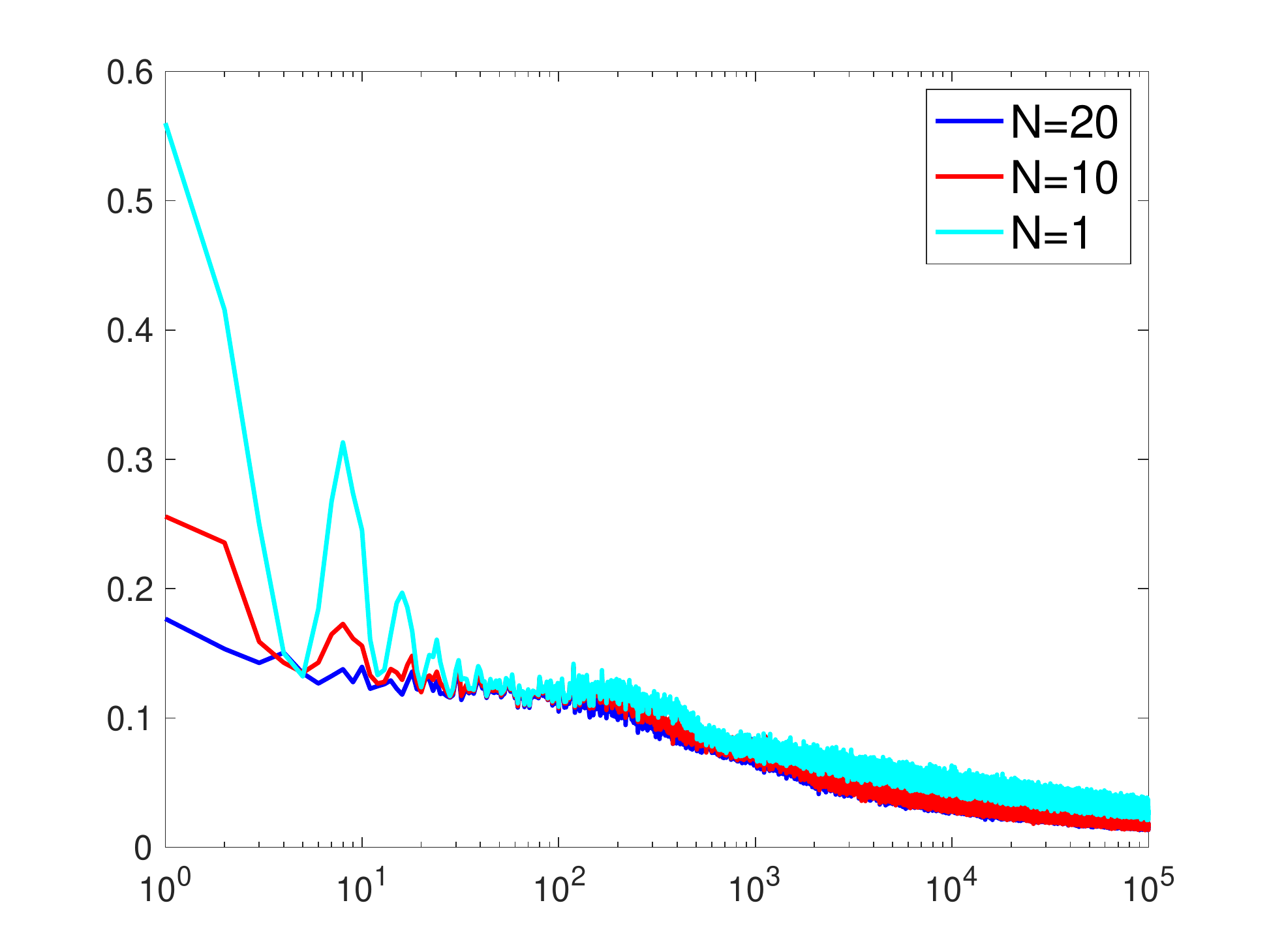}
\end{subfigure}
\begin{subfigure}{.32\textwidth}
     \includegraphics[width=2.2in]{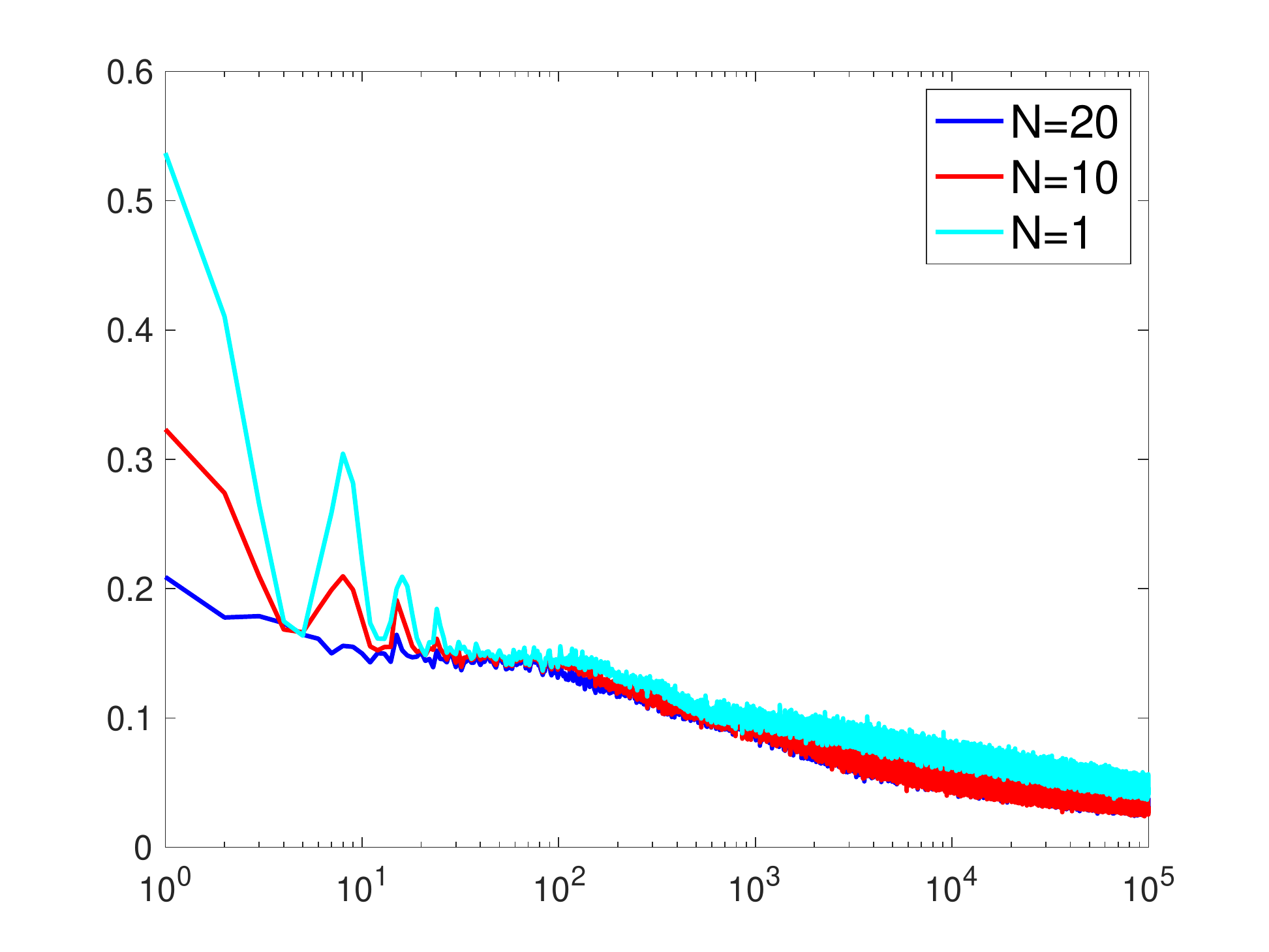}
\end{subfigure}
\begin{subfigure}{.32\textwidth}
     \includegraphics[width=2.2in]{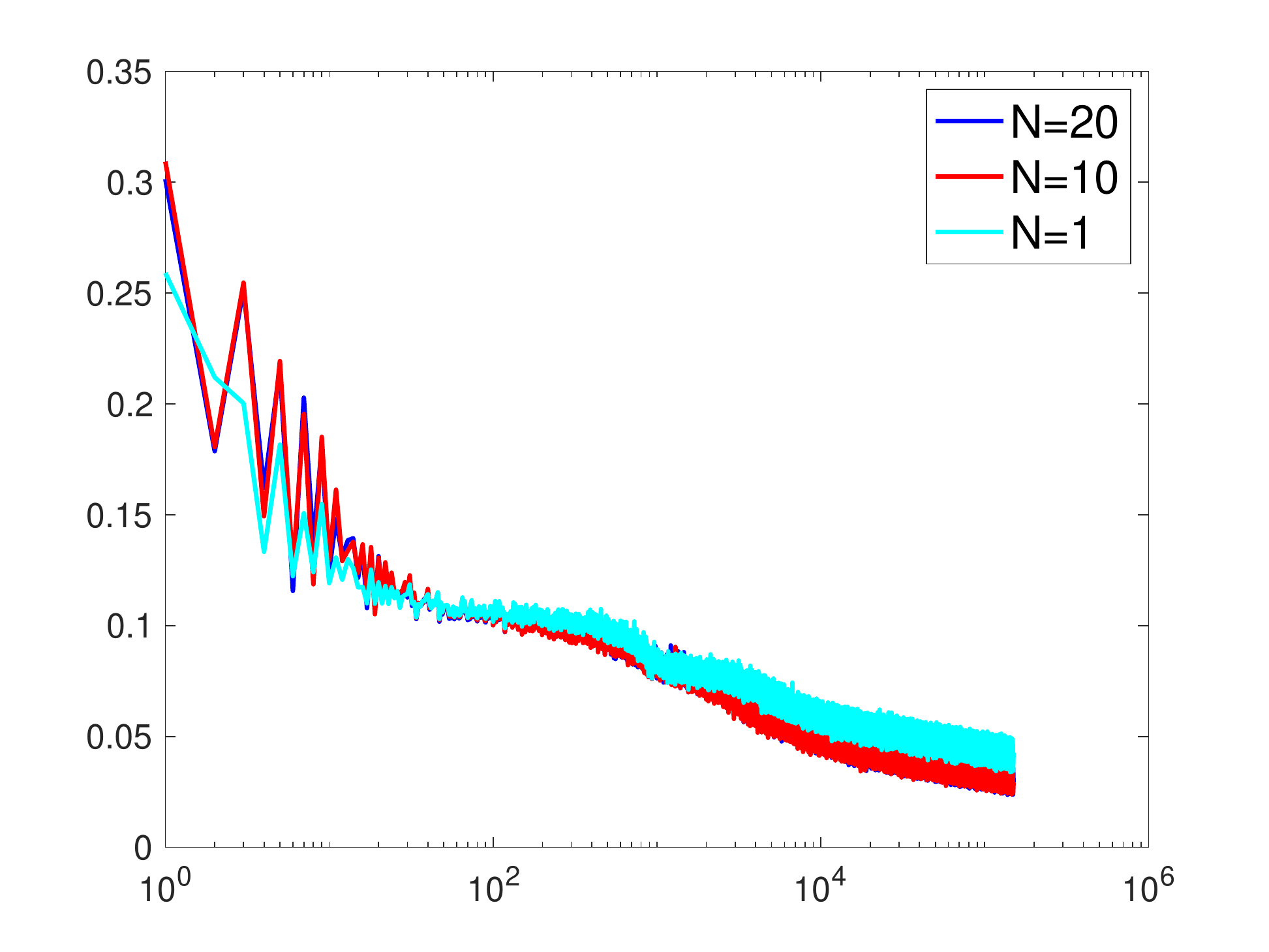}
\end{subfigure}
\caption{Evolution of MSE loss values versus training iterations of the CNN-DDSM (from left to right) for the \textbf{Scenario 1}(3 circles), \textbf{Scenario 2}(5 circles) and \textbf{Scenario 3}(4 ellipses)}
  \label{fig:train_process} 
\end{figure}

\begin{figure}[htbp]
\begin{tabular}{ >{\centering\arraybackslash}m{0.9in} >{\centering\arraybackslash}m{0.9in} >{\centering\arraybackslash}m{0.9in}  >{\centering\arraybackslash}m{0.9in}  >{\centering\arraybackslash}m{0.9in}  >{\centering\arraybackslash}m{0.9in} }
\centering
True coefficients &
N=1, $\delta=0$&
N=10, $\delta=0$&
N=20, $\delta=0$&
N=20, $\delta=10\%$ &
N=20, $\delta=20\%$ \\
\includegraphics[width=1.1in]{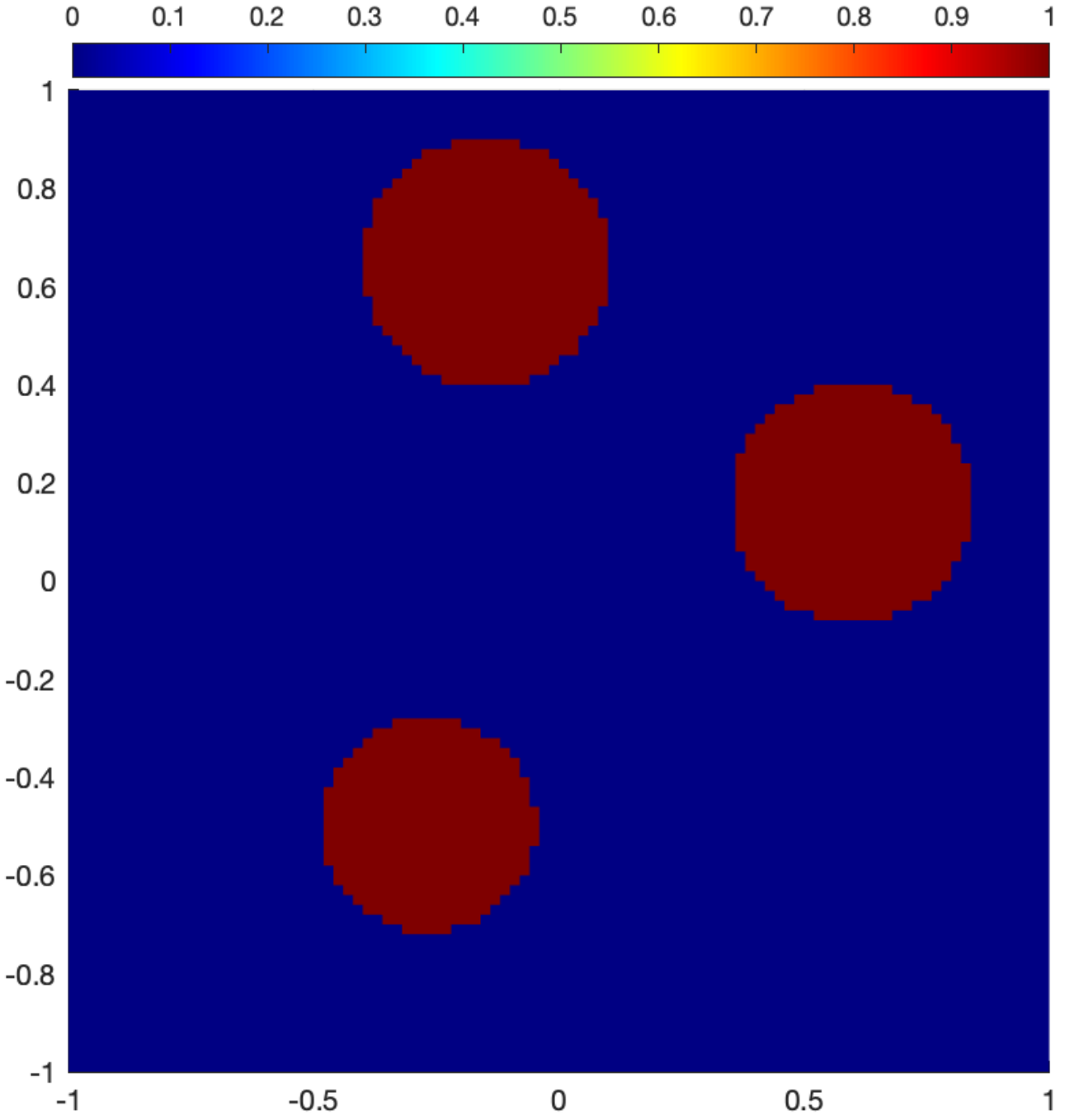}&
\includegraphics[width=1.1in]{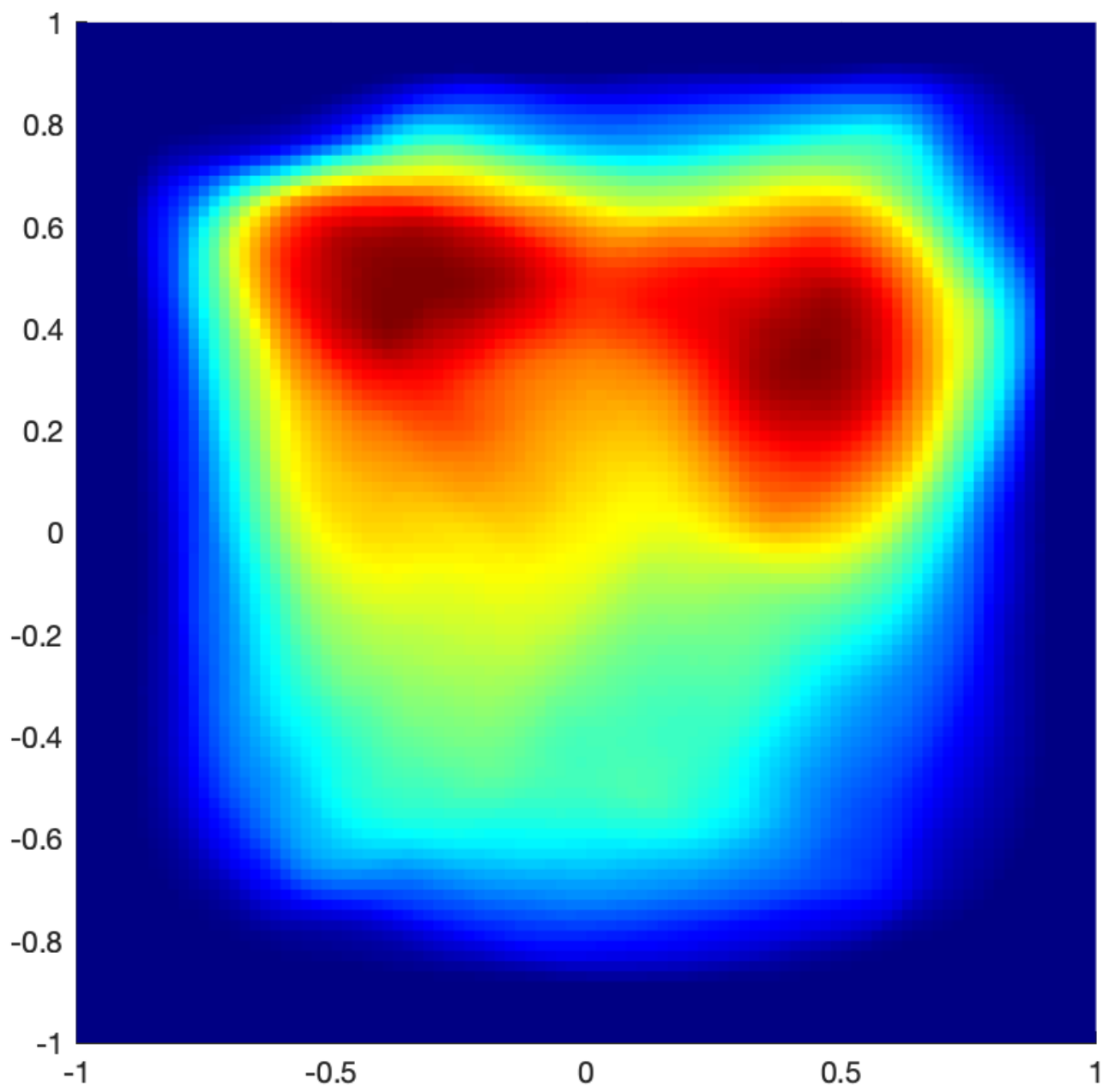}&
\includegraphics[width=1.1in]{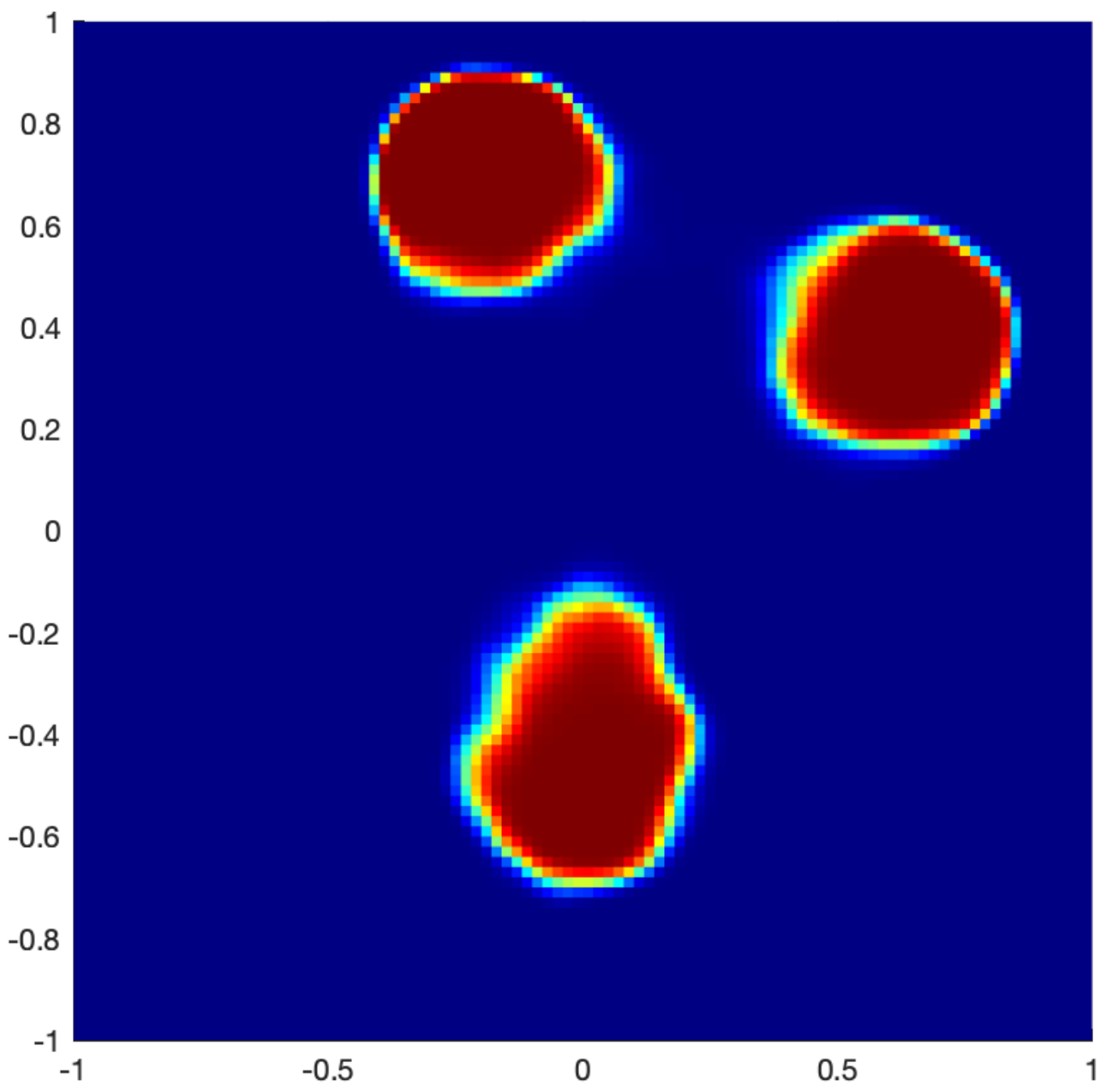}&
\includegraphics[width=1.1in]{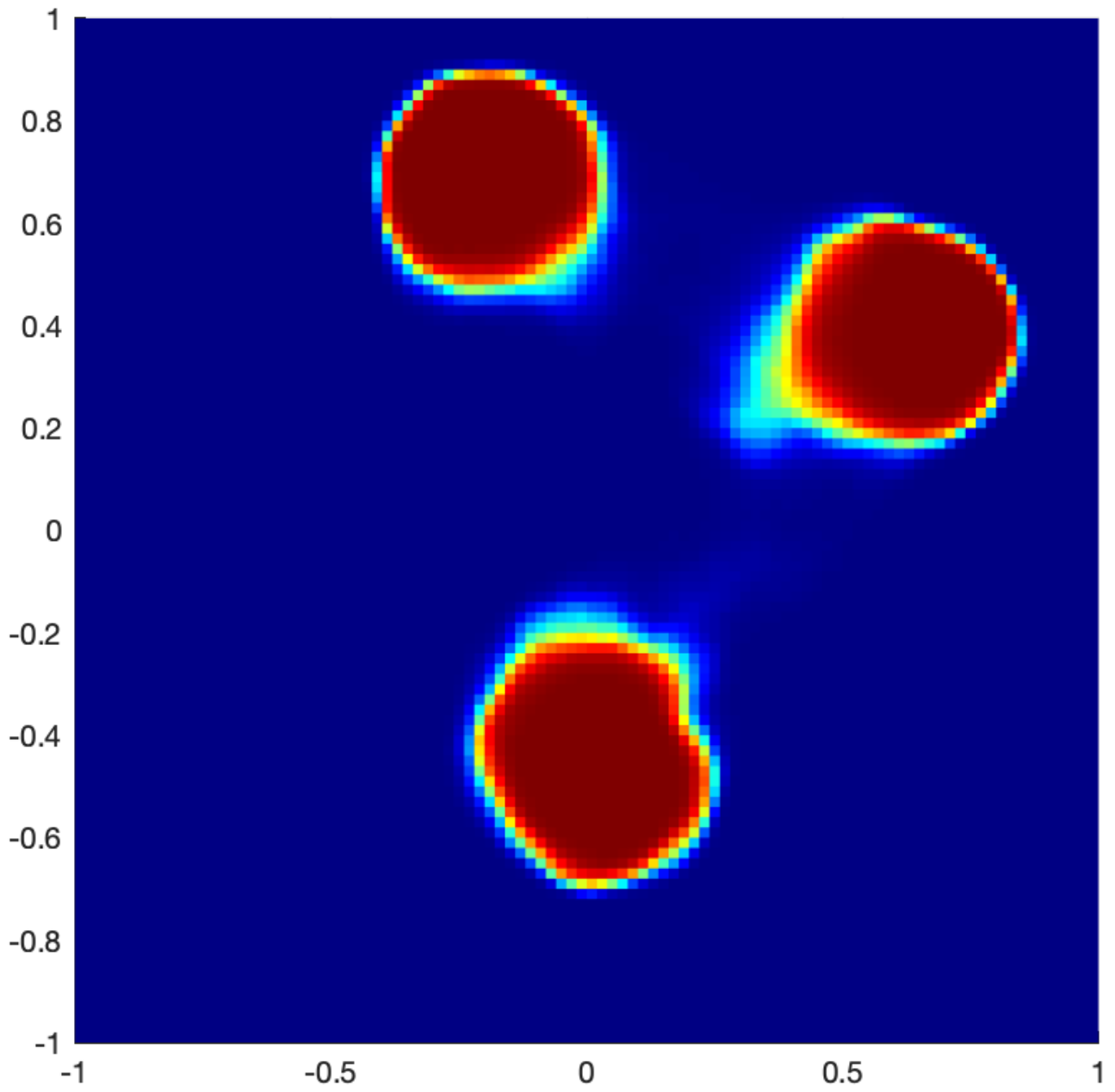}&
\includegraphics[width=1.1in]{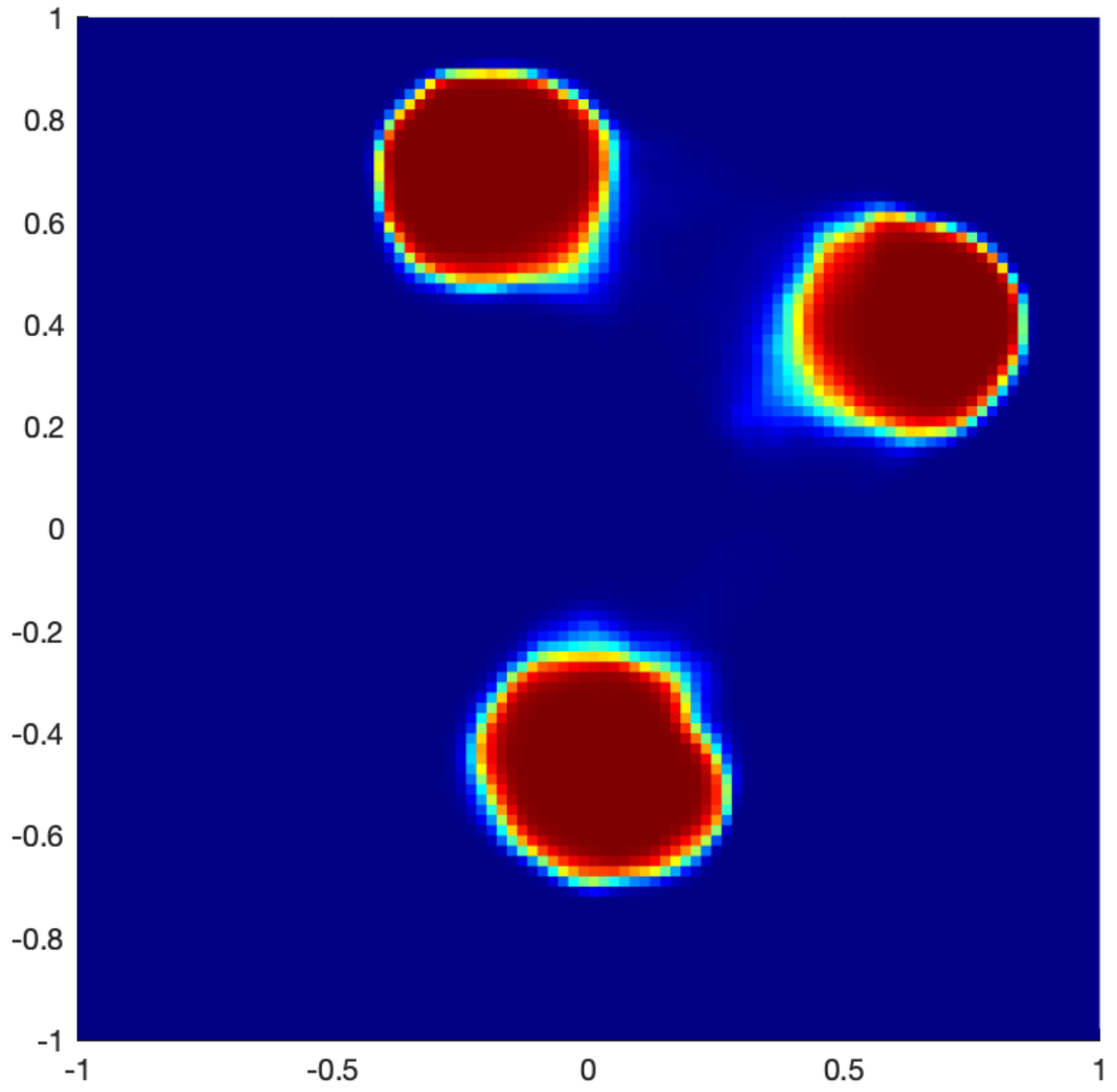}&
\includegraphics[width=1.1in]{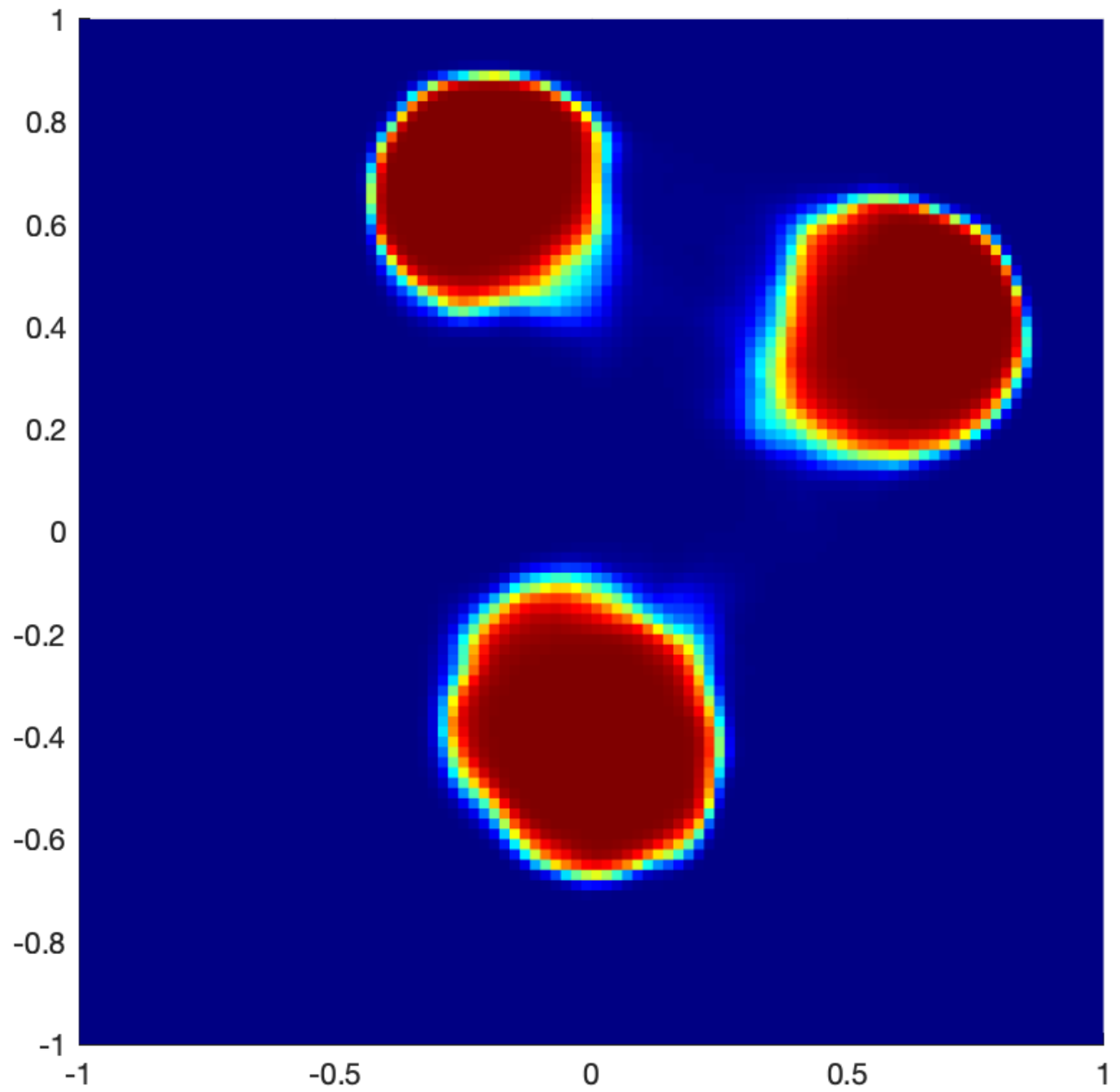}\\
\includegraphics[width=1.1in]{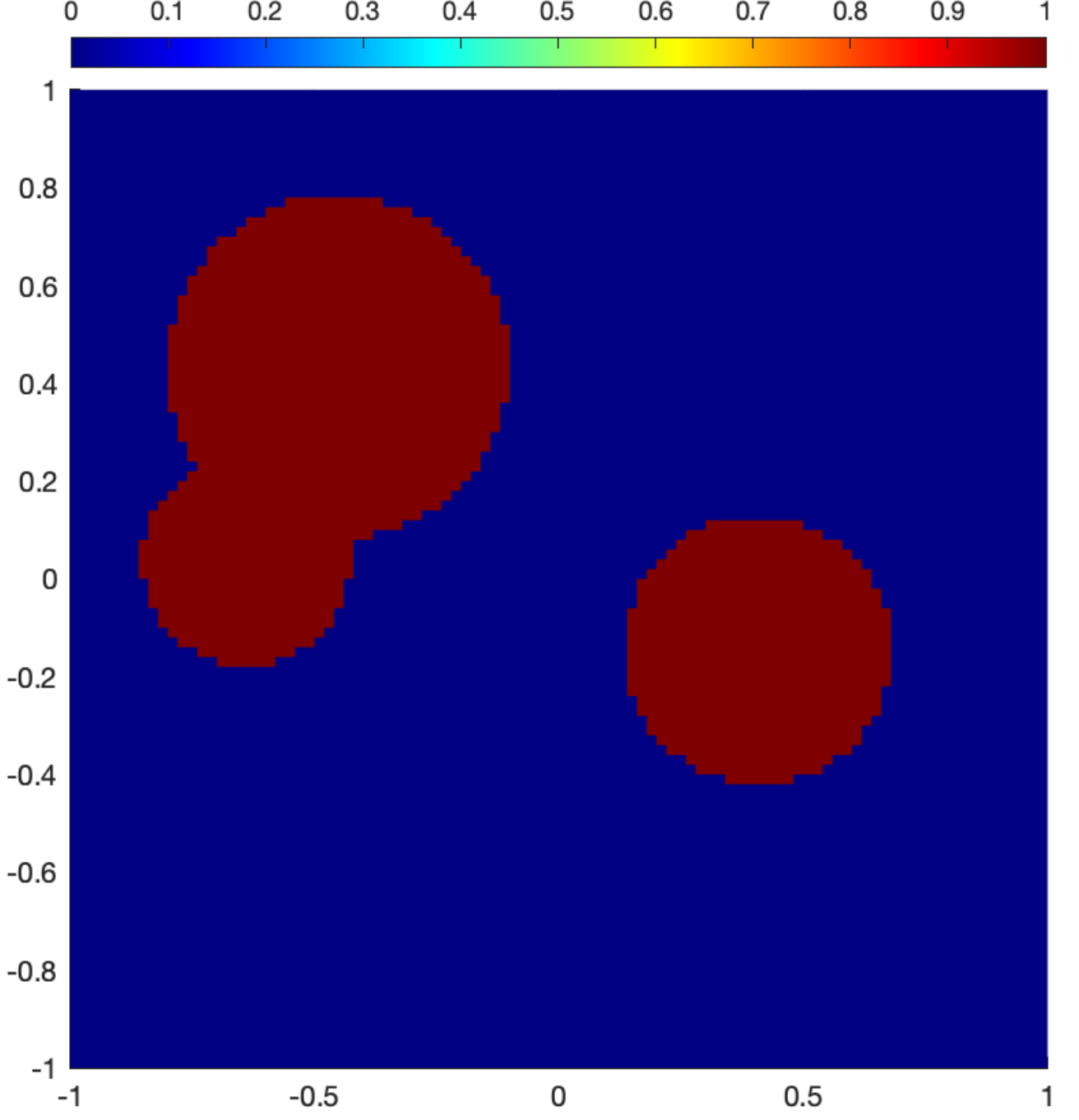}&
\includegraphics[width=1.1in]{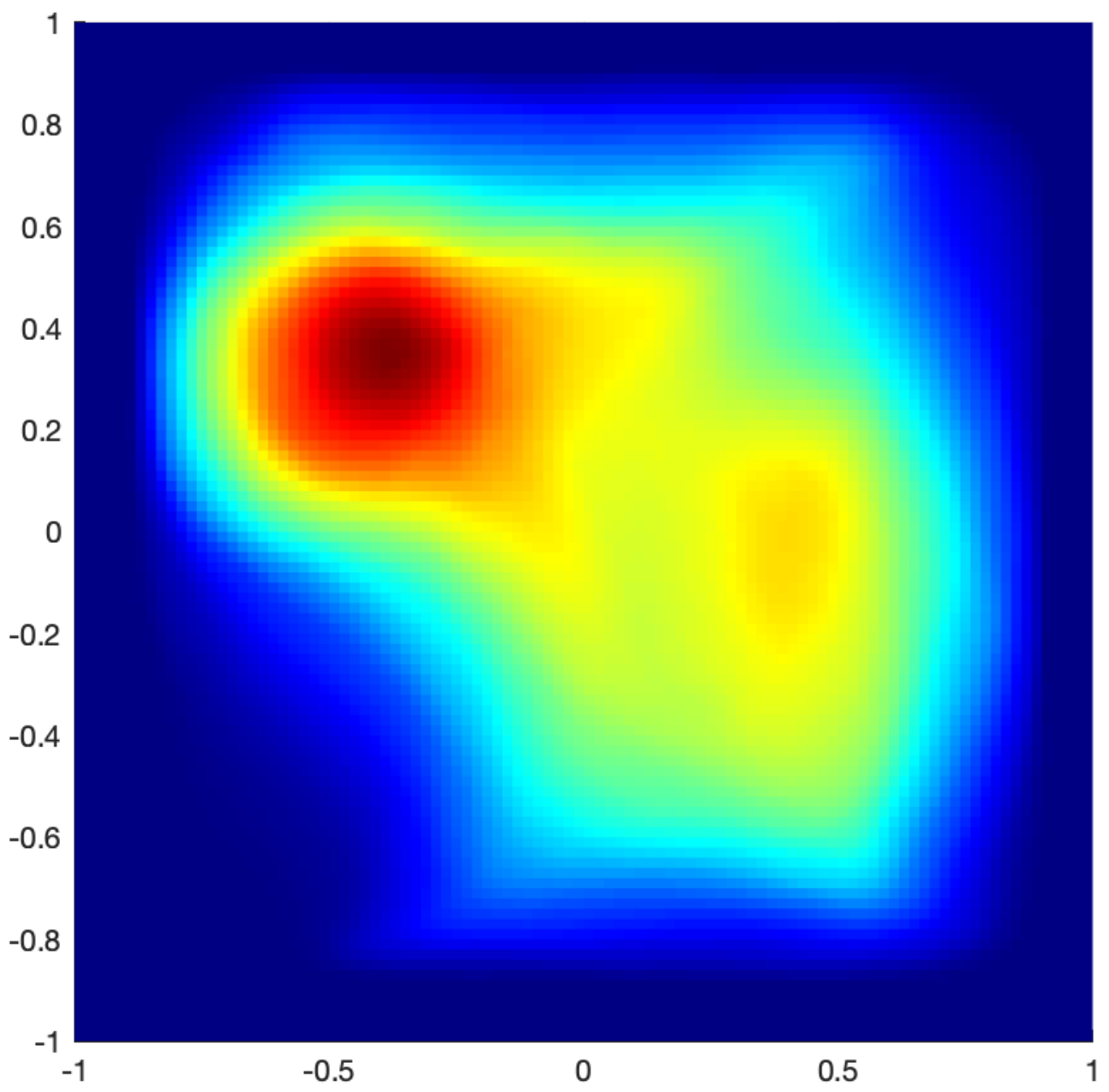}&
\includegraphics[width=1.1in]{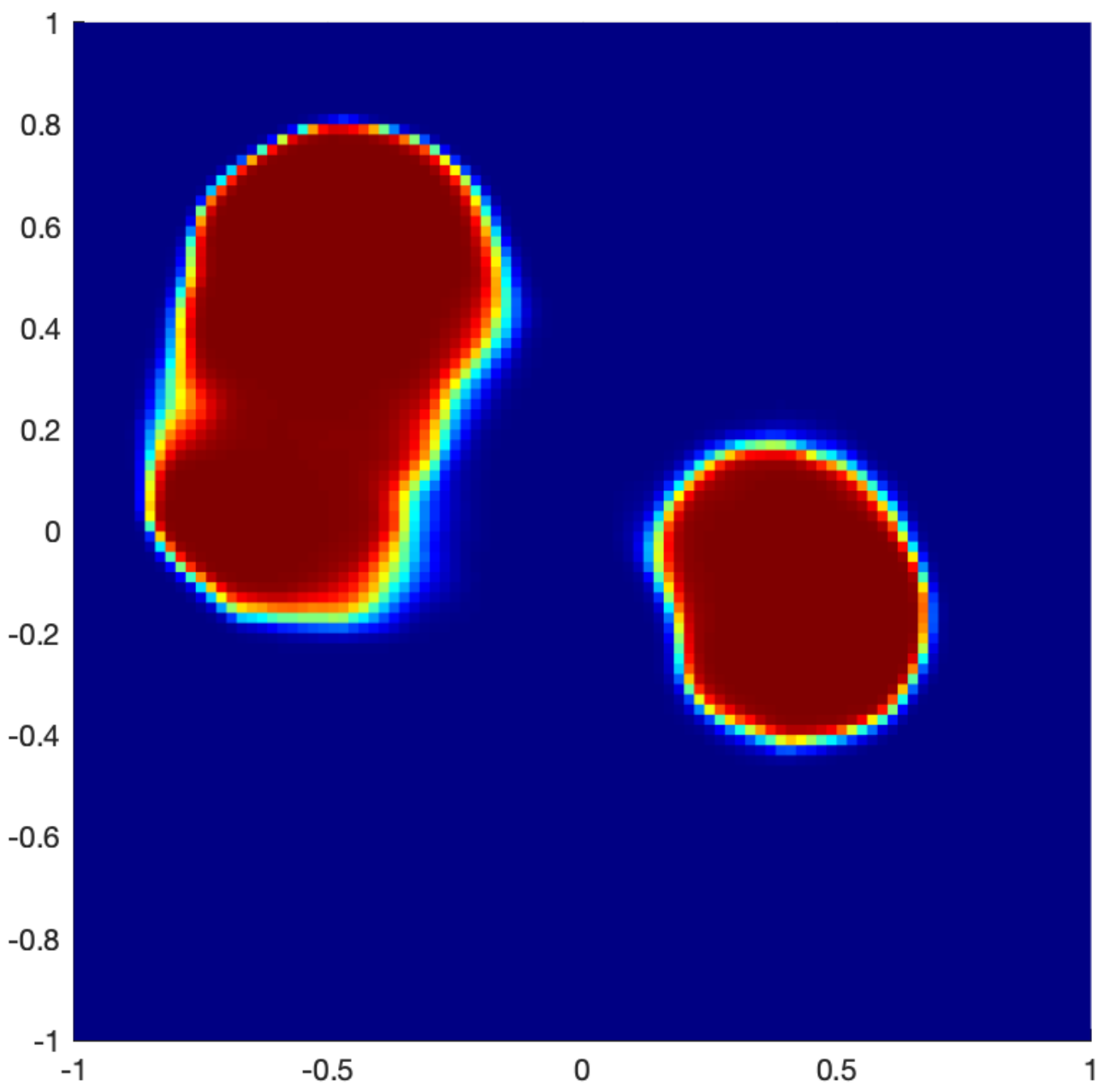}&
\includegraphics[width=1.1in]{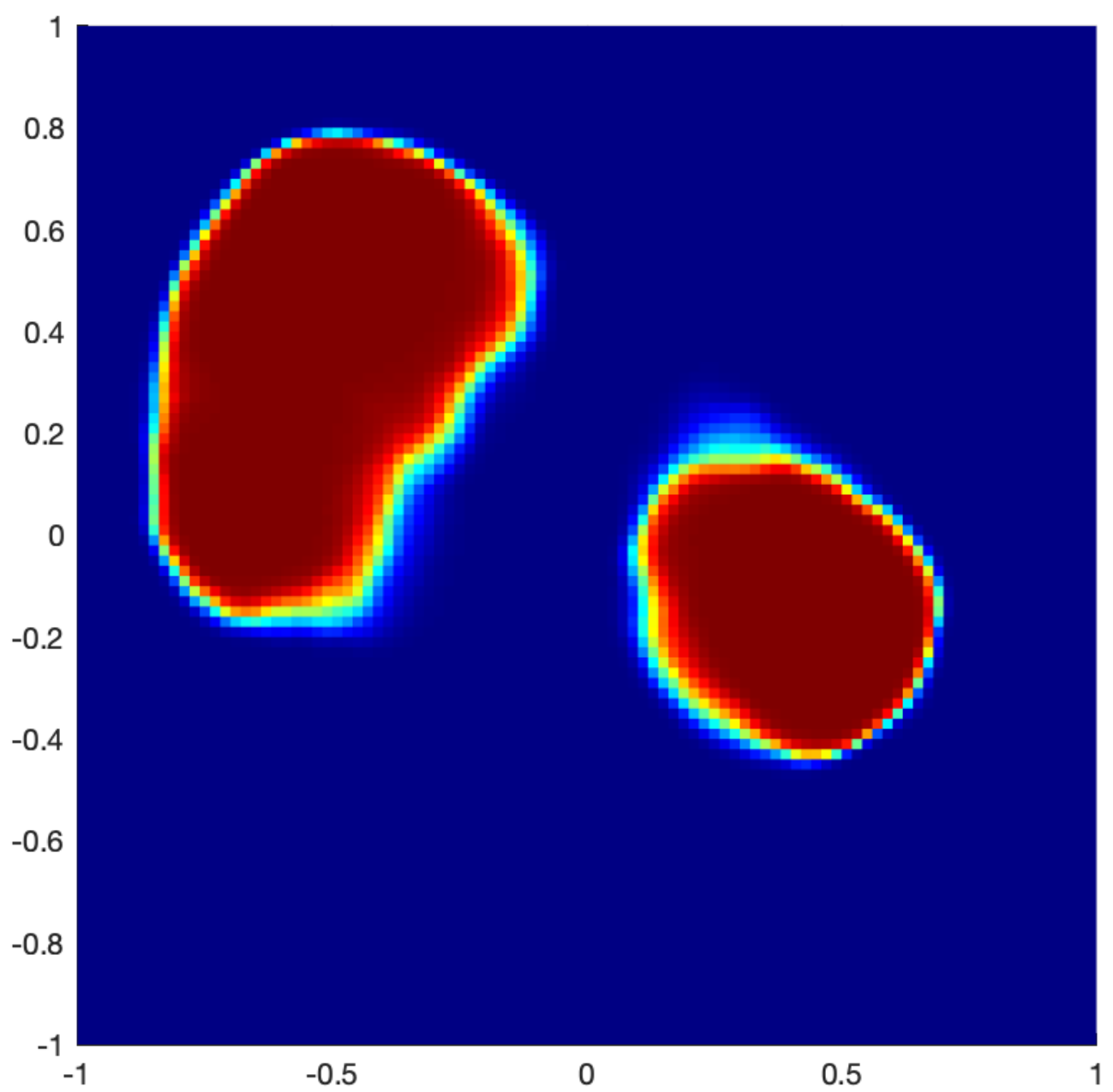}&
\includegraphics[width=1.1in]{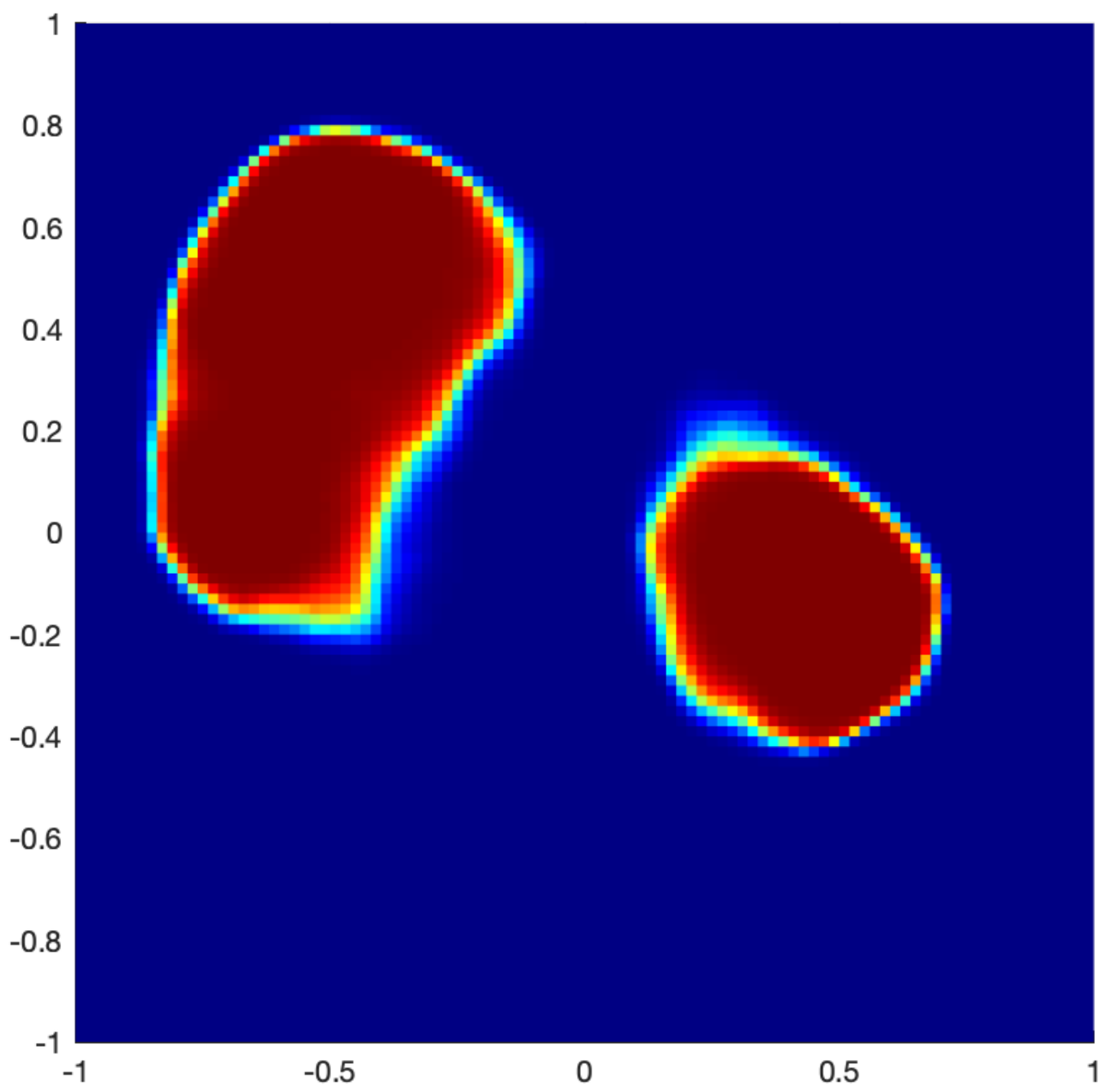}&
\includegraphics[width=1.1in]{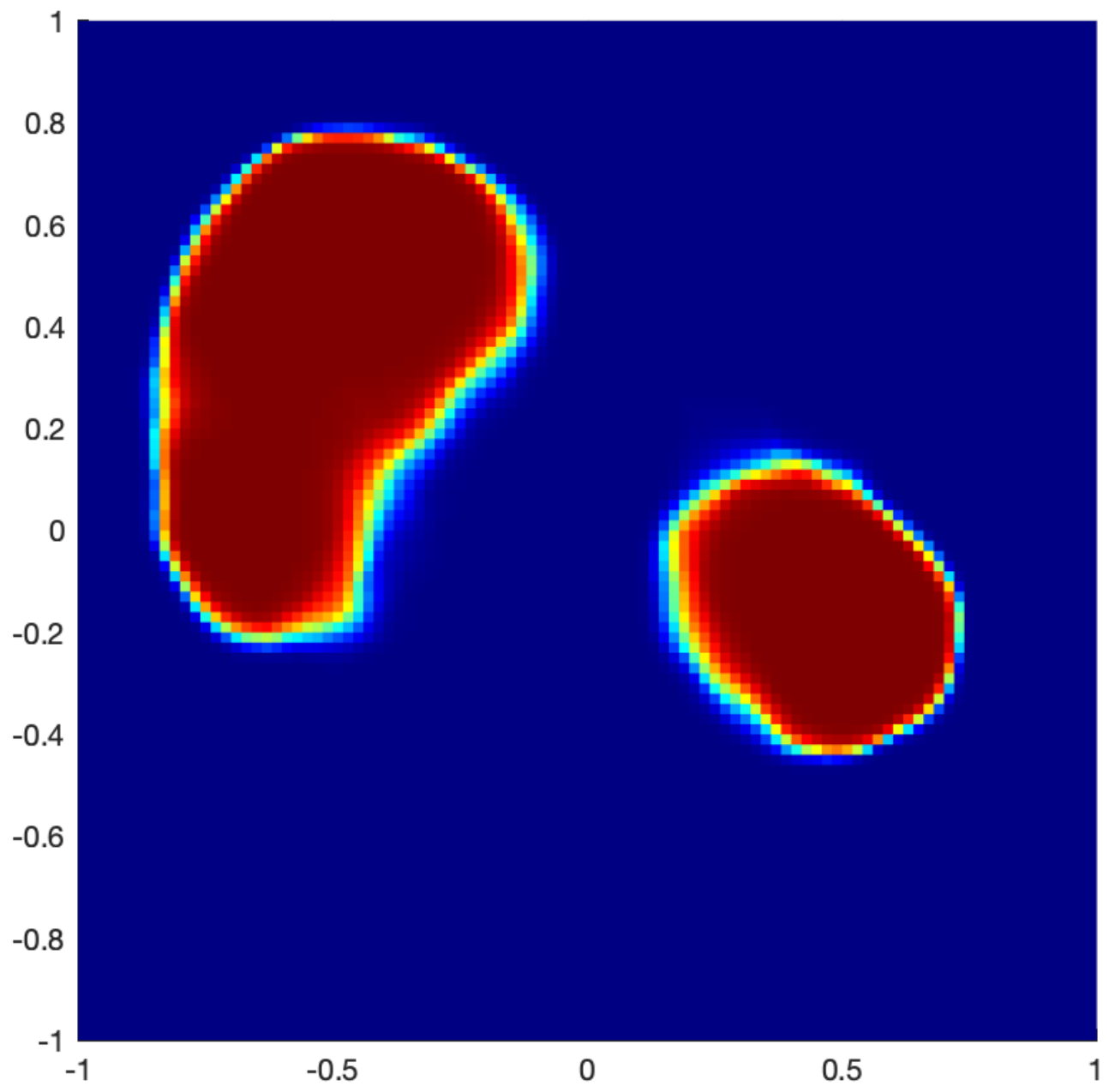}\\
\includegraphics[width=1.1in]{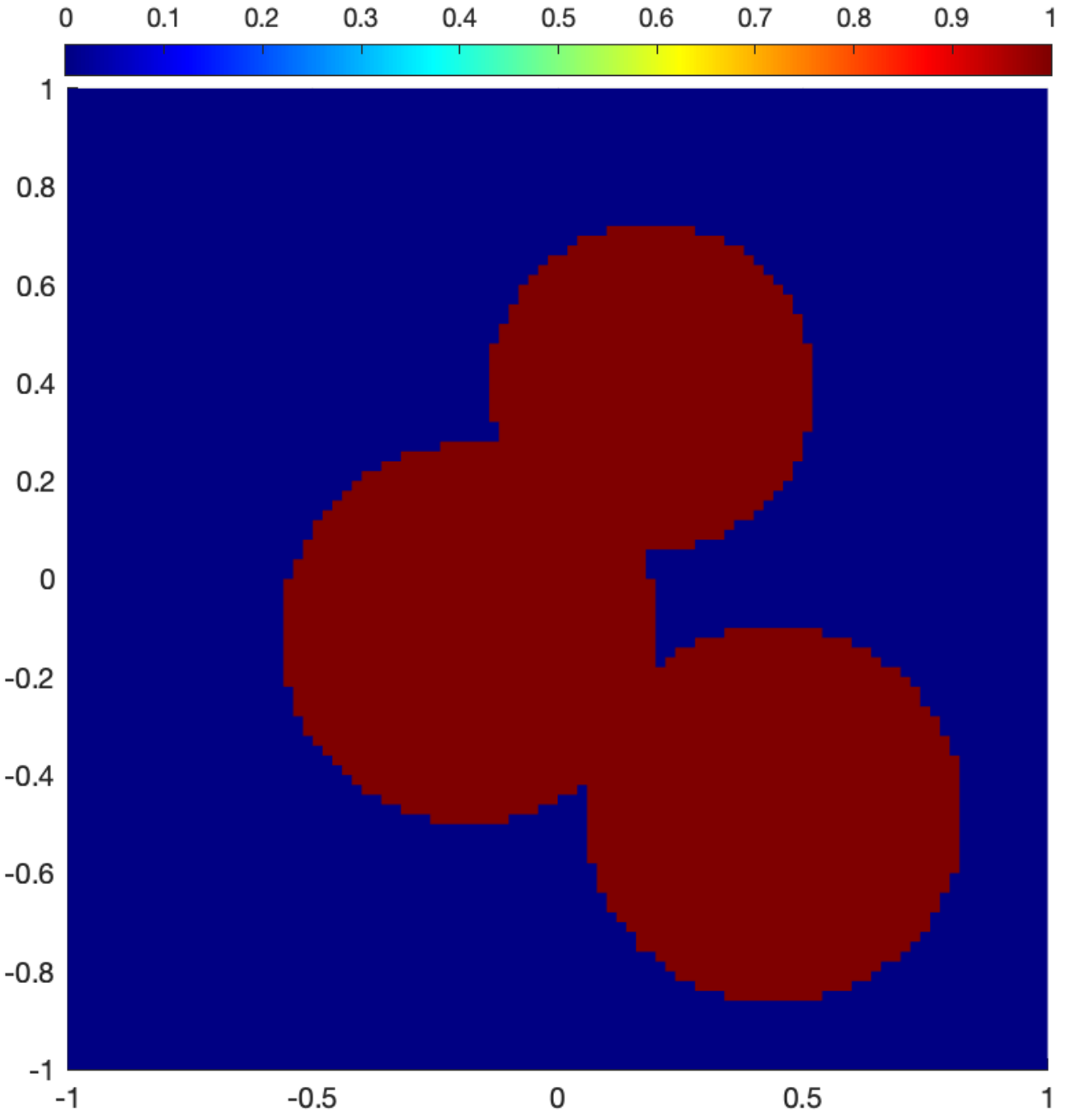}&
\includegraphics[width=1.1in]{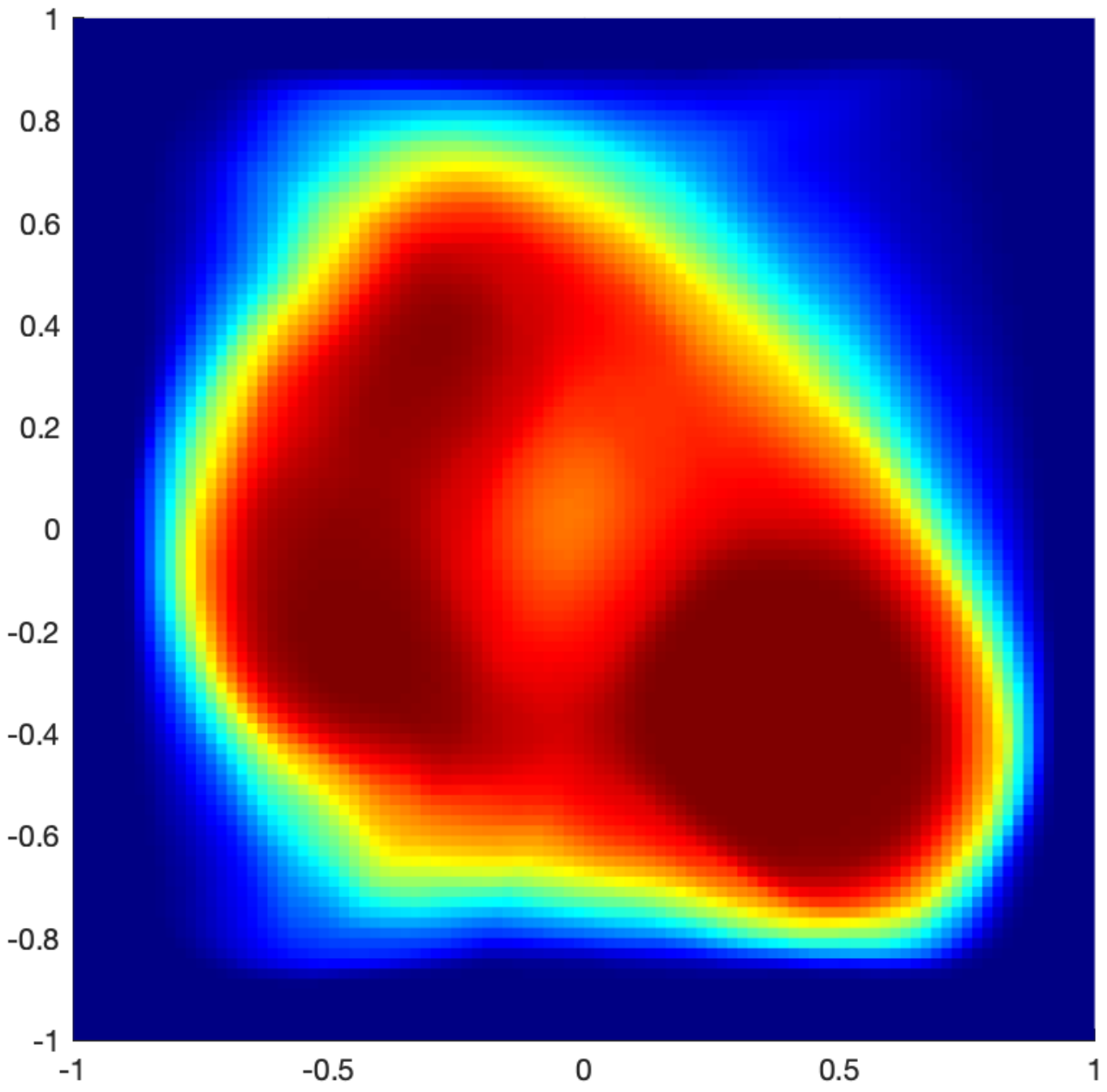}&
\includegraphics[width=1.1in]{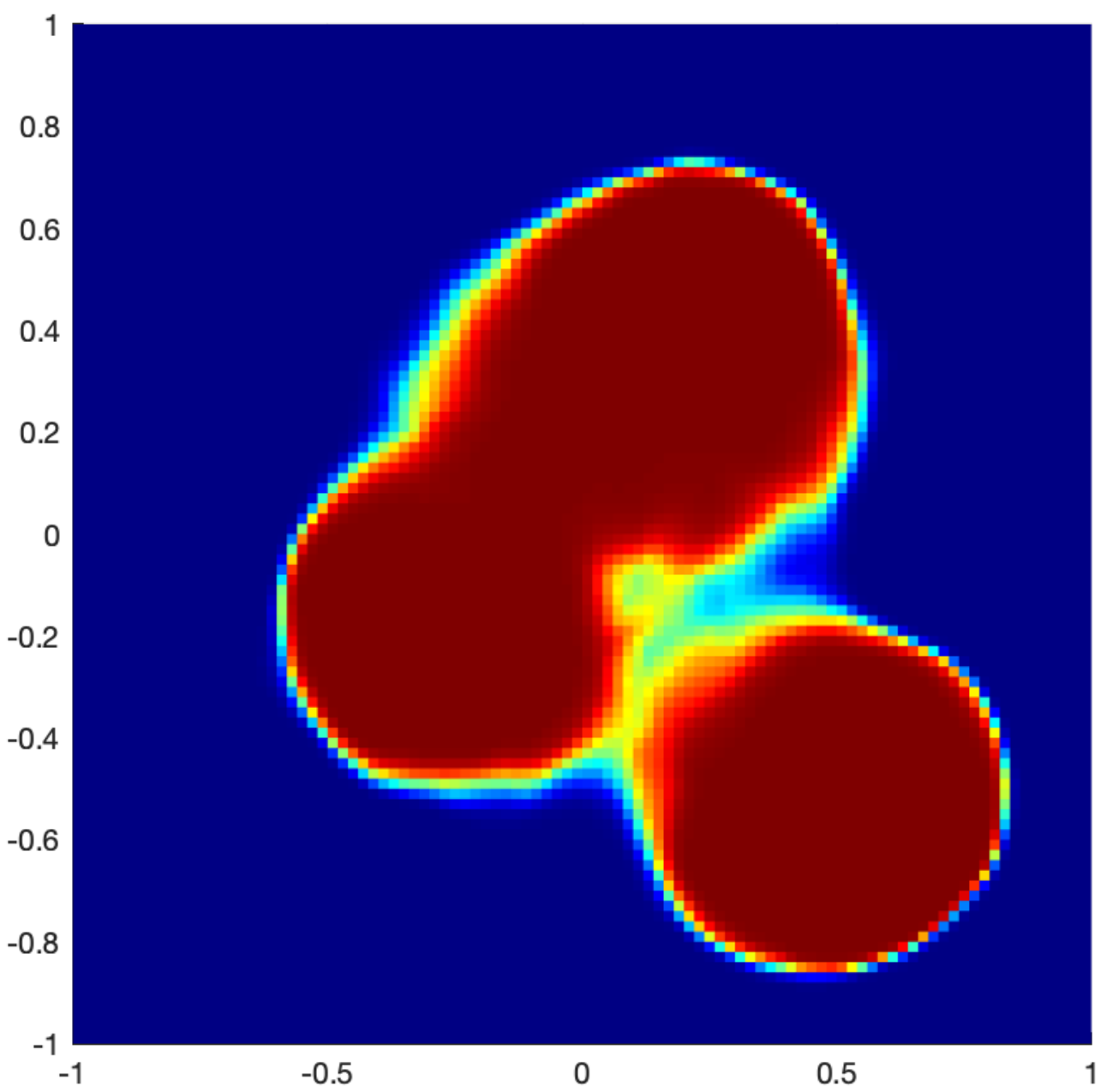}&
\includegraphics[width=1.1in]{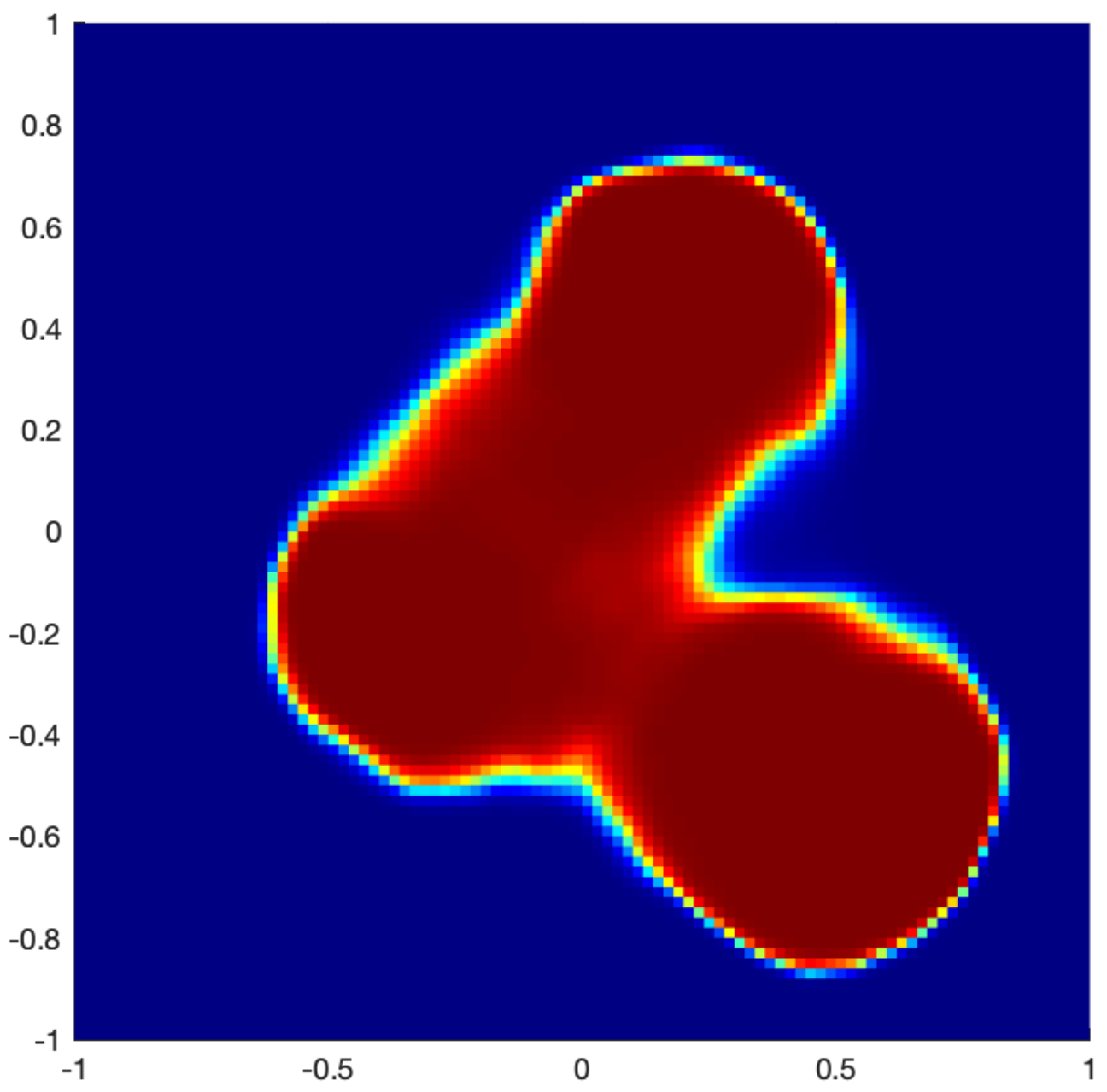}&
\includegraphics[width=1.1in]{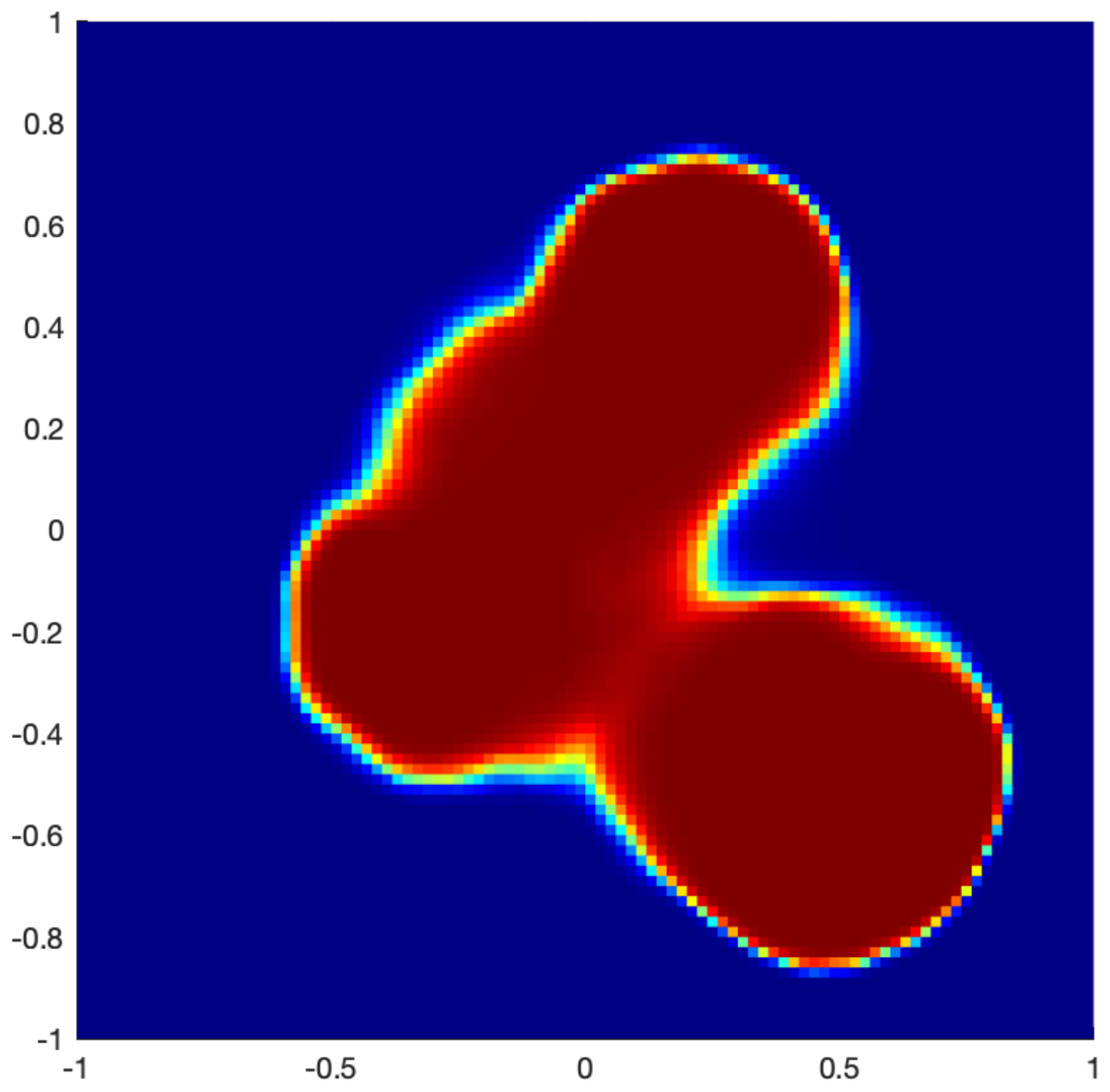}&
\includegraphics[width=1.1in]{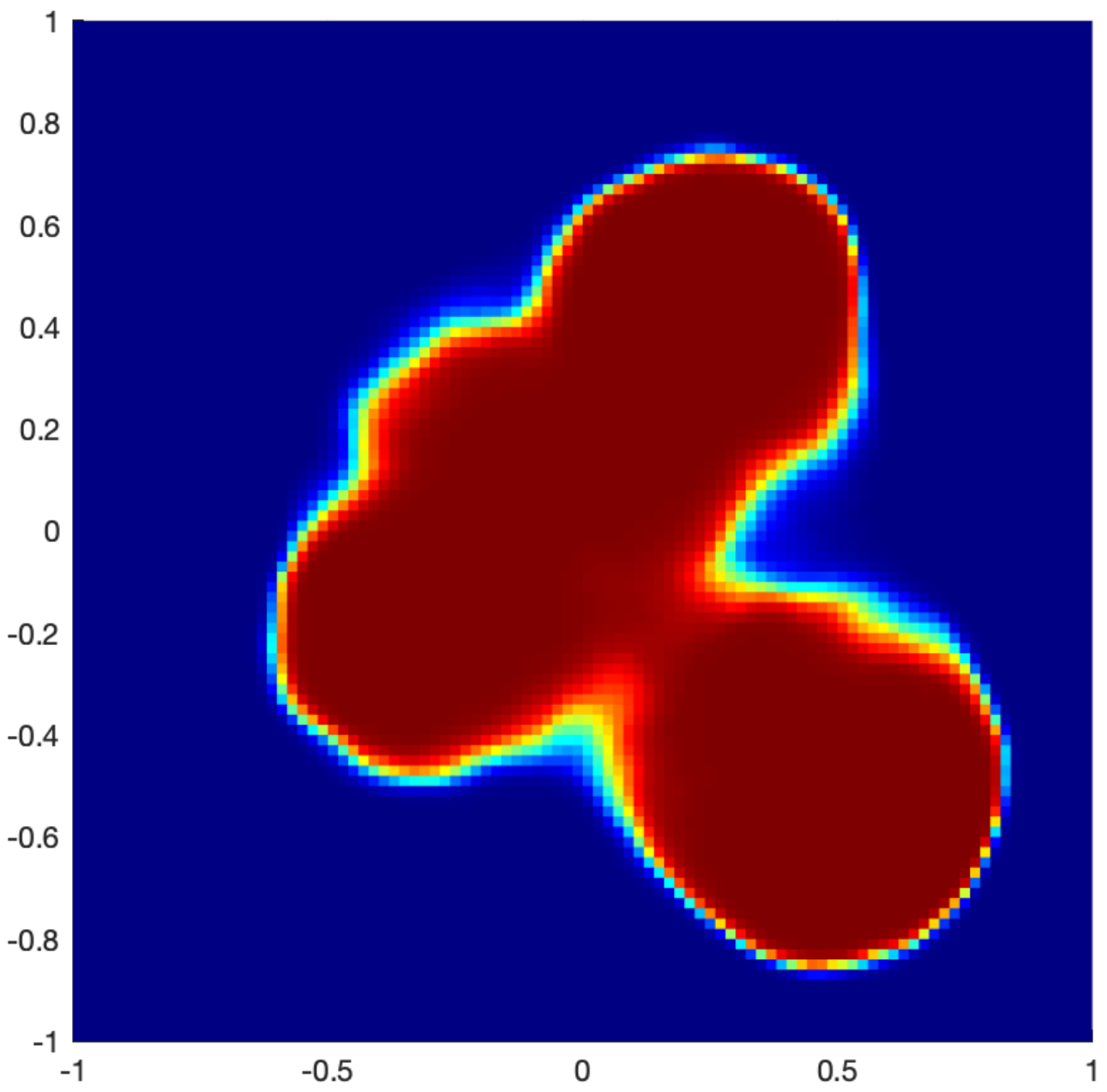}\\
\end{tabular}
  \caption{FNN-DDSM reconstruction for 3 cases in \textbf{Scenario 1} (3 circles) with different Cauchy data number and noise level: Case 1(top), Case 2(middle) and Case 3(bottom) } 
  \label{tab_FN_3cir}
\end{figure}

\begin{figure}[htbp]
\begin{tabular}{ >{\centering\arraybackslash}m{0.9in} >{\centering\arraybackslash}m{0.9in} >{\centering\arraybackslash}m{0.9in}  >{\centering\arraybackslash}m{0.9in}  >{\centering\arraybackslash}m{0.9in}  >{\centering\arraybackslash}m{0.9in} }
\centering
True coefficients &
N=1, $\delta=0$&
N=10, $\delta=0$&
N=20, $\delta=0$&
N=20, $\delta=10\%$ &
N=20, $\delta=20\%$ \\
\includegraphics[width=1.1in]{cir_num3_case1_0-eps-converted-to.pdf}&
\includegraphics[width=1.1in]{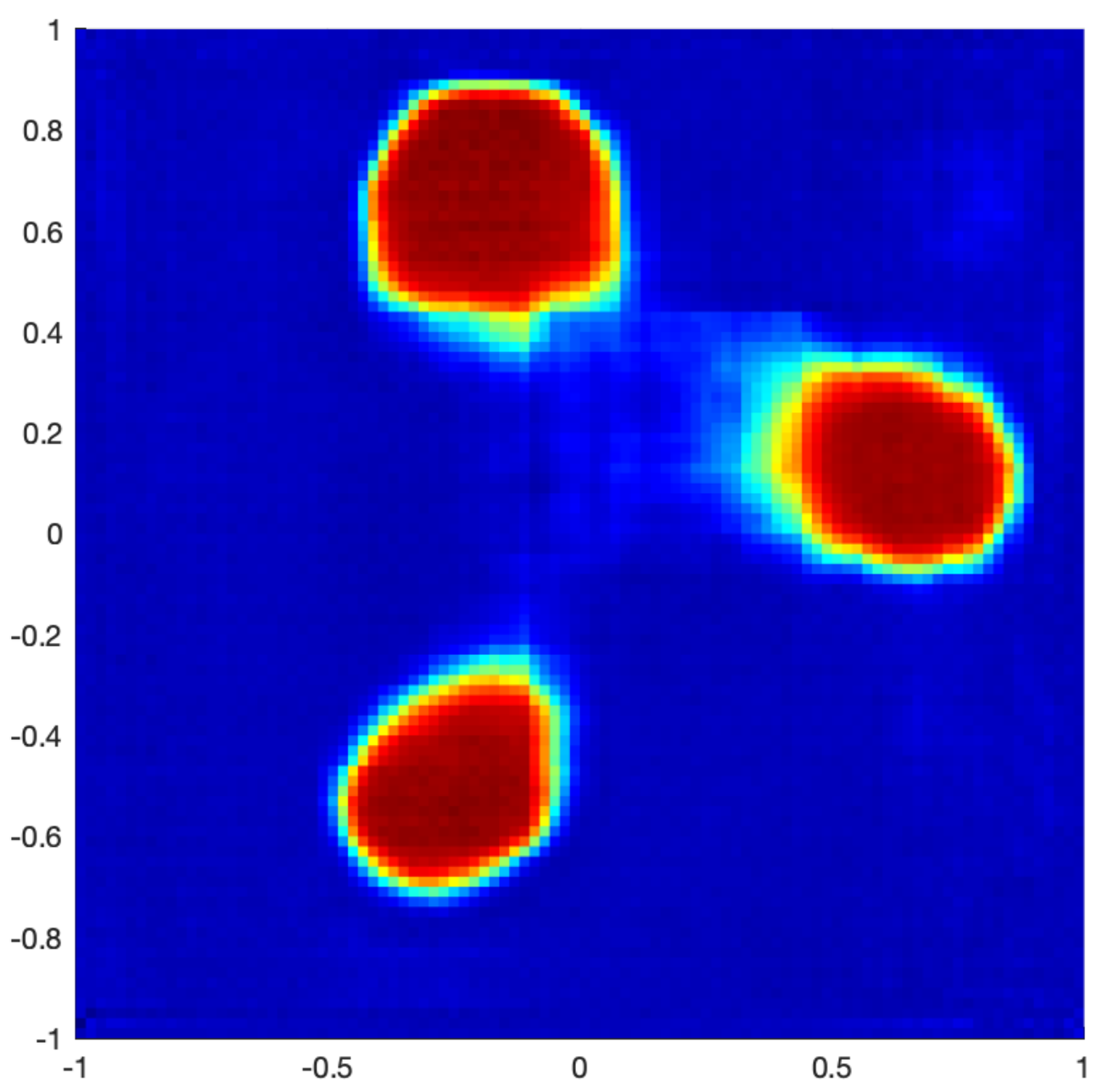}&
\includegraphics[width=1.1in]{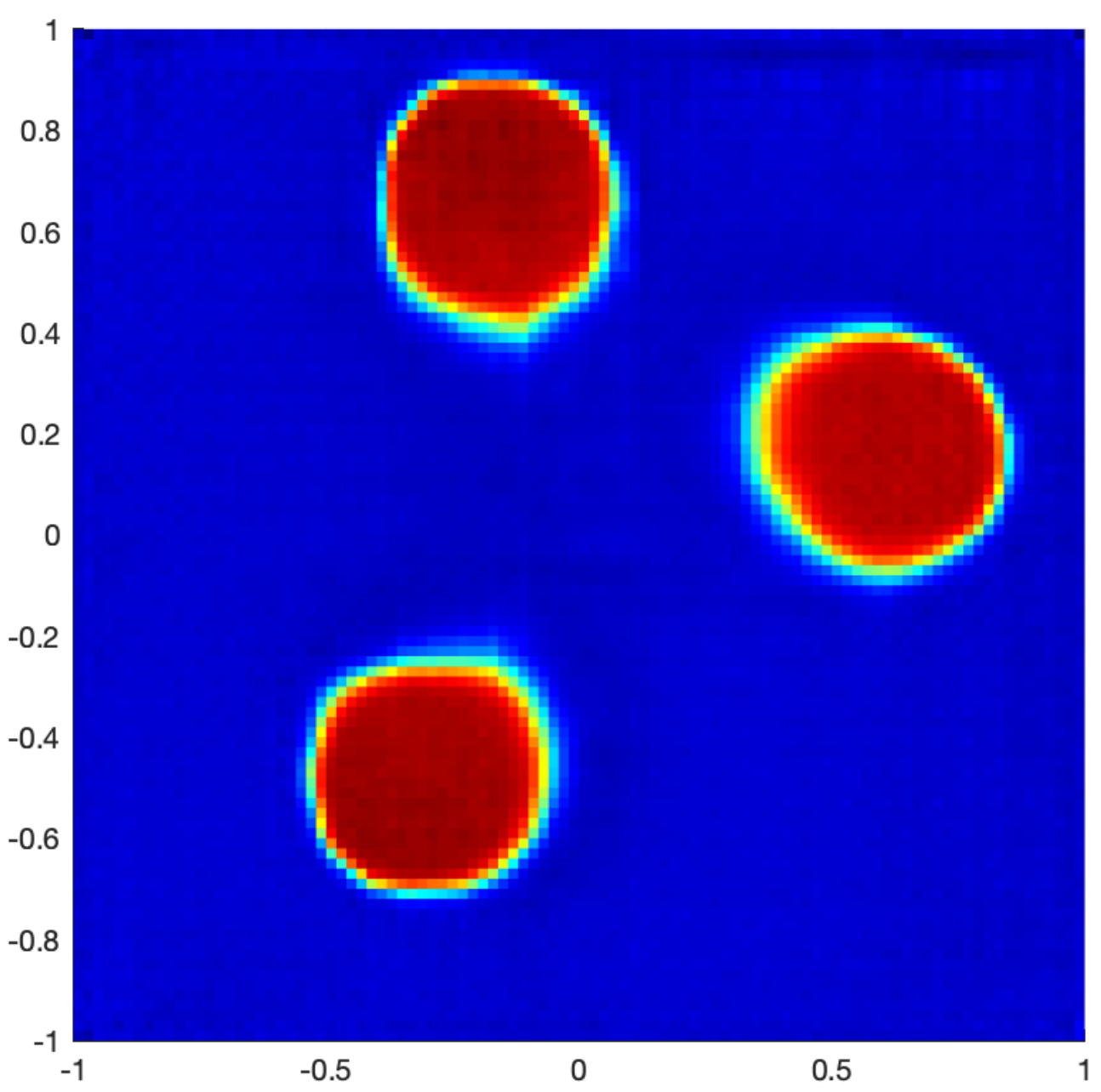}&
\includegraphics[width=1.1in]{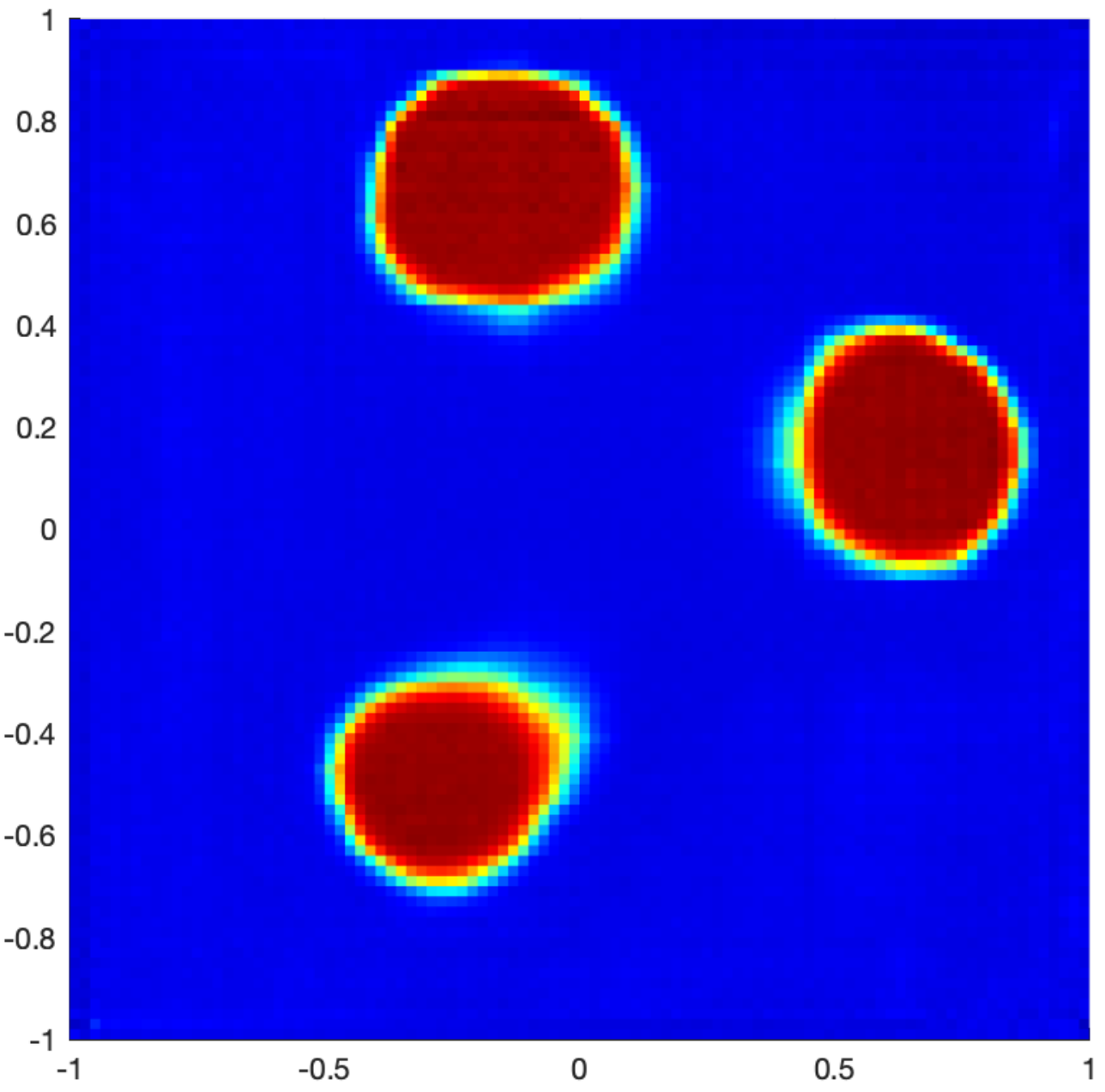}&
\includegraphics[width=1.1in]{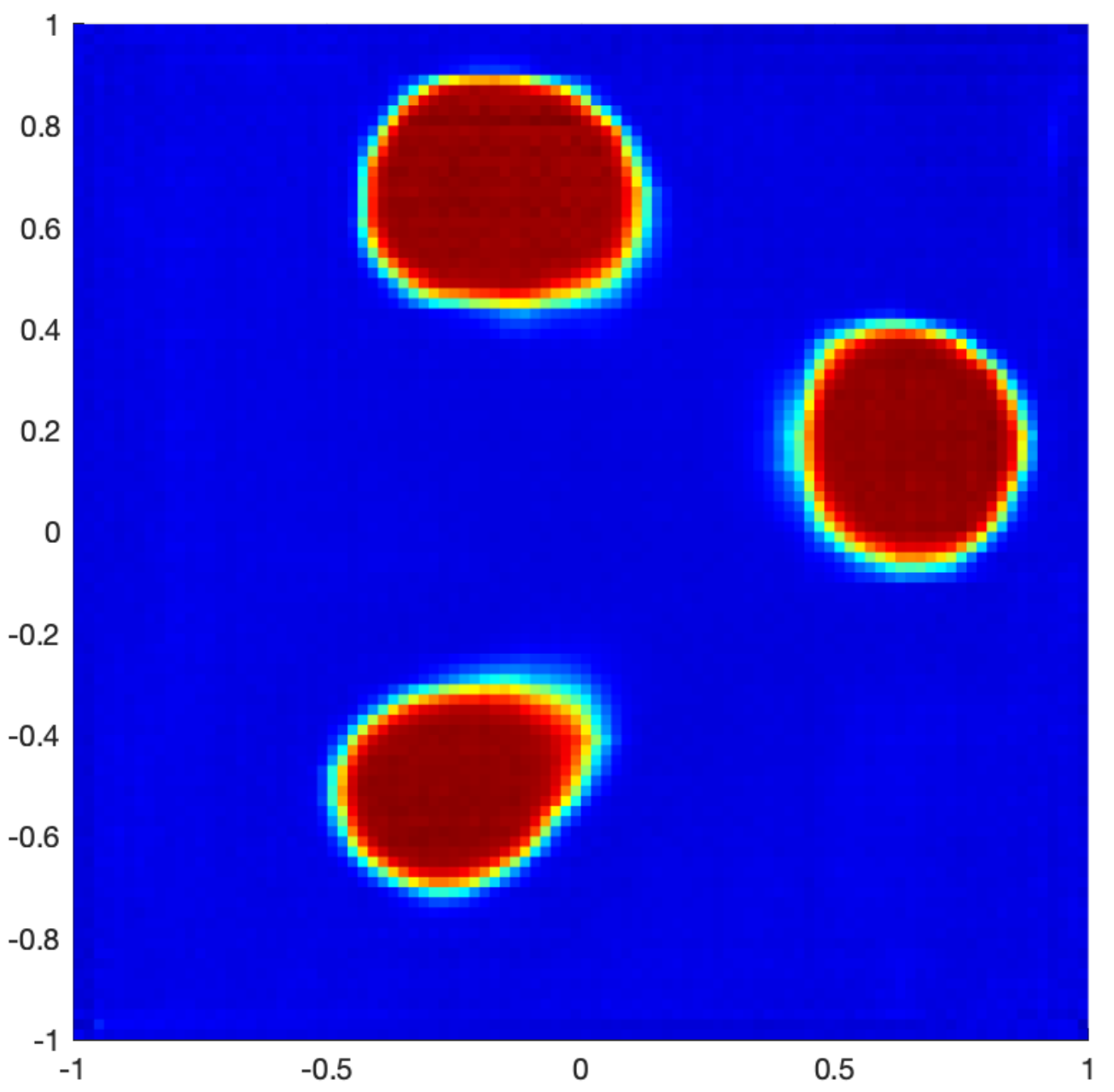}&
\includegraphics[width=1.1in]{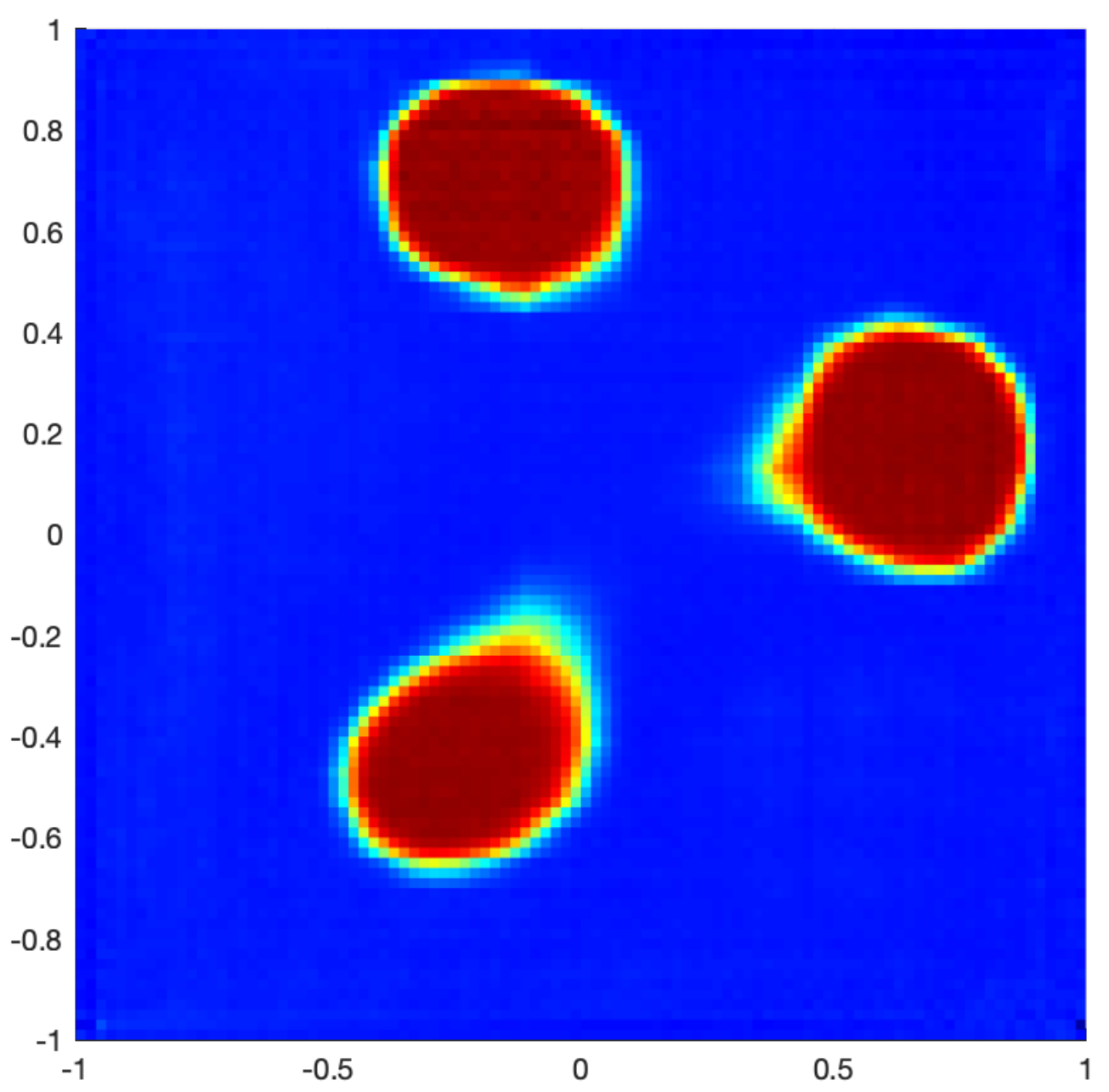}\\
\includegraphics[width=1.1in]{cir_num3_case2_0-eps-converted-to.pdf}&
\includegraphics[width=1.1in]{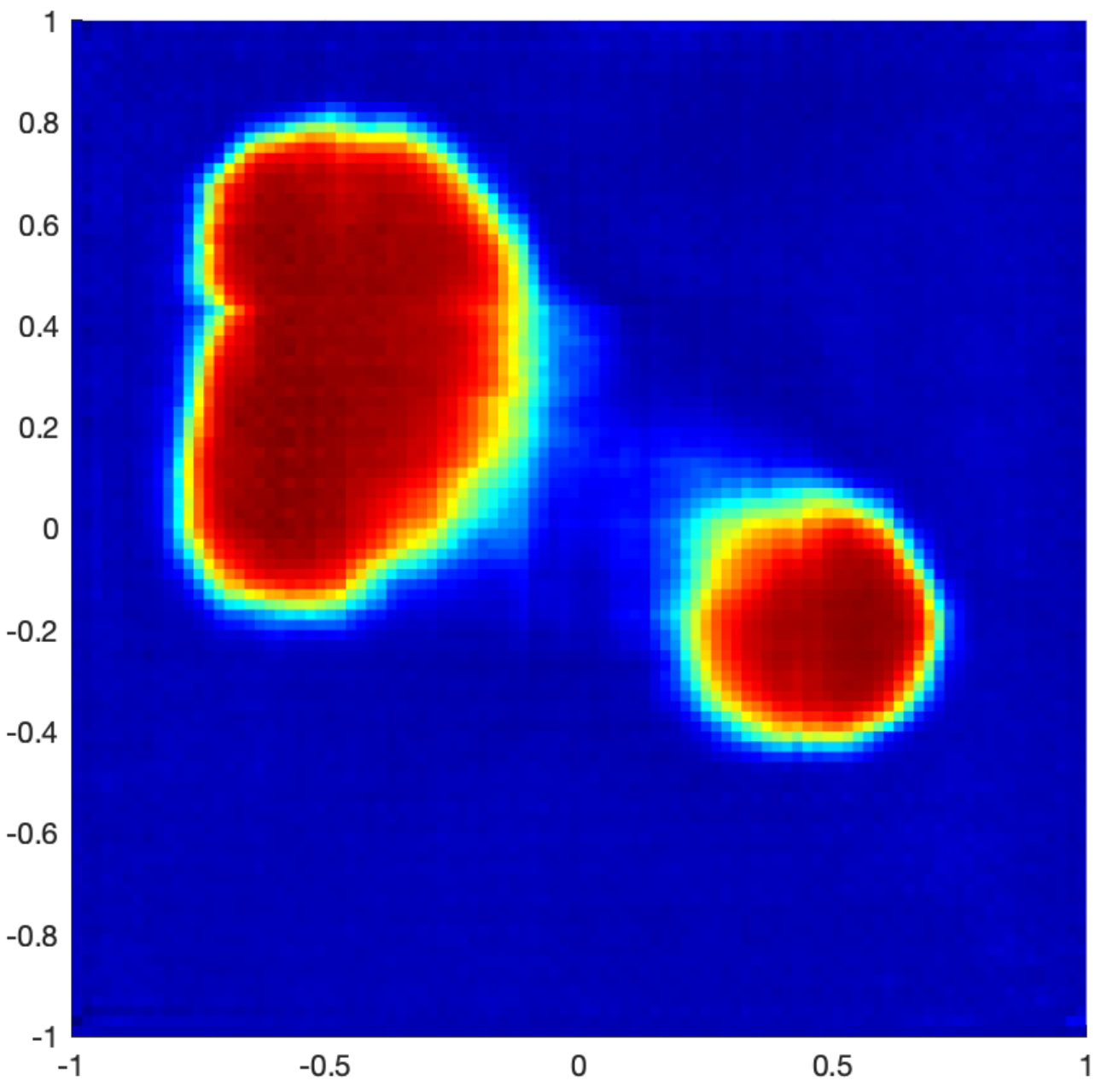}&
\includegraphics[width=1.1in]{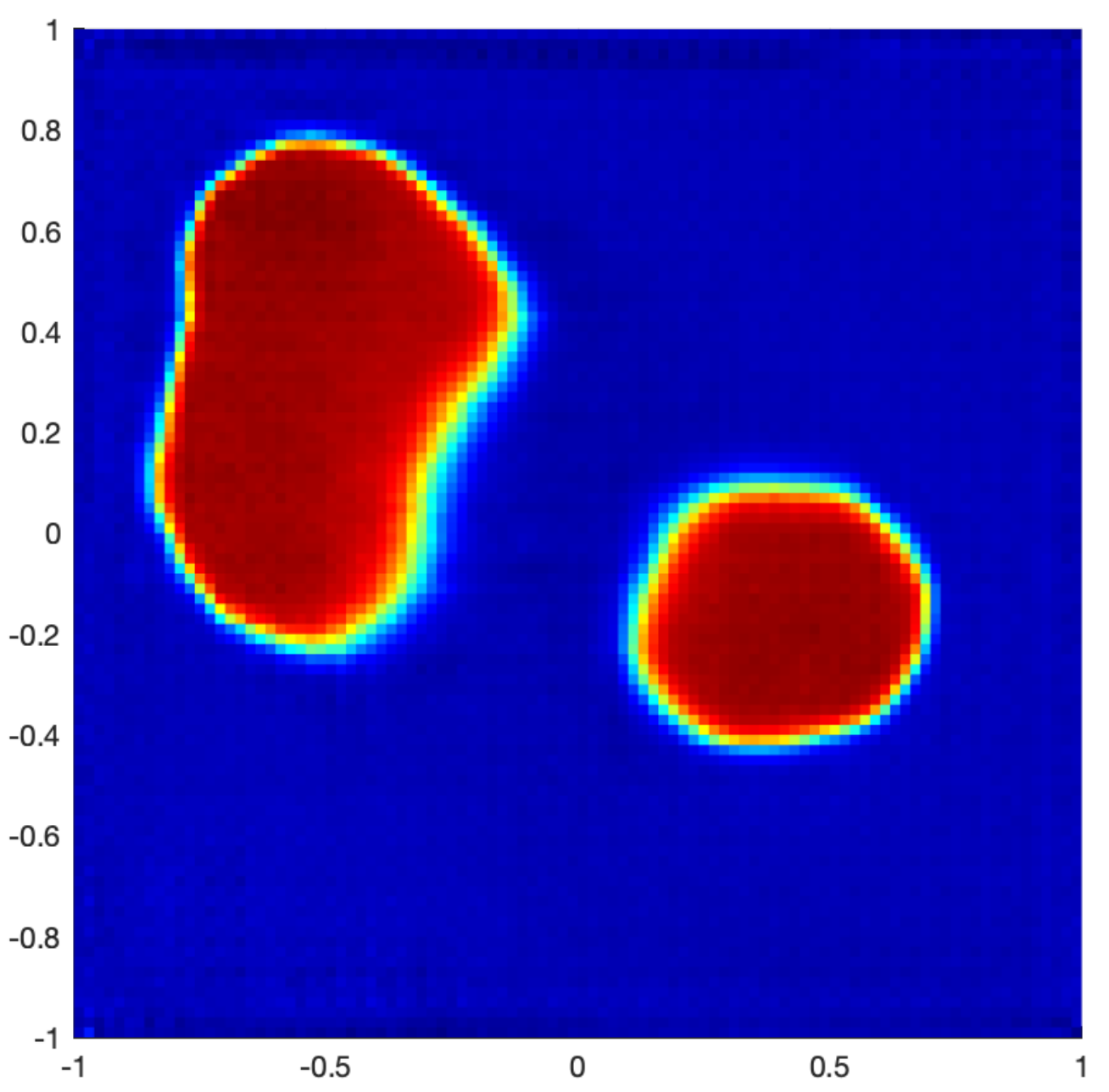}&
\includegraphics[width=1.1in]{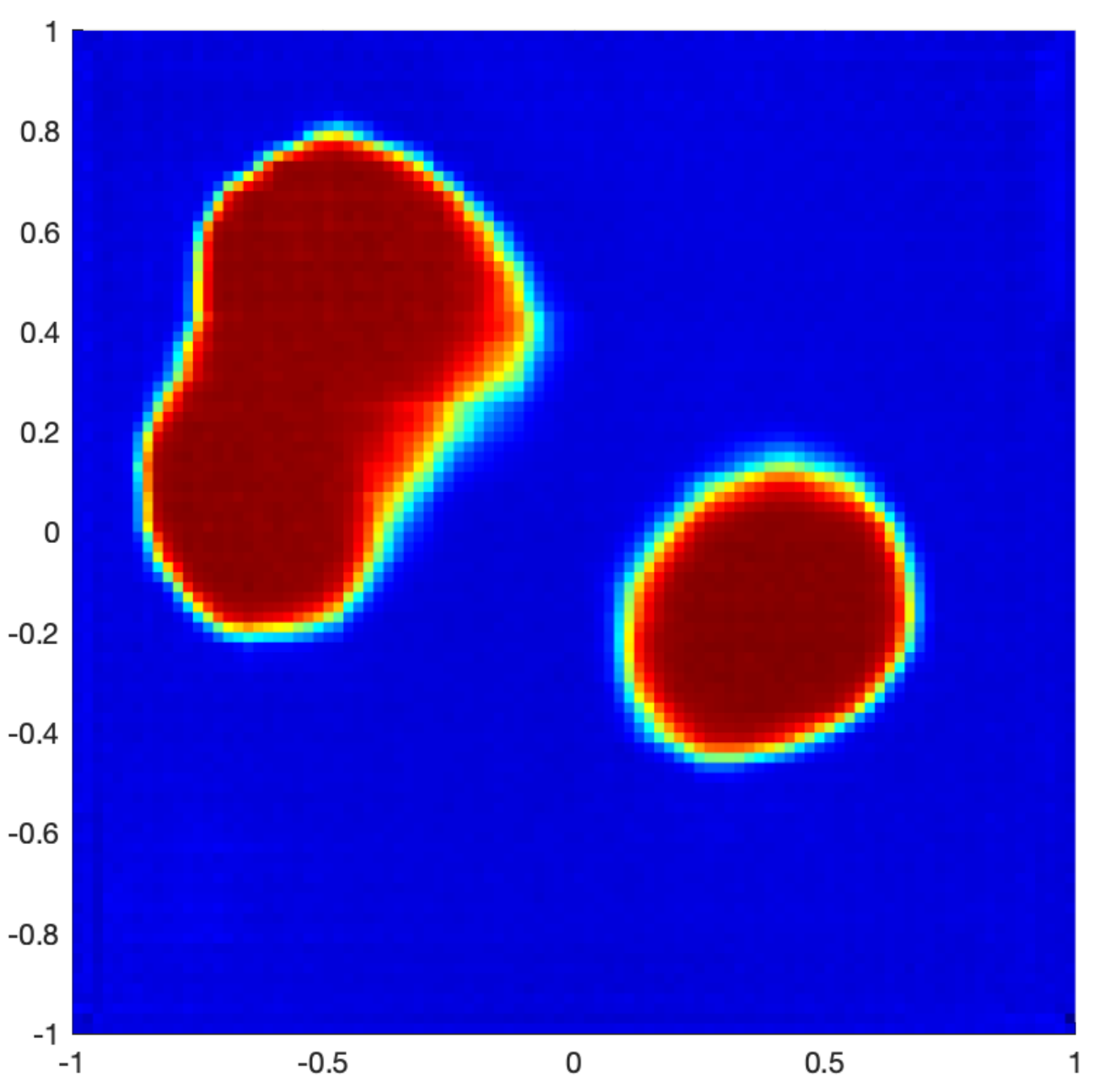}&
\includegraphics[width=1.1in]{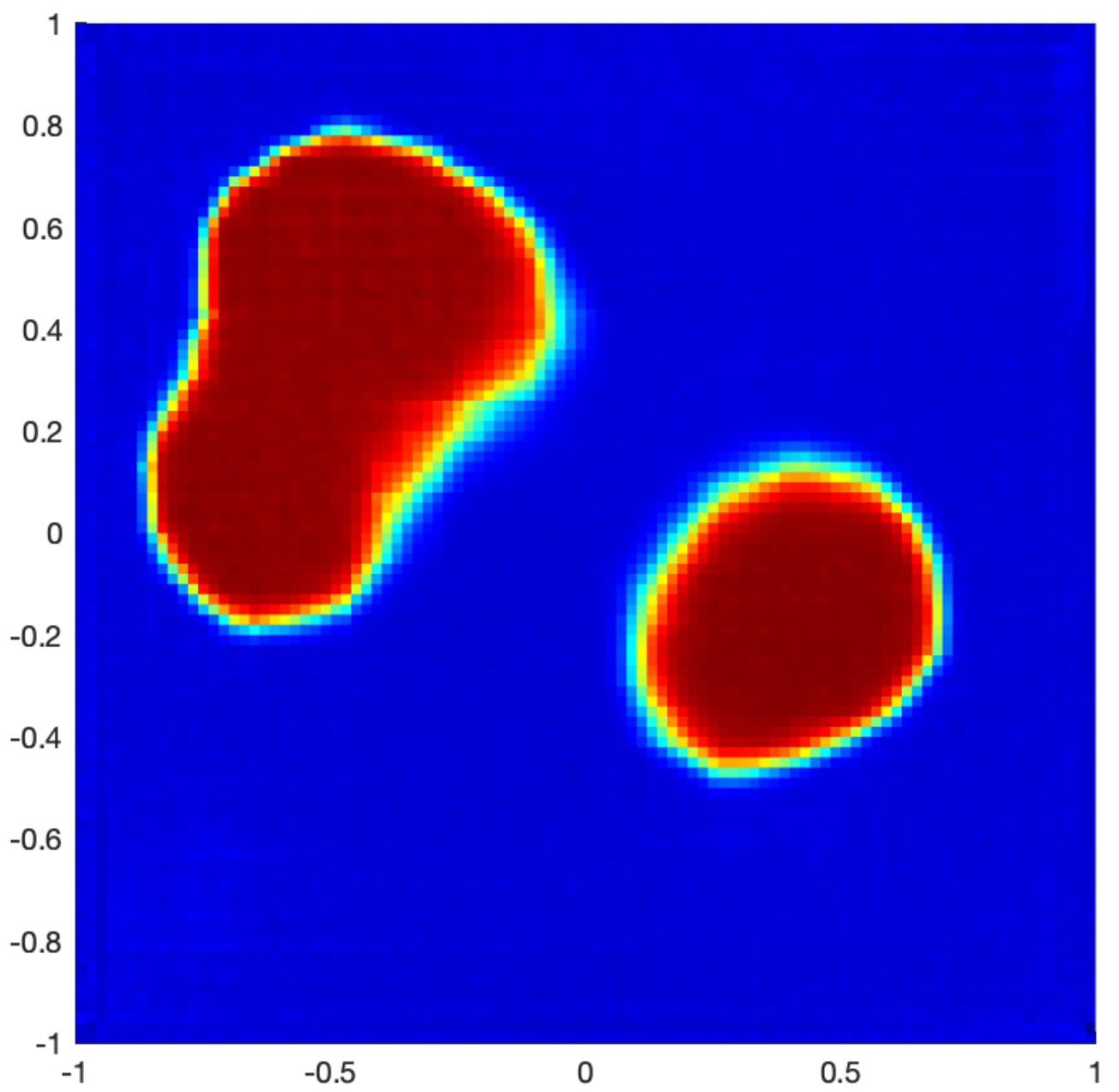}&
\includegraphics[width=1.1in]{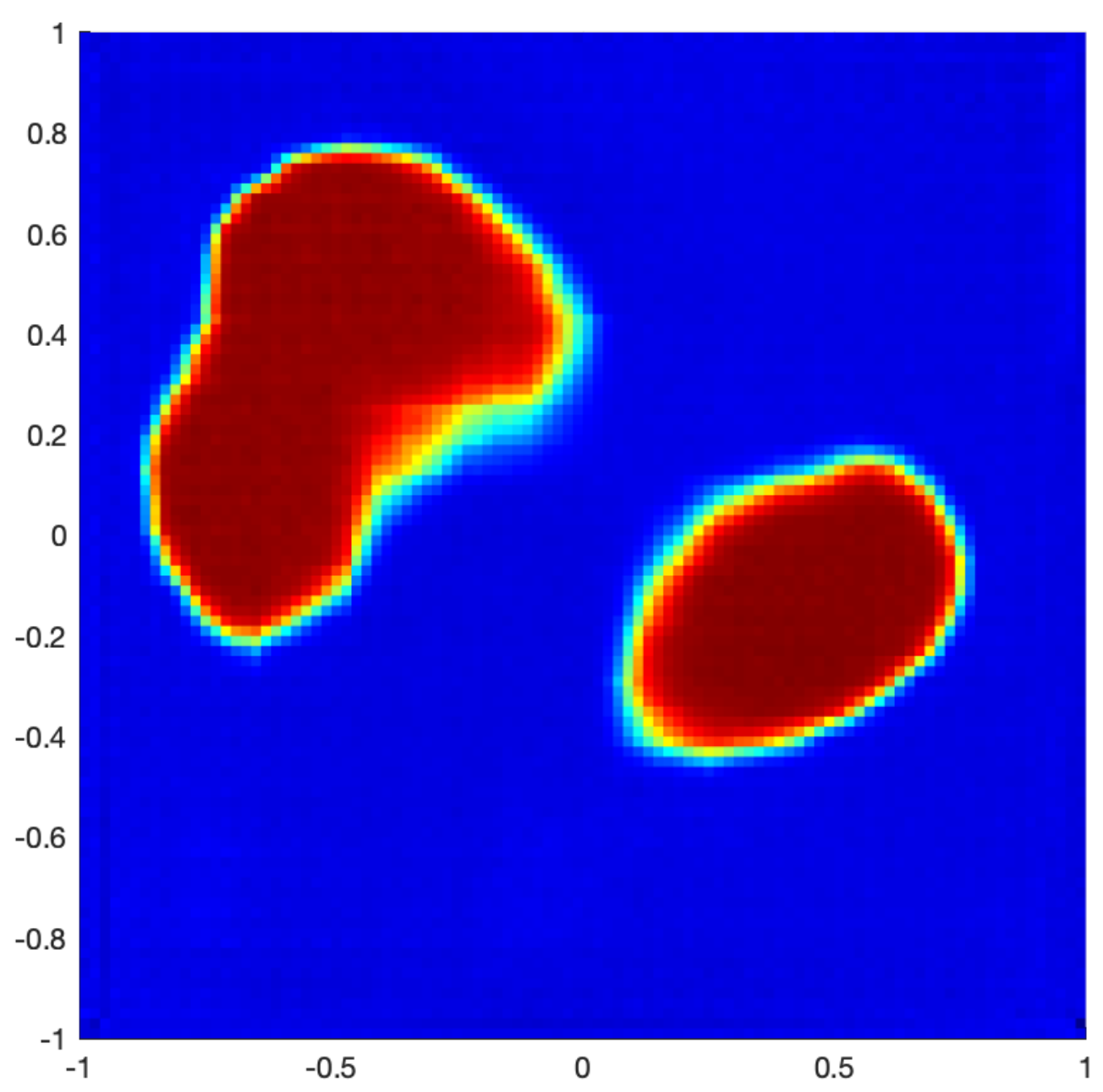}\\
\includegraphics[width=1.1in]{cir_num3_case3_0-eps-converted-to.pdf}&
\includegraphics[width=1.1in]{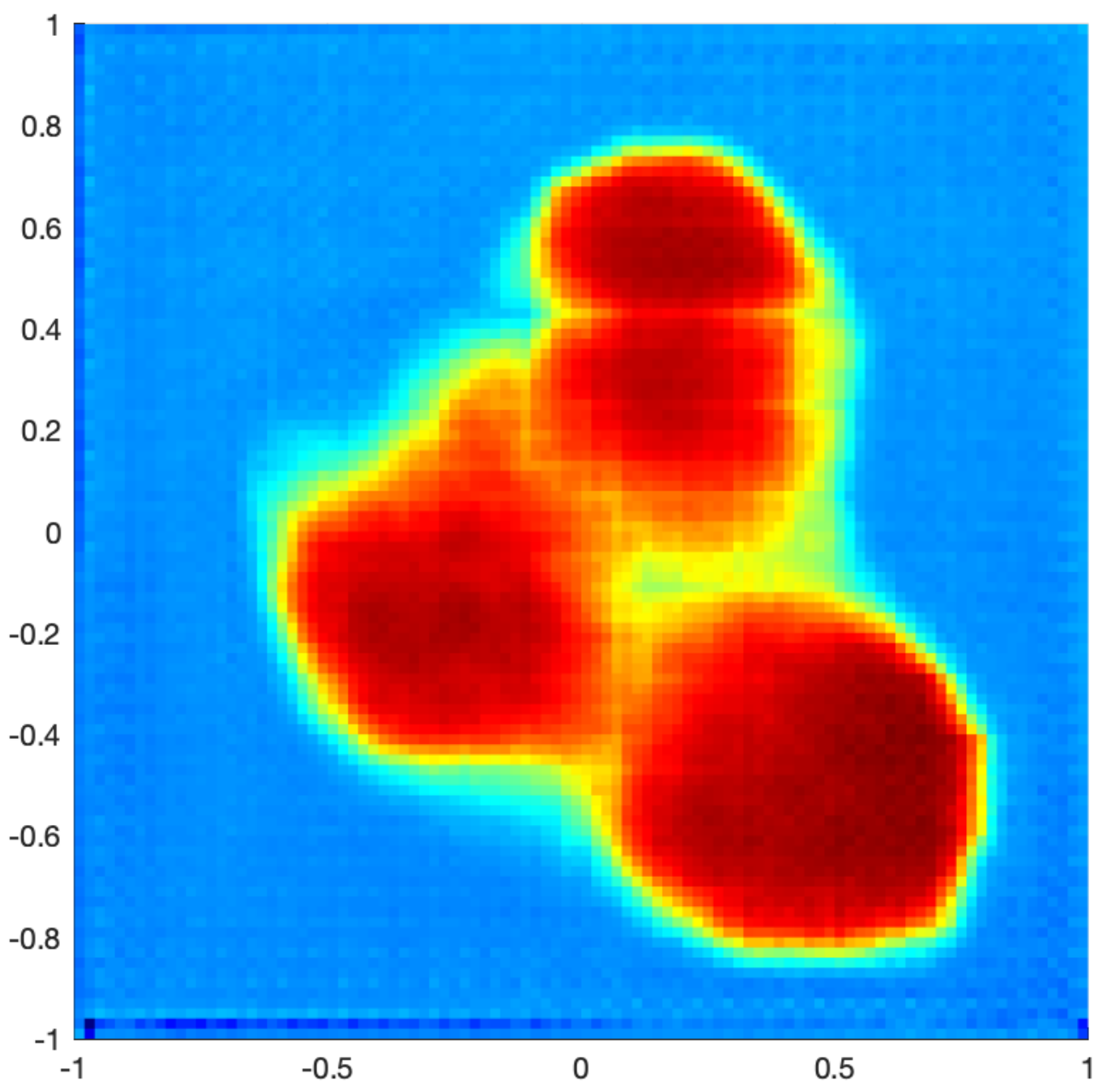}&
\includegraphics[width=1.1in]{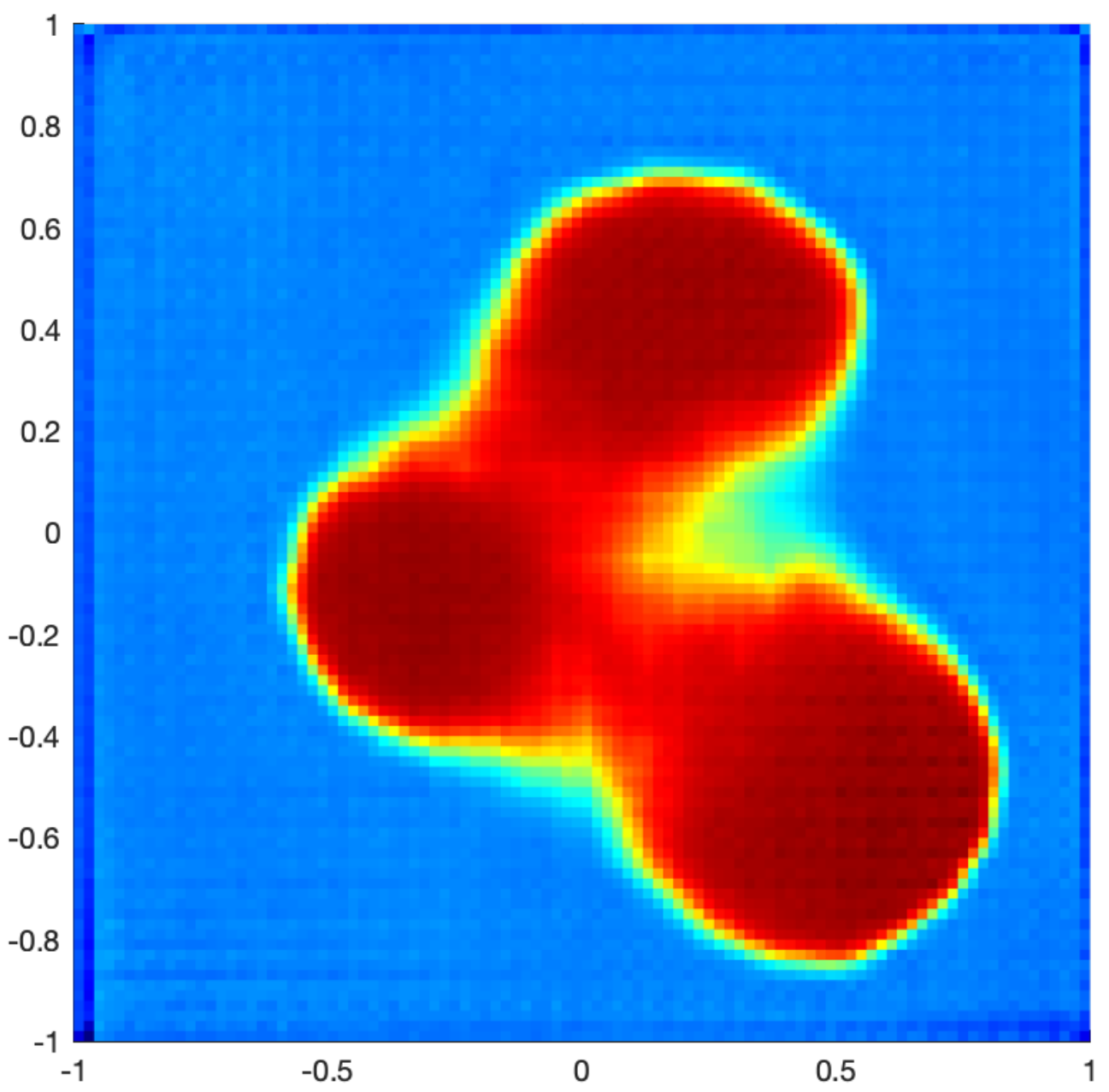}&
\includegraphics[width=1.1in]{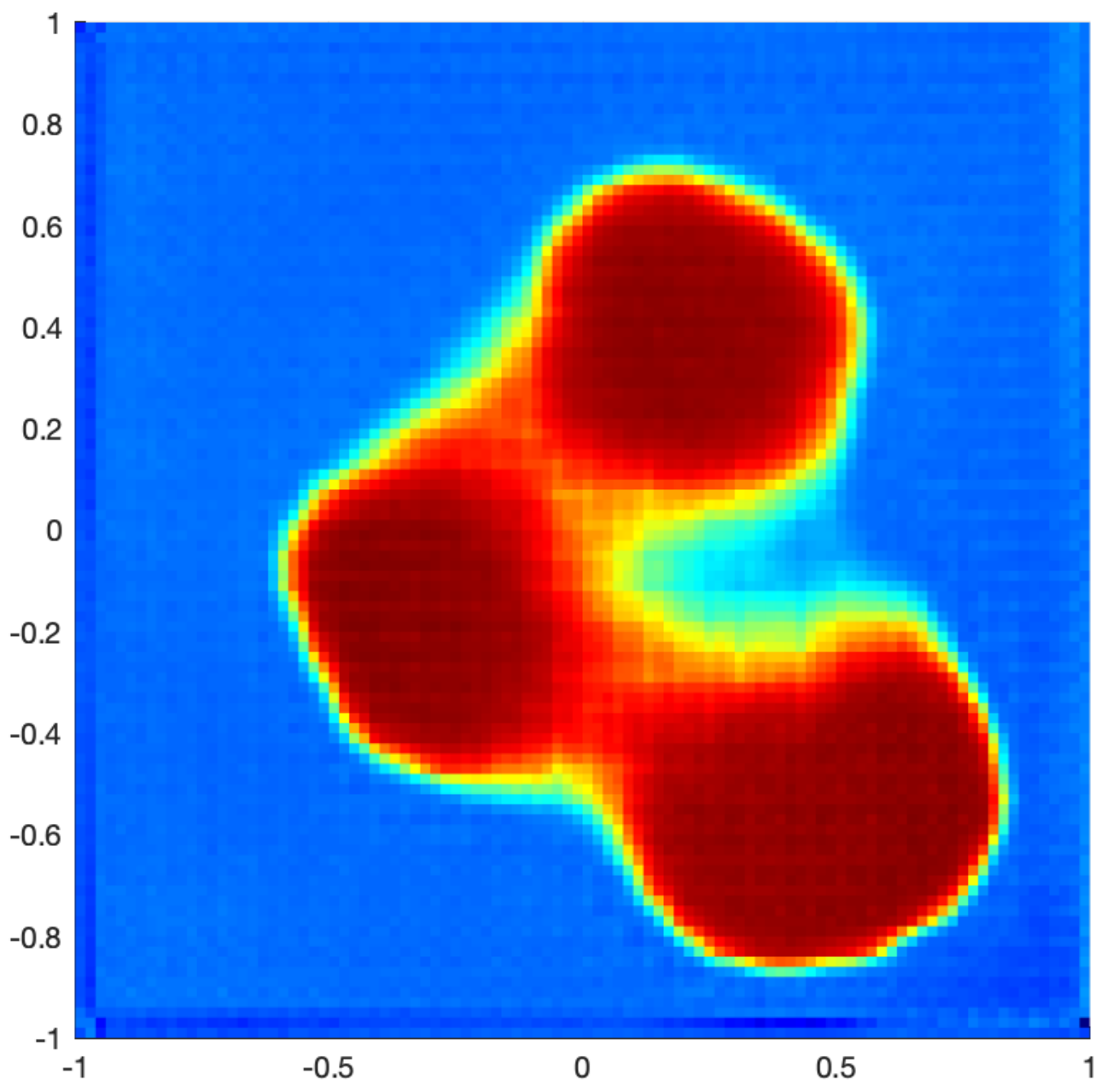}&
\includegraphics[width=1.1in]{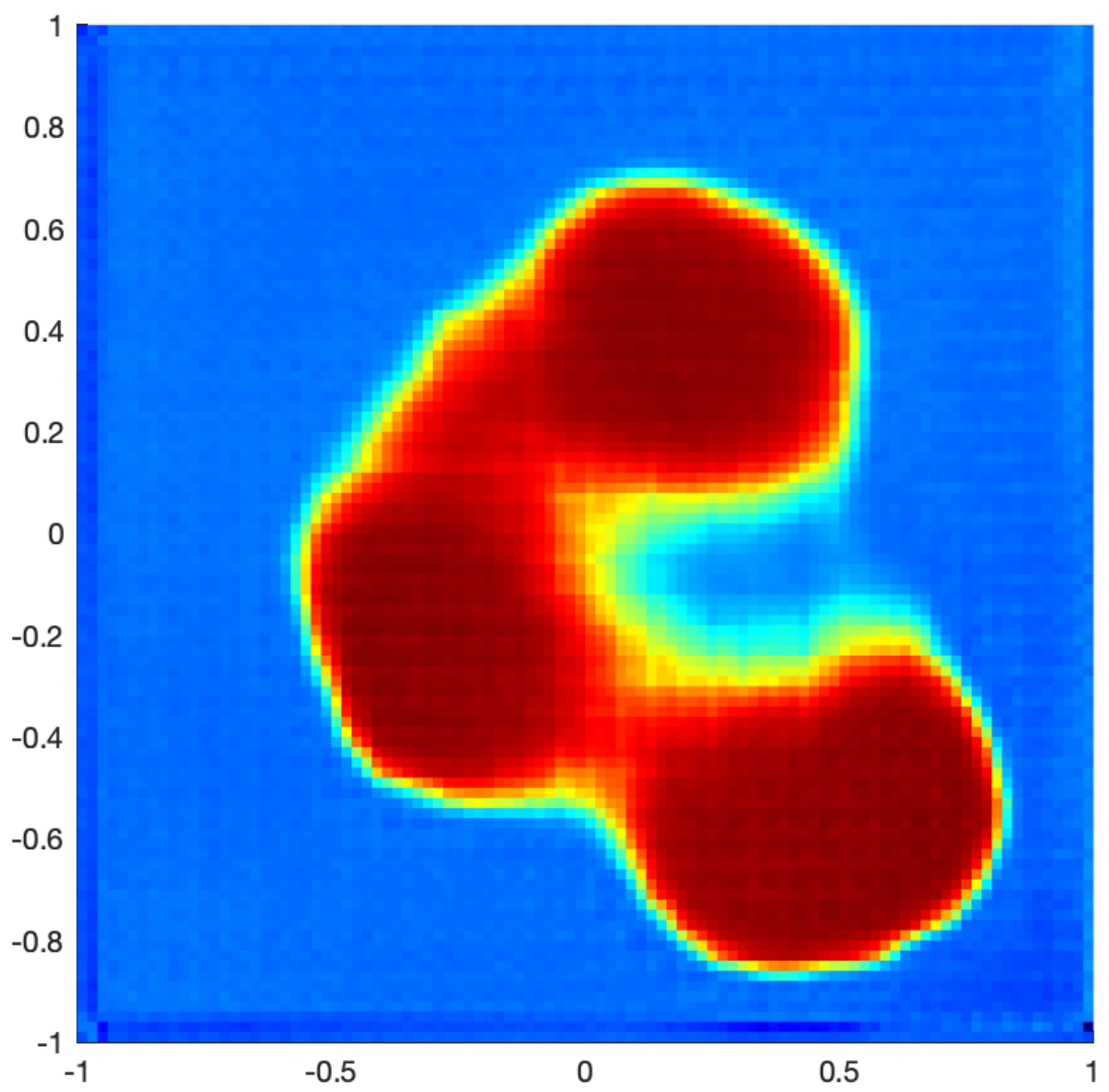}&
\includegraphics[width=1.1in]{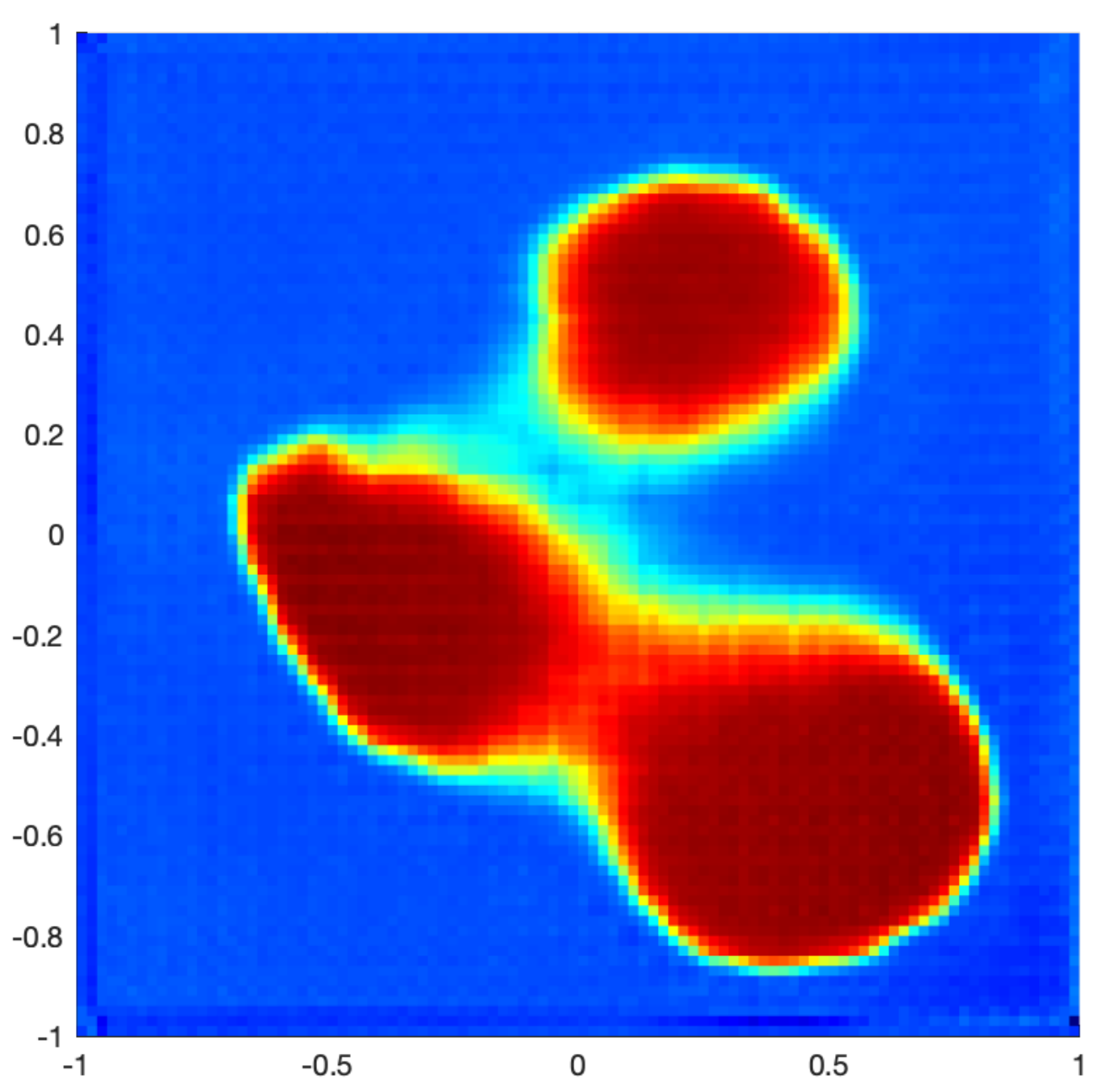}\\
\end{tabular}
  \caption{CNN-DDSM reconstruction for 3 cases in \textbf{Scenario 1} (3 circles) with different Cauchy data number and noise level: Case 1(top), Case 2(middle) and Case 3(bottom) } 
  \label{tab_3cir}
\end{figure}

\begin{figure}[htbp]
\begin{tabular}{ >{\centering\arraybackslash}m{0.9in} >{\centering\arraybackslash}m{0.9in} >{\centering\arraybackslash}m{0.9in}  >{\centering\arraybackslash}m{0.9in}  >{\centering\arraybackslash}m{0.9in}  >{\centering\arraybackslash}m{0.9in} }
\centering
True coefficients &
N=1, $\delta=0$&
N=10, $\delta=0$&
N=20, $\delta=0$&
N=20, $\delta=10\%$ &
N=20, $\delta=20\%$ \\
\includegraphics[width=1.1in]{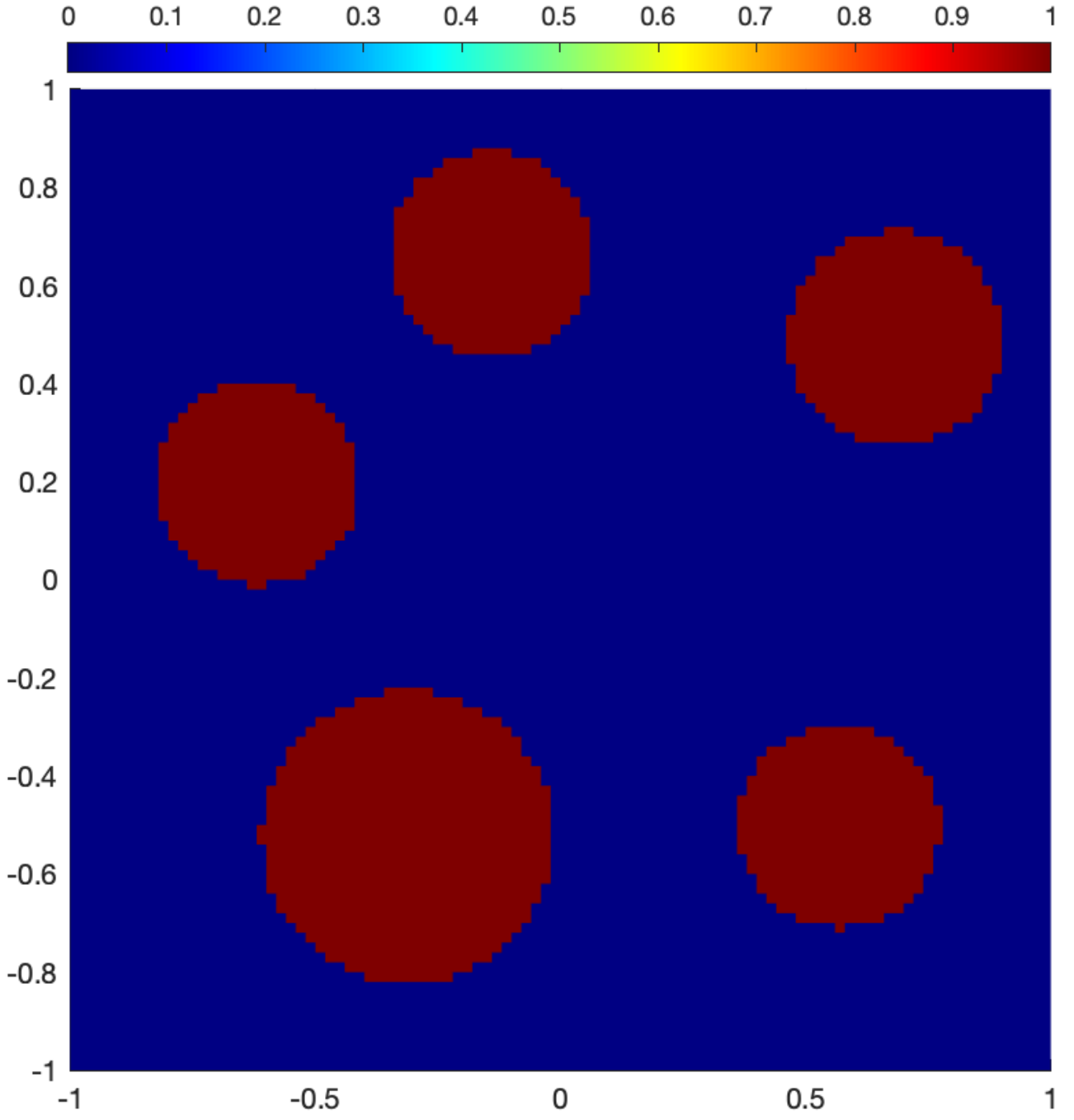}&
\includegraphics[width=1.1in]{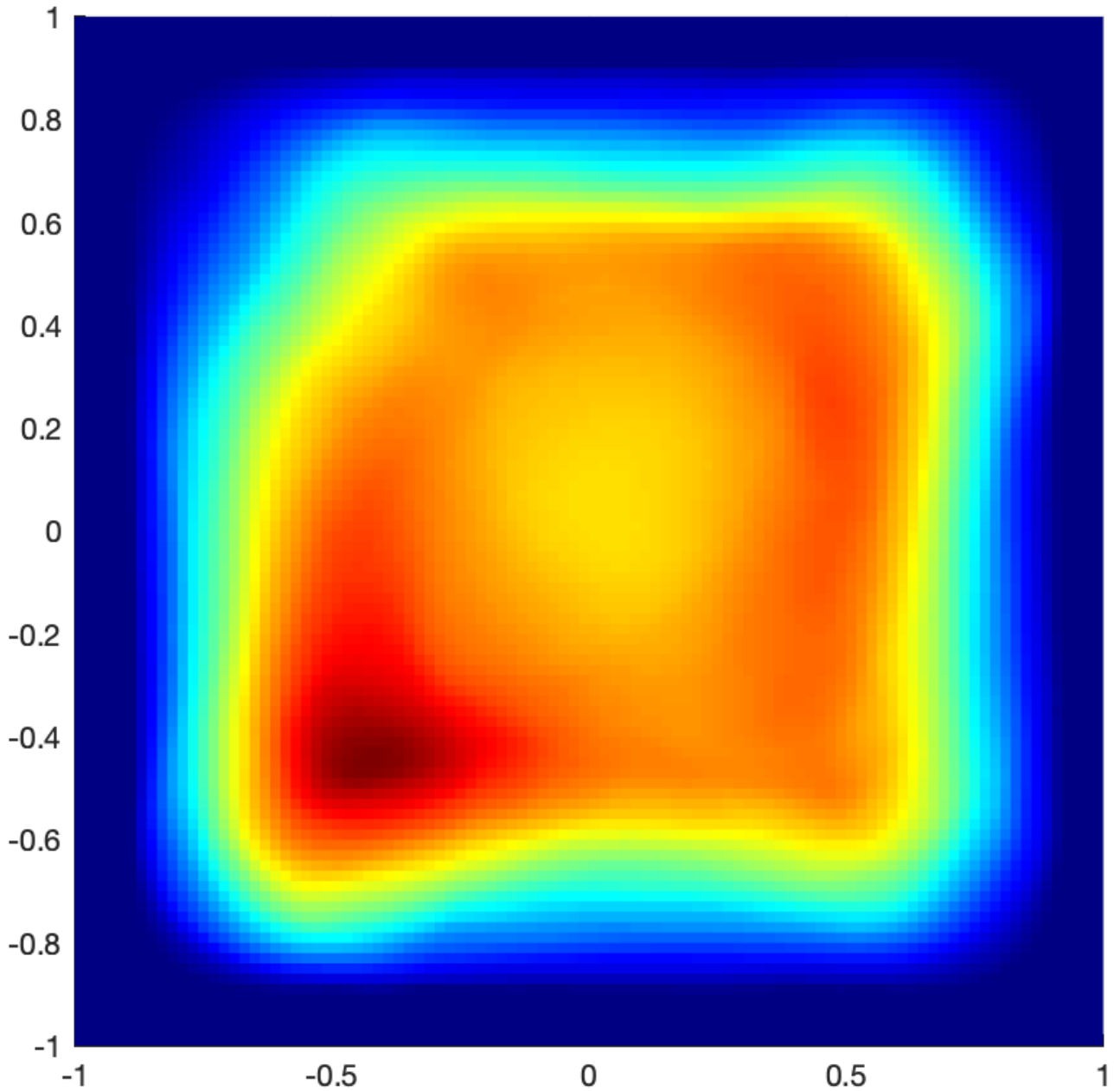}&
\includegraphics[width=1.1in]{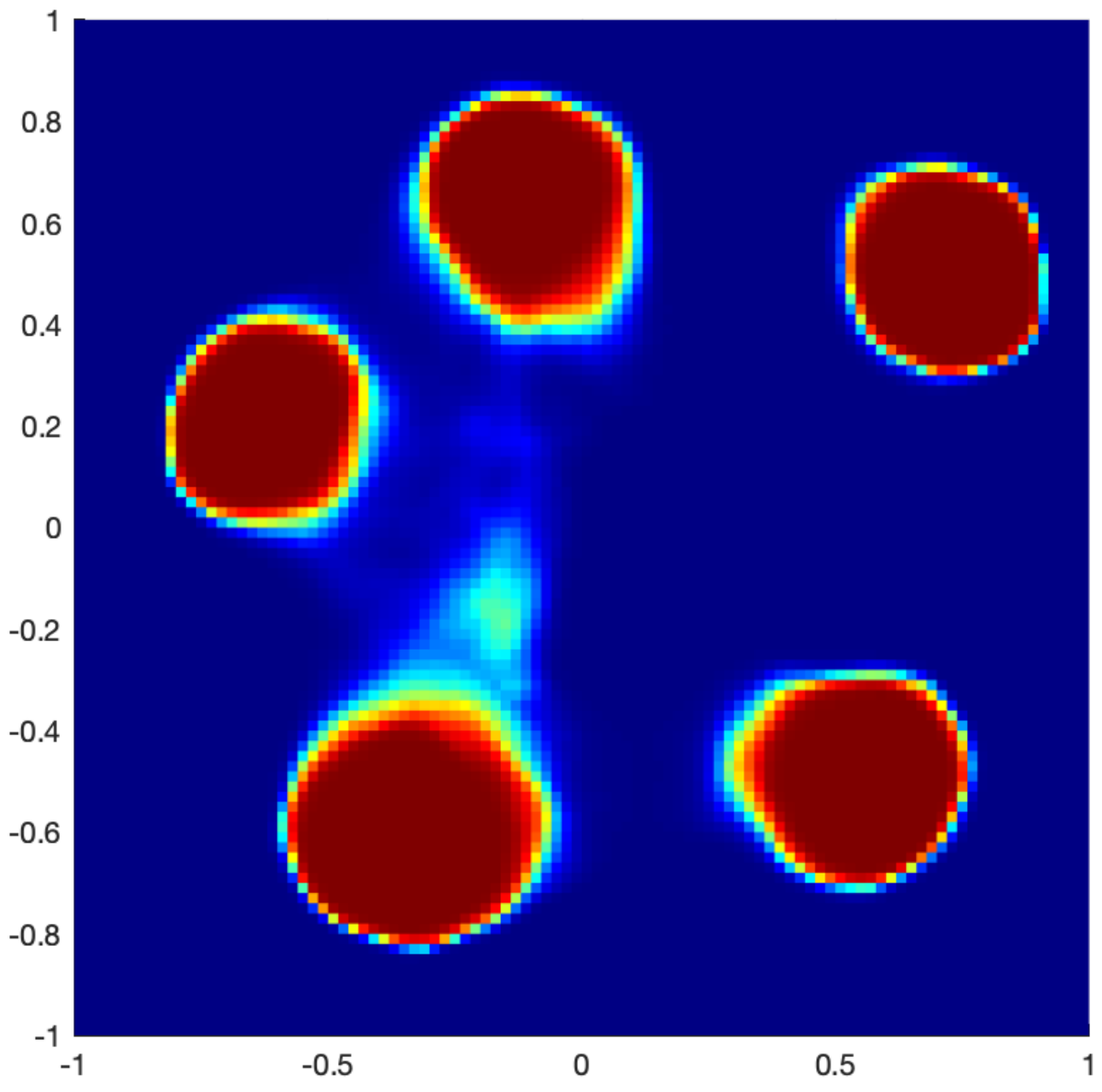}&
\includegraphics[width=1.1in]{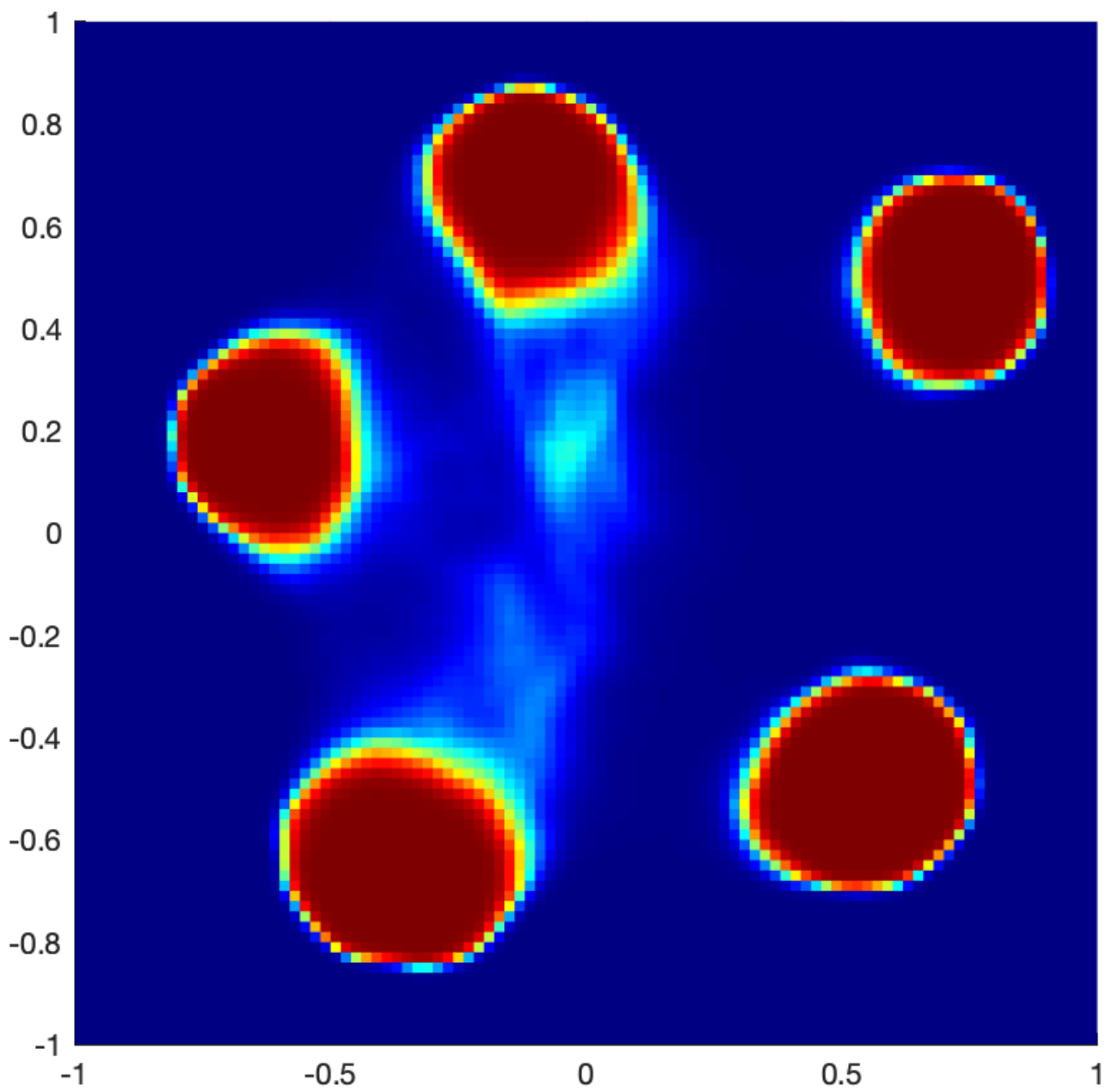}&
\includegraphics[width=1.1in]{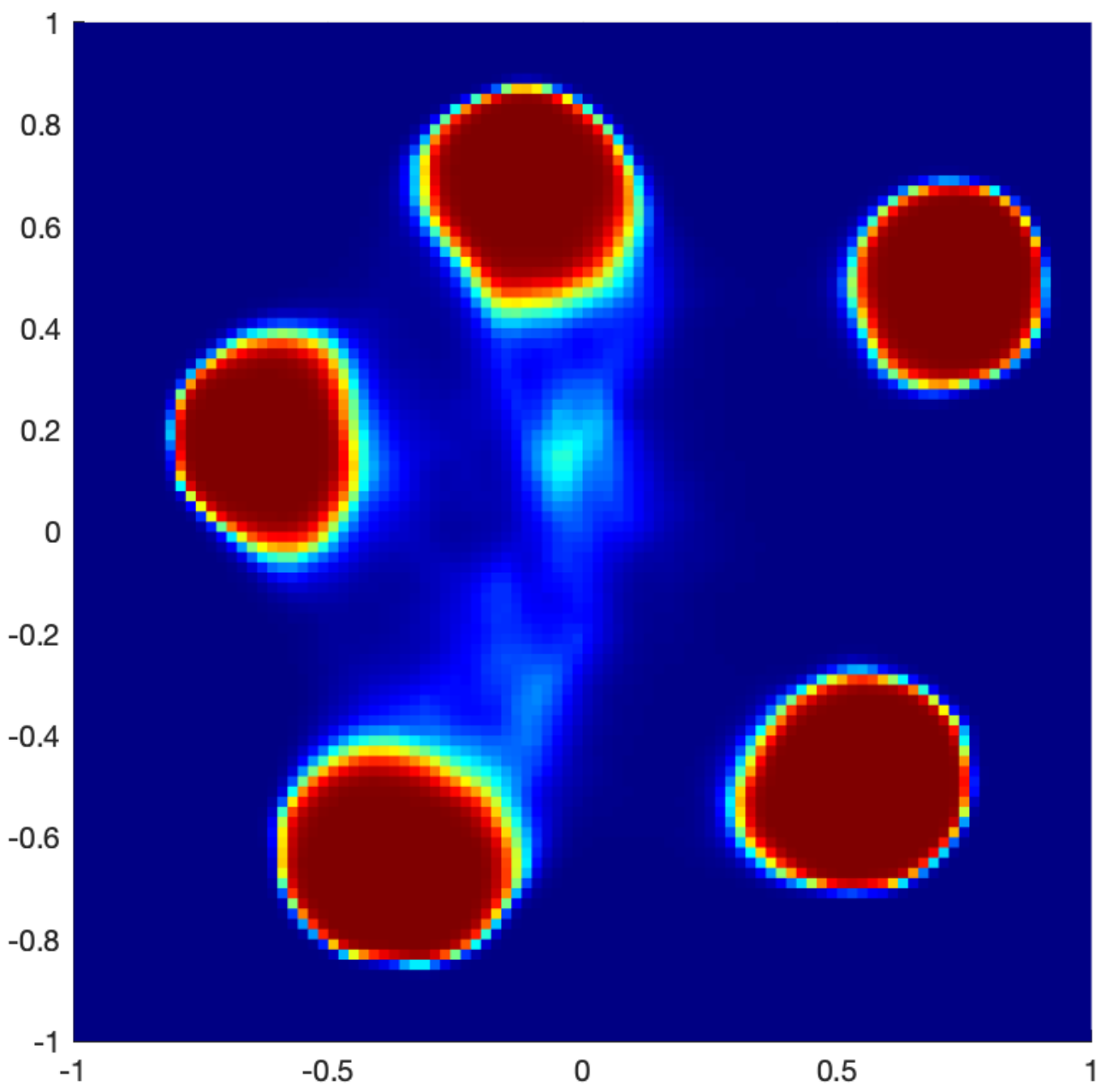}&
\includegraphics[width=1.1in]{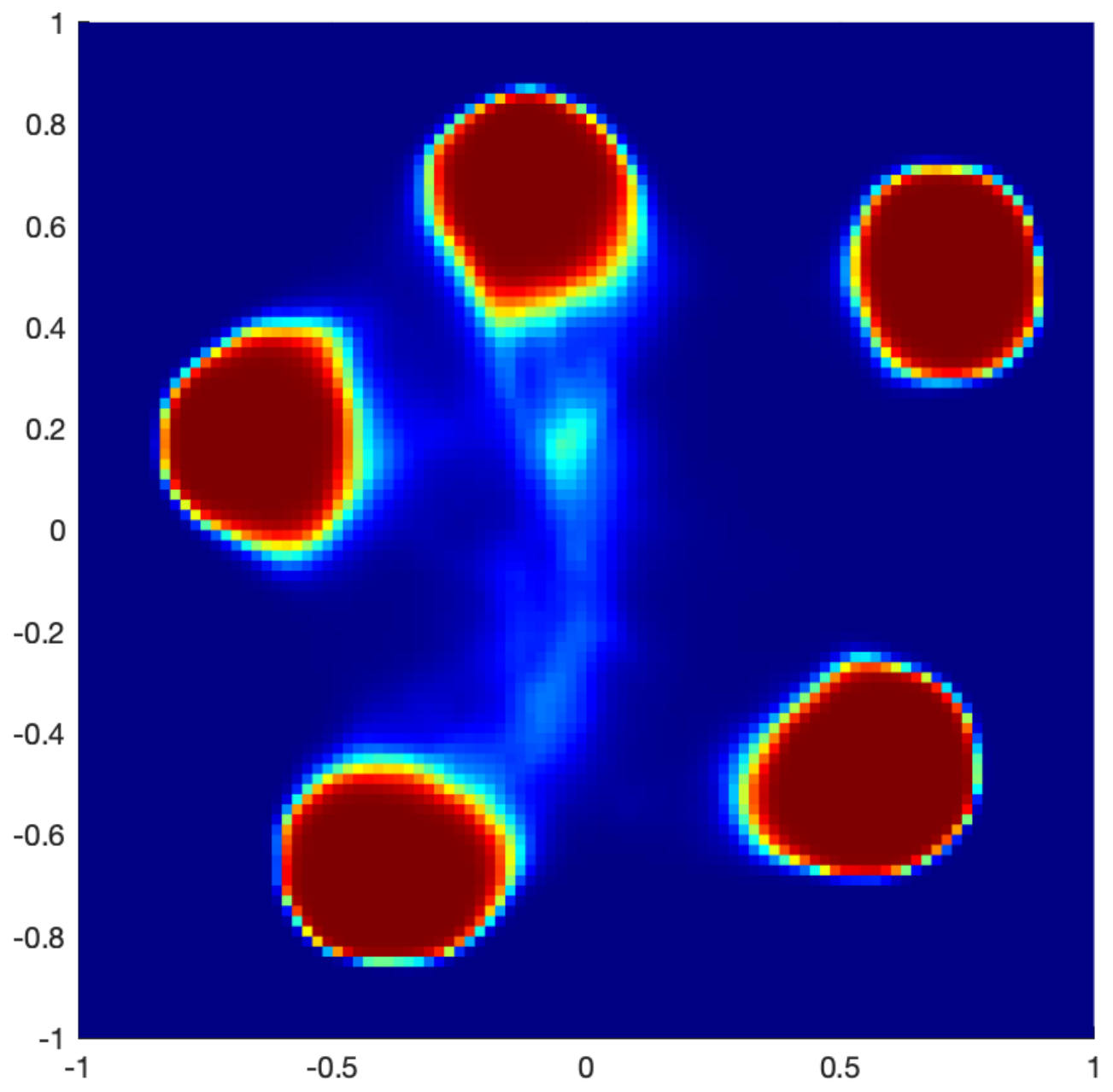}\\
\includegraphics[width=1.1in]{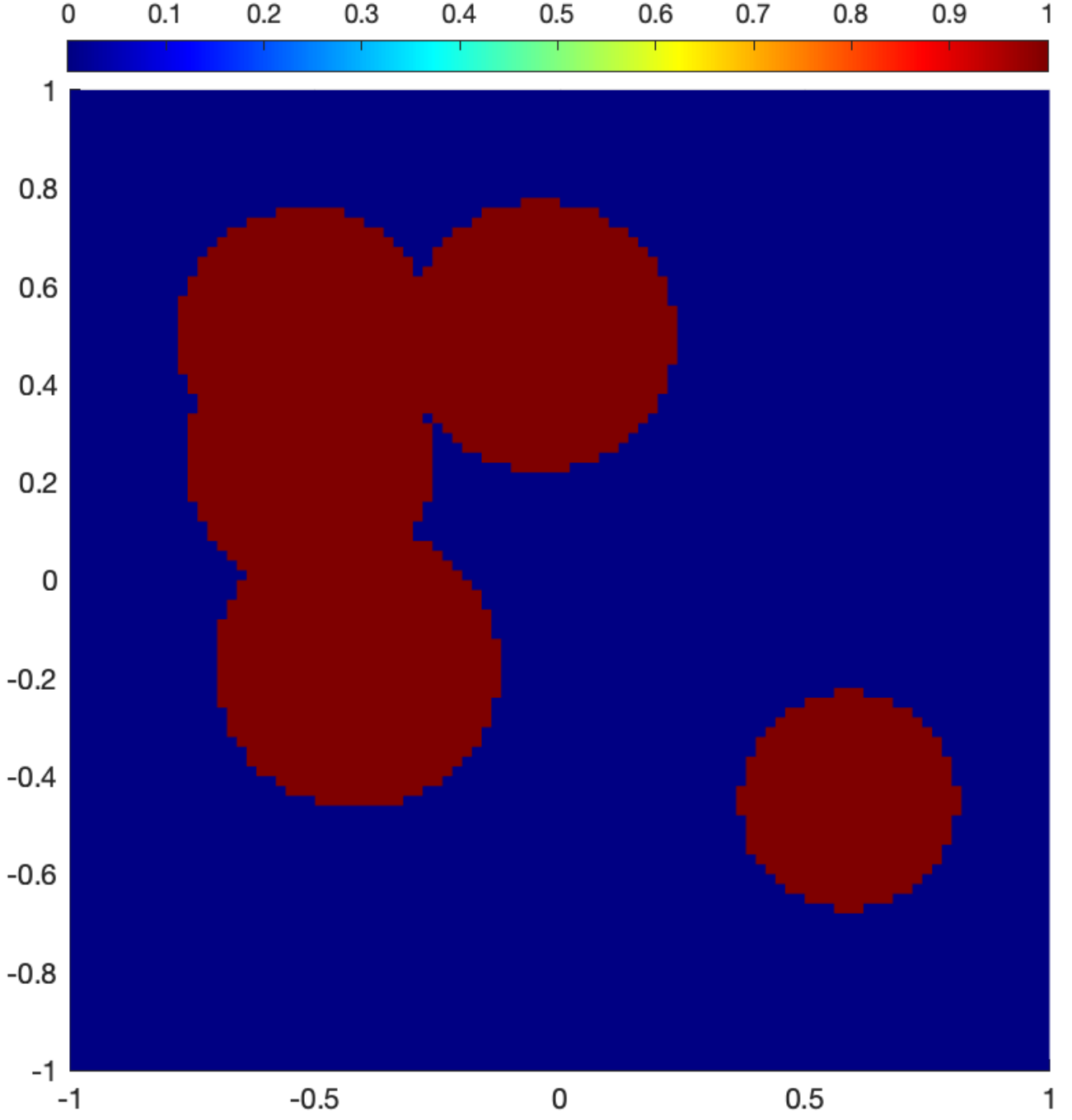}&
\includegraphics[width=1.1in]{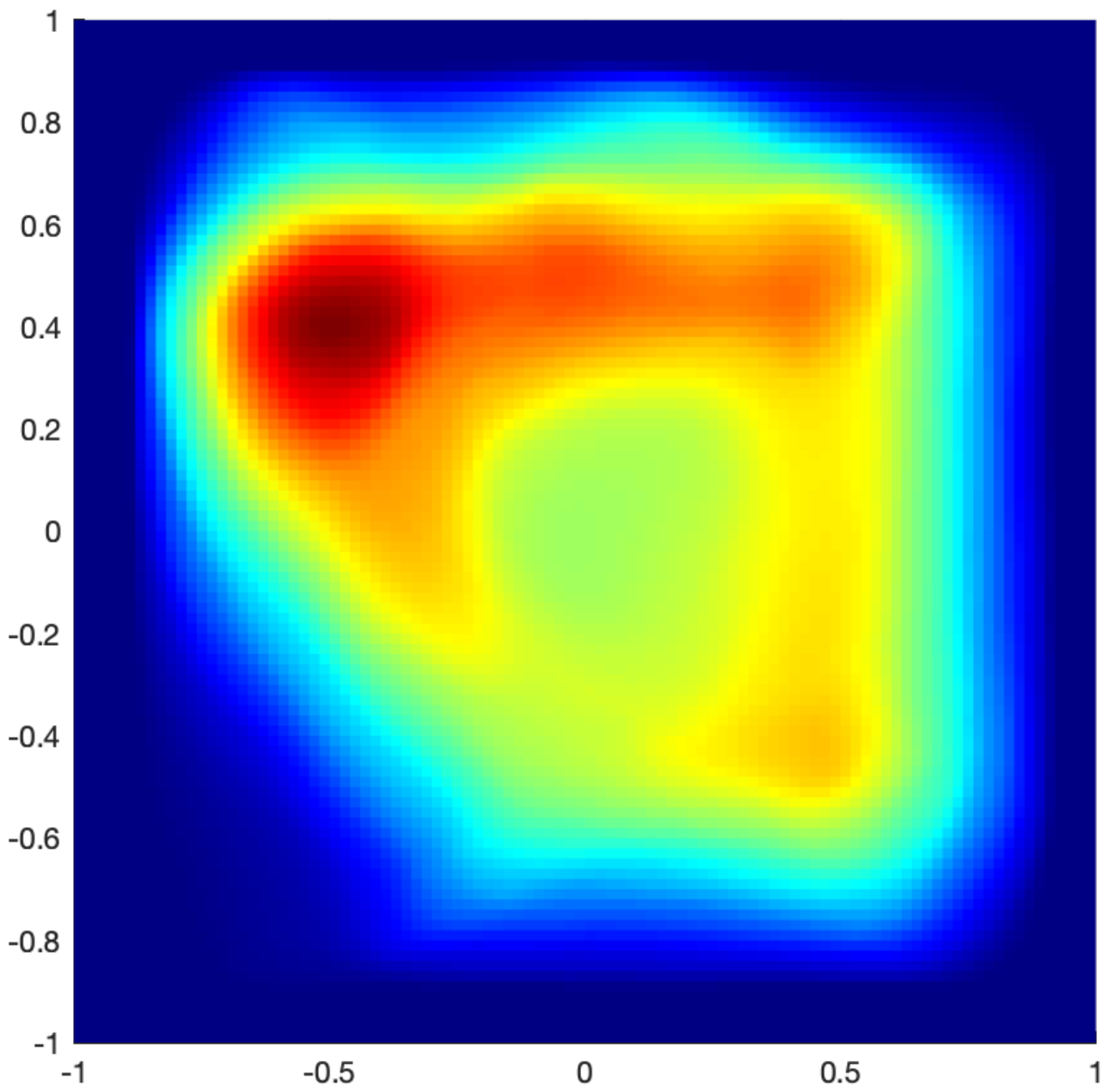}&
\includegraphics[width=1.1in]{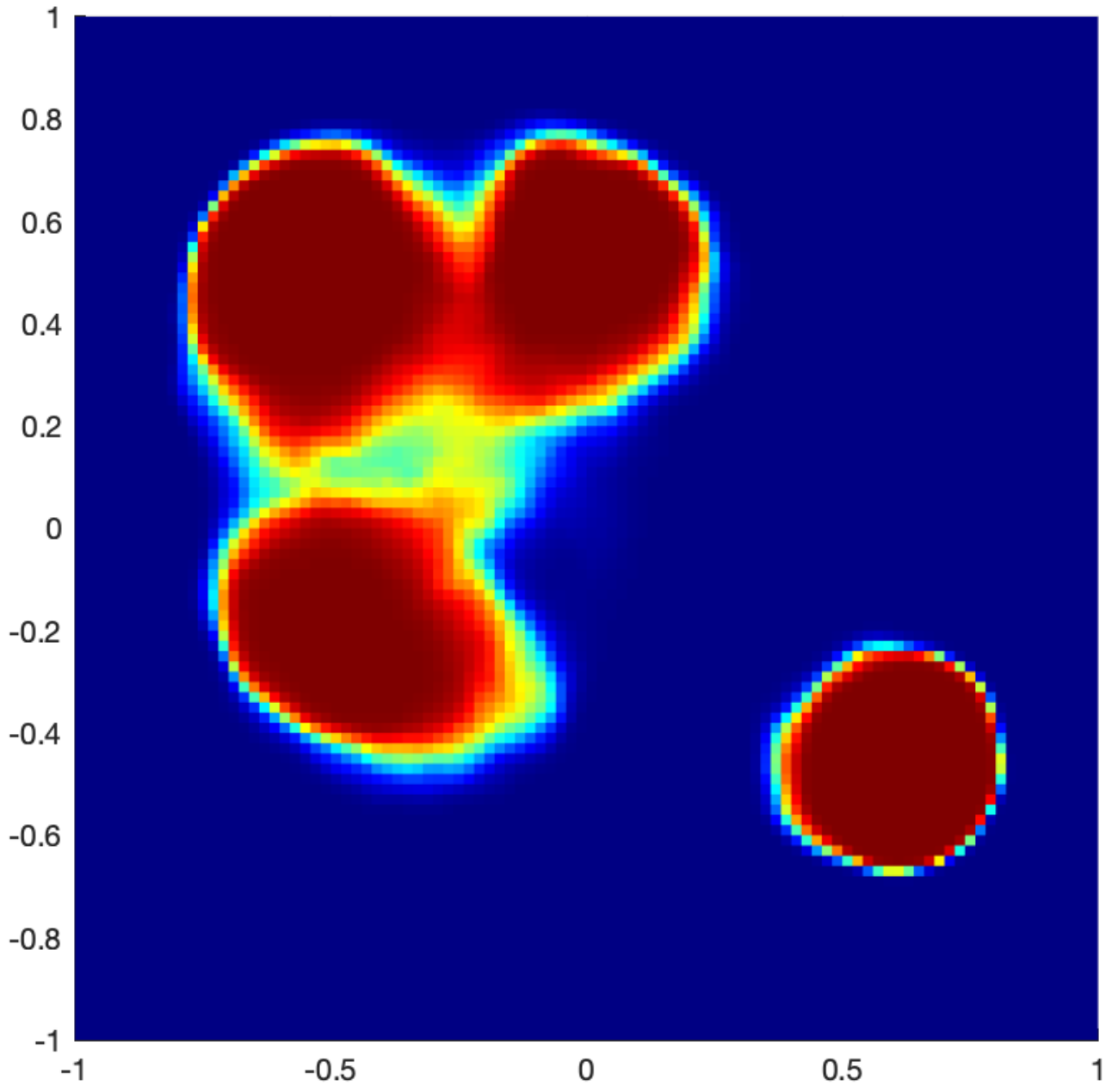}&
\includegraphics[width=1.1in]{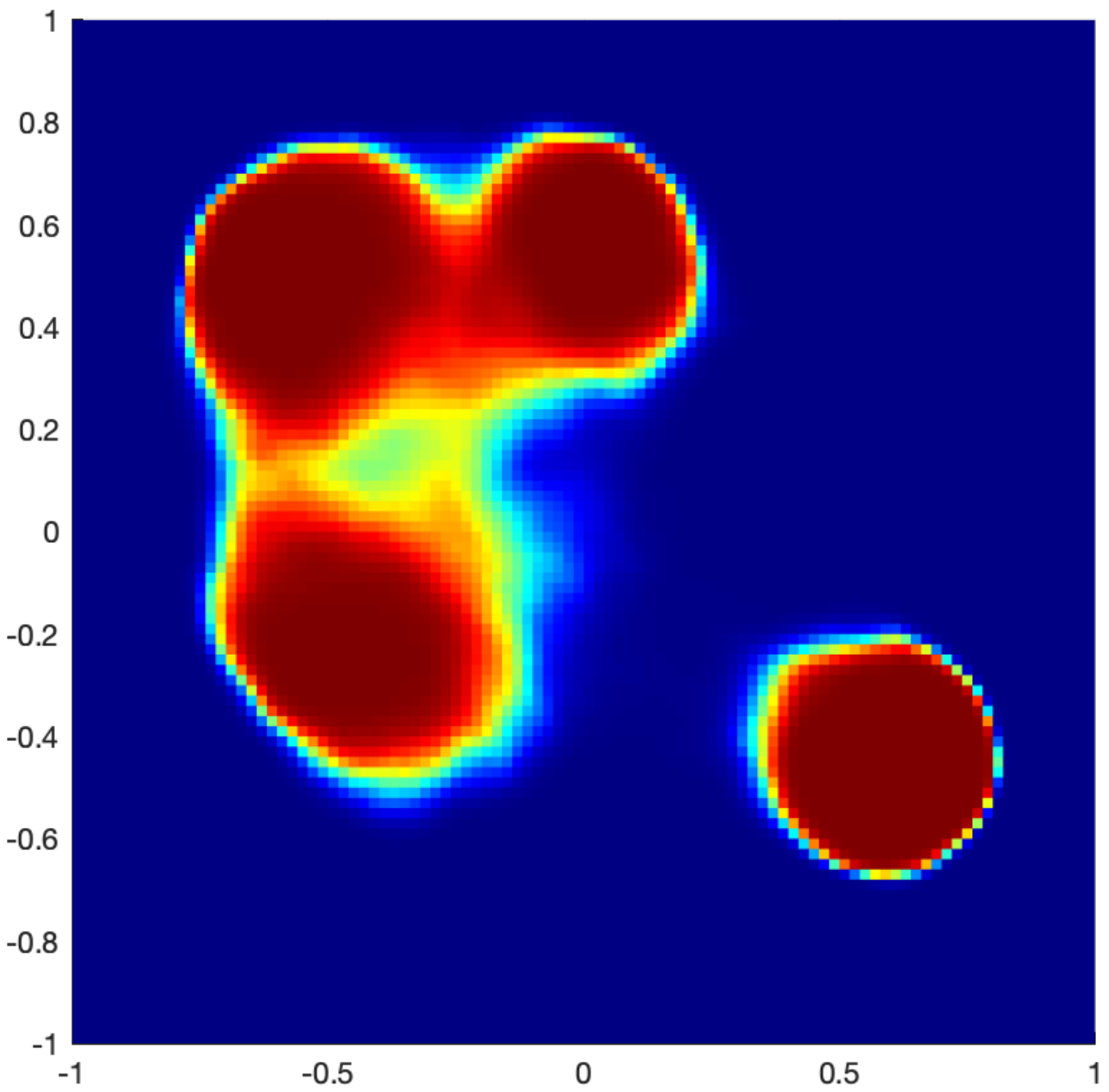}&
\includegraphics[width=1.1in]{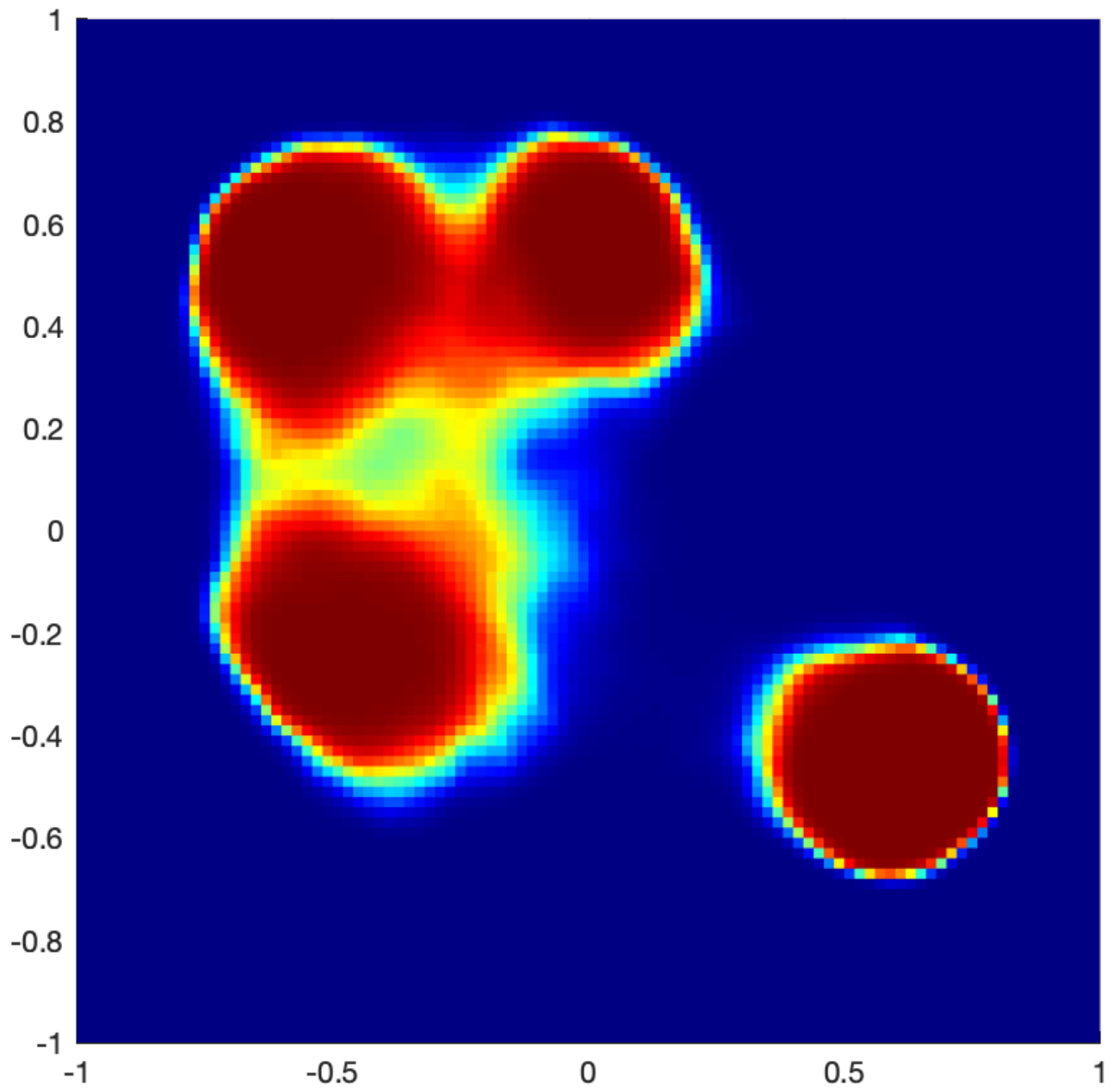}&
\includegraphics[width=1.1in]{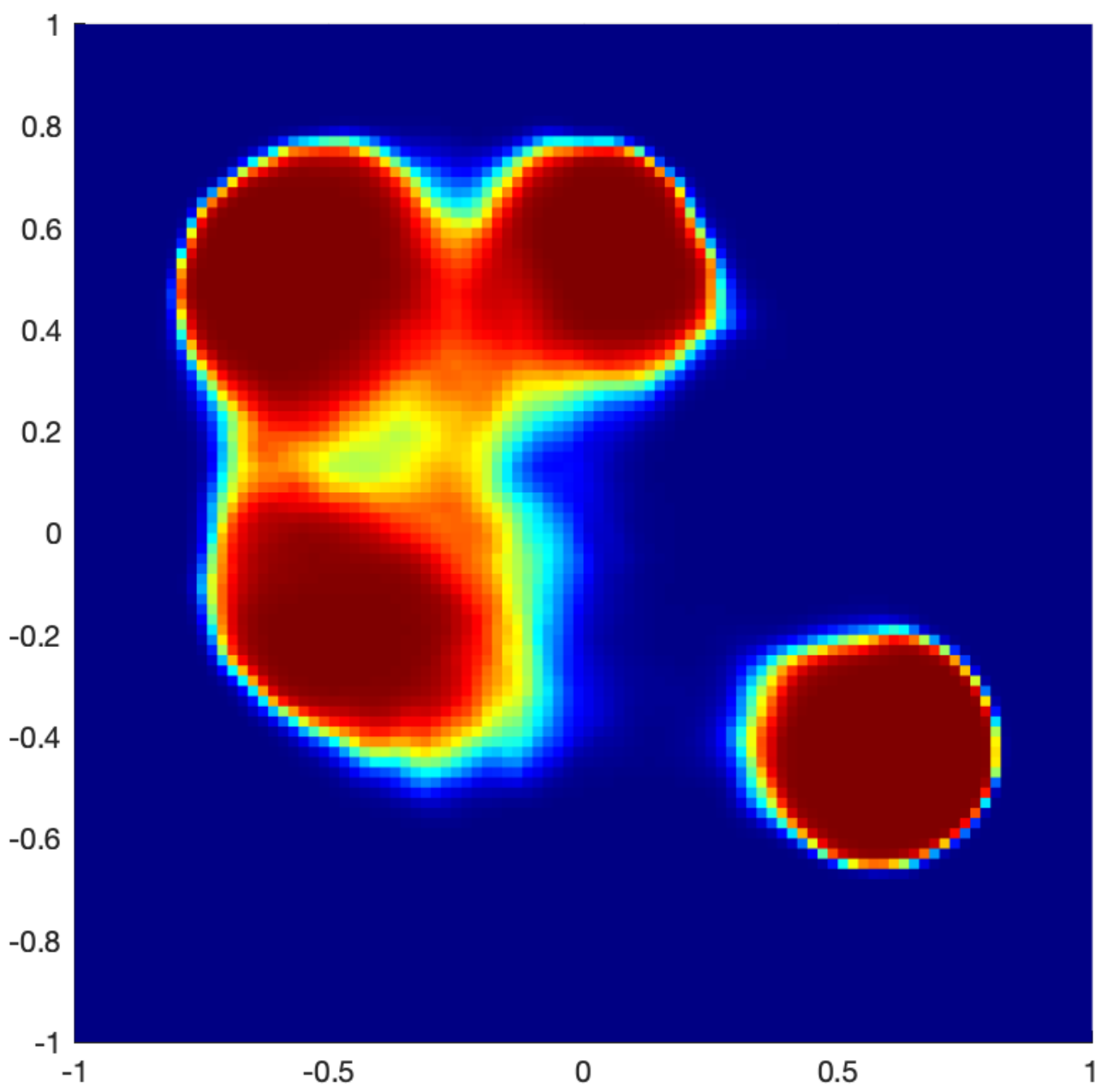}\\
\includegraphics[width=1.1in]{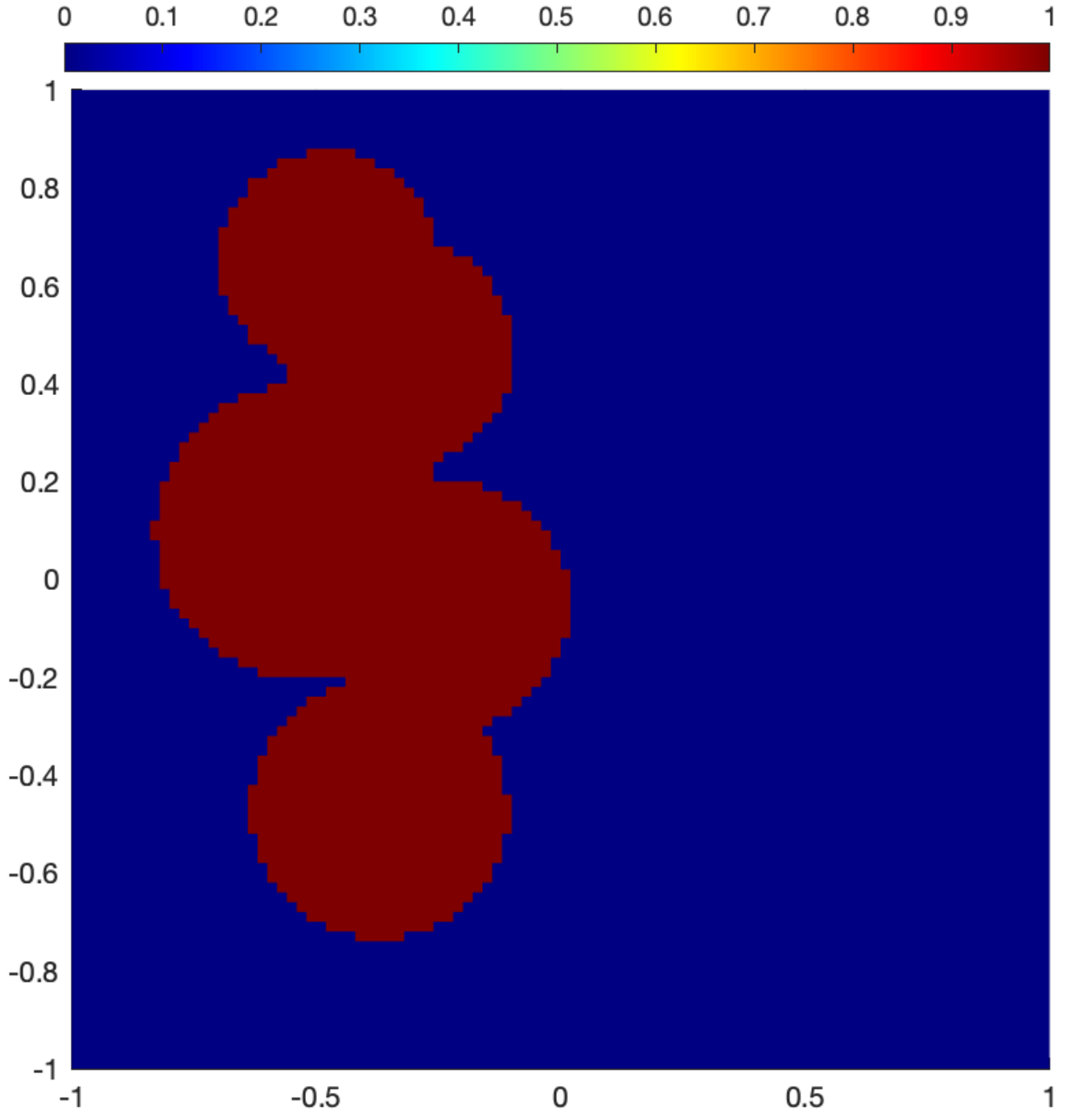}&
\includegraphics[width=1.1in]{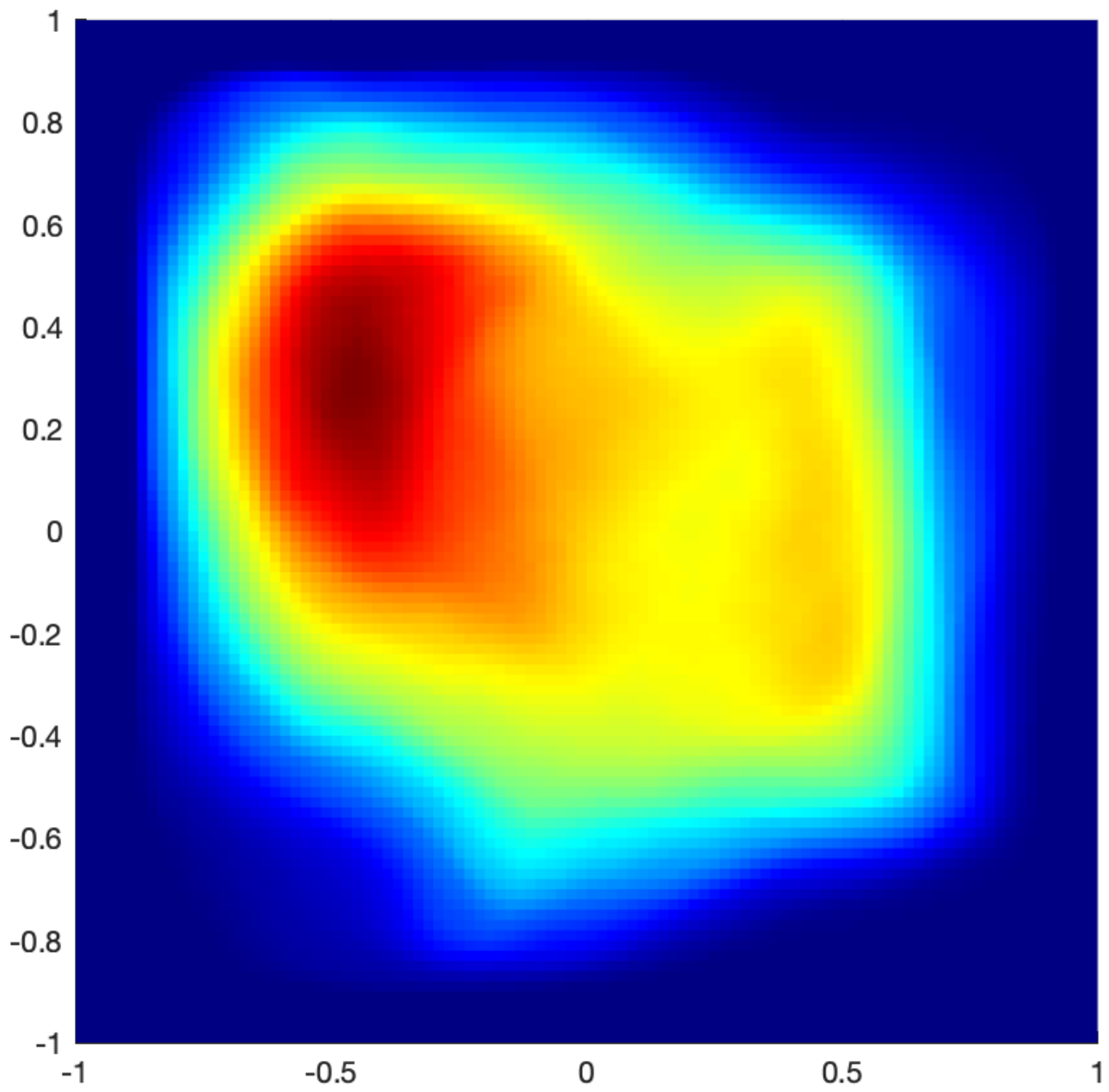}&
\includegraphics[width=1.1in]{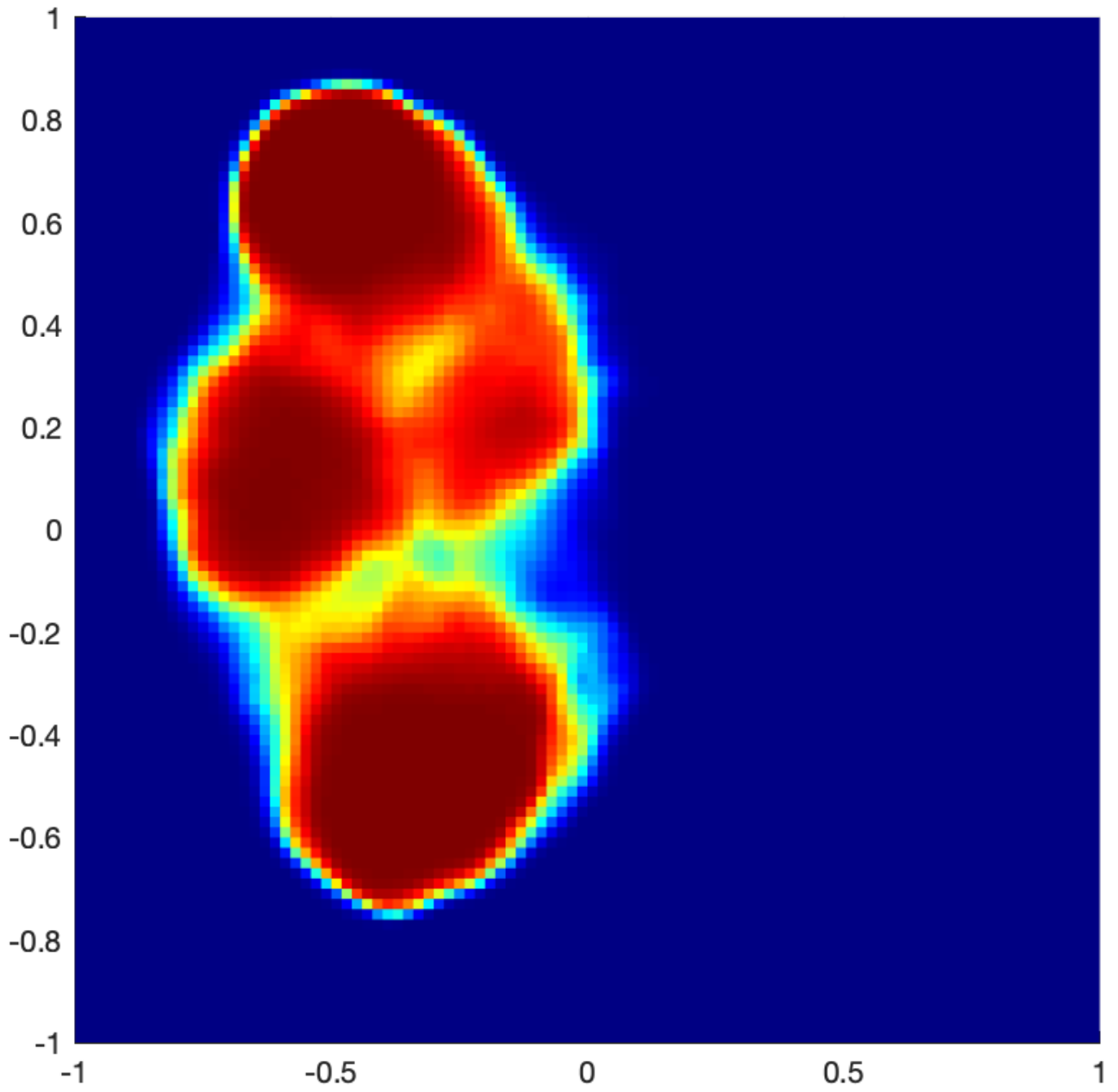}&
\includegraphics[width=1.1in]{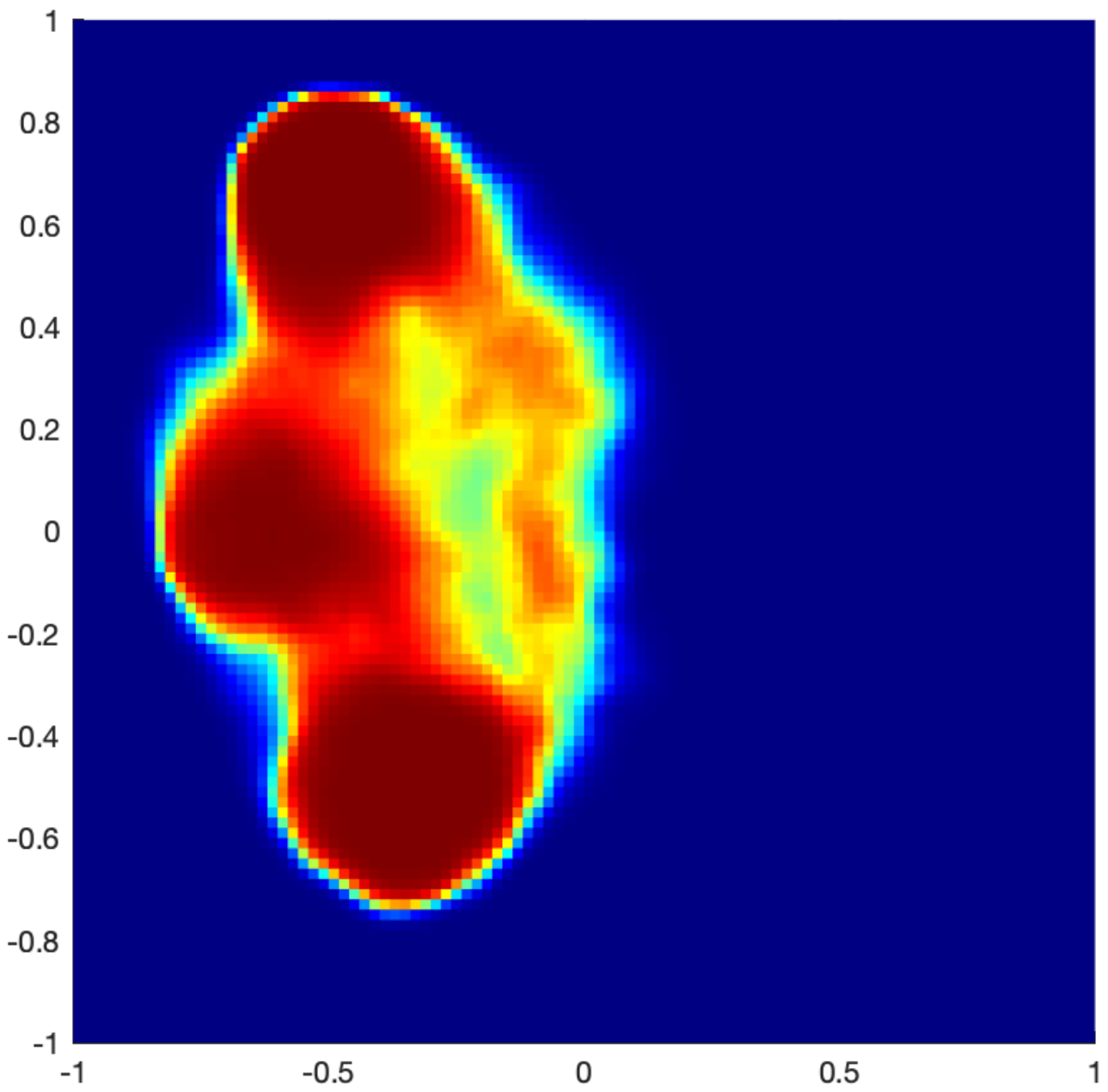}&
\includegraphics[width=1.1in]{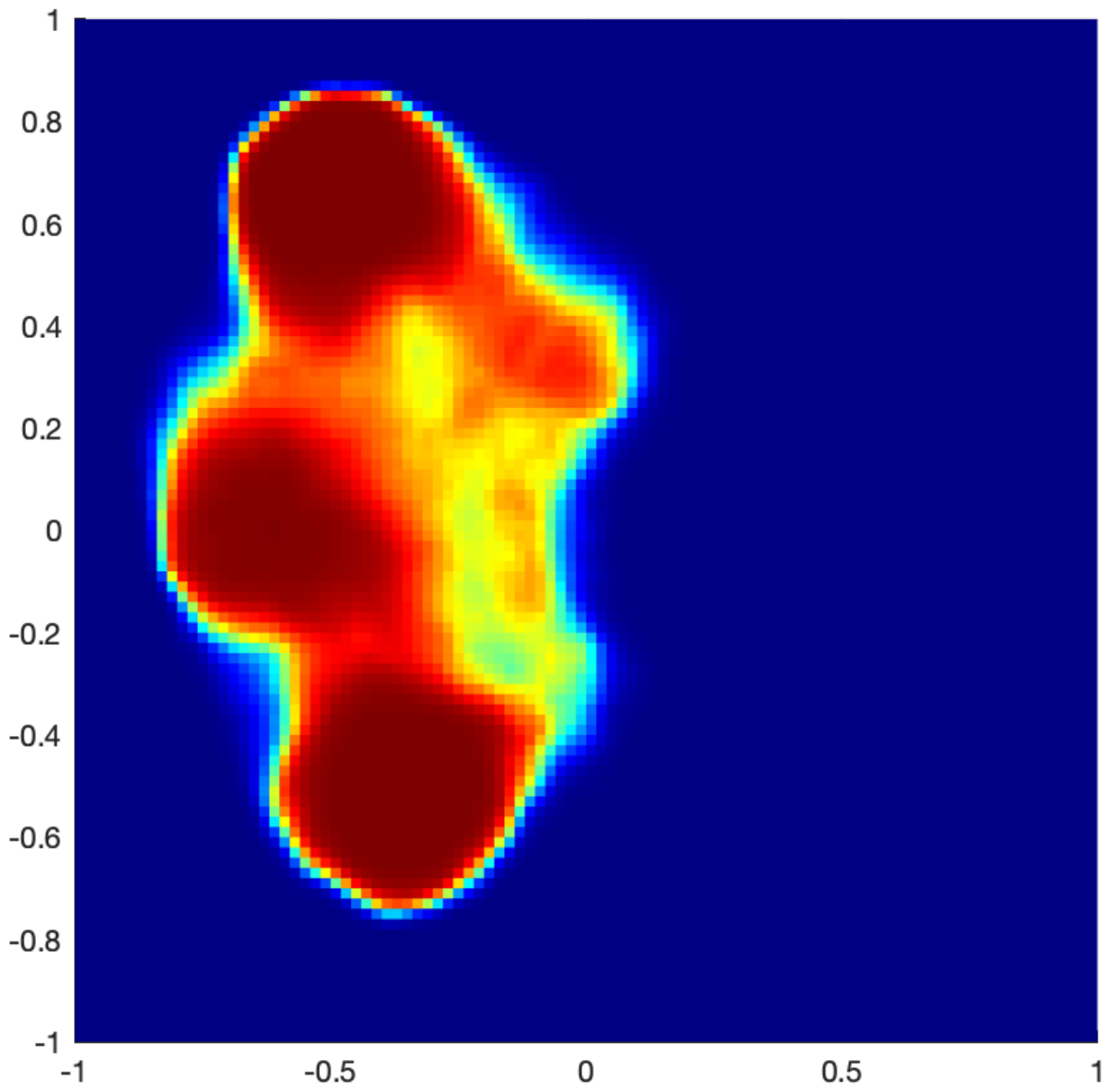}&
\includegraphics[width=1.1in]{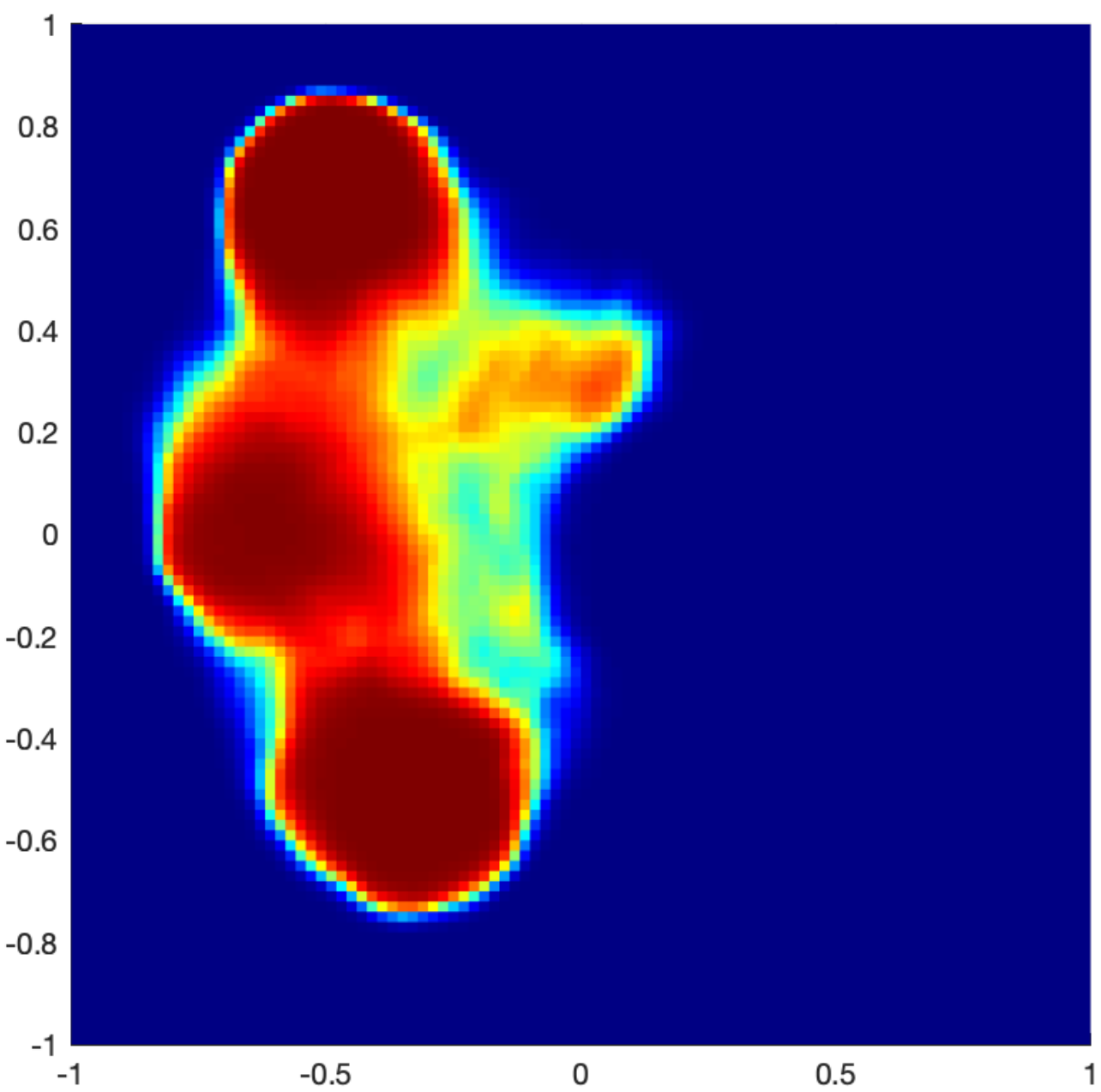}\\
\end{tabular}
  \caption{FNN-DDSM reconstruction for 3 cases in \textbf{Scenario 2} (5 circles) with different Cauchy data number and noise level: Case 1(top), Case 2(middle) and Case 3(bottom) } 
  \label{tab_FN_5cir}
\end{figure}

\begin{figure}[htbp]
\begin{tabular}{ >{\centering\arraybackslash}m{0.9in} >{\centering\arraybackslash}m{0.9in} >{\centering\arraybackslash}m{0.9in}  >{\centering\arraybackslash}m{0.9in}  >{\centering\arraybackslash}m{0.9in}  >{\centering\arraybackslash}m{0.9in} }
\centering
True coefficients &
N=1, $\delta=0$&
N=10, $\delta=0$&
N=20, $\delta=0$&
N=20, $\delta=10\%$ &
N=20, $\delta=20\%$ \\
\includegraphics[width=1.1in]{cir_num5_case1_0-eps-converted-to.pdf}&
\includegraphics[width=1.1in]{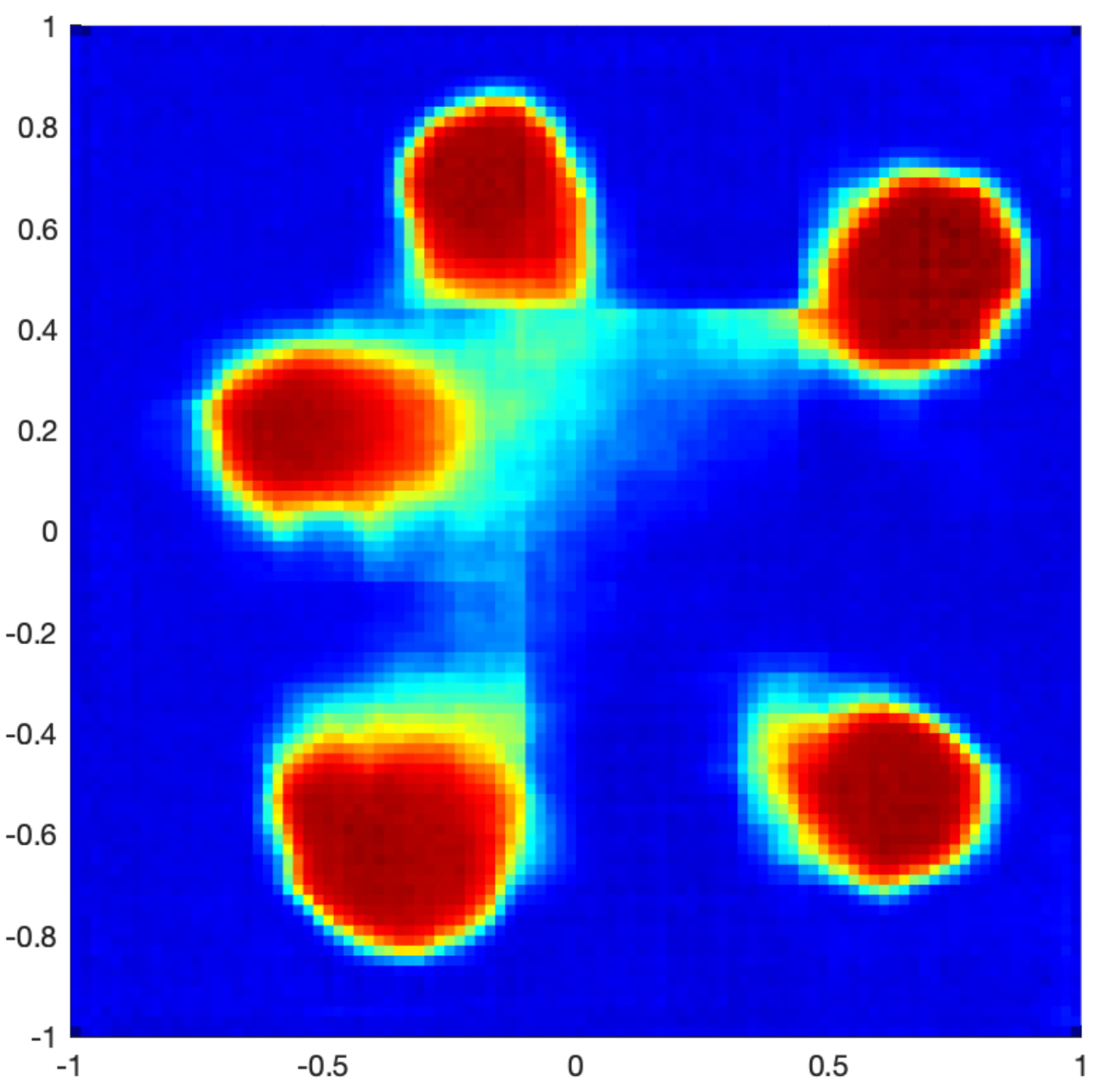}&
\includegraphics[width=1.1in]{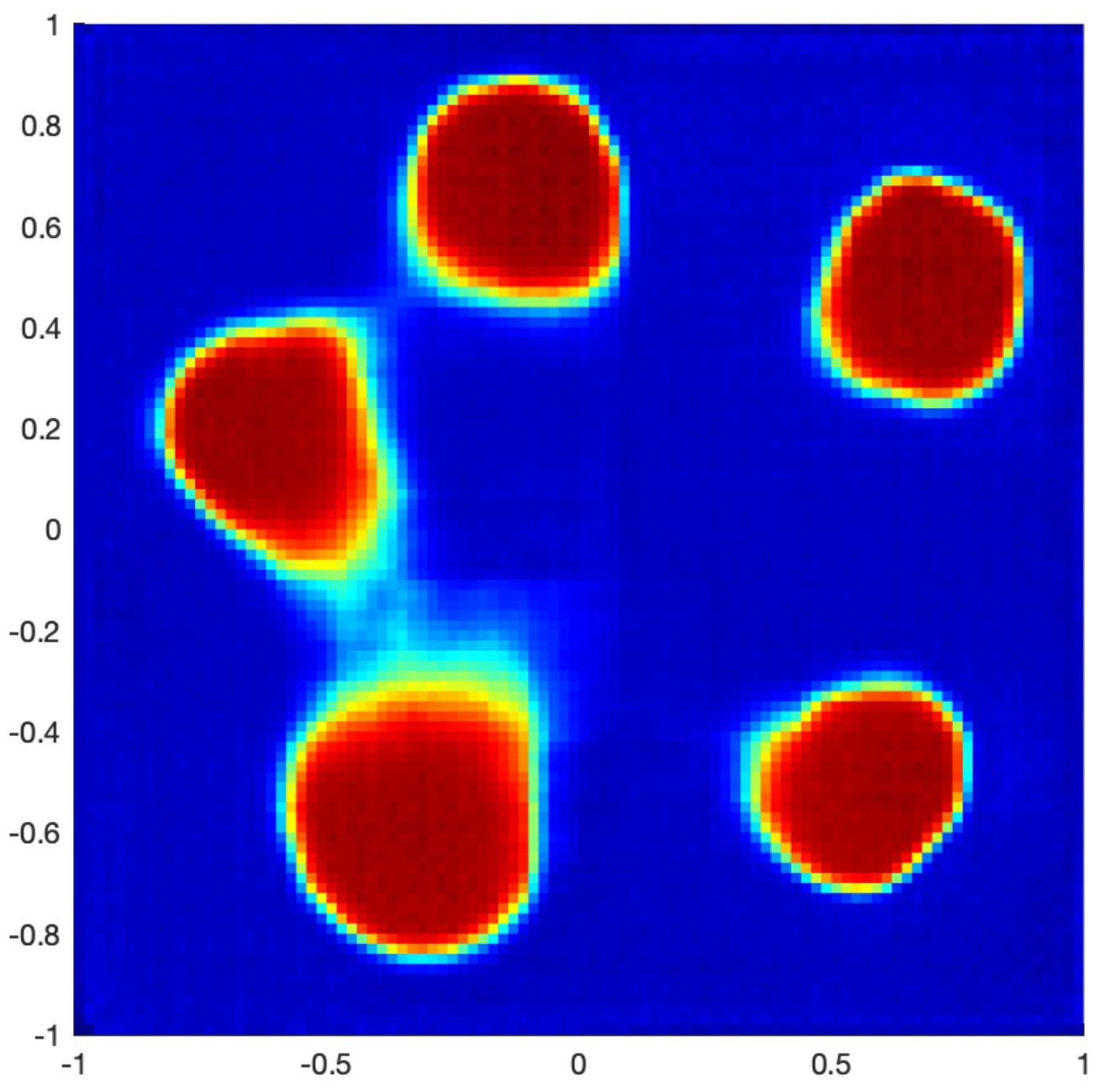}&
\includegraphics[width=1.1in]{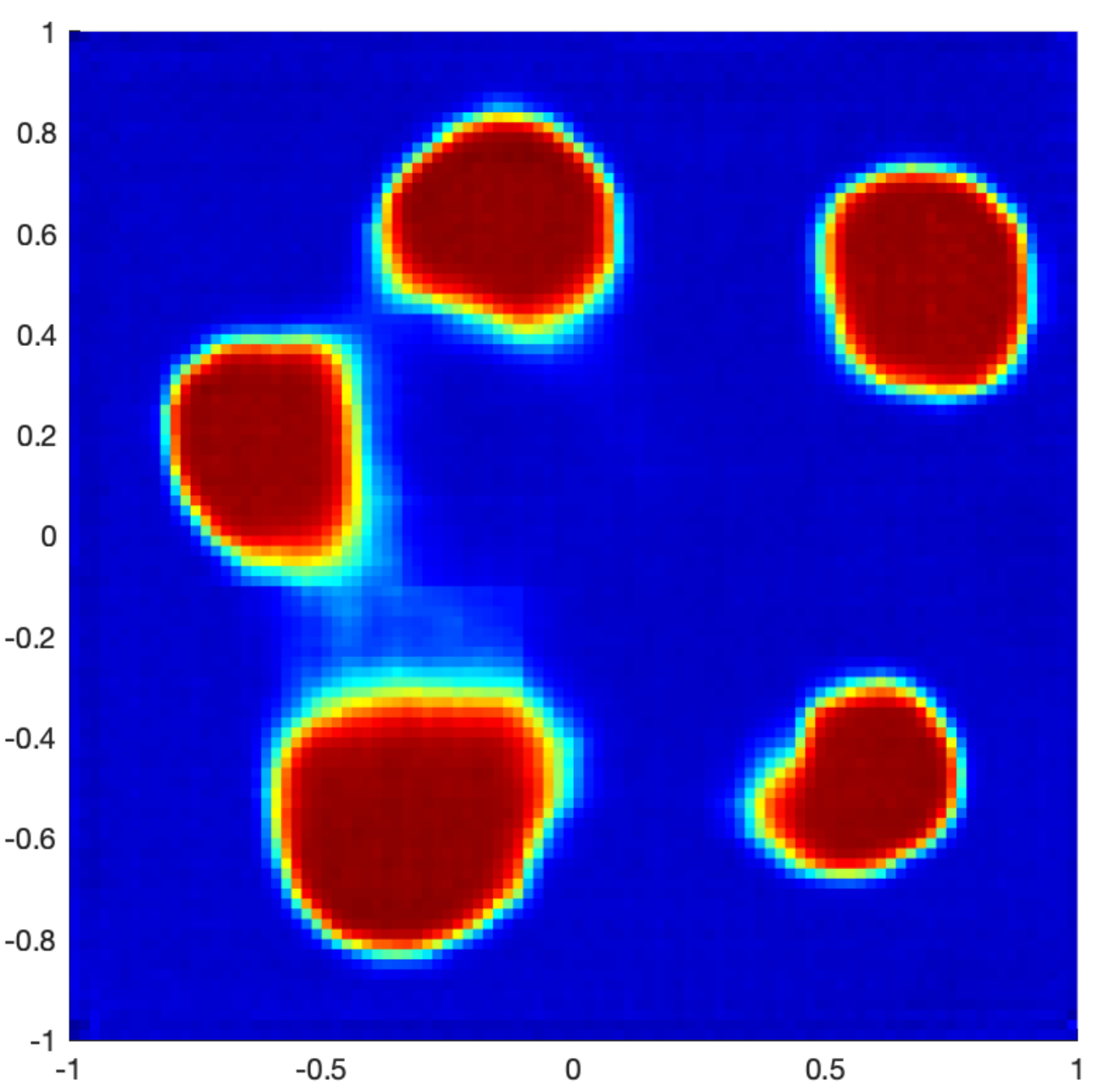}&
\includegraphics[width=1.1in]{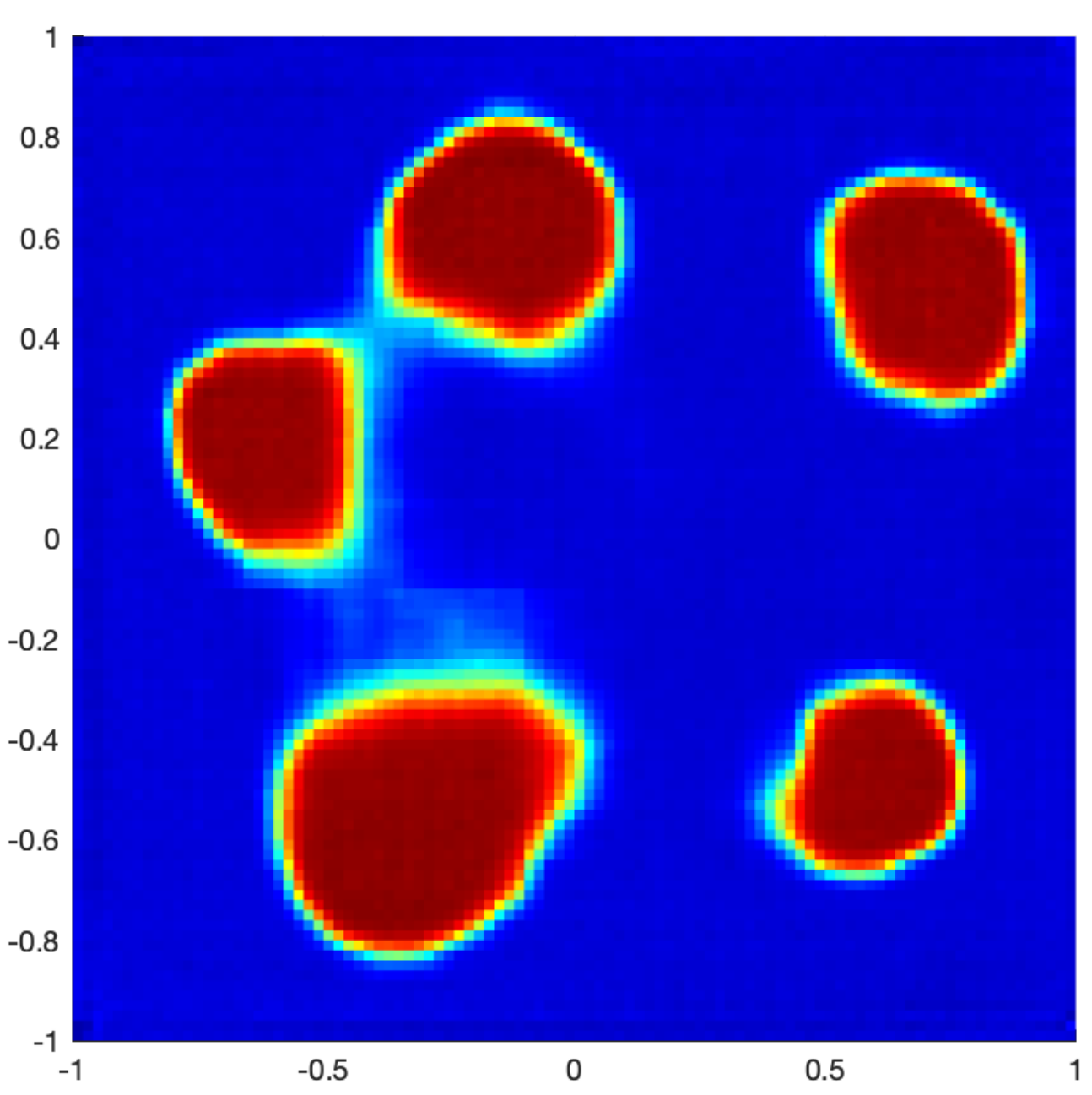}&
\includegraphics[width=1.1in]{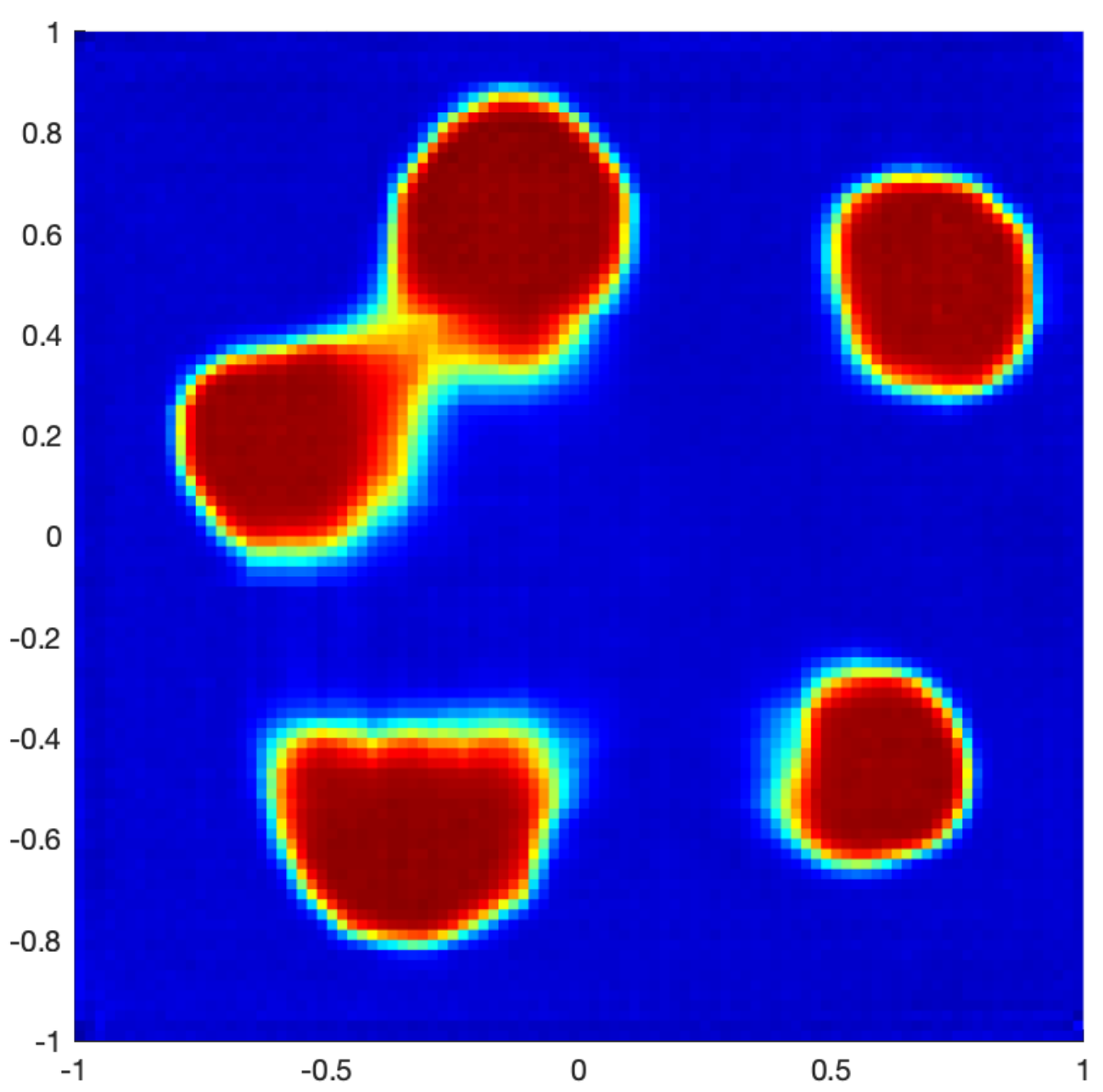}\\
\includegraphics[width=1.1in]{cir_num5_case2_0-eps-converted-to.pdf}&
\includegraphics[width=1.1in]{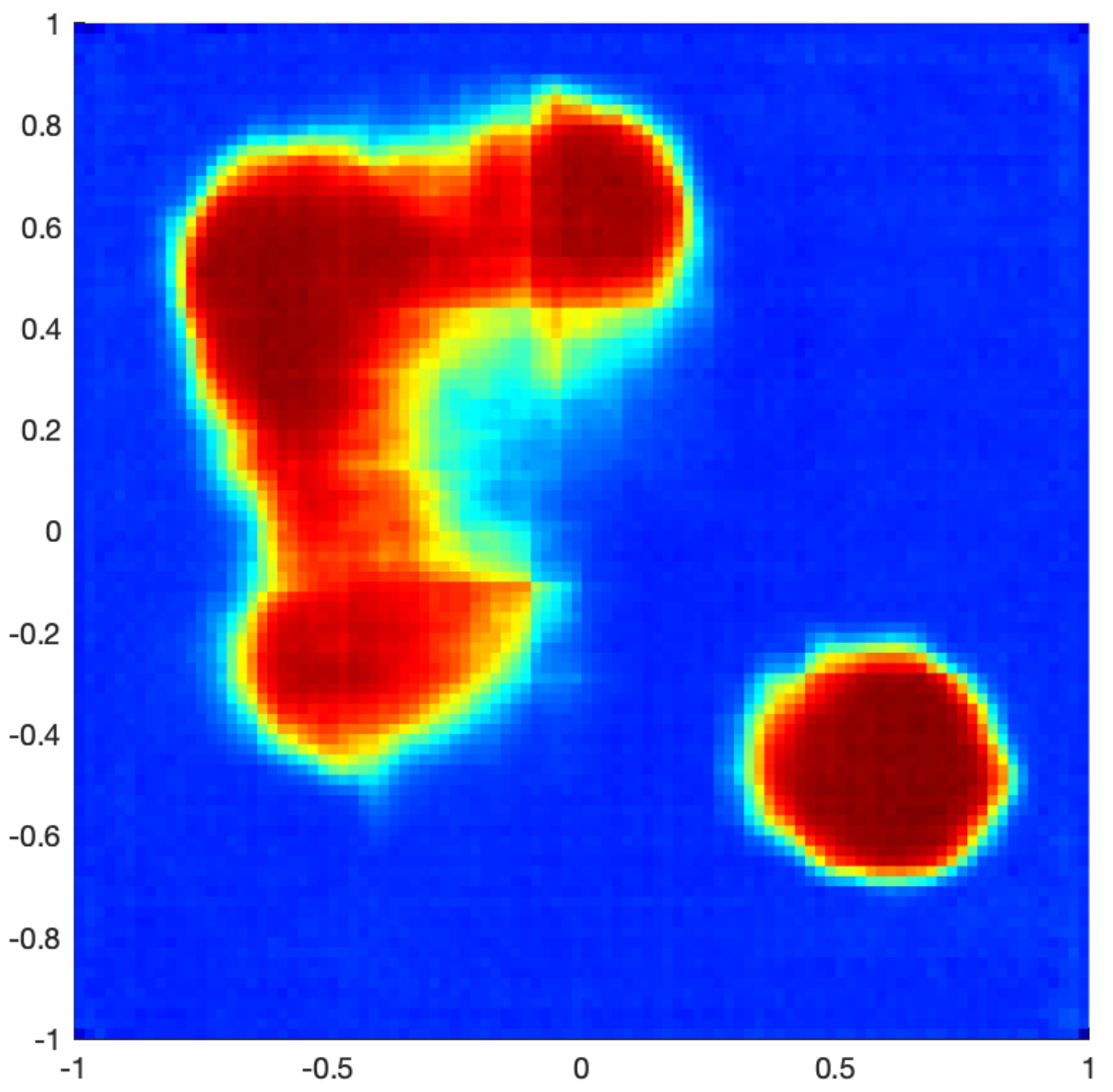}&
\includegraphics[width=1.1in]{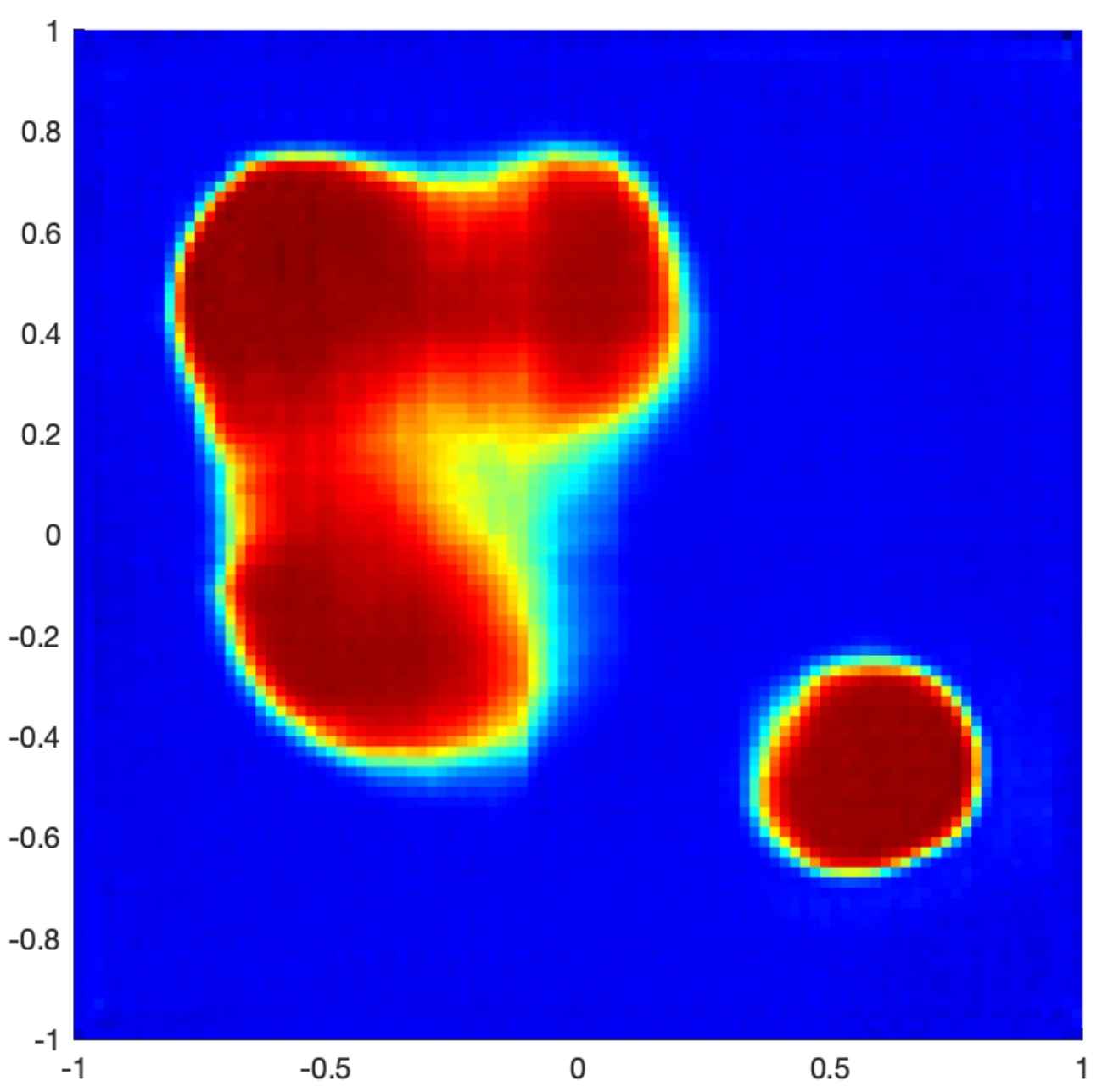}&
\includegraphics[width=1.1in]{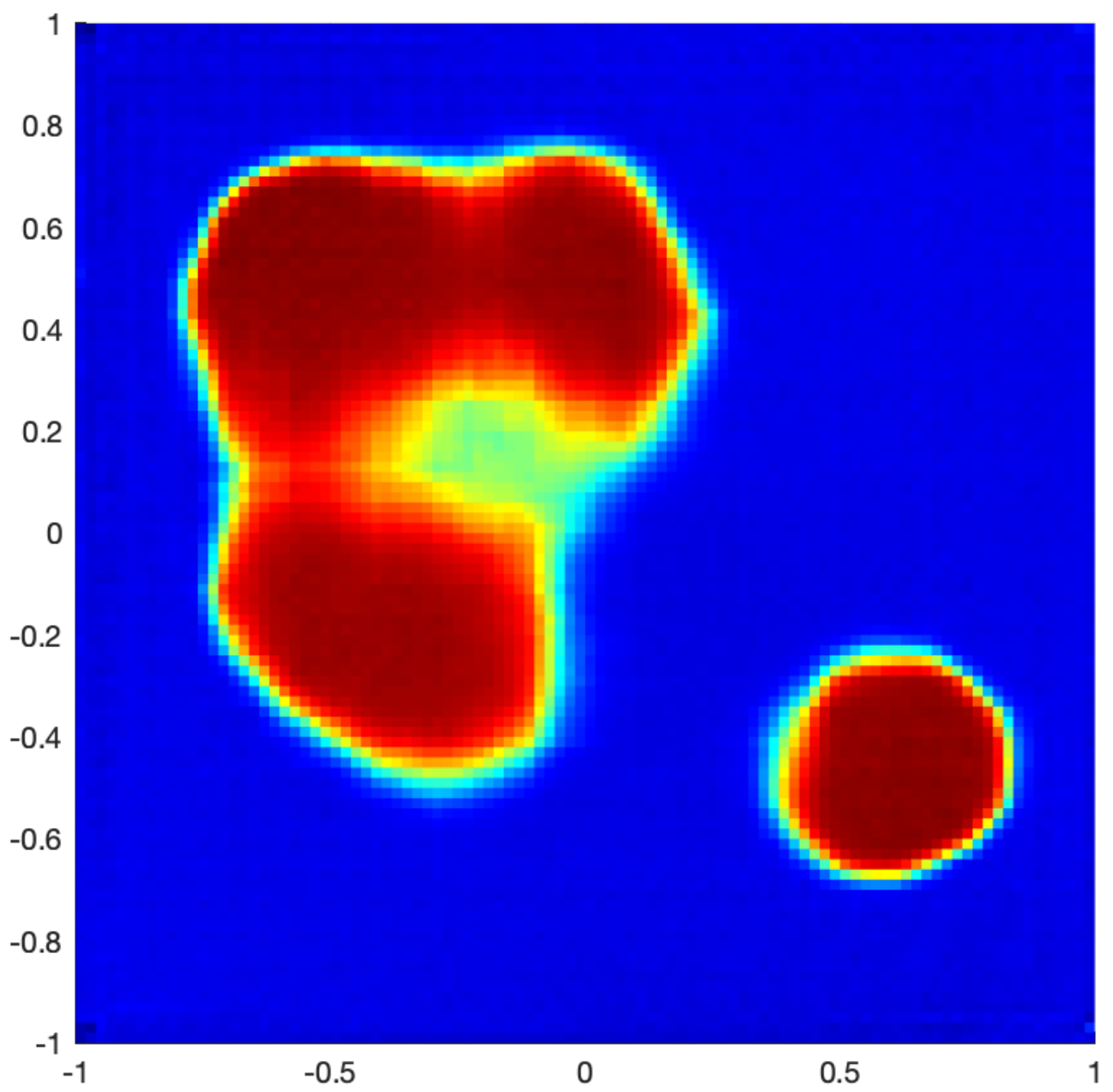}&
\includegraphics[width=1.1in]{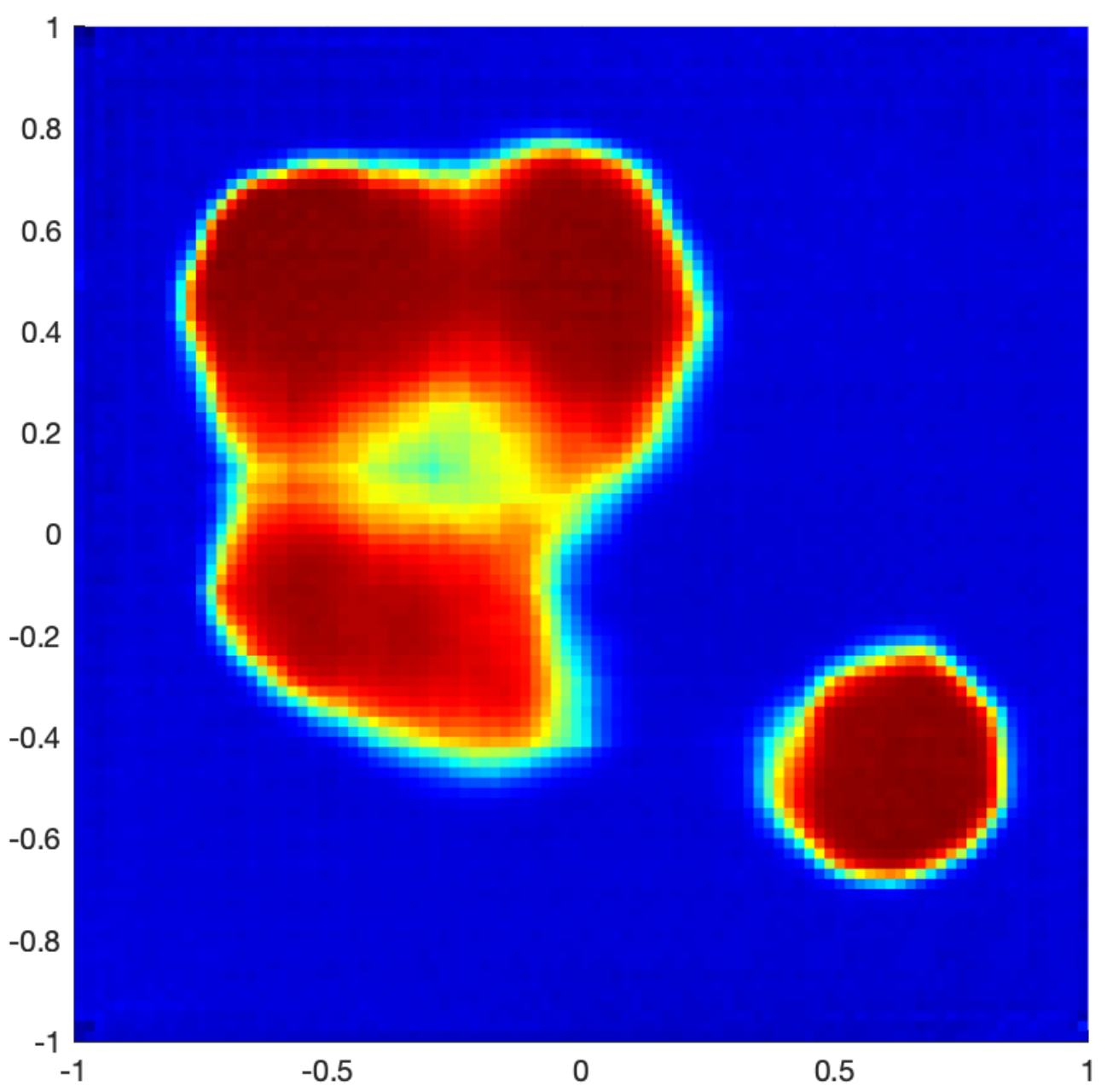}&
\includegraphics[width=1.1in]{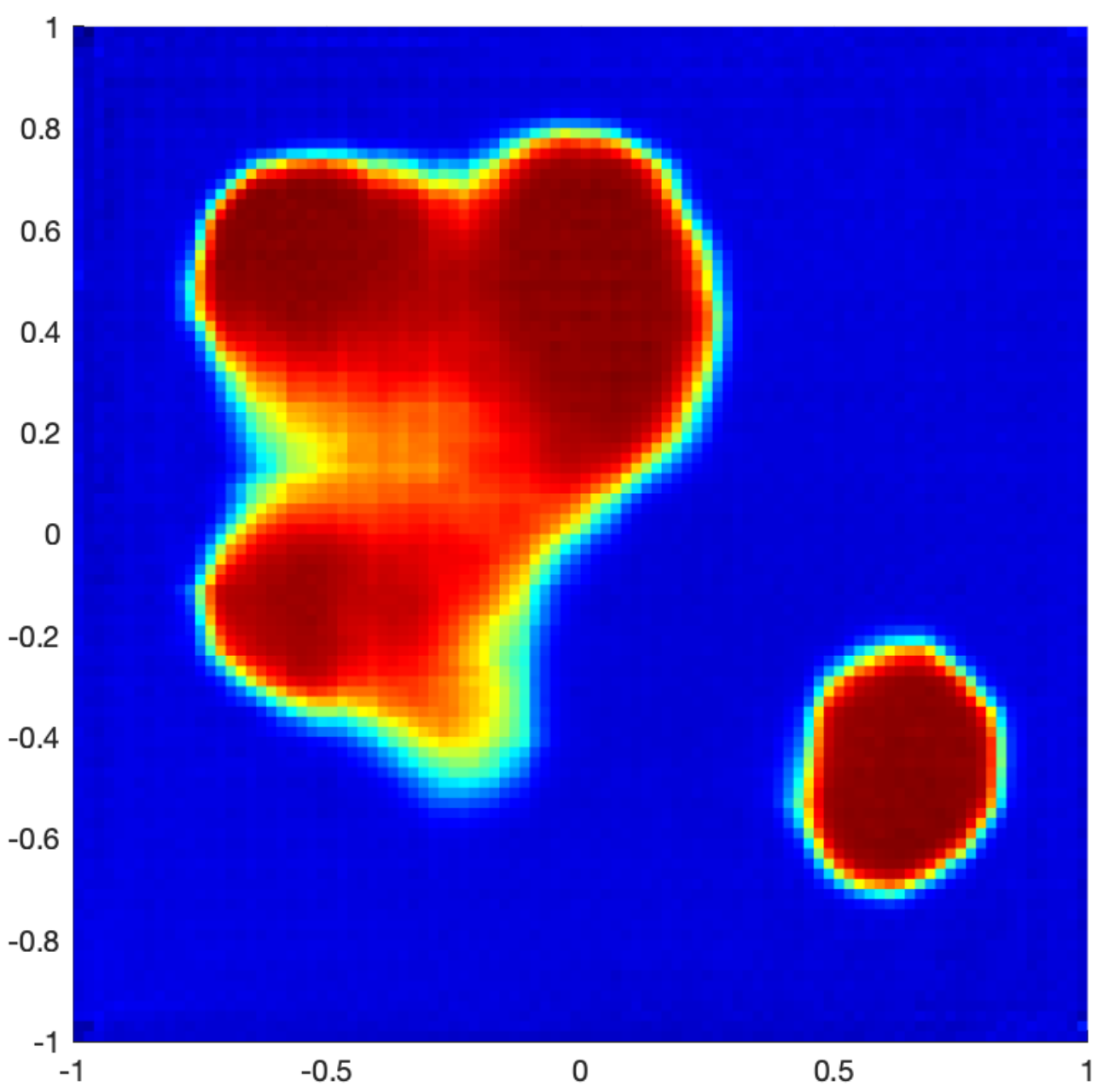}\\
\includegraphics[width=1.1in]{cir_num5_case3_0-eps-converted-to.pdf}&
\includegraphics[width=1.1in]{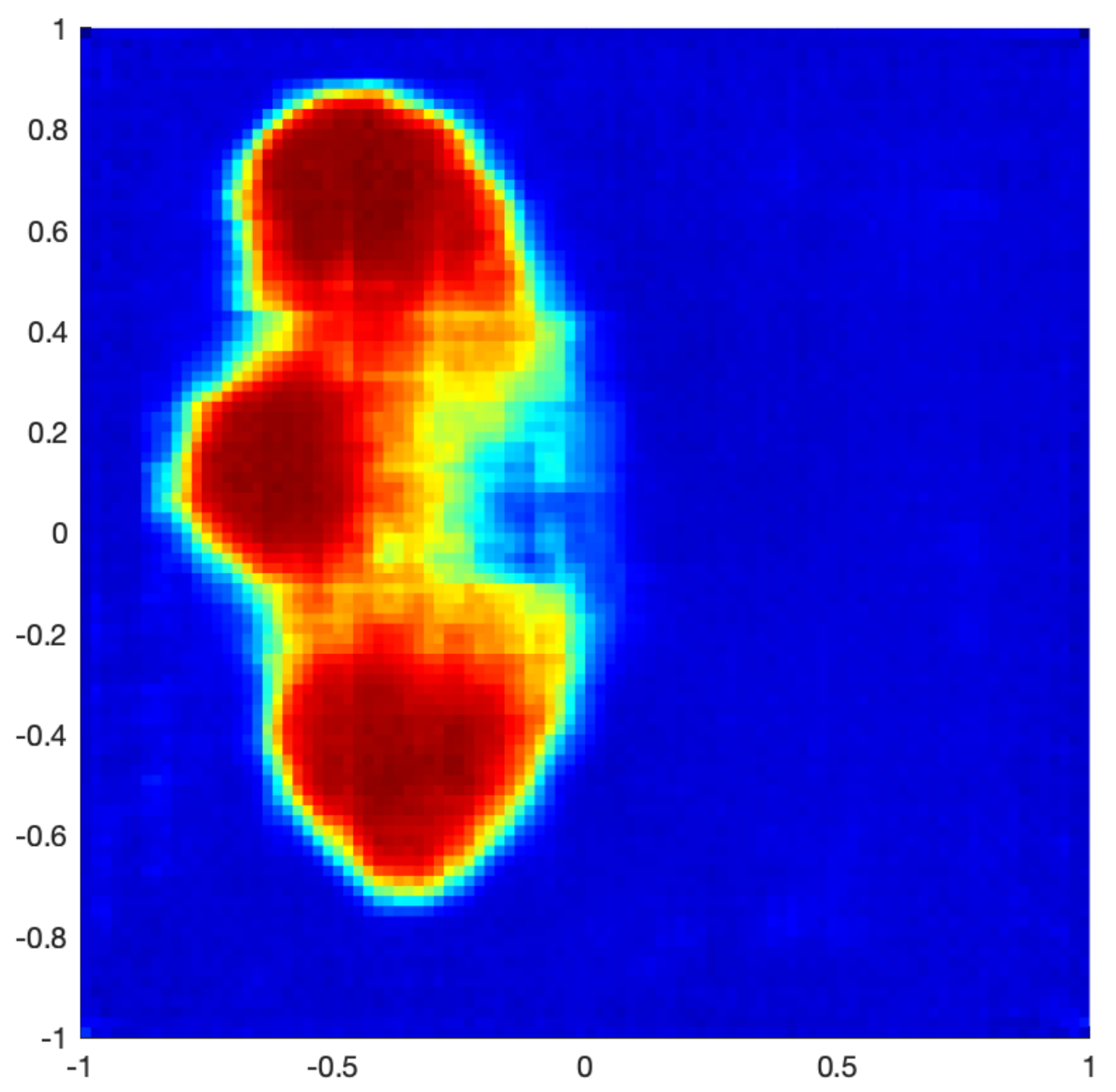}&
\includegraphics[width=1.1in]{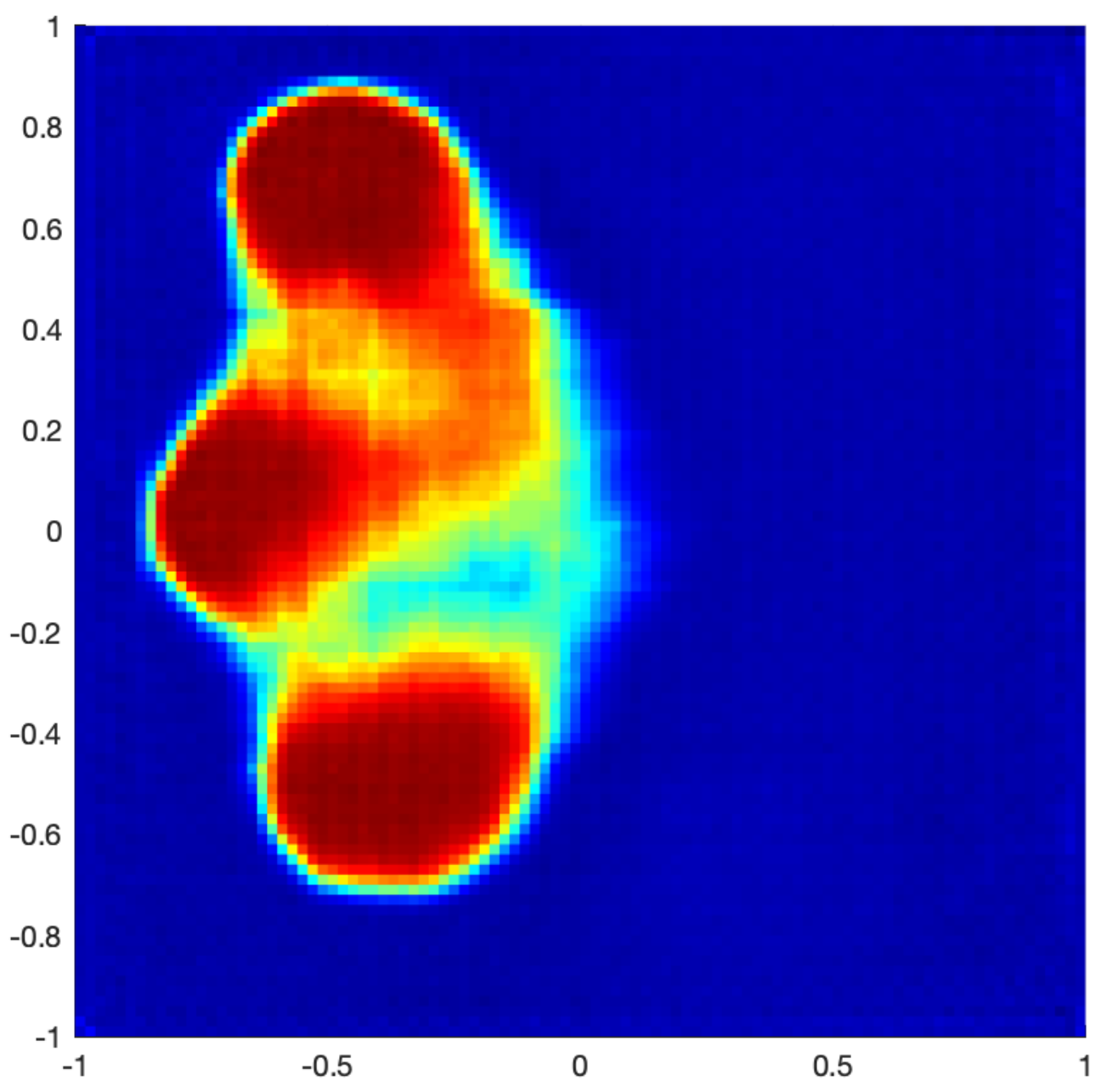}&
\includegraphics[width=1.1in]{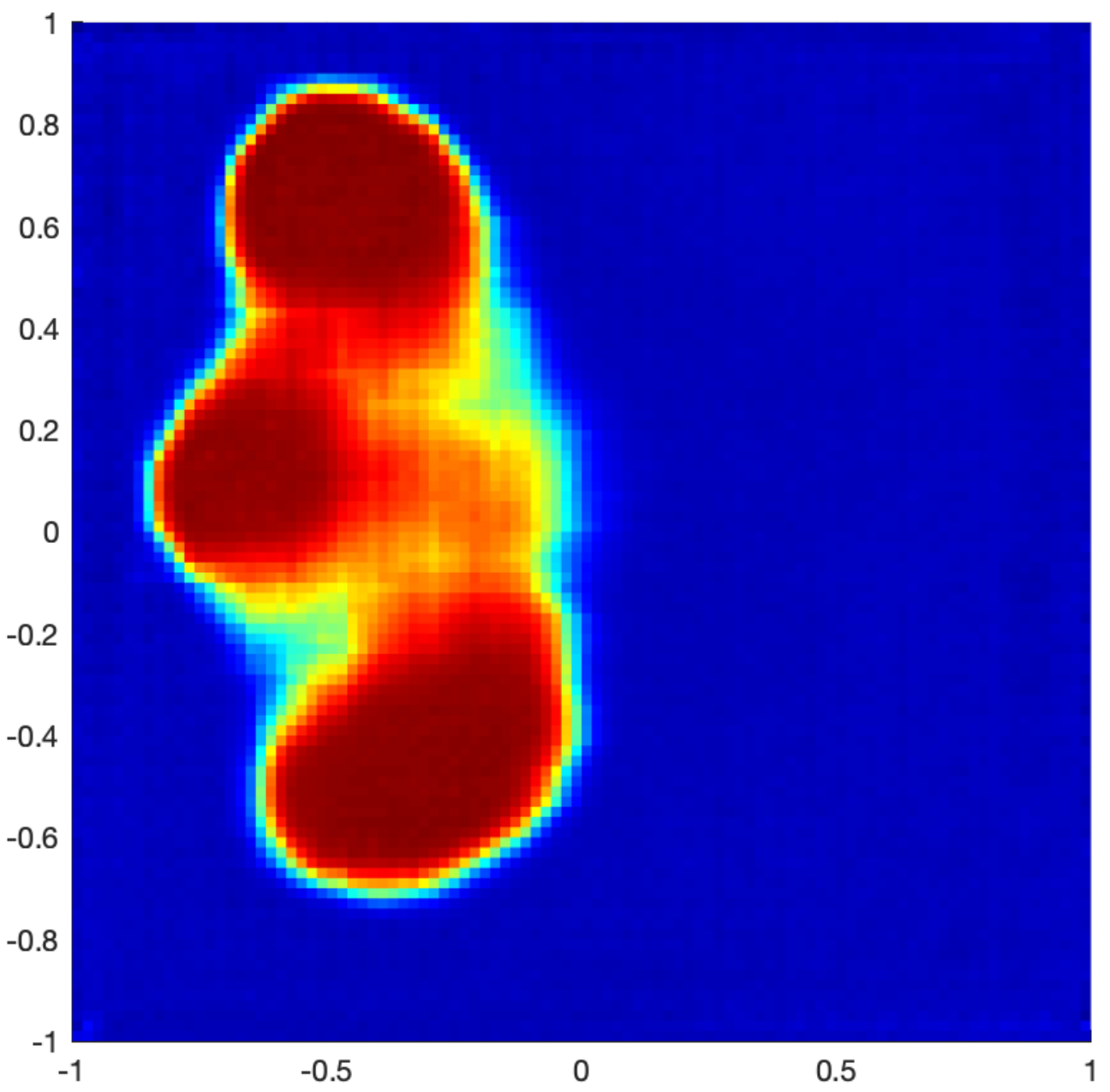}&
\includegraphics[width=1.1in]{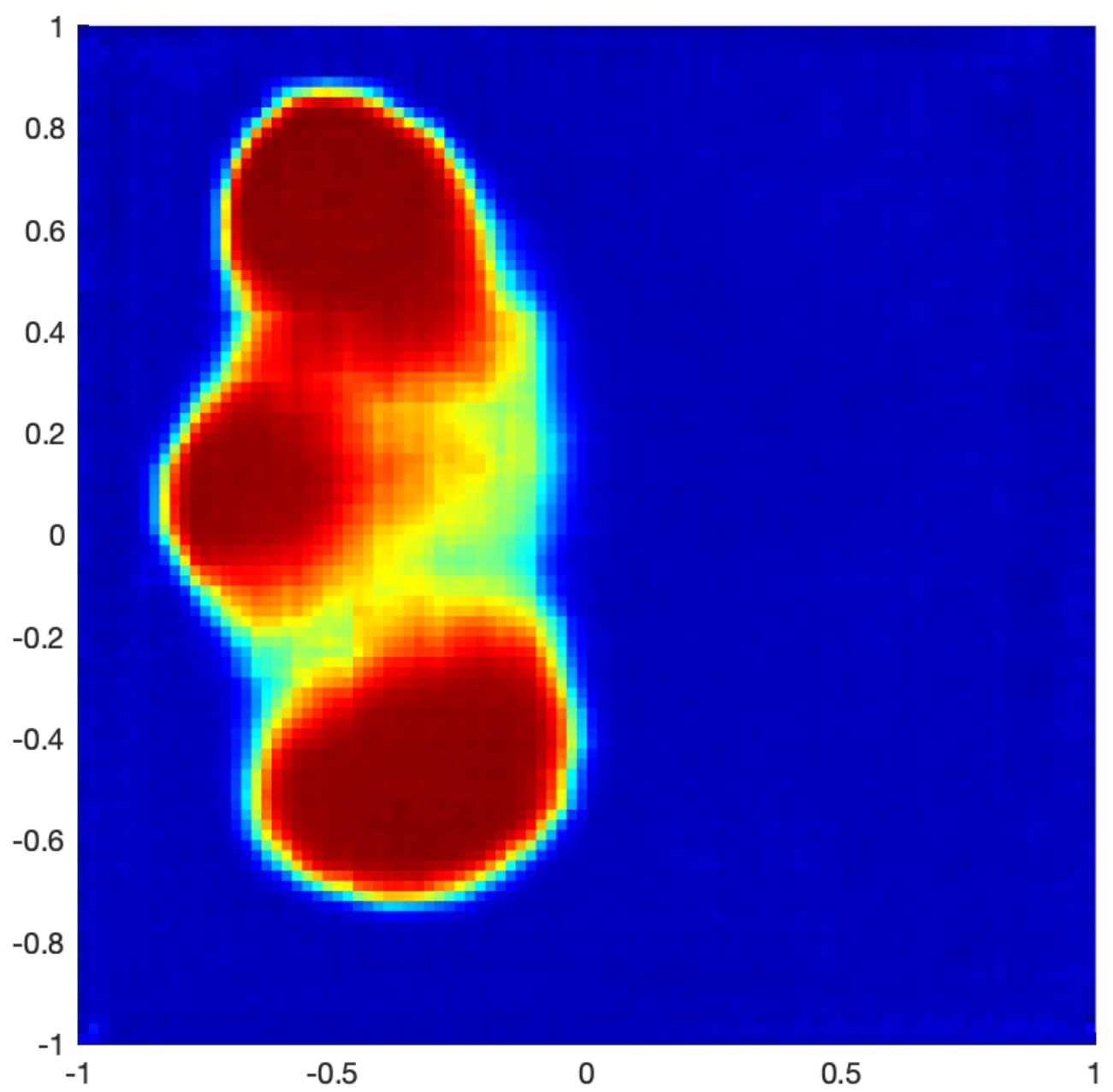}&
\includegraphics[width=1.1in]{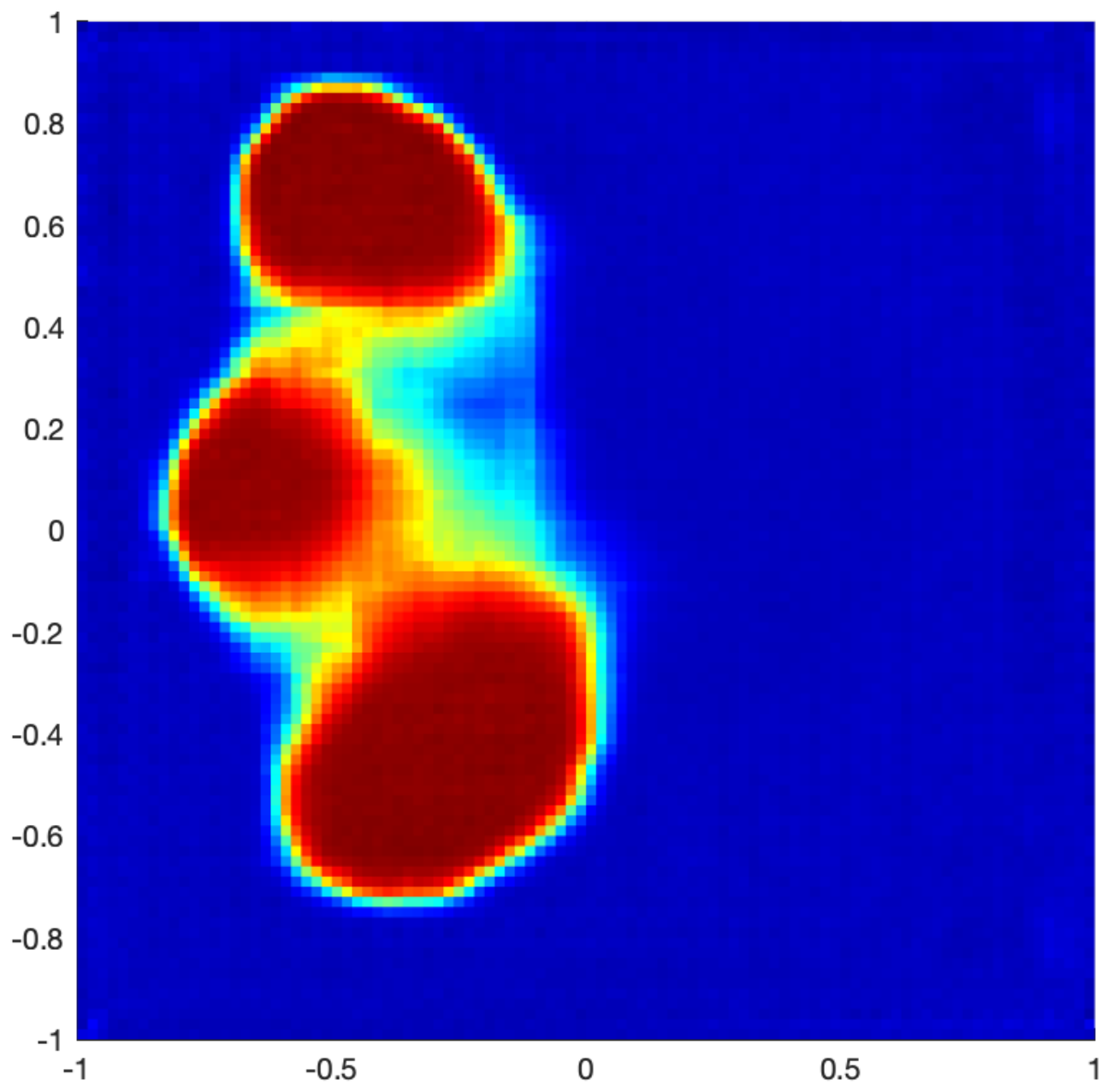}\\
\end{tabular}
  \caption{CNN-DDSM reconstruction for 3 cases in \textbf{Scenario 2} (5 circles) with different Cauchy data number and noise level: Case 1(top), Case 2(middle) and Case 3(bottom) } 
  \label{tab_5cir}
\end{figure}

\begin{figure}[htbp]
\begin{tabular}{ >{\centering\arraybackslash}m{0.9in} >{\centering\arraybackslash}m{0.9in} >{\centering\arraybackslash}m{0.9in}  >{\centering\arraybackslash}m{0.9in}  >{\centering\arraybackslash}m{0.9in}  >{\centering\arraybackslash}m{0.9in} }
\centering
True coefficients &
N=1, $\delta=0$&
N=10, $\delta=0$&
N=20, $\delta=0$&
N=20, $\delta=10\%$ &
N=20, $\delta=20\%$ \\
\includegraphics[width=1.1in]{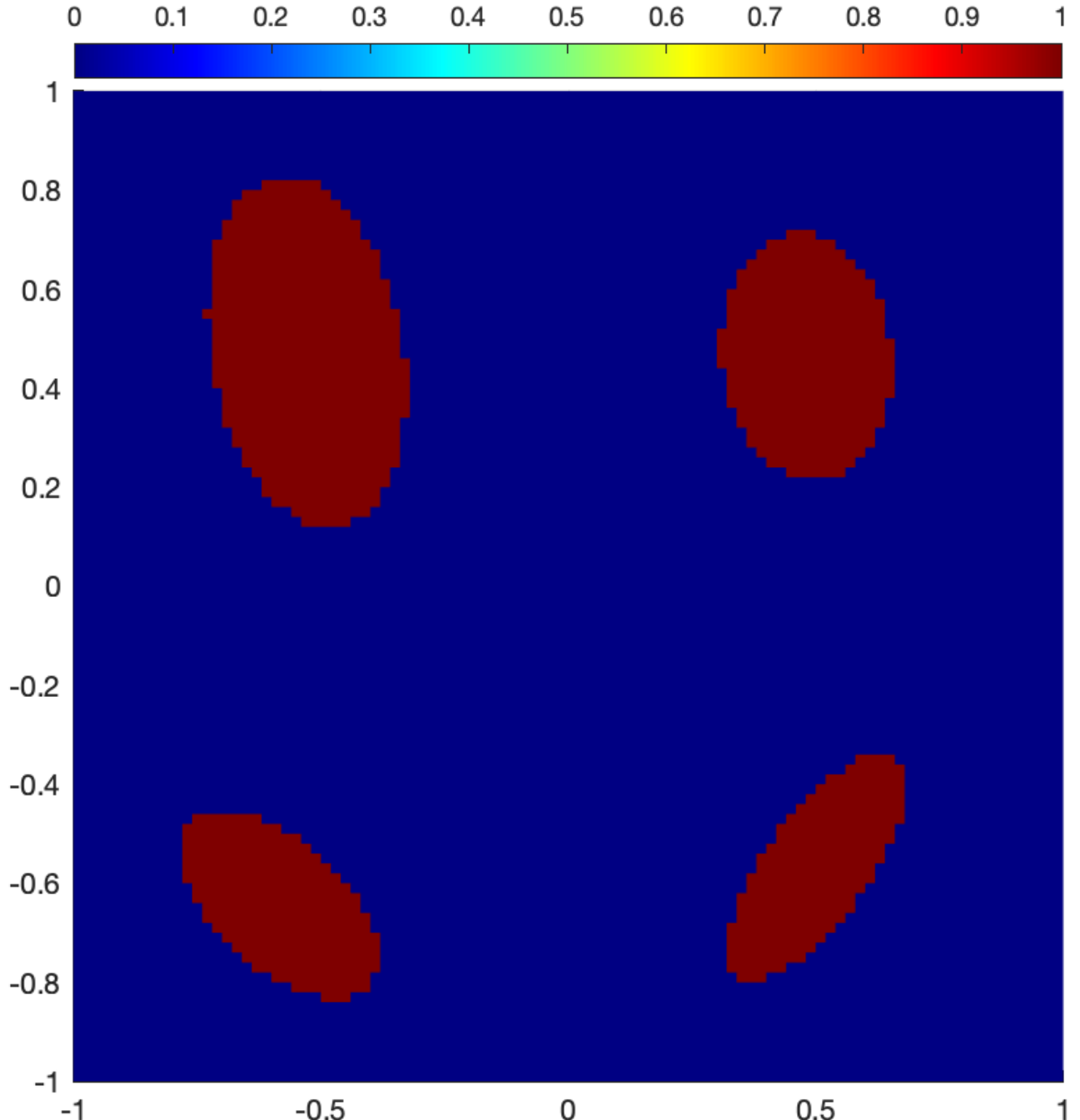}&
\includegraphics[width=1.1in]{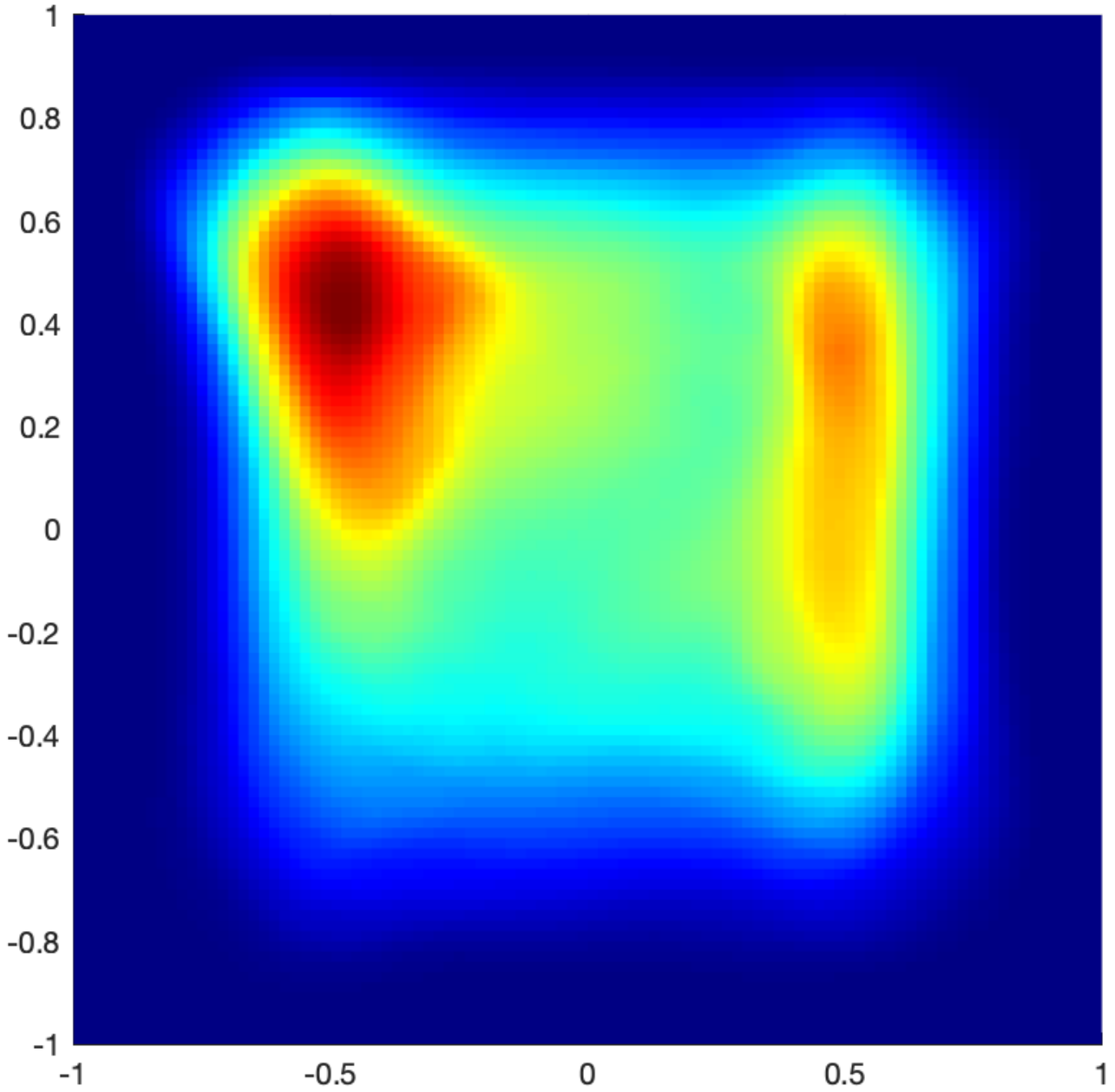}&
\includegraphics[width=1.1in]{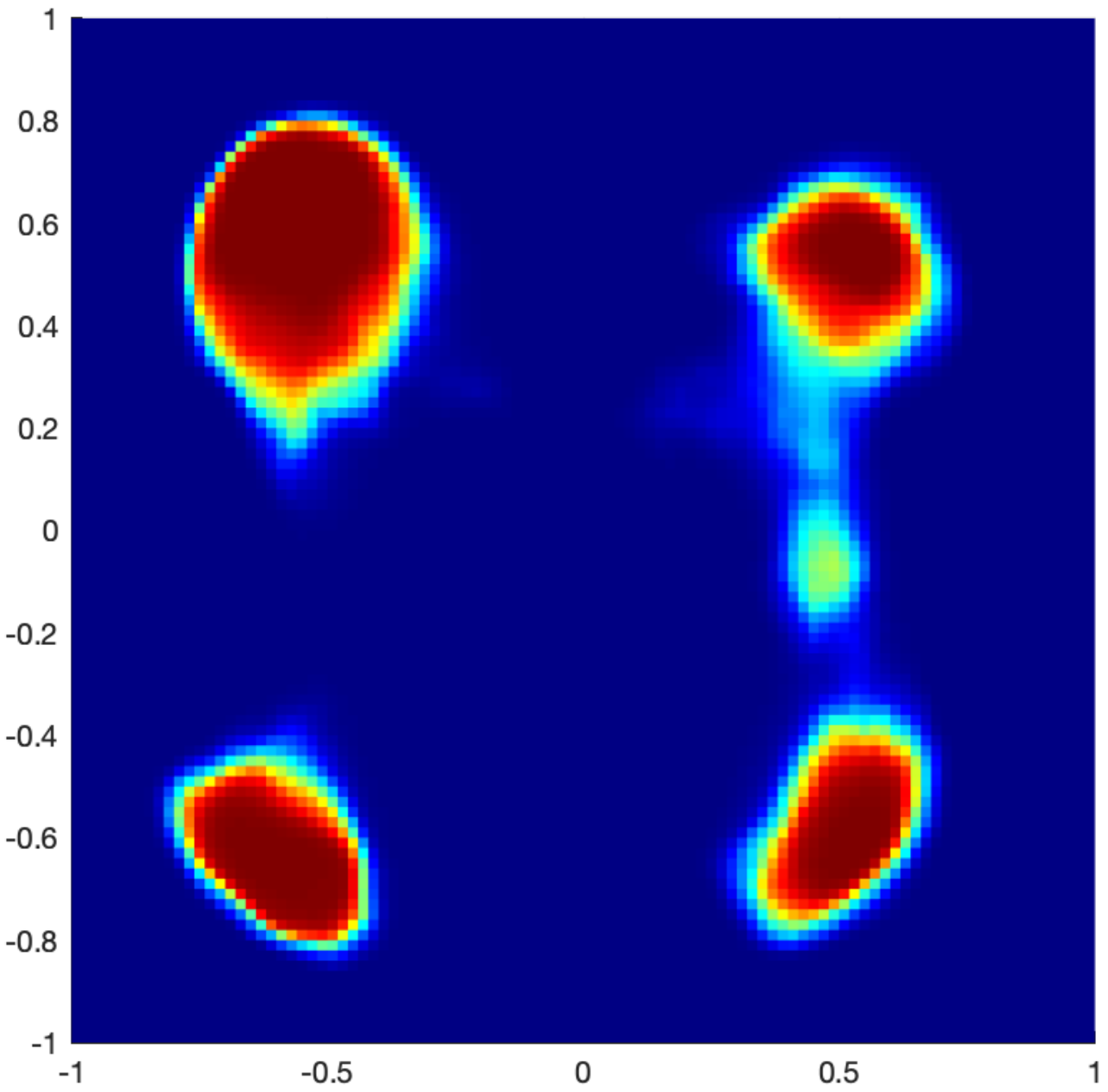}&
\includegraphics[width=1.1in]{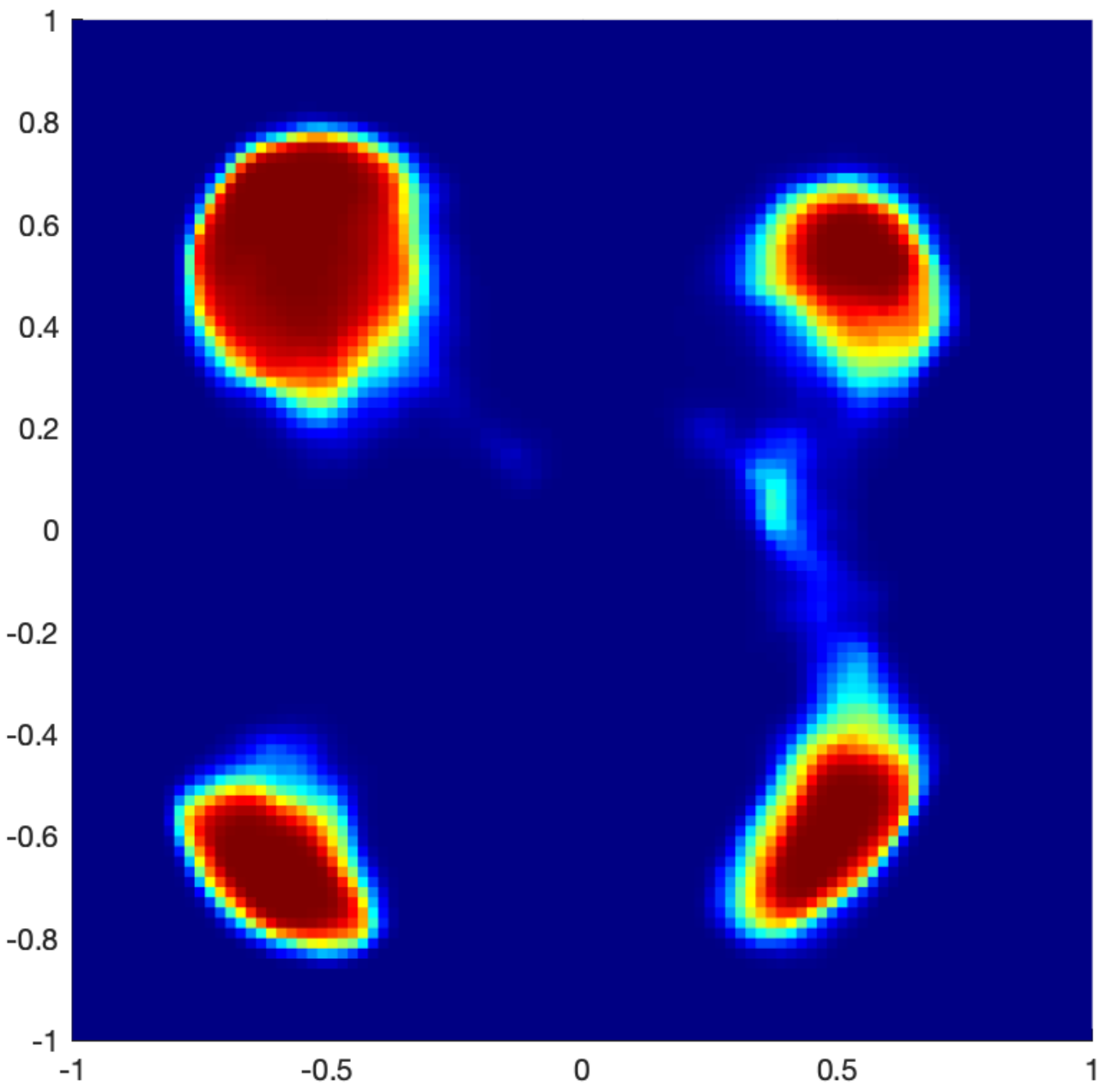}&
\includegraphics[width=1.1in]{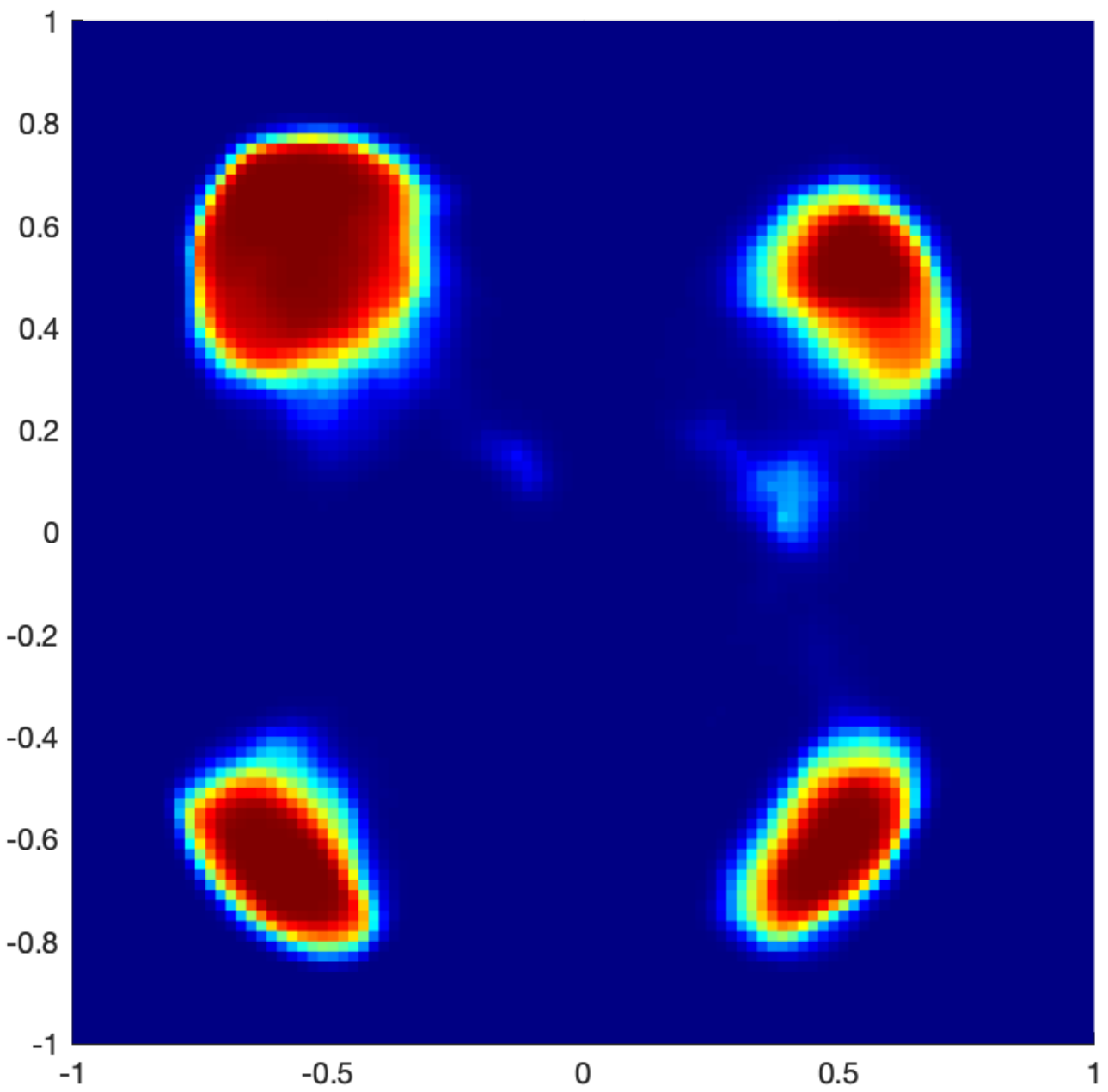}&
\includegraphics[width=1.1in]{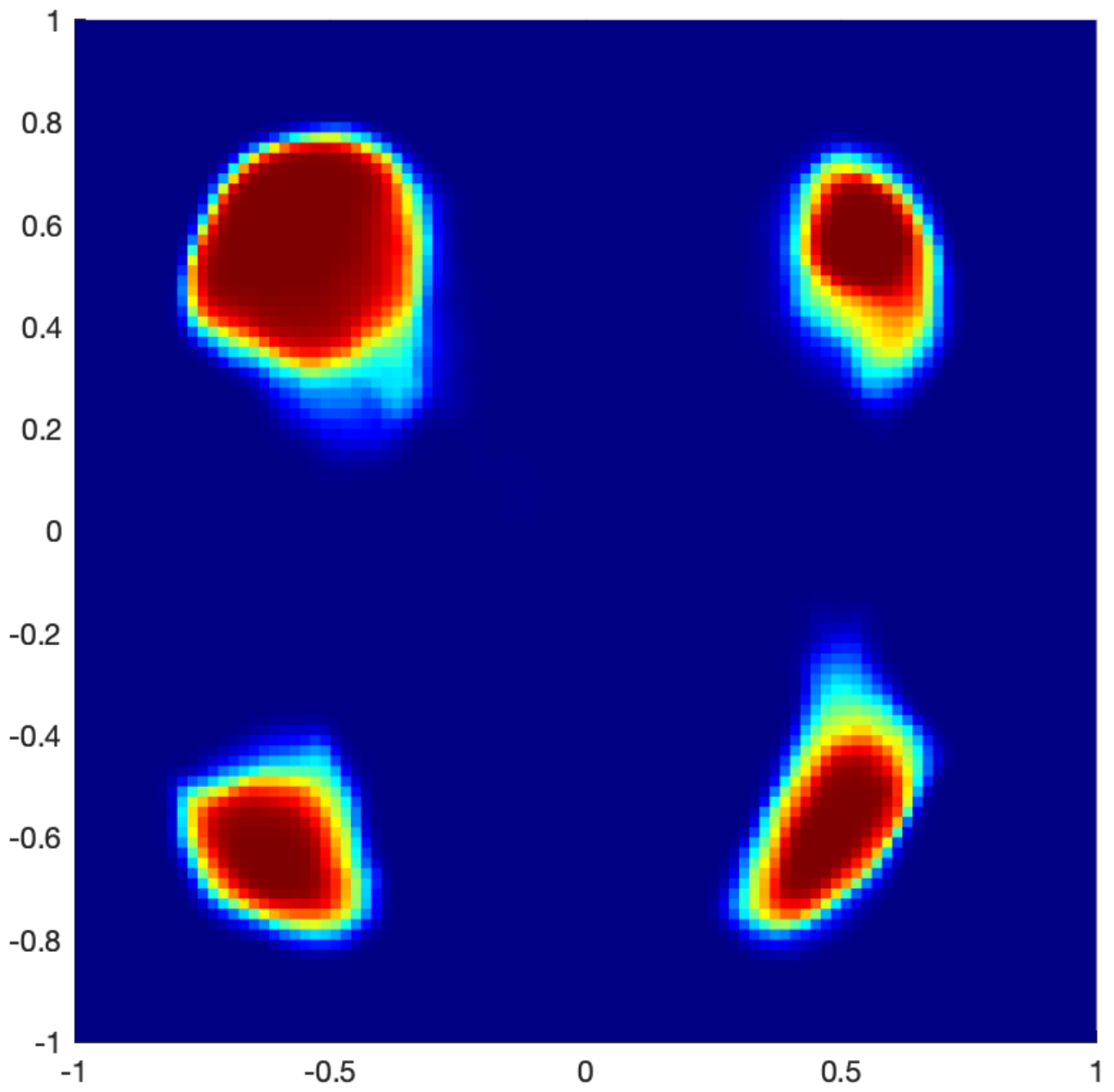}\\
\includegraphics[width=1.1in]{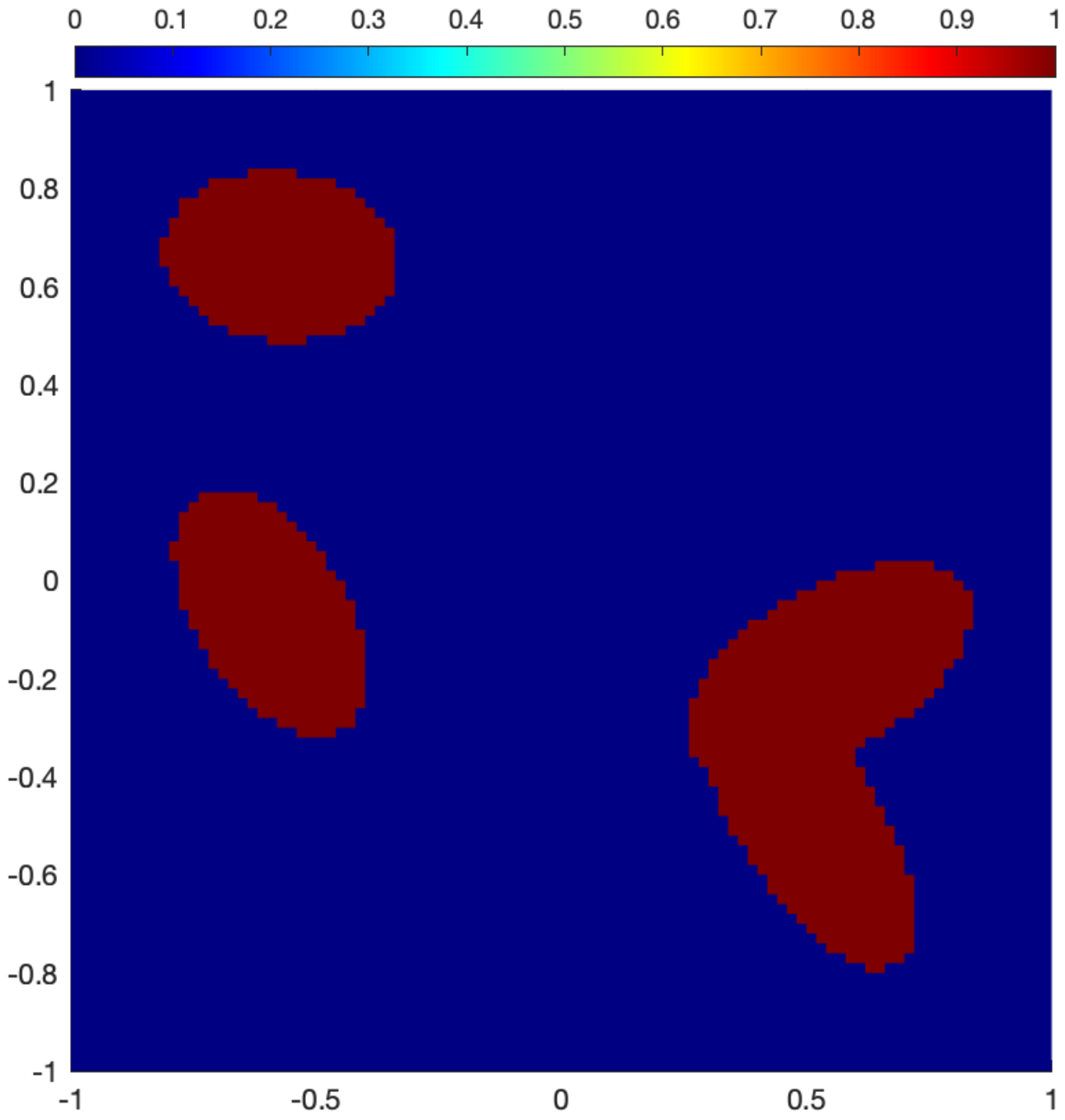}&
\includegraphics[width=1.1in]{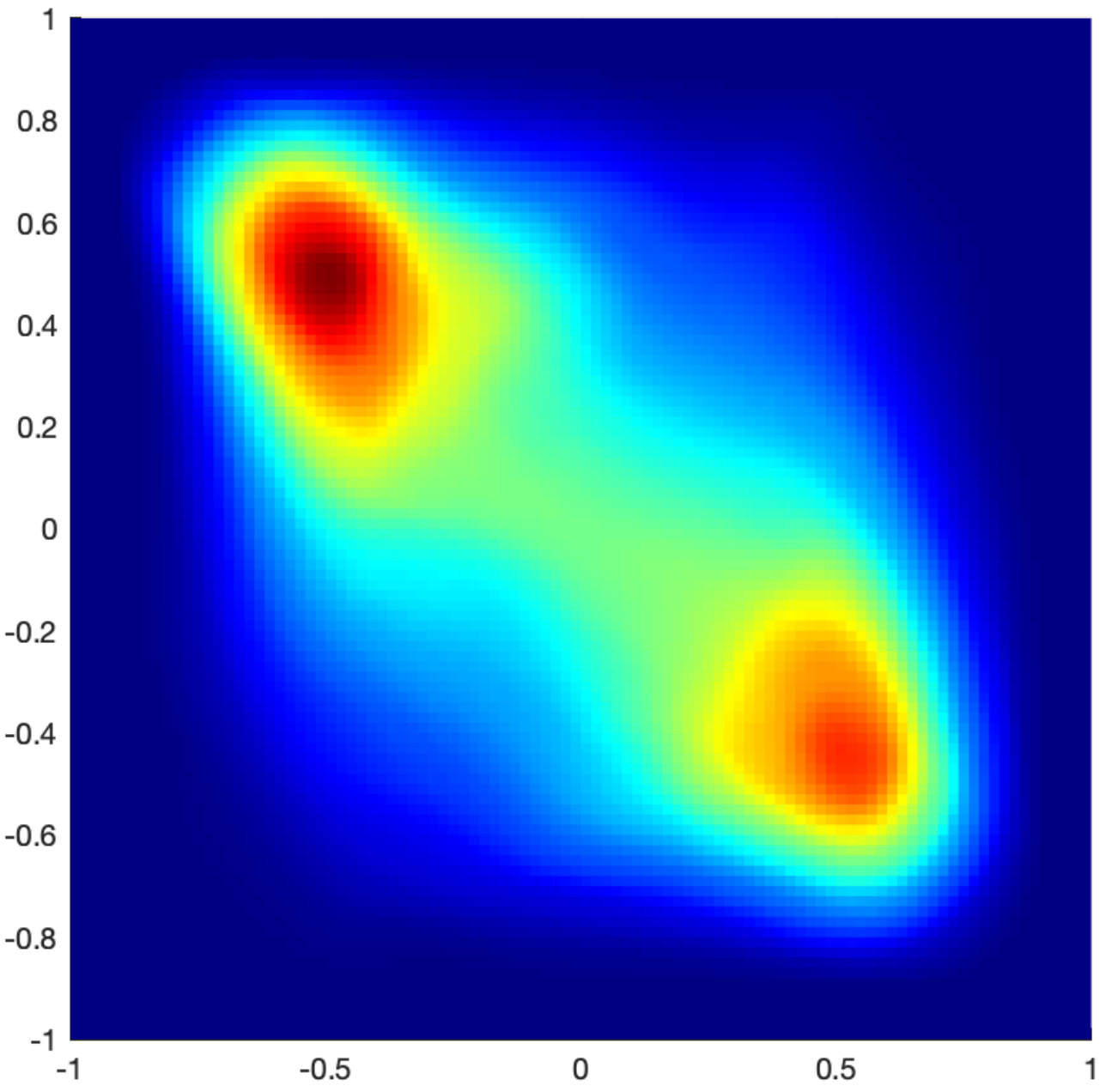}&
\includegraphics[width=1.1in]{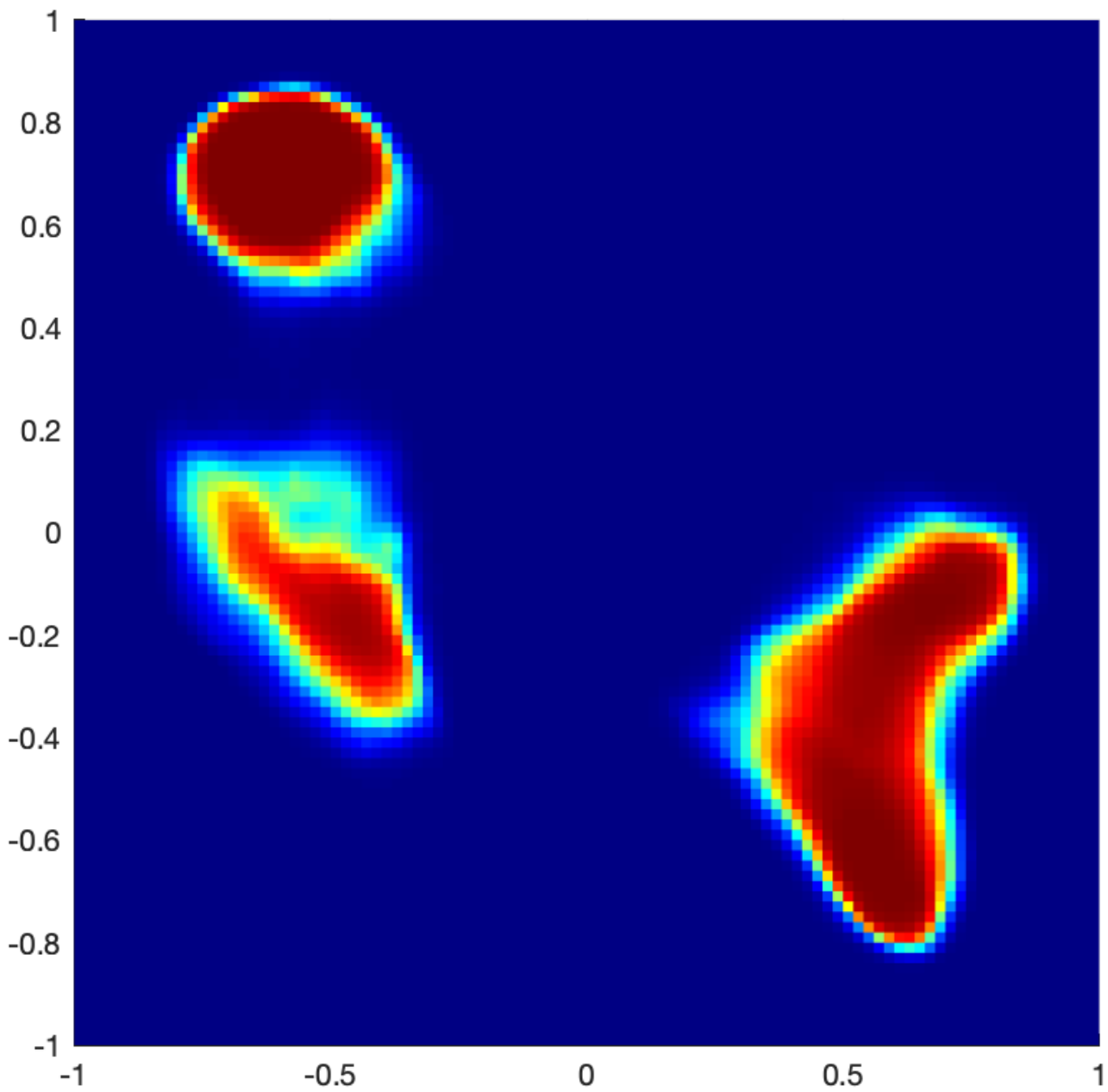}&
\includegraphics[width=1.1in]{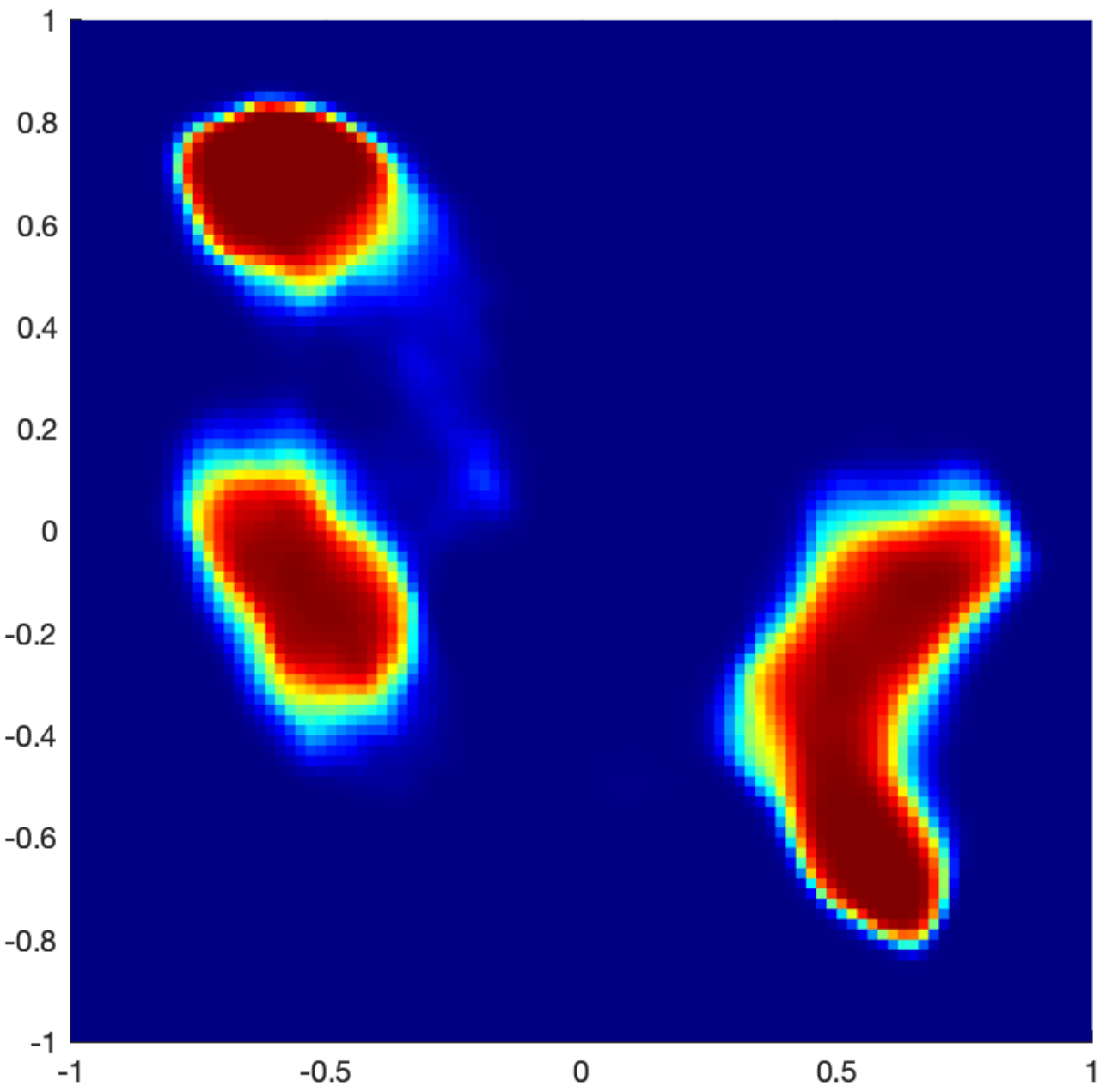}&
\includegraphics[width=1.1in]{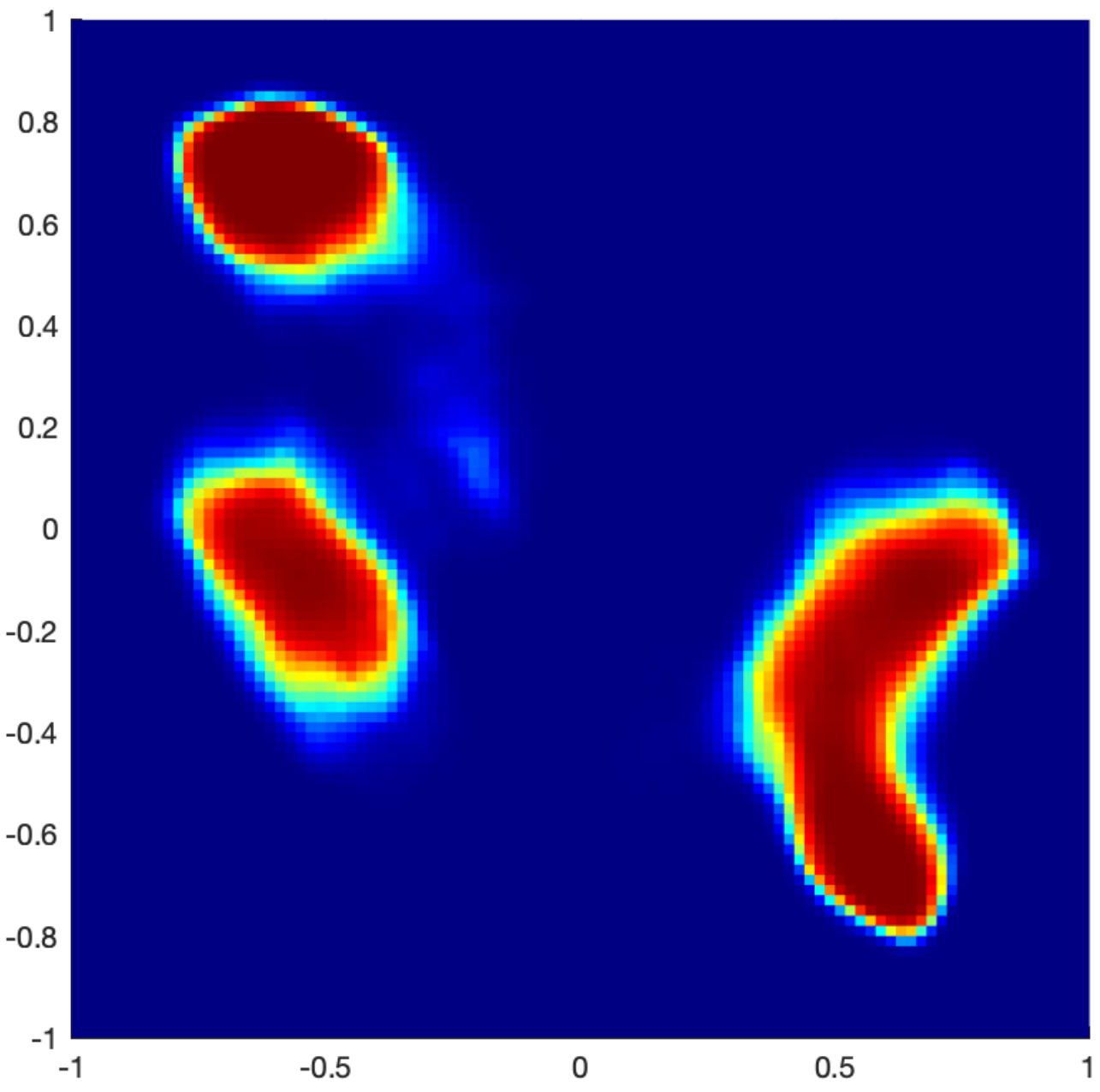}&
\includegraphics[width=1.1in]{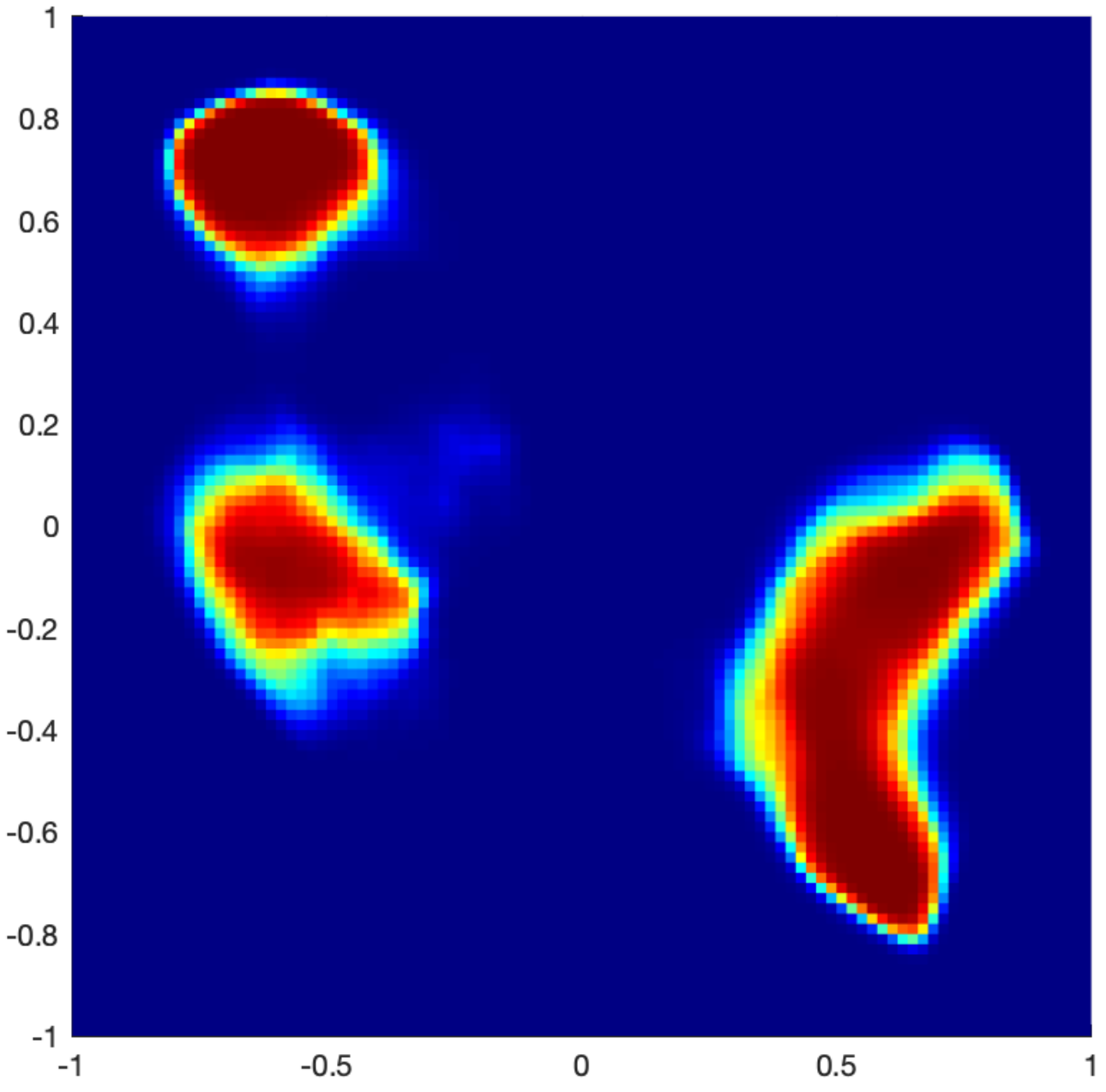}\\
\includegraphics[width=1.1in]{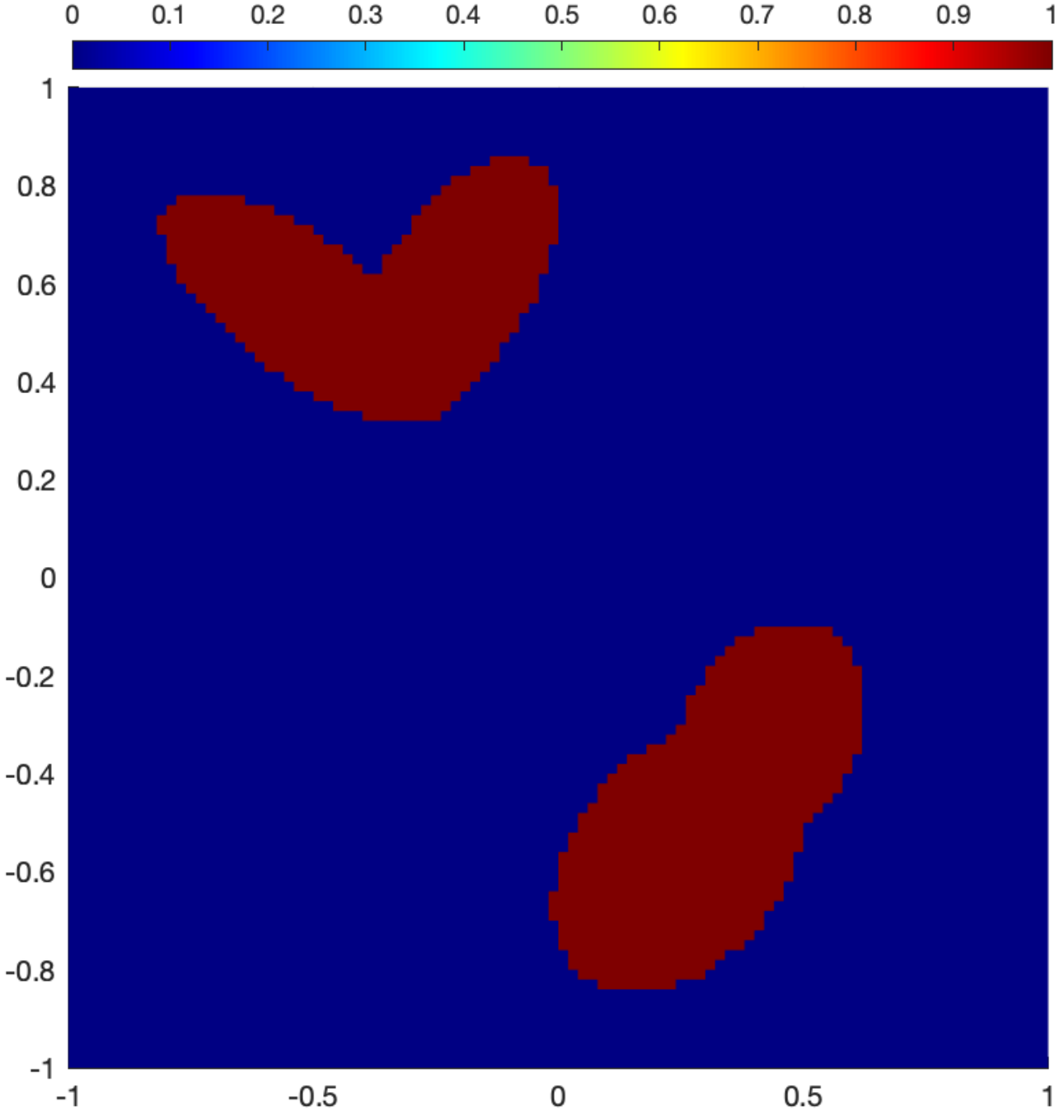}&
\includegraphics[width=1.1in]{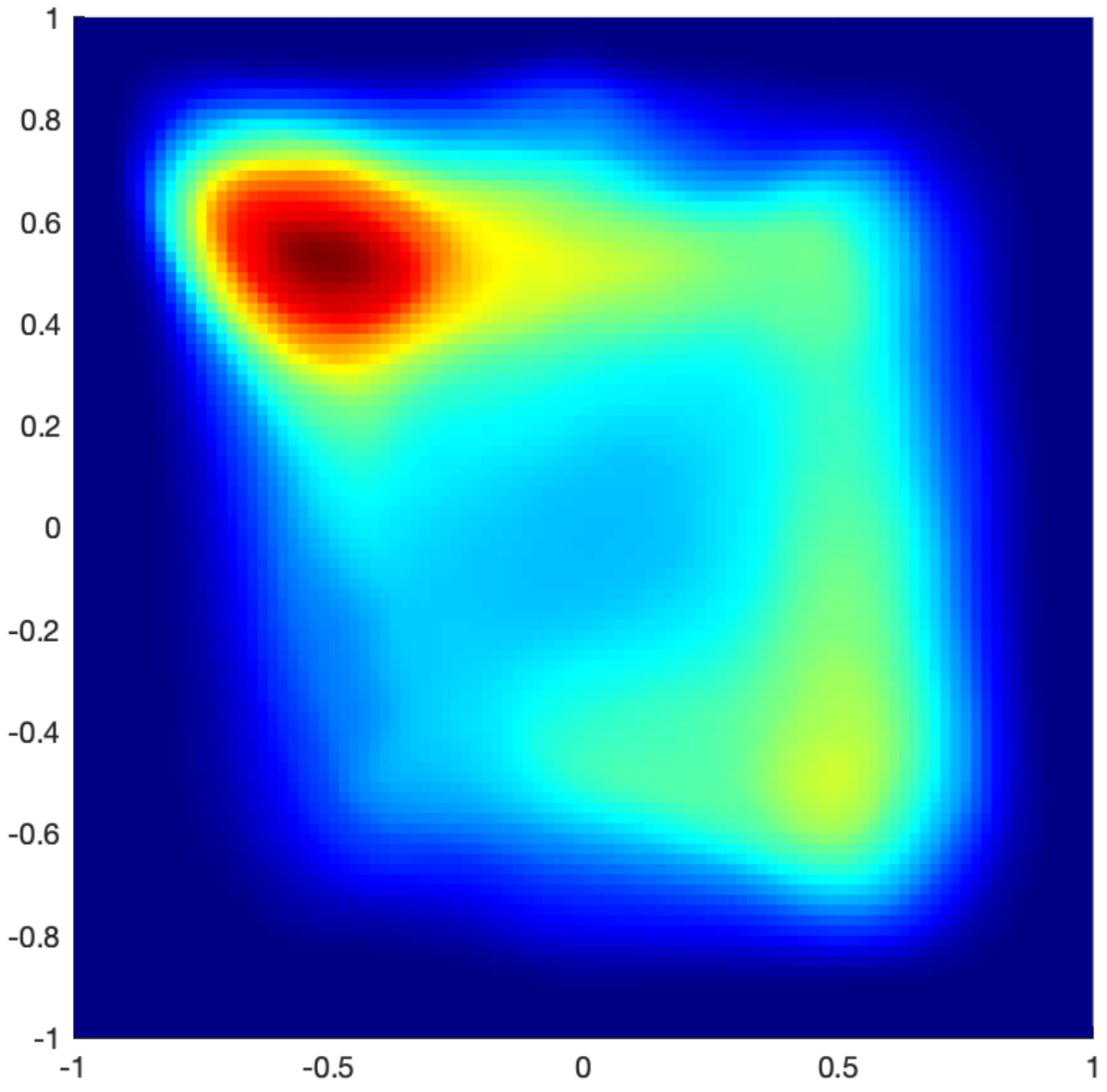}&
\includegraphics[width=1.1in]{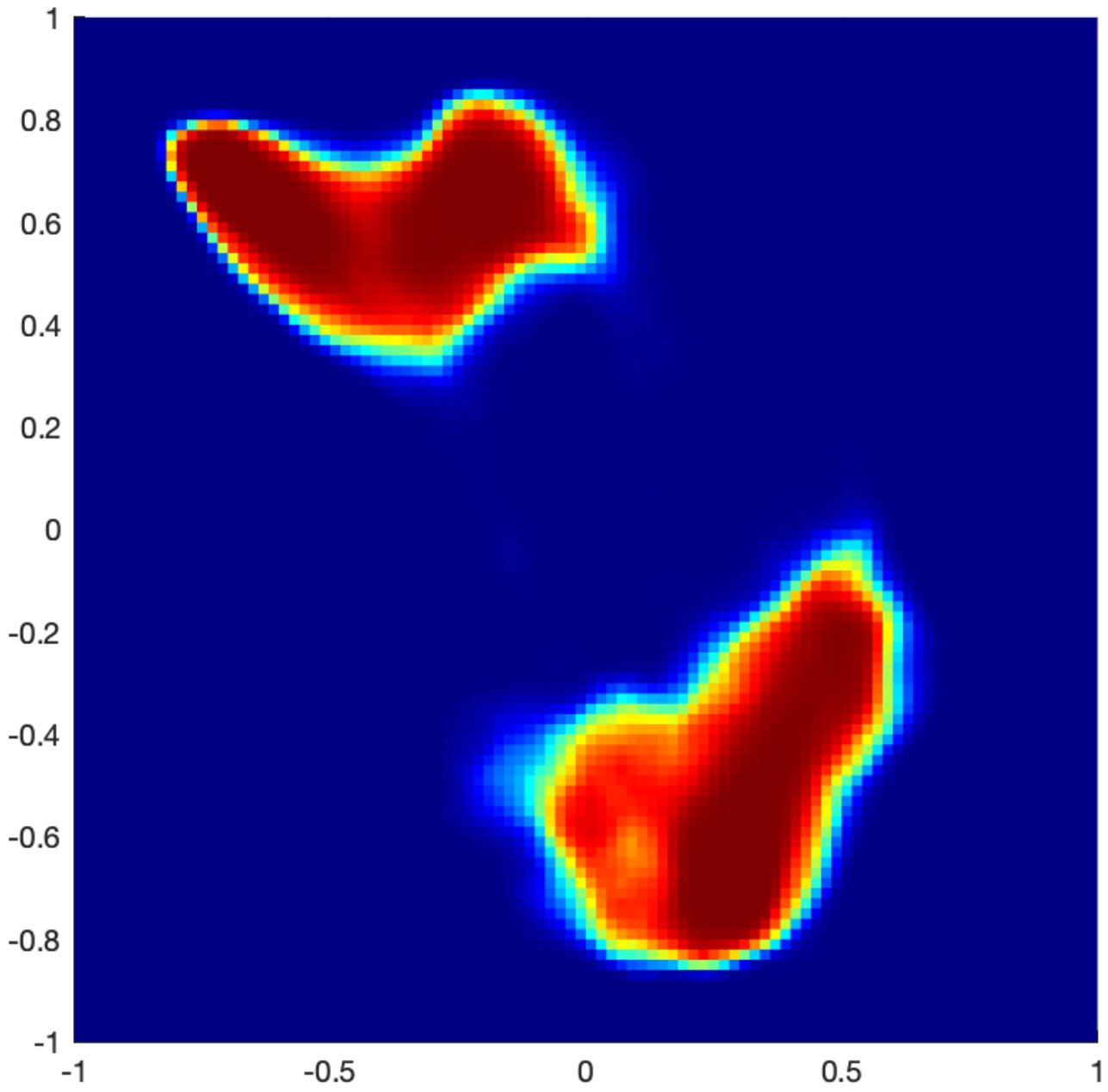}&
\includegraphics[width=1.1in]{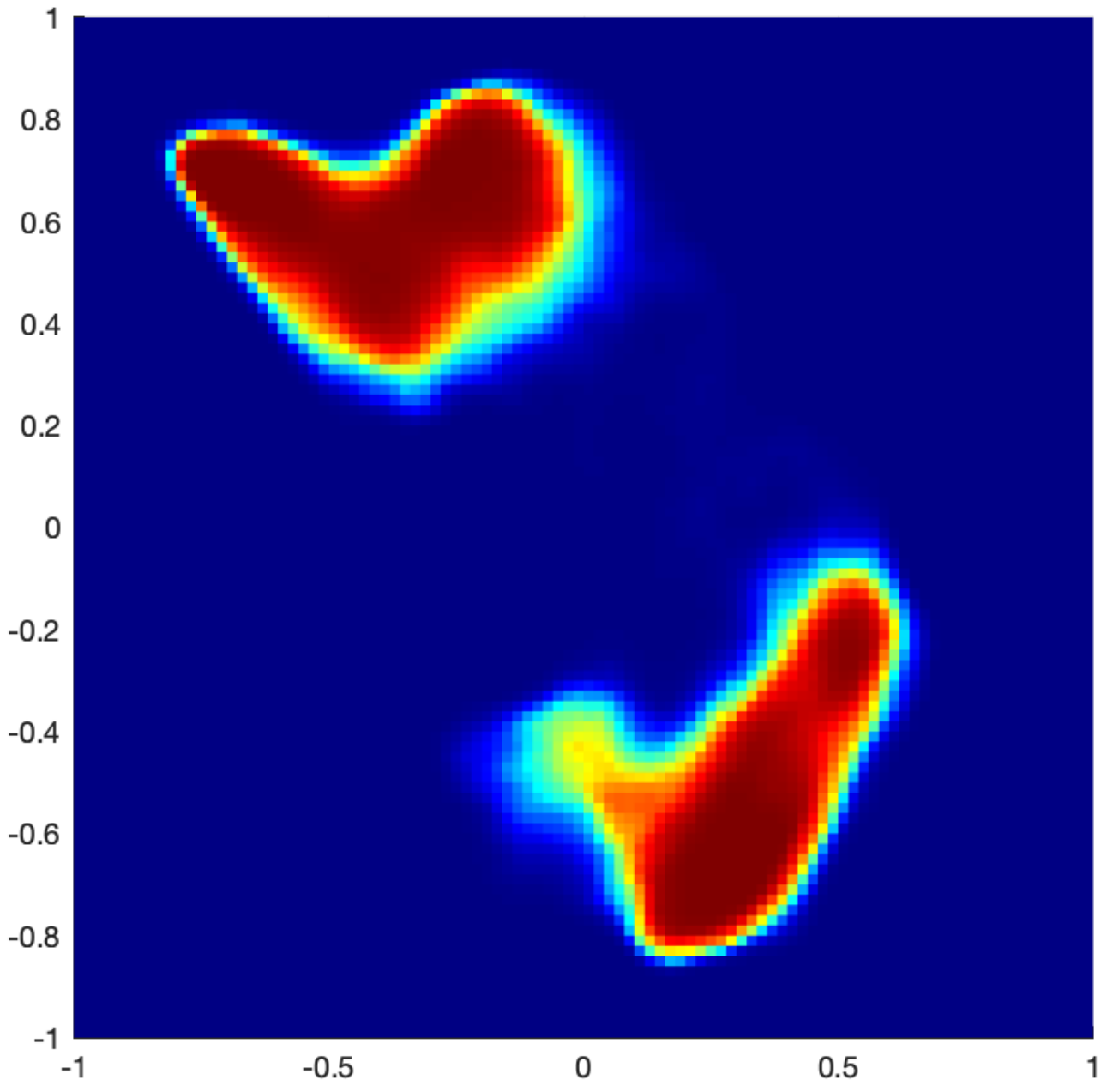}&
\includegraphics[width=1.1in]{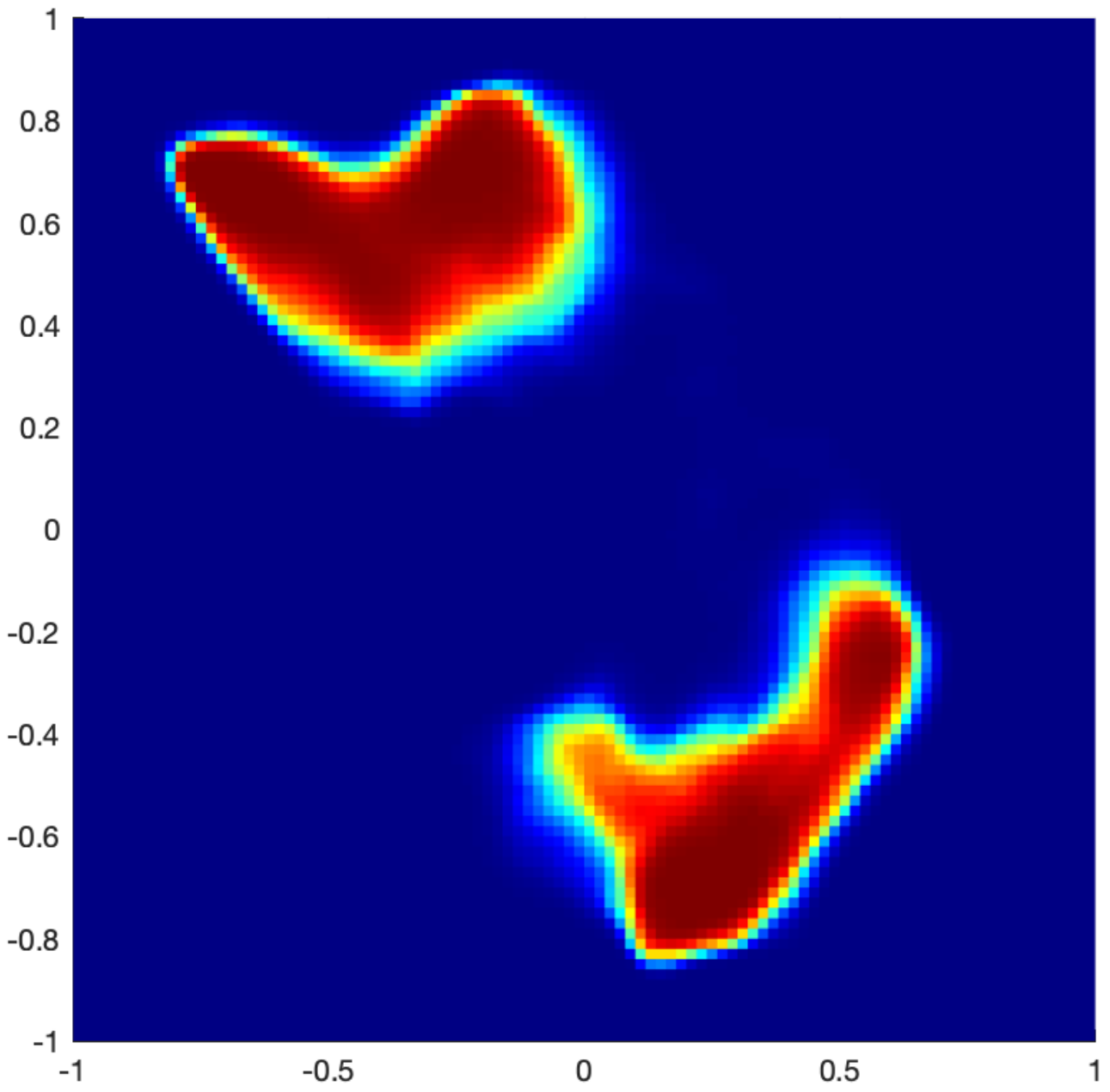}&
\includegraphics[width=1.1in]{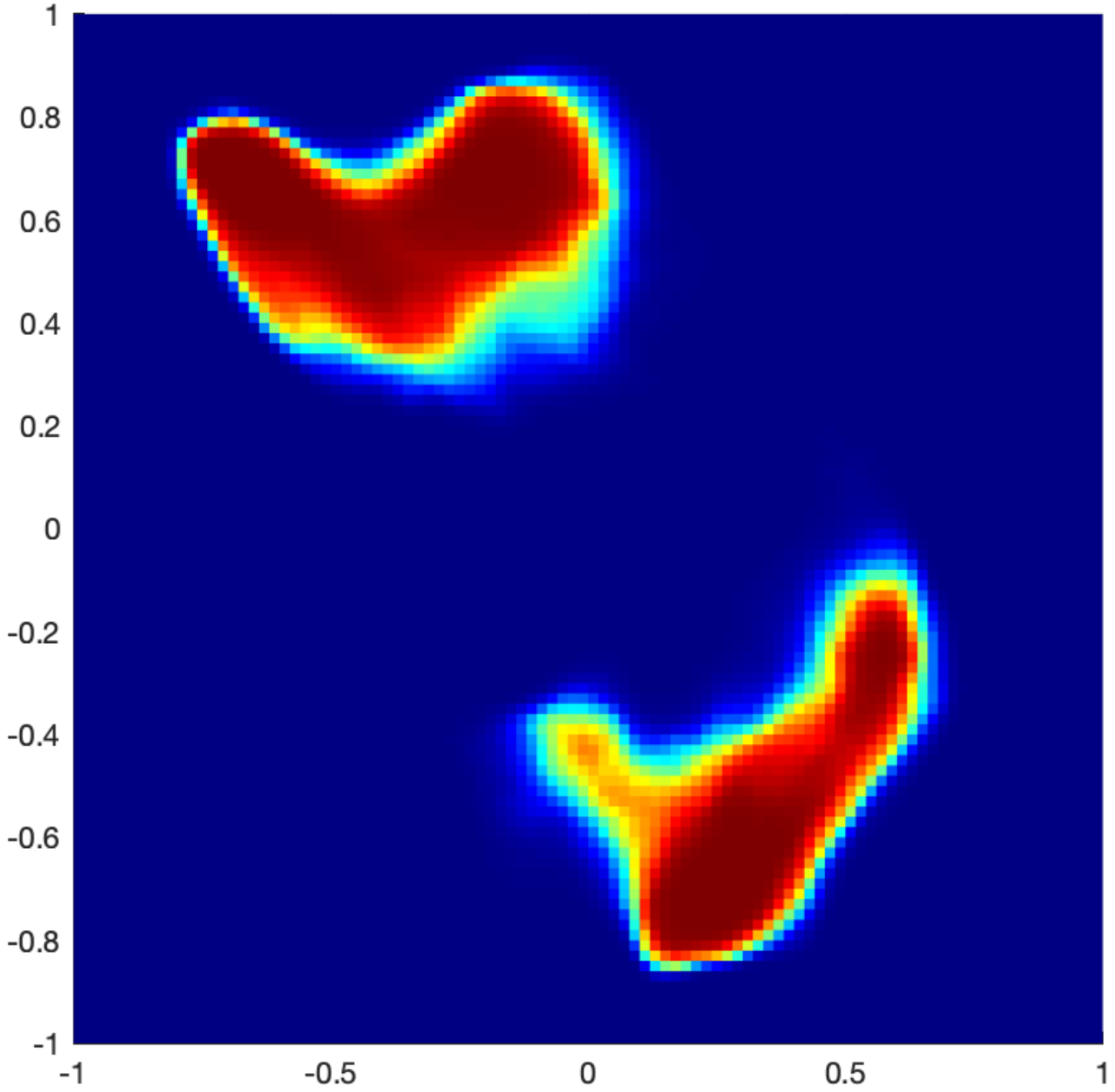}\\
\end{tabular}
  \caption{FNN-DDSM reconstruction for 3 cases in \textbf{Scenario 3} (4 ellipses) with different Cauchy data number and noise level: Case 1(top), Case 2(middle) and Case 3(bottom) } 
  \label{tab_FN_ell}
\end{figure}

\begin{figure}[htbp]
\begin{tabular}{ >{\centering\arraybackslash}m{0.9in} >{\centering\arraybackslash}m{0.9in} >{\centering\arraybackslash}m{0.9in}  >{\centering\arraybackslash}m{0.9in}  >{\centering\arraybackslash}m{0.9in}  >{\centering\arraybackslash}m{0.9in} }
\centering
True coefficients &
N=1, $\delta=0$&
N=10, $\delta=0$&
N=20, $\delta=0$&
N=20, $\delta=10\%$ &
N=20, $\delta=20\%$ \\
\includegraphics[width=1.1in]{ell_case1_0-eps-converted-to.pdf}&
\includegraphics[width=1.1in]{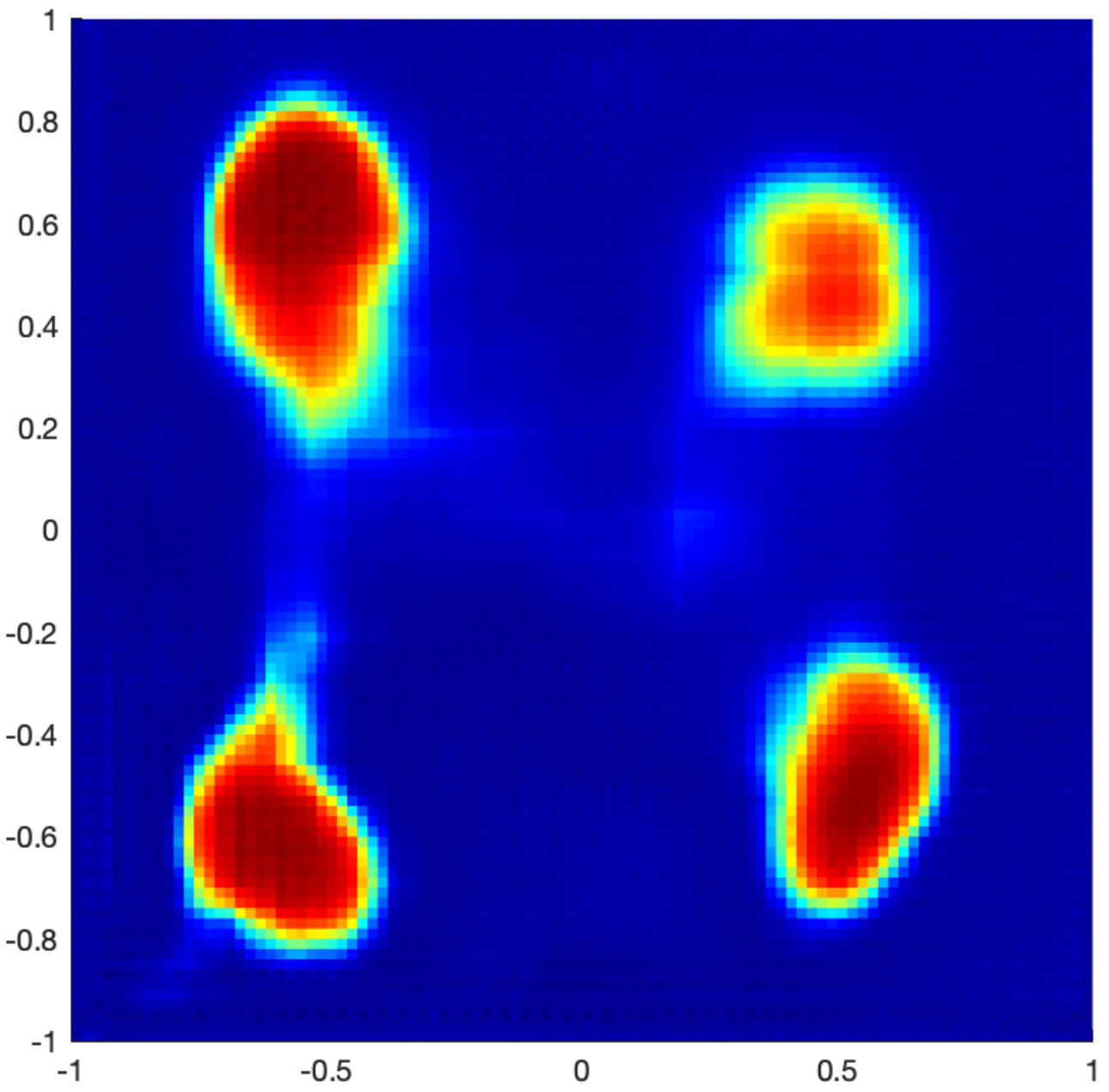}&
\includegraphics[width=1.1in]{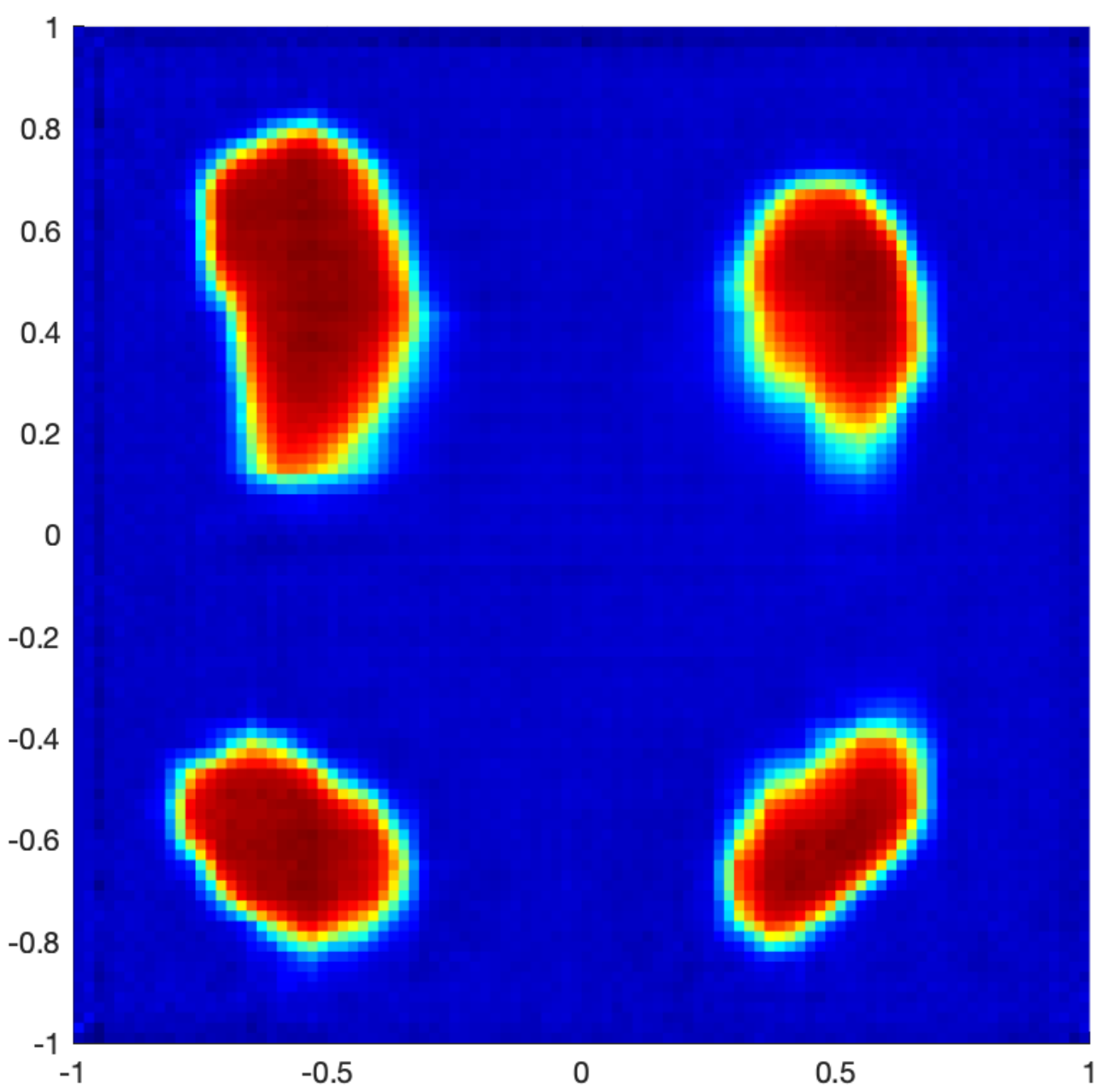}&
\includegraphics[width=1.1in]{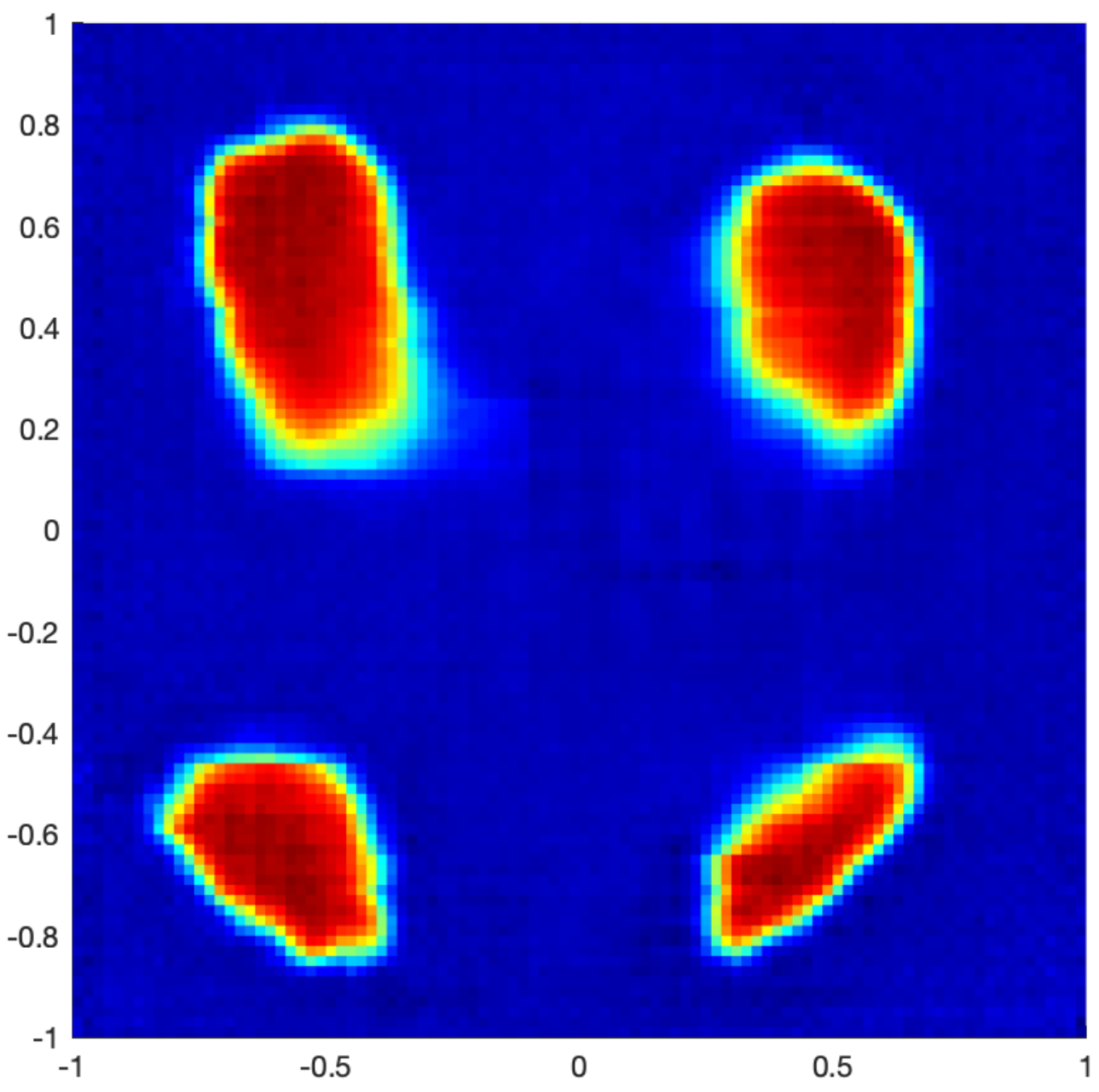}&
\includegraphics[width=1.1in]{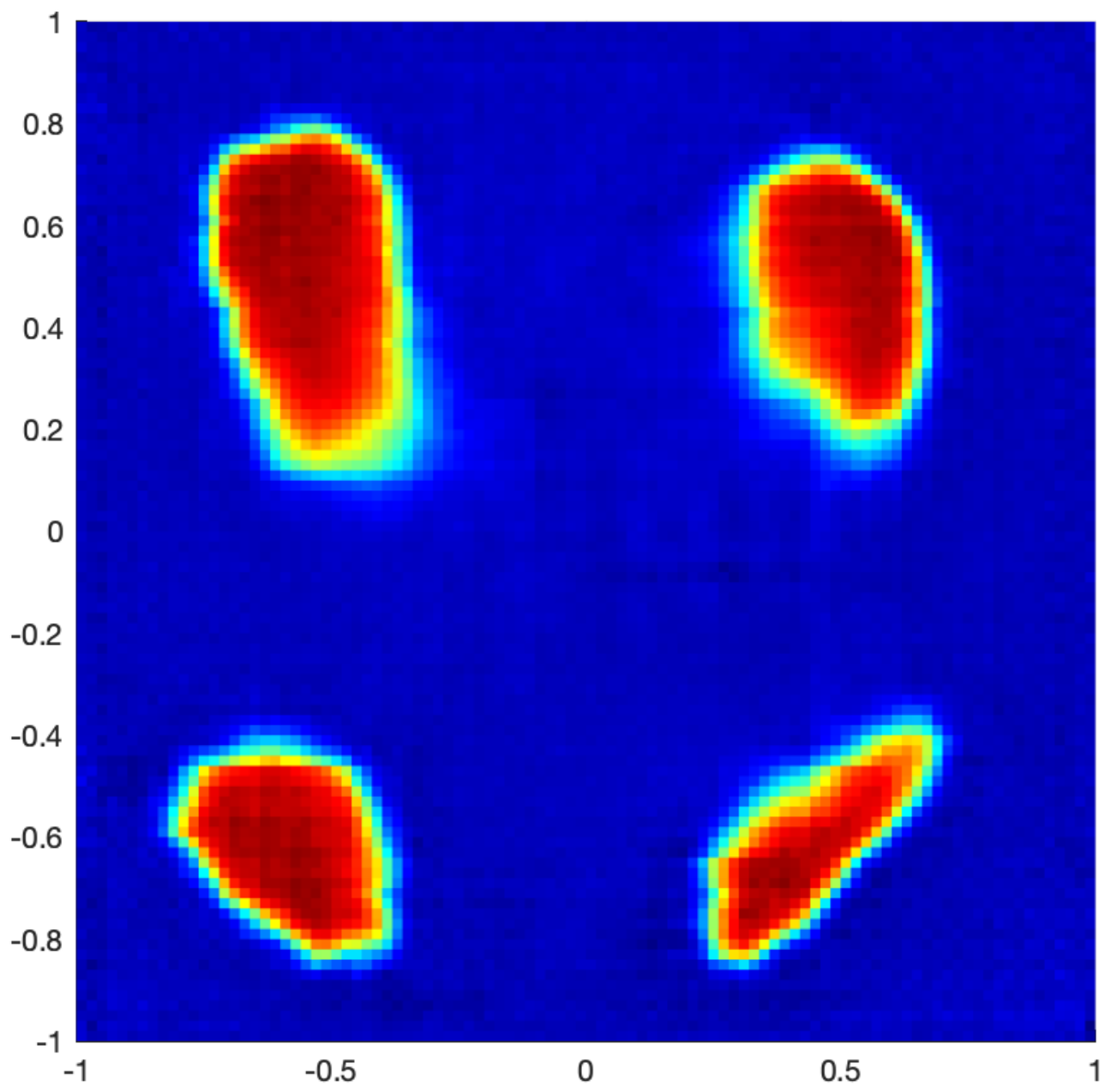}&
\includegraphics[width=1.1in]{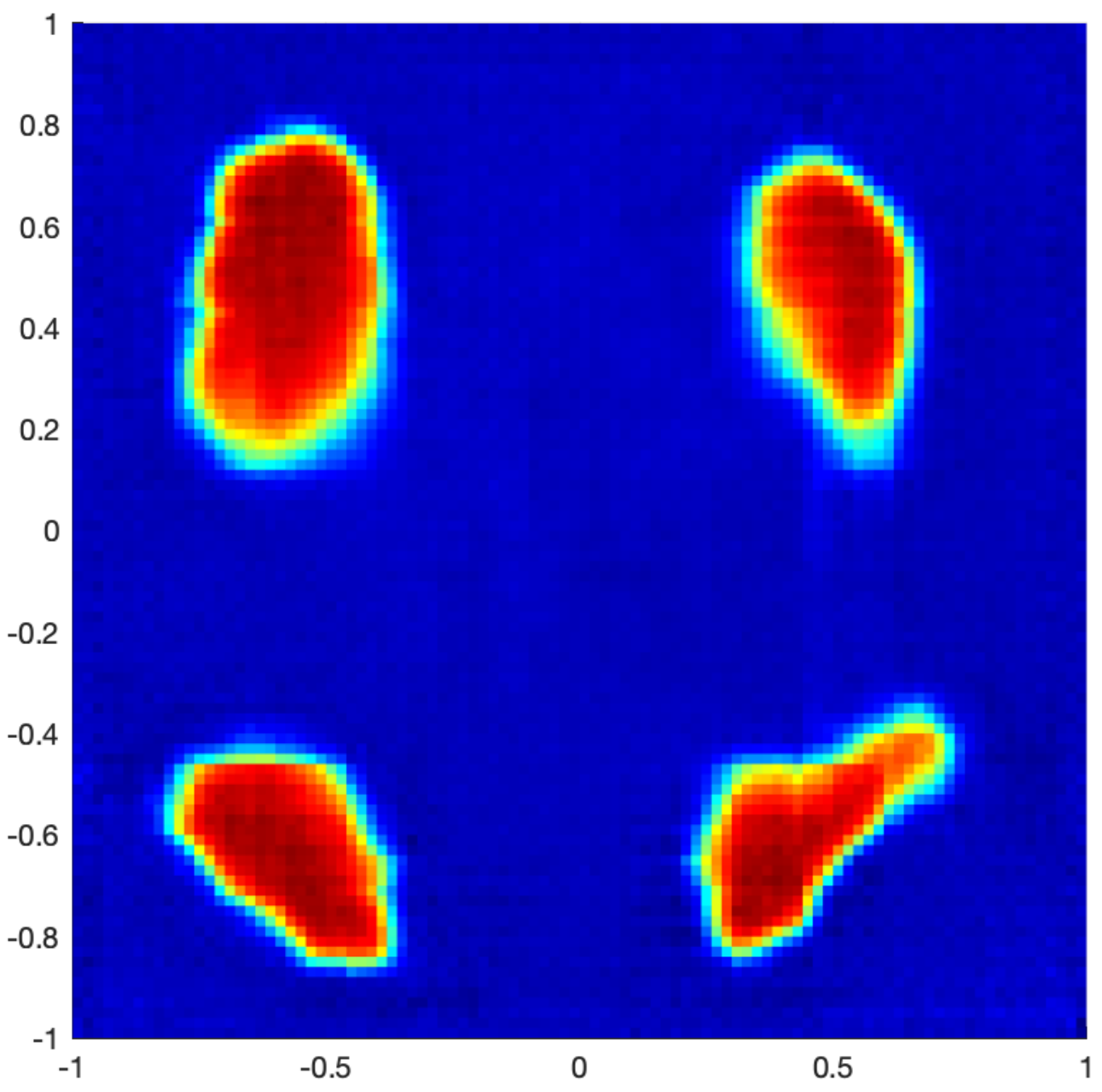}\\
\includegraphics[width=1.1in]{ell_case2_0-eps-converted-to.pdf}&
\includegraphics[width=1.1in]{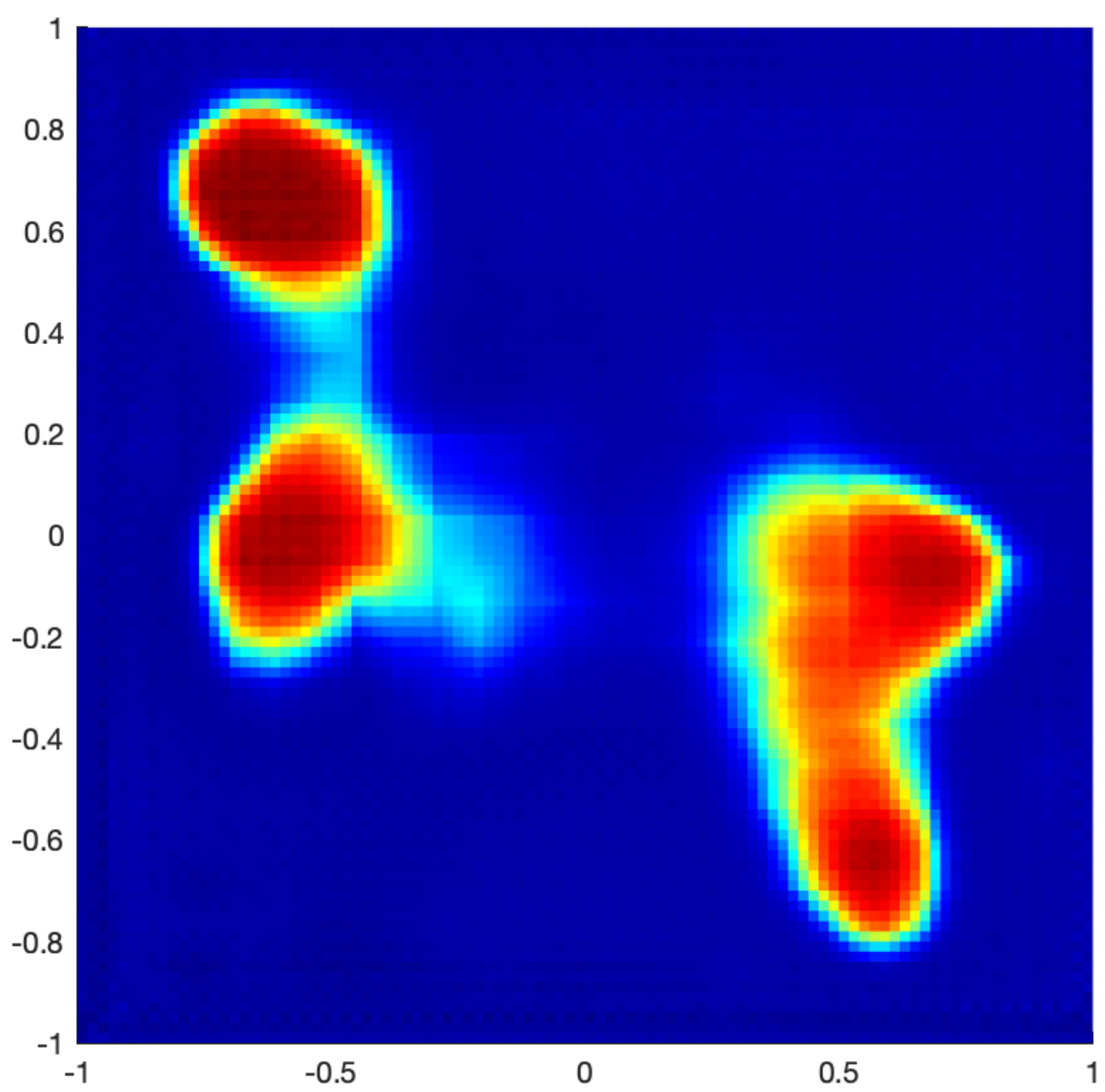}&
\includegraphics[width=1.1in]{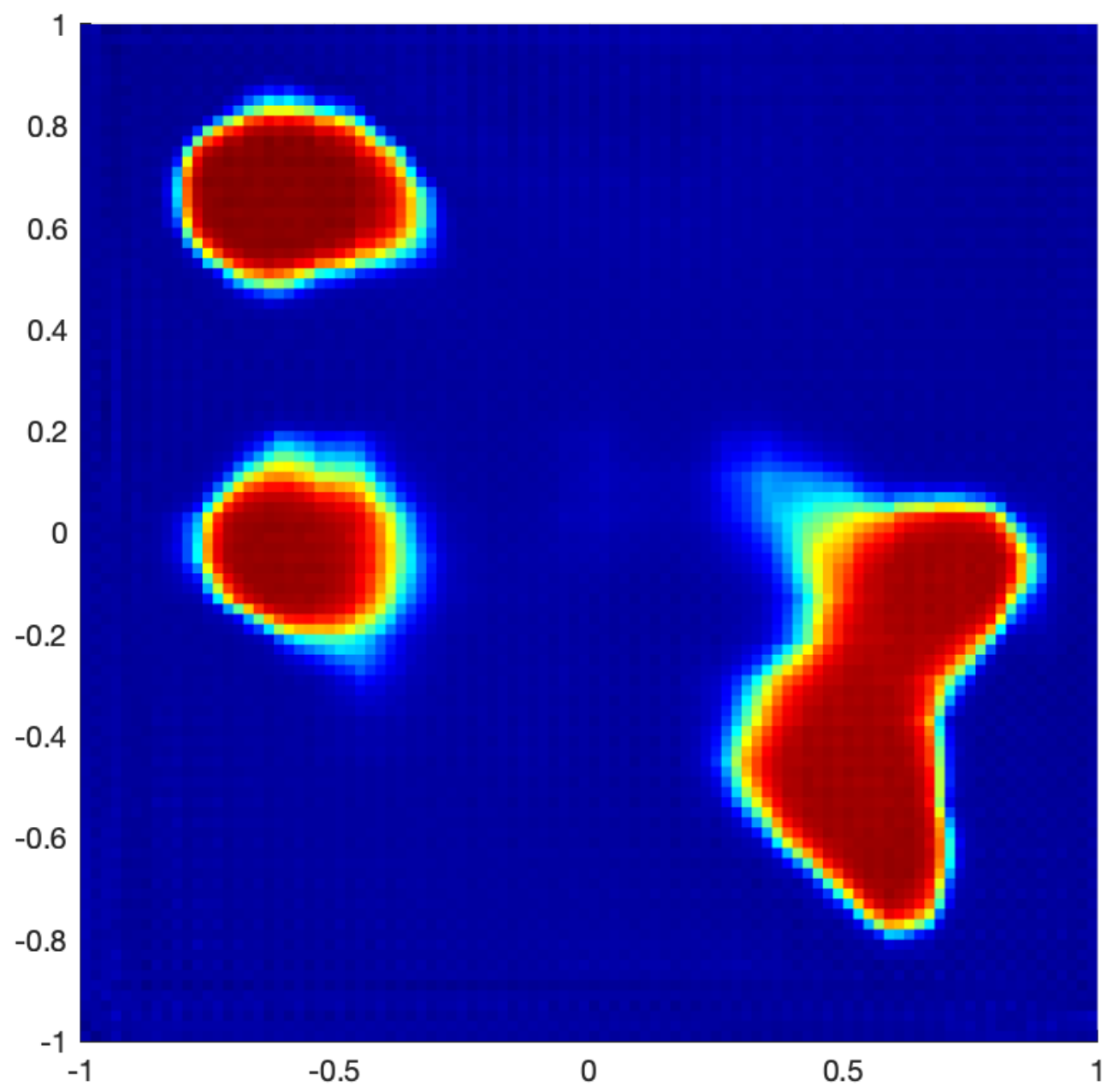}&
\includegraphics[width=1.1in]{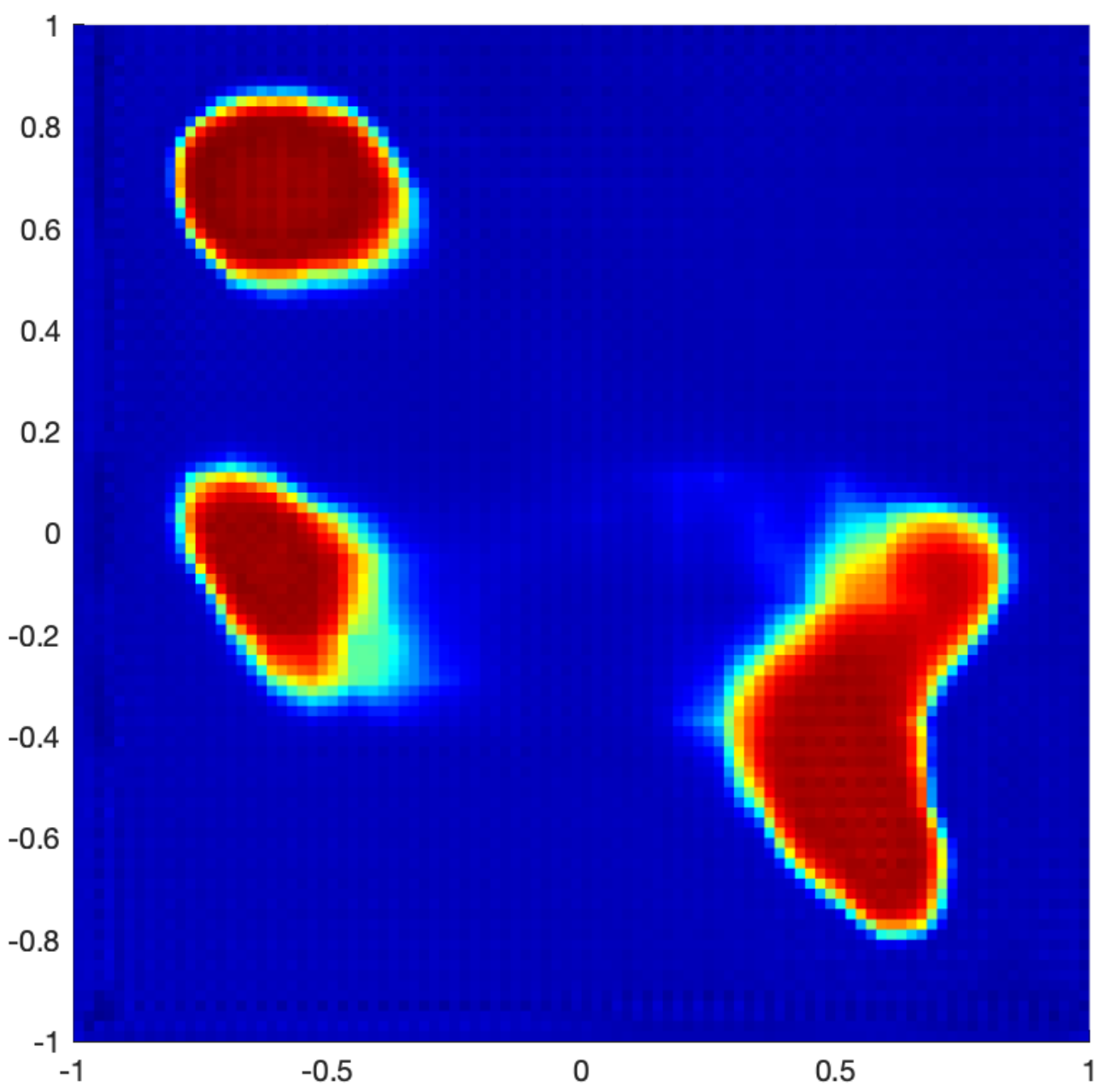}&
\includegraphics[width=1.1in]{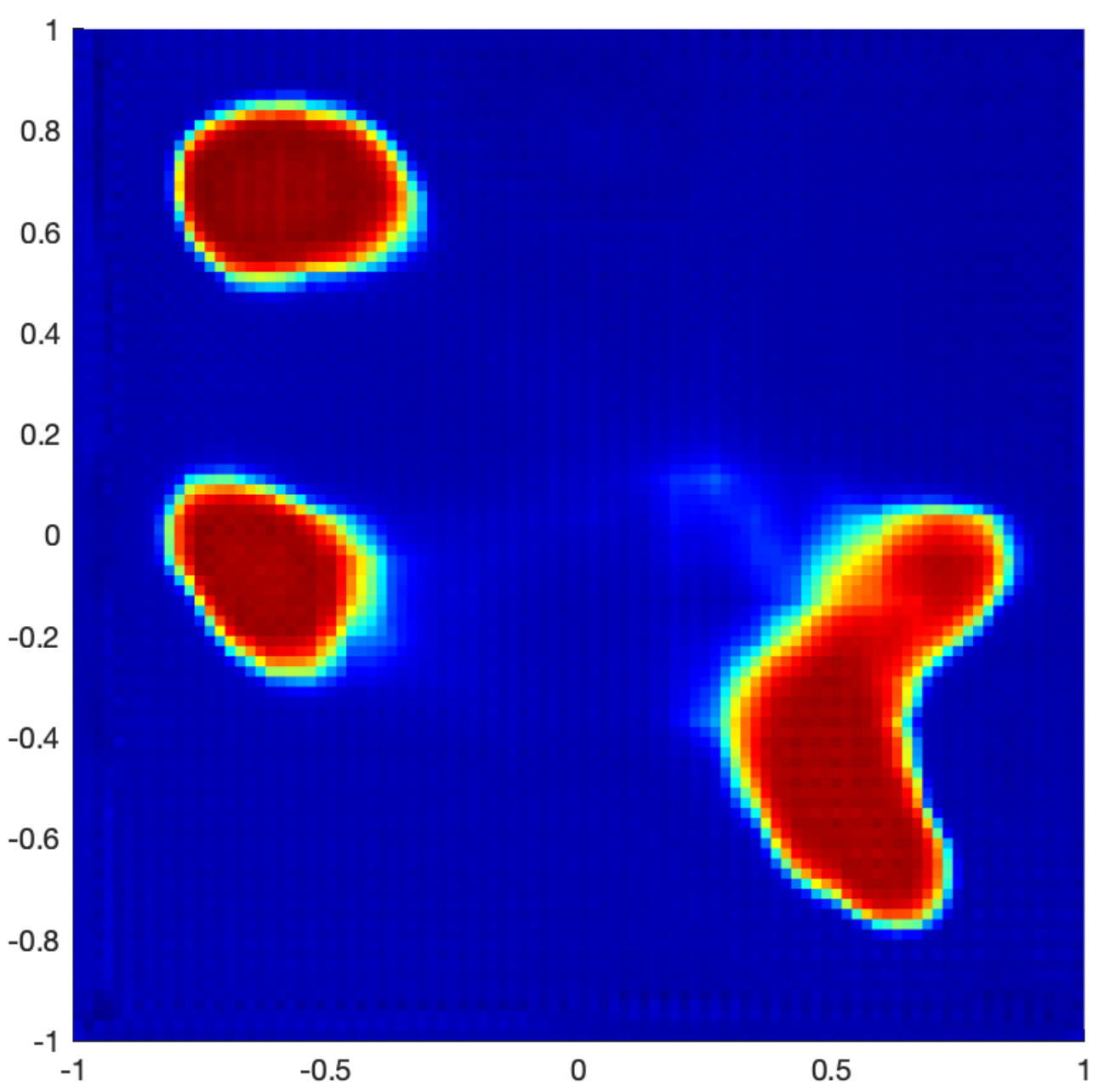}&
\includegraphics[width=1.1in]{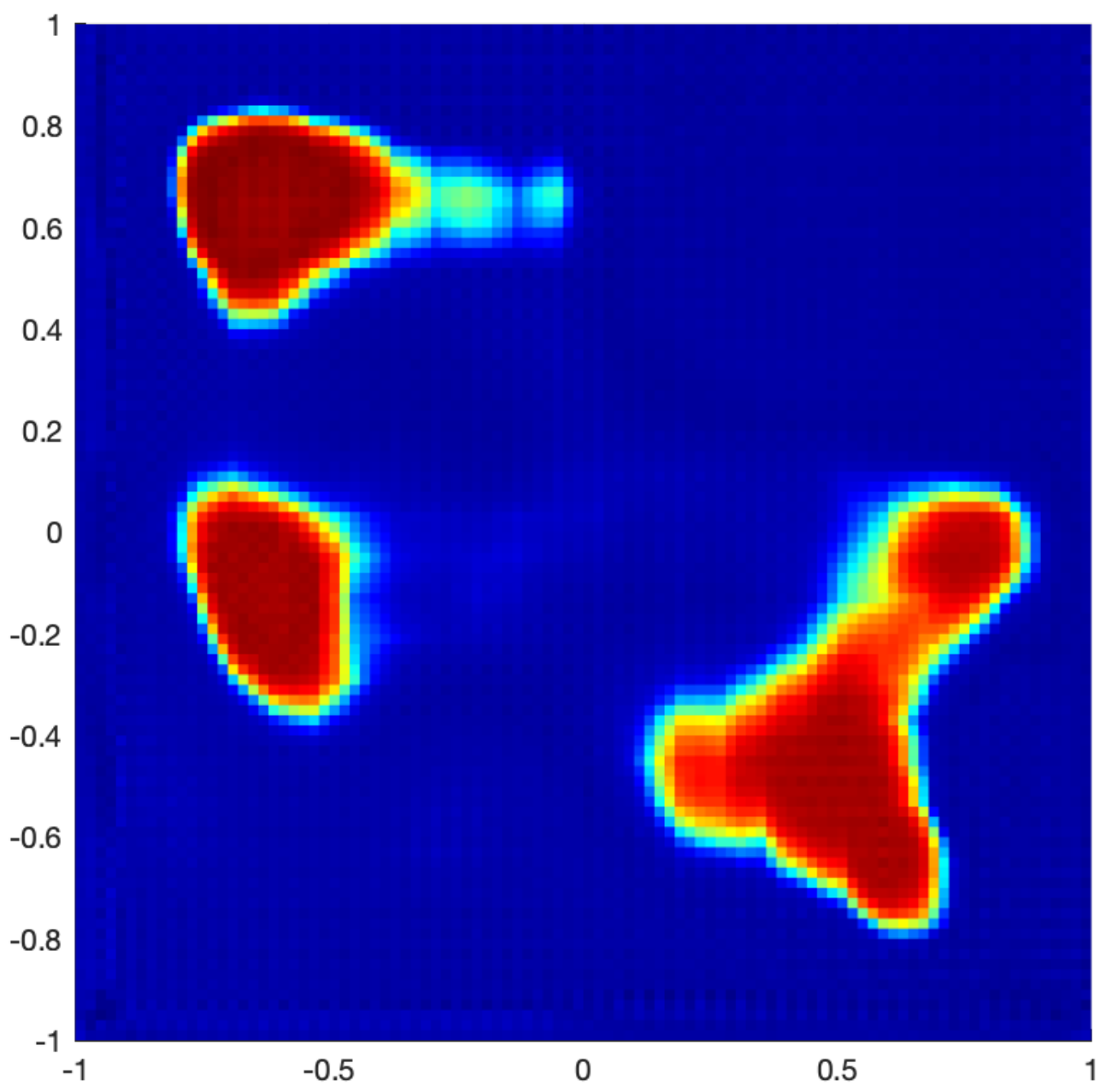}\\
\includegraphics[width=1.1in]{ell_case3_0-eps-converted-to.pdf}&
\includegraphics[width=1.1in]{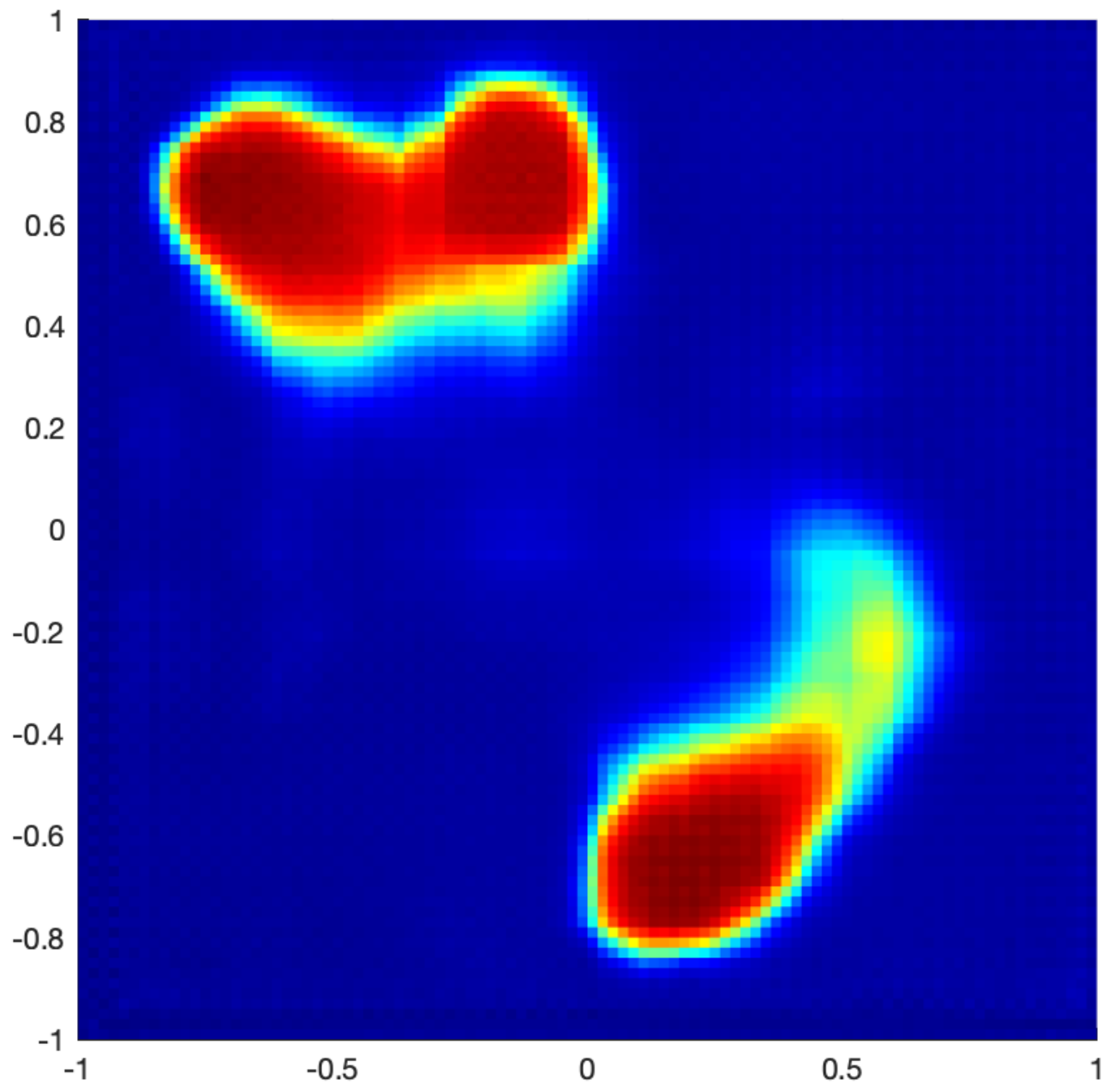}&
\includegraphics[width=1.1in]{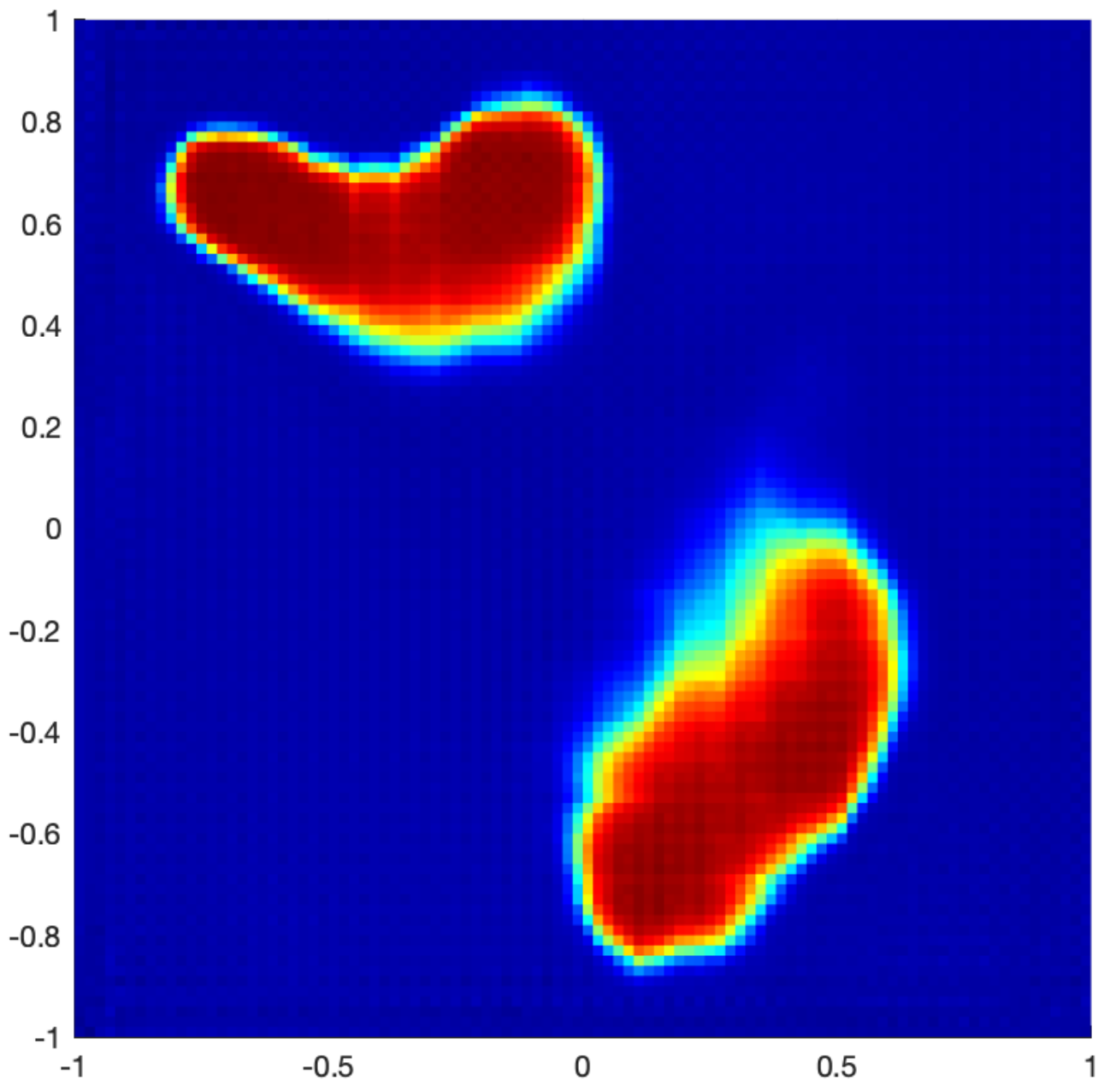}&
\includegraphics[width=1.1in]{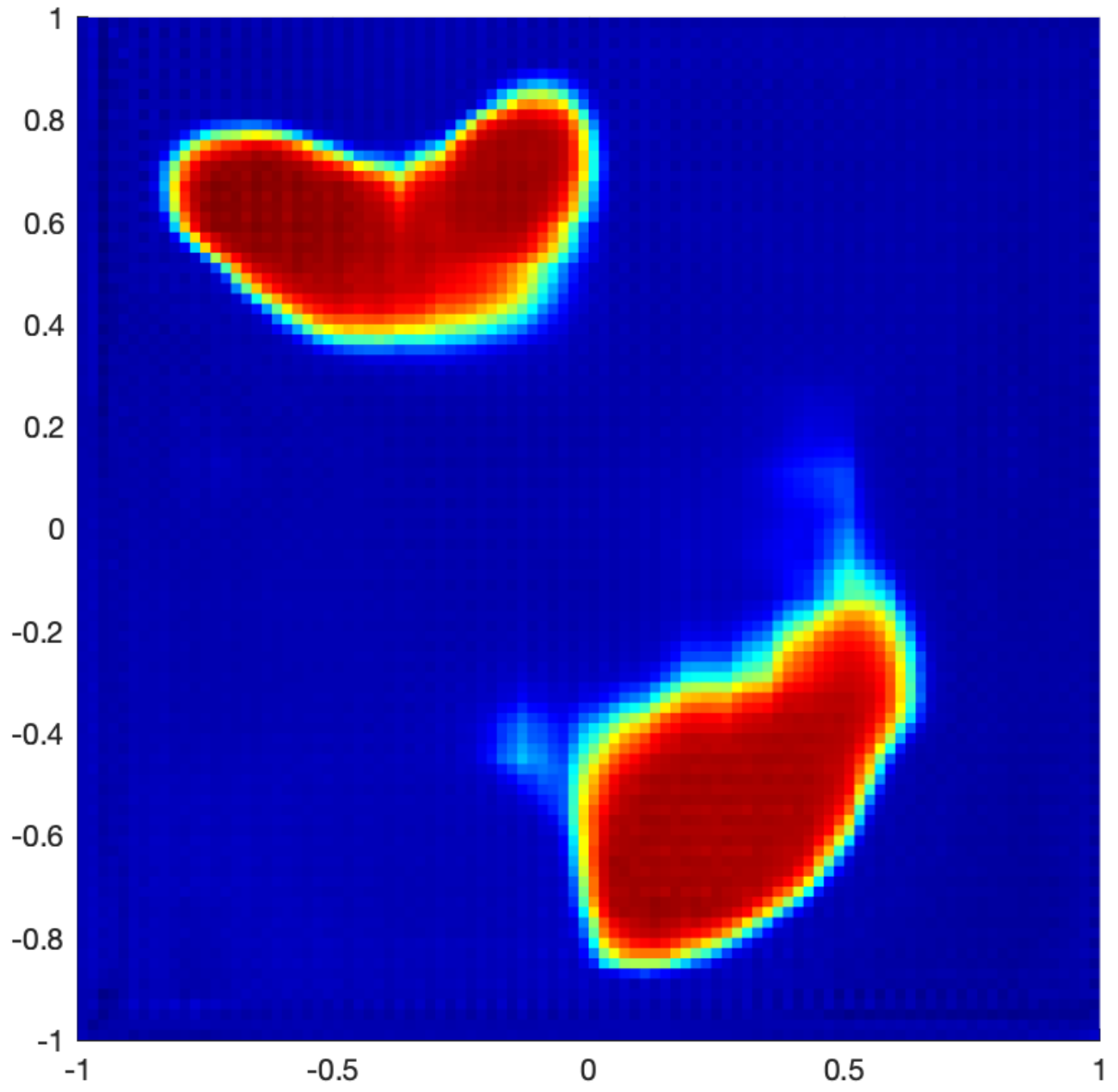}&
\includegraphics[width=1.1in]{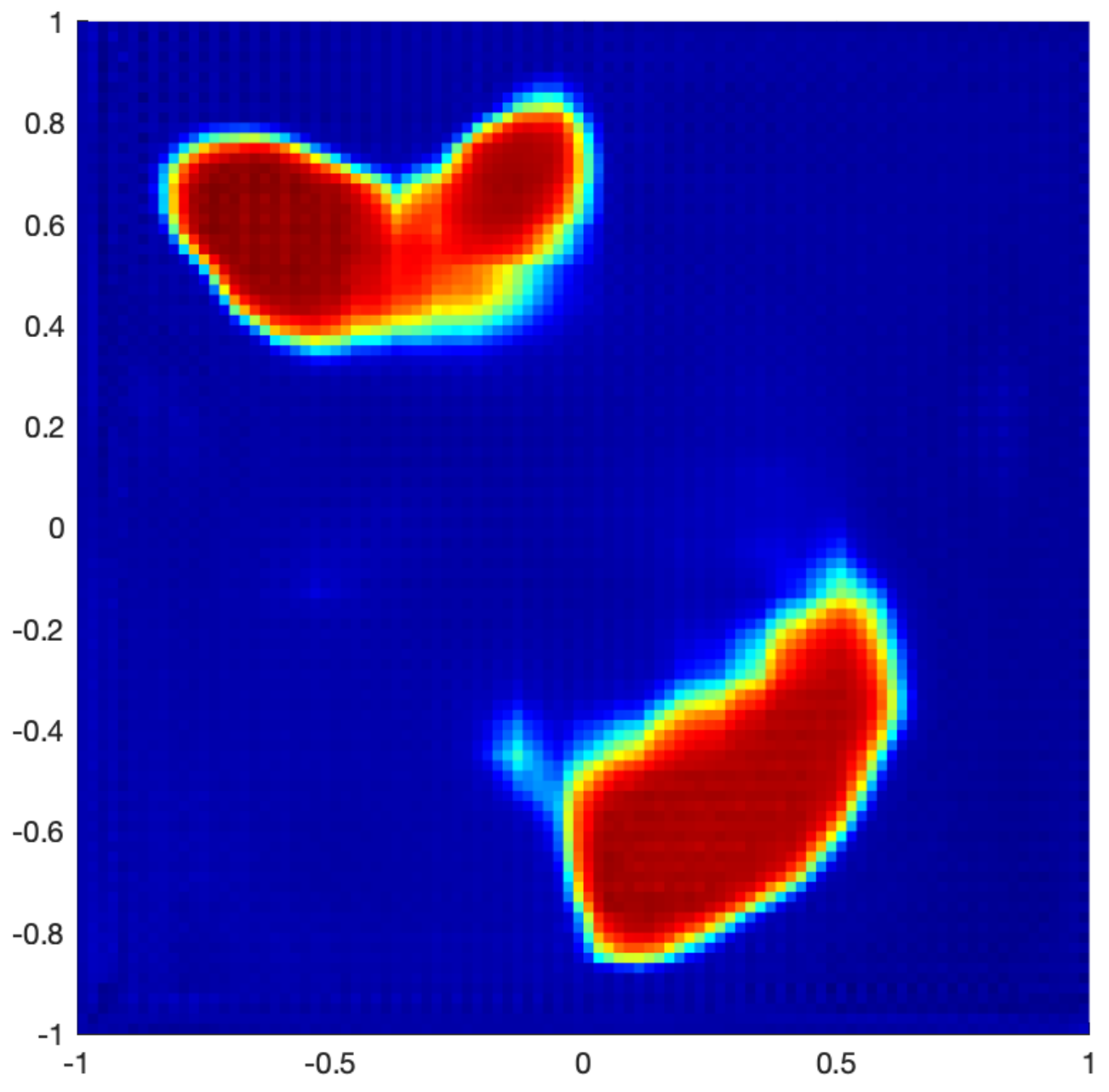}&
\includegraphics[width=1.1in]{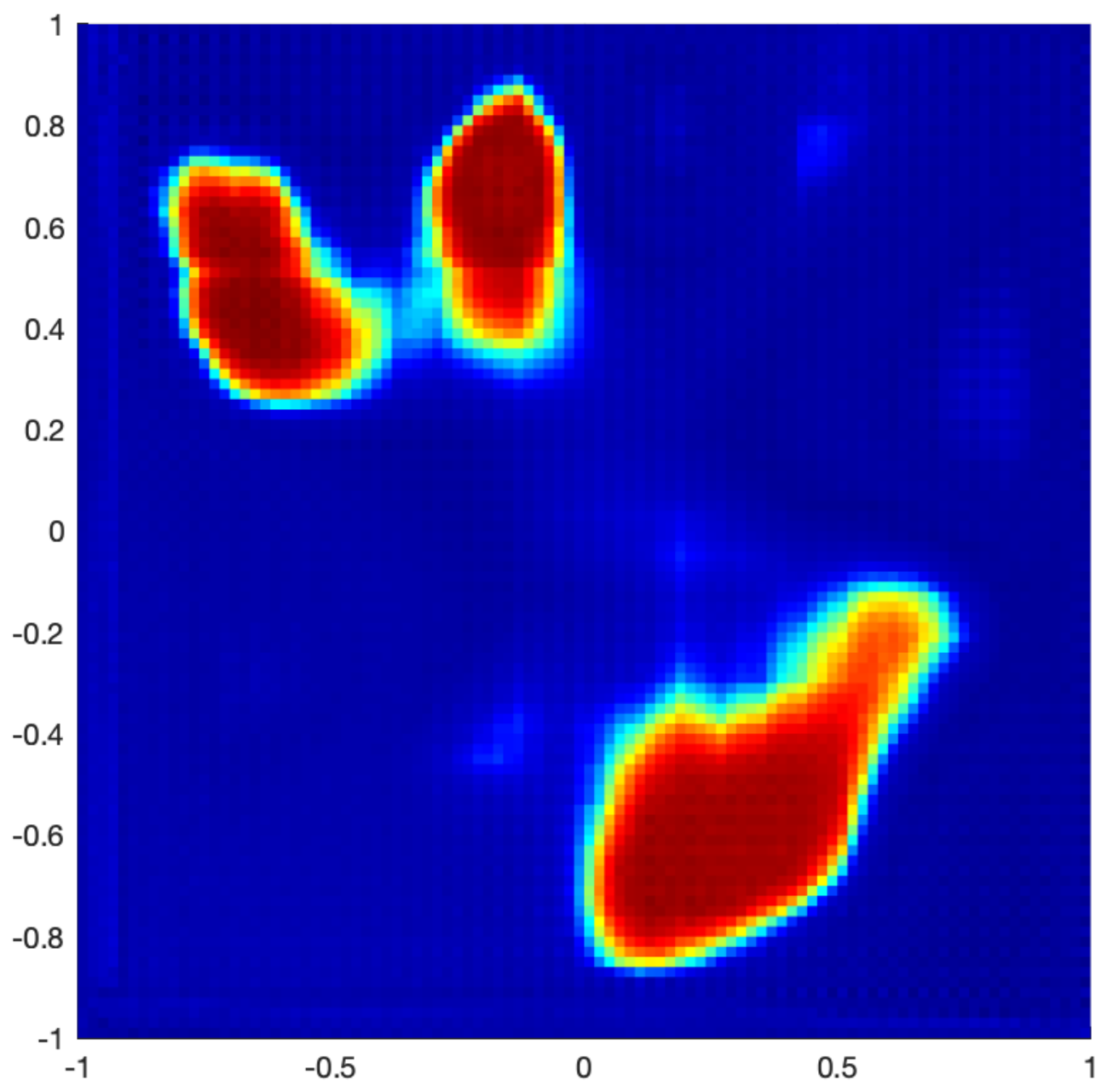}\\
\end{tabular}
  \caption{CNN-DDSM reconstruction for 3 cases in \textbf{Scenario 3} (4 ellipses) with different Cauchy data number and noise level: Case 1(top), Case 2(middle) and Case 3(bottom) } 
  \label{tab_ell}
\end{figure}

In this subsection, we present and discuss the training and test results of both the FNN-DDSM and CNN-DDSM for all three scenario mentioned above. 
We first show the evolution of loss function values versus training iterations in Figures \ref{fig:FN_train_process} and \ref{fig:train_process} of which the $x$-axis is in log-scale to capture the very quick convergence at the beginning. 
For both DNNs, the training errors of the single Cauchy data pair are larger than those of $10$ and $20$ pairs, while the errors of $20$ pairs are almost comparable with those of $10$ pairs. For FNN-DDSM, the error of the single pair stagnate after the first thousand of iterations. {\grc The gap of CNN-DDSM between the error of single pair and the errors of multiple pairs is clearly much smaller than the one of FNN-DDSM.} This also reflects that with single pair of Cauchy data, the CNN-DDSM performs better than the FNN-DDSM, which can be noticed in the second column of Figures \ref{tab_FN_3cir} - \ref{tab_ell} and Figures \ref{tab_FN_cir_comp} -\ref{tab_ell_comp}.
Although a single pair of Cauchy data is sufficient for CNN-DDSM to generate a satisfying reconstruction for some simple shaped inclusions, including multiple measurements are needed for more complicated cases.

\begin{figure}[htbp]
\begin{tabular}{ >{\centering\arraybackslash}m{1.5in}  >{\centering\arraybackslash}m{1.5in}  >{\centering\arraybackslash}m{1.5in}  >{\centering\arraybackslash}m{1.5in} }
\centering
\includegraphics[width=1.3in]{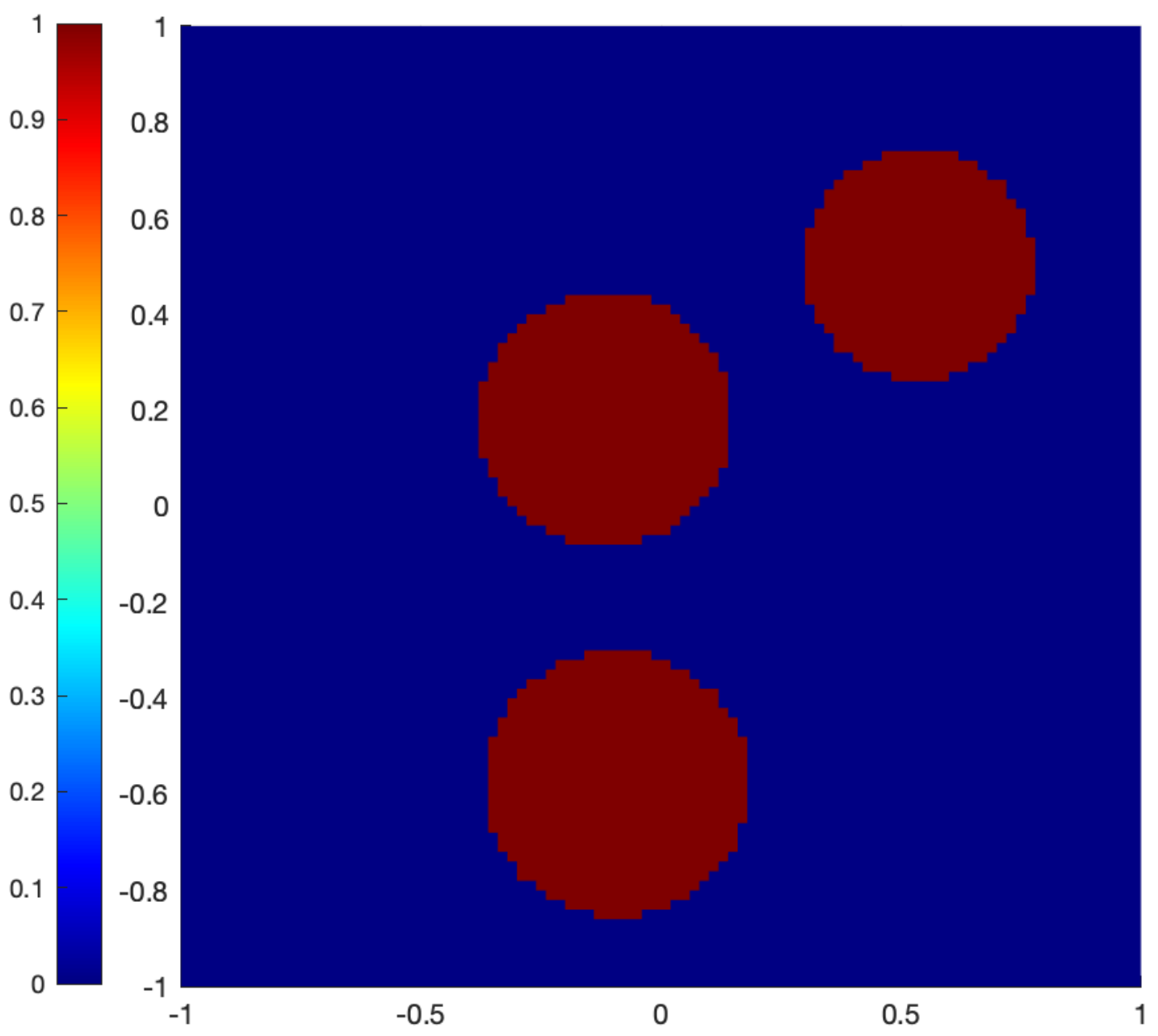}&
\includegraphics[width=1.2in]{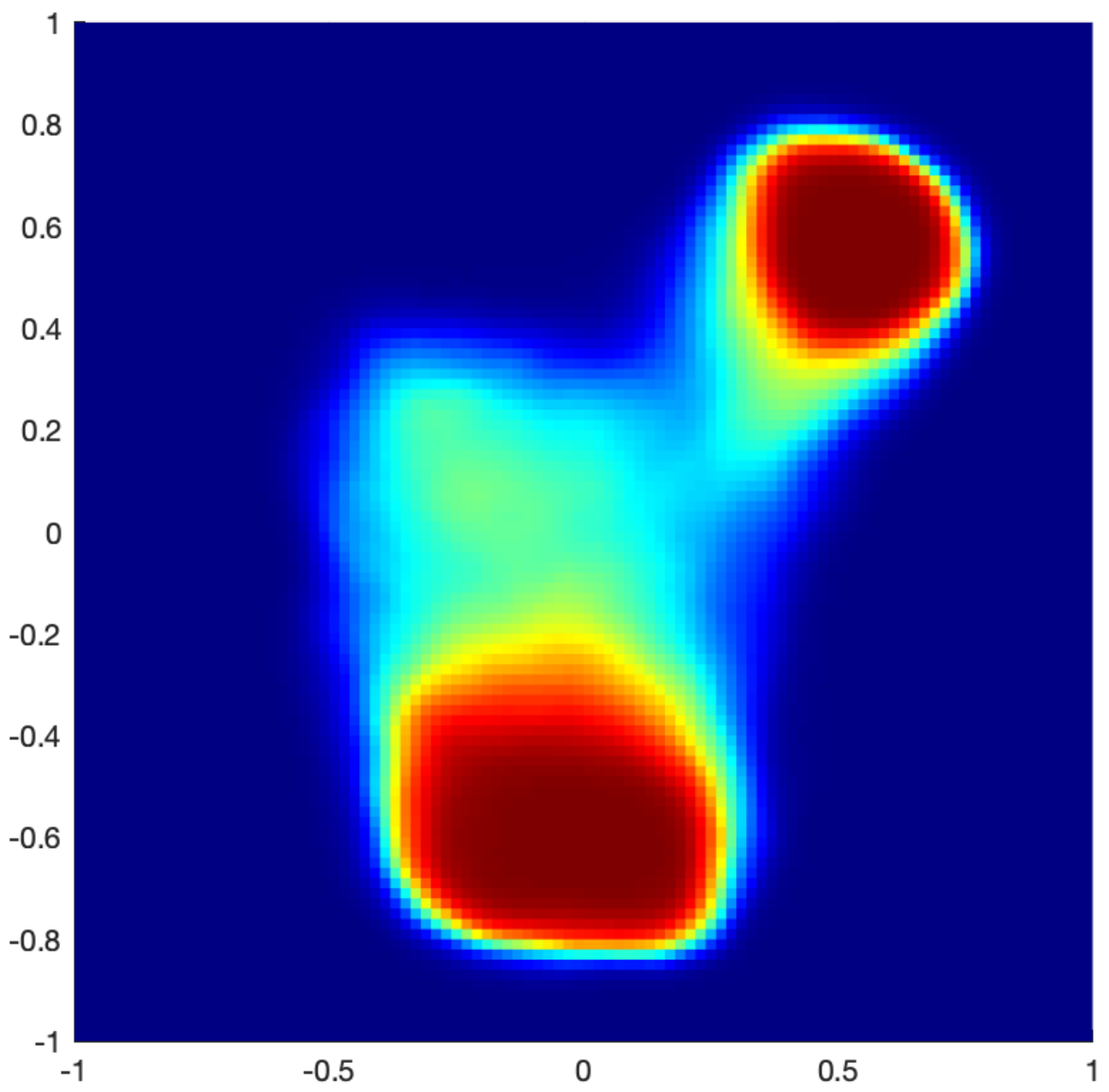}&
\includegraphics[width=1.2in]{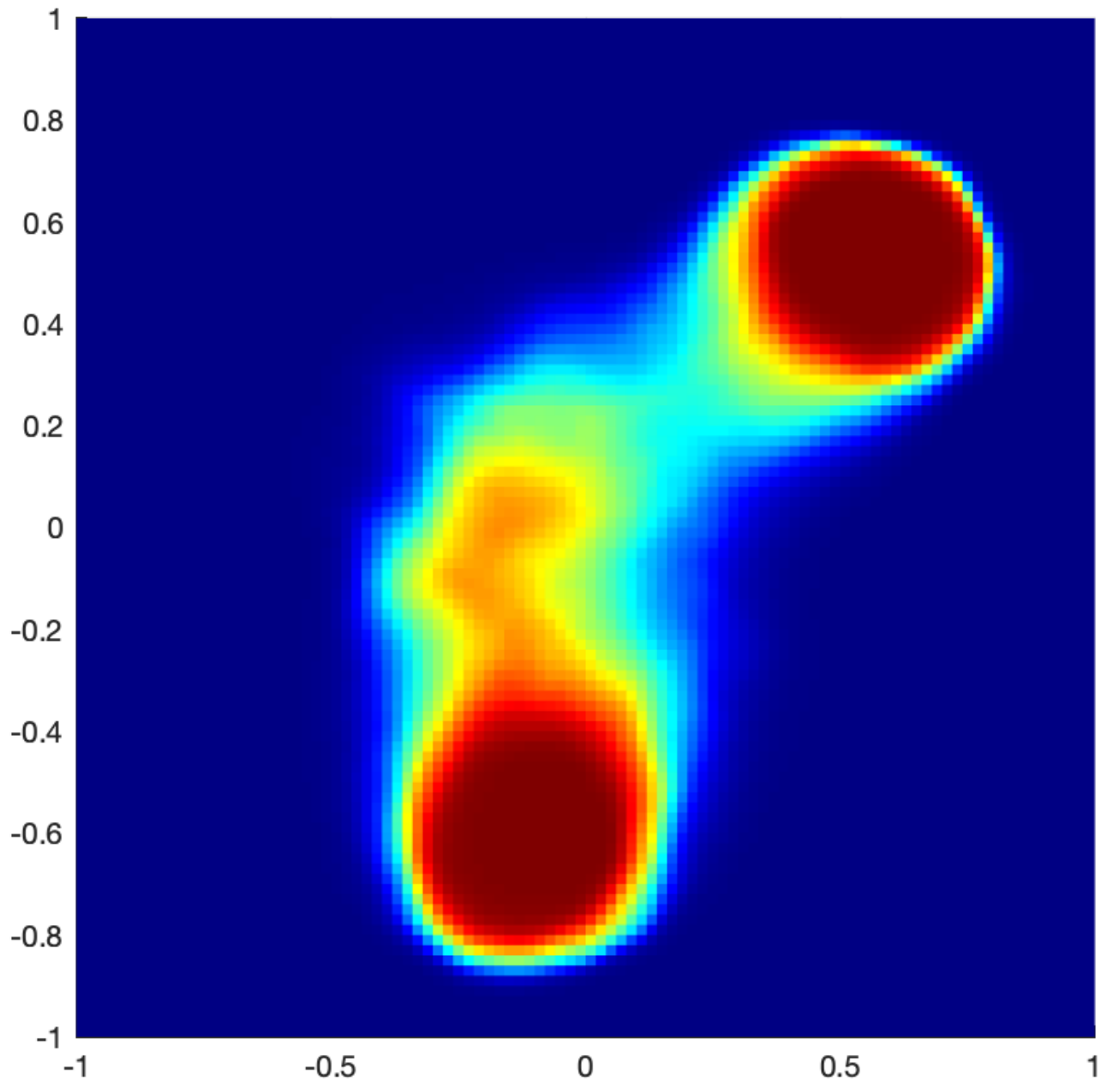}&
\includegraphics[width=1.2in]{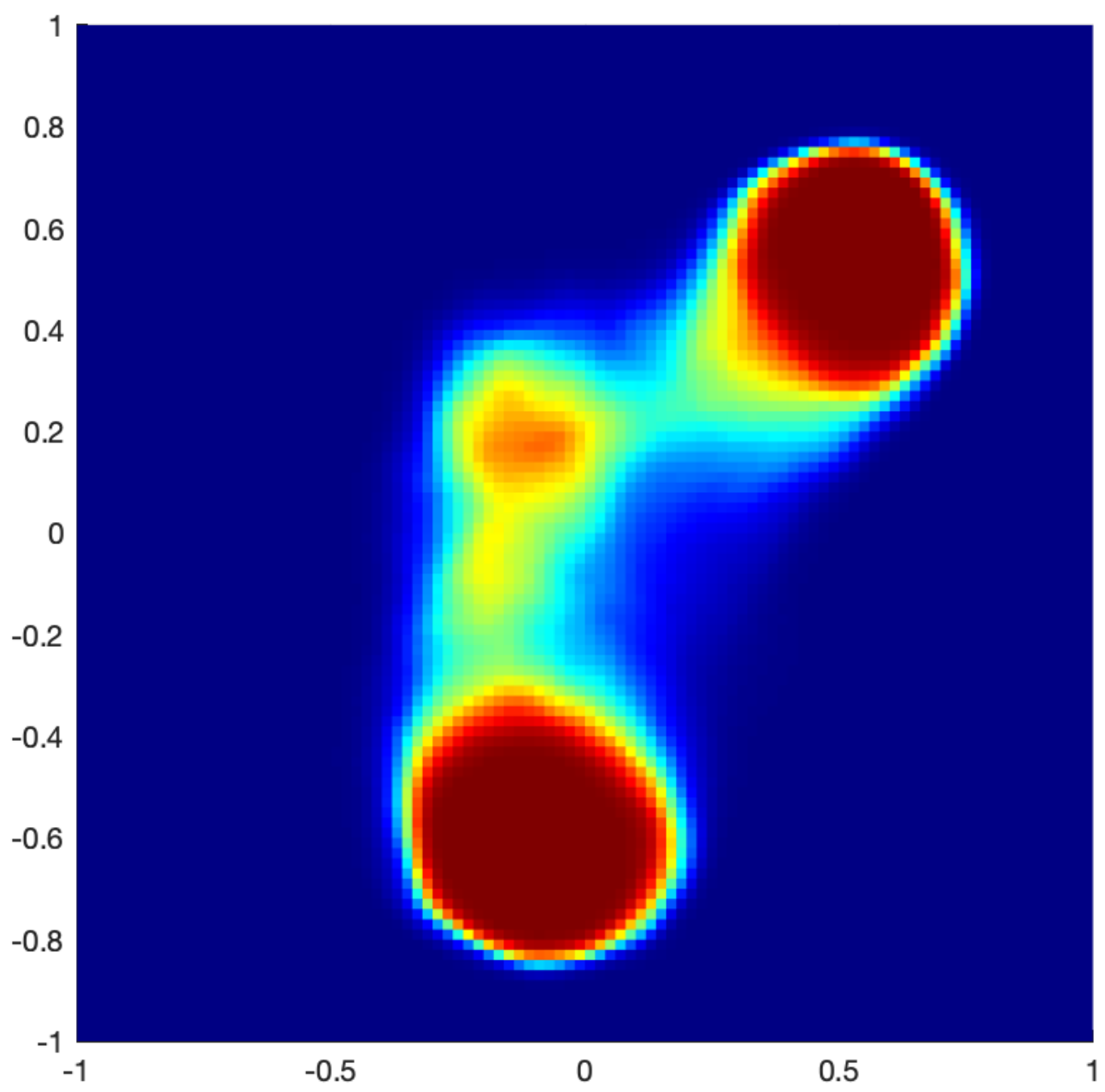}\\
\includegraphics[width=1.2in]{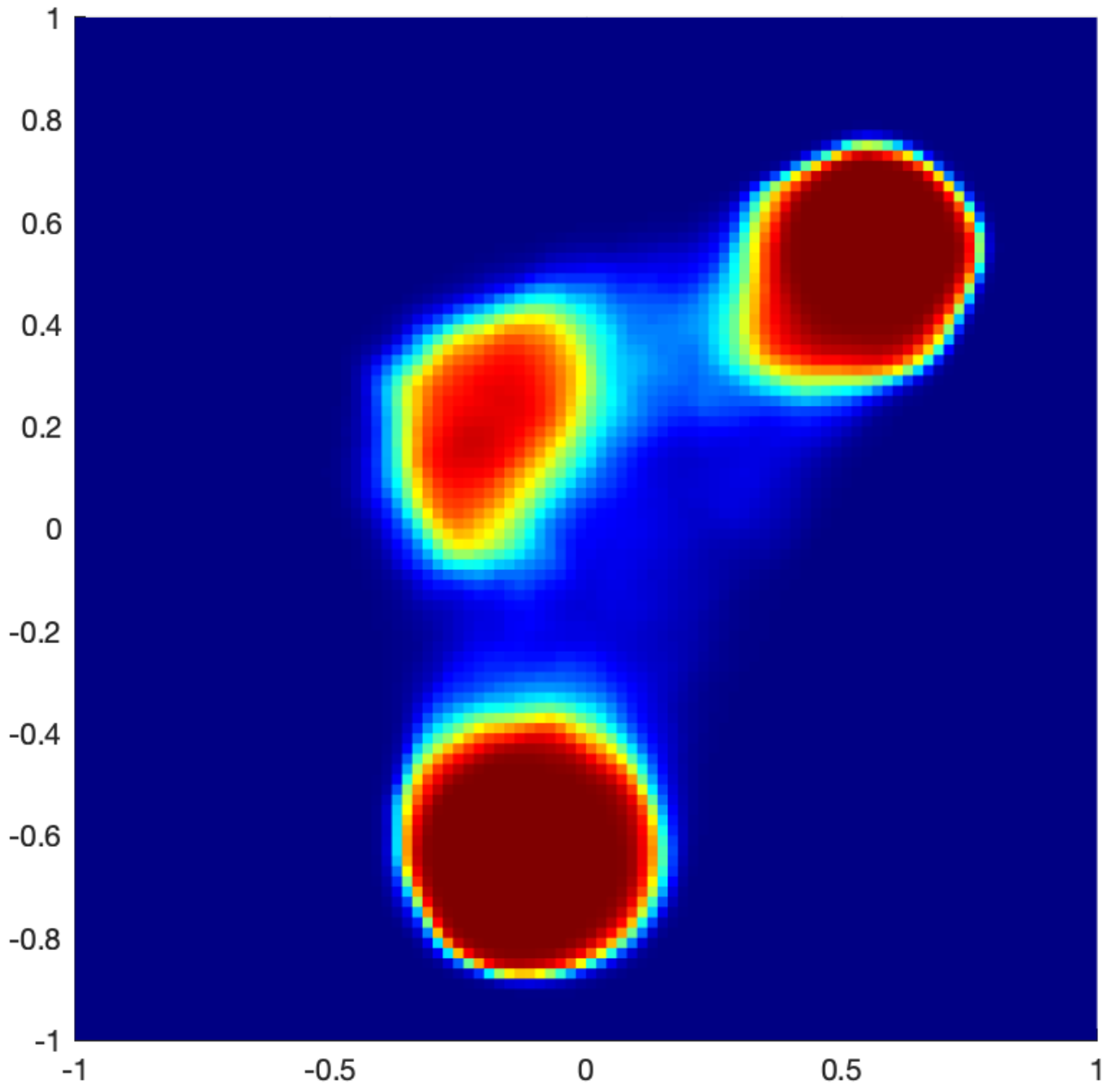}&
\includegraphics[width=1.2in]{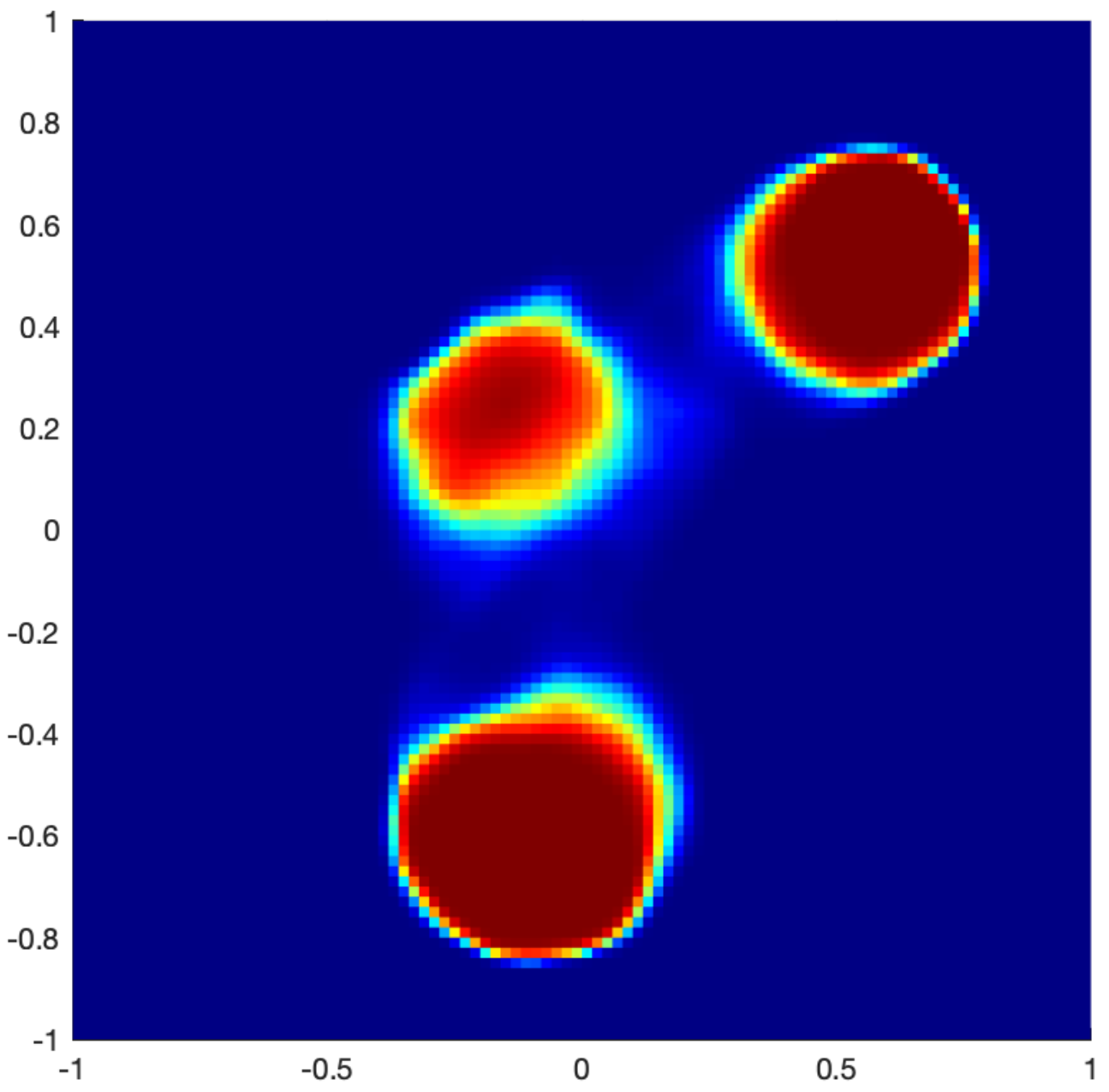}&
\includegraphics[width=1.2in]{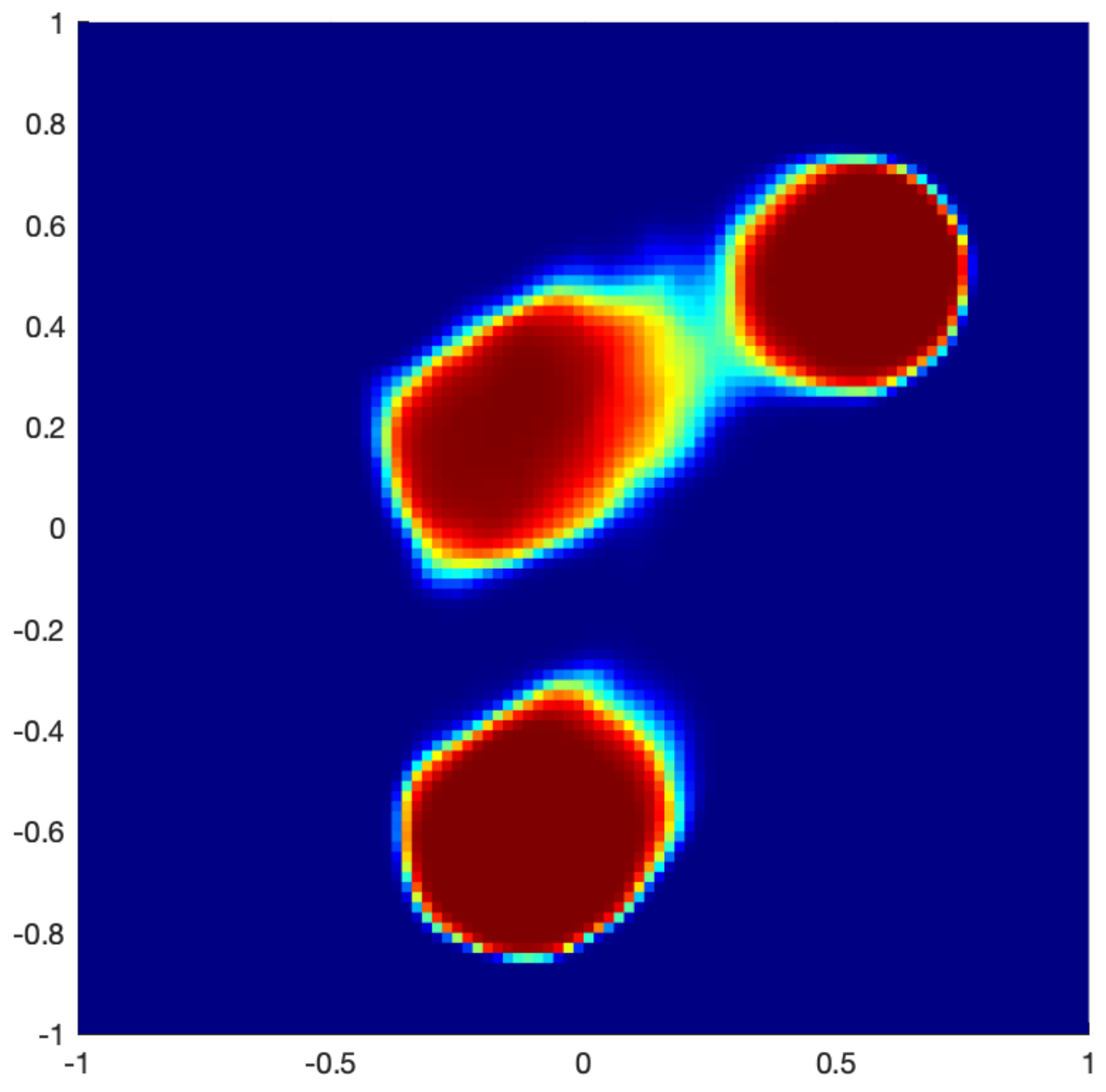}&
\includegraphics[width=1.2in]{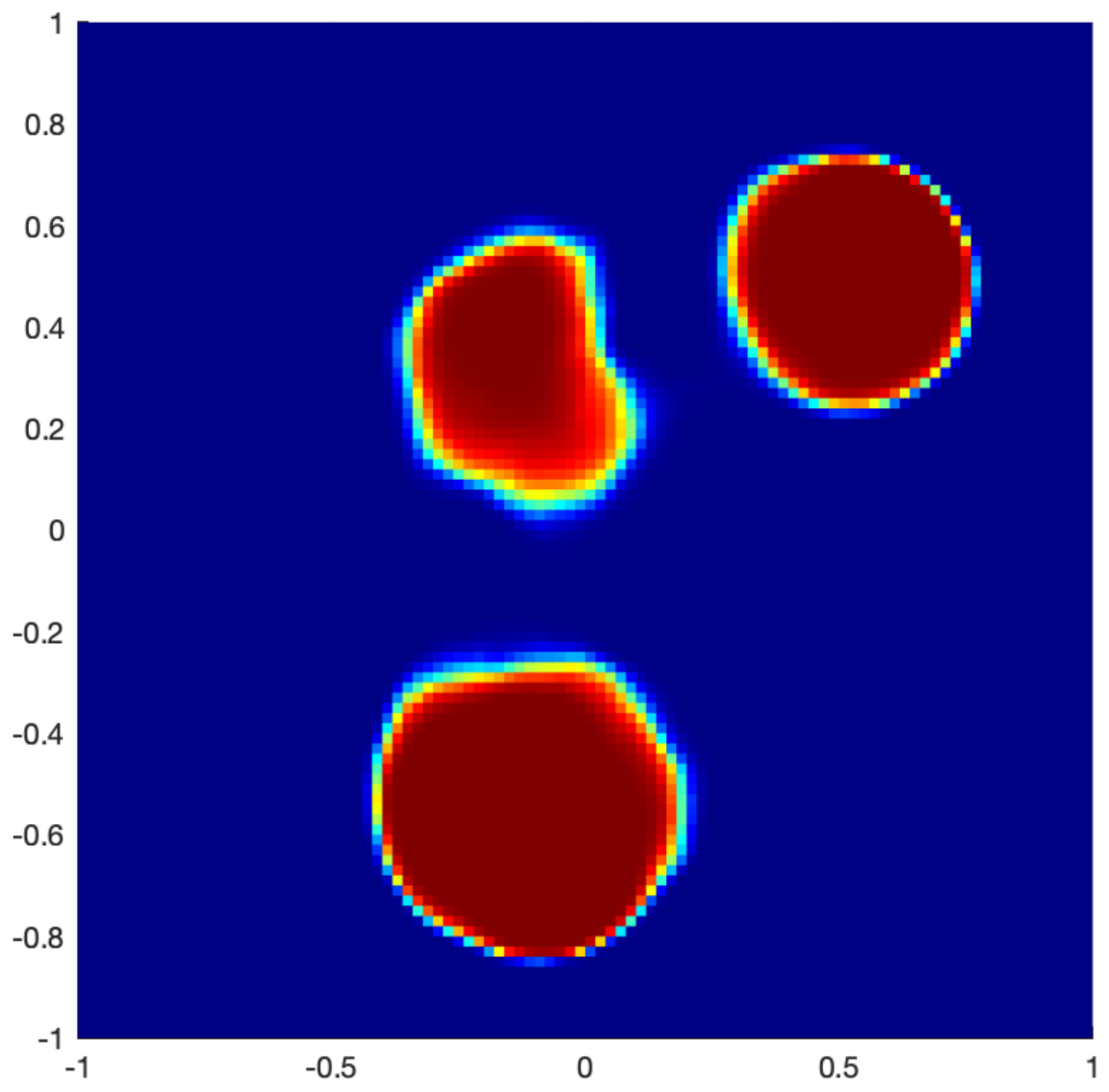}\\
\end{tabular}
  \caption{Some intermediate results in the training progress of the FNN-DDSM for  \textbf{Scenario 1} at iterations: $2000$, $5000$, $10000$, $25000$, $50000$, $75000$ and $100000$ (ordered from left to right and from top to bottom).} 
  \label{tab_FN_cir_train}
\end{figure}

\begin{figure}[htbp]
\begin{tabular}{ >{\centering\arraybackslash}m{1.5in}  >{\centering\arraybackslash}m{1.5in}  >{\centering\arraybackslash}m{1.5in}  >{\centering\arraybackslash}m{1.5in} }
\centering
\includegraphics[width=1.3in]{cir3_train0-eps-converted-to.pdf}&
\includegraphics[width=1.2in]{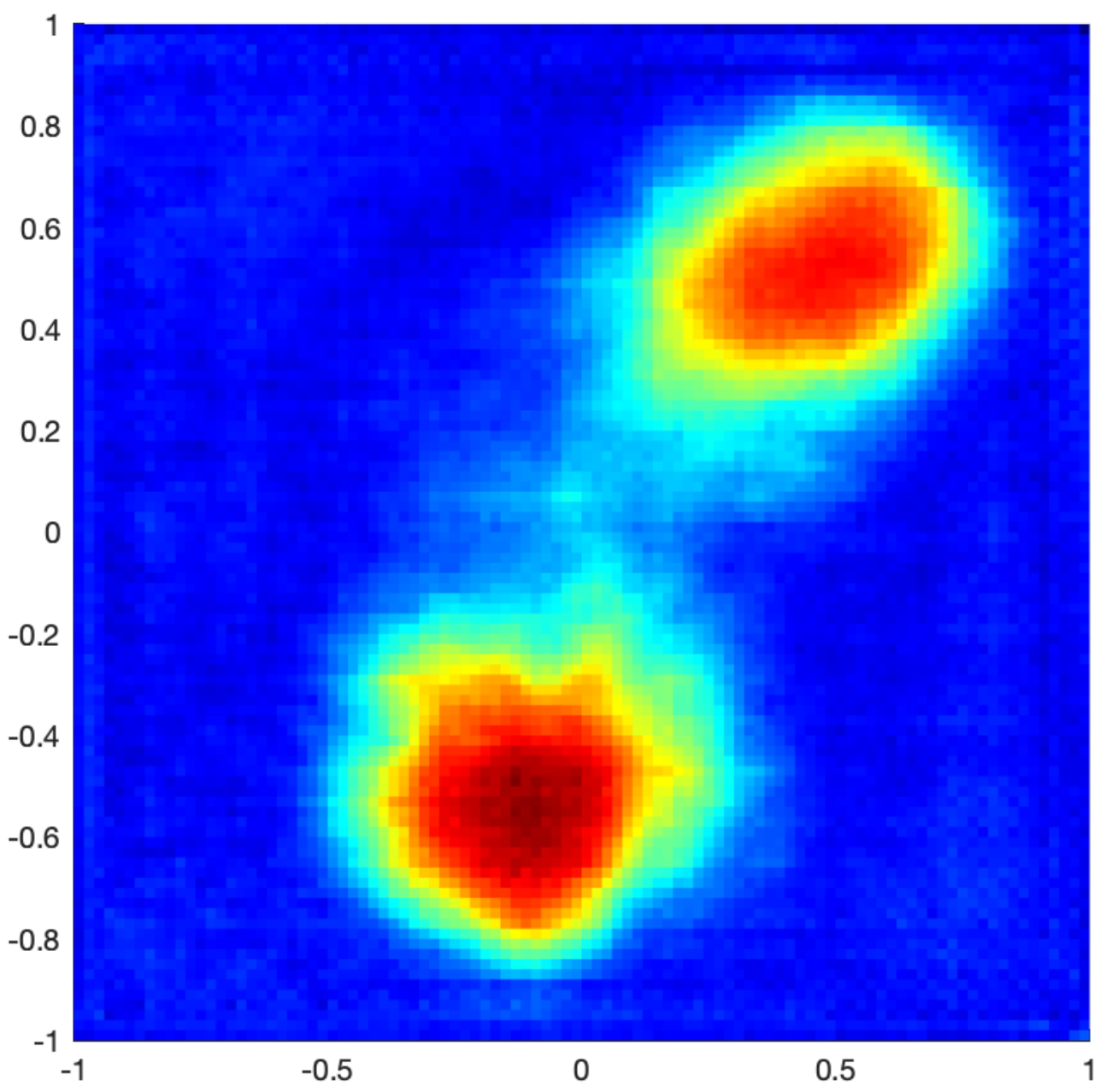}&
\includegraphics[width=1.2in]{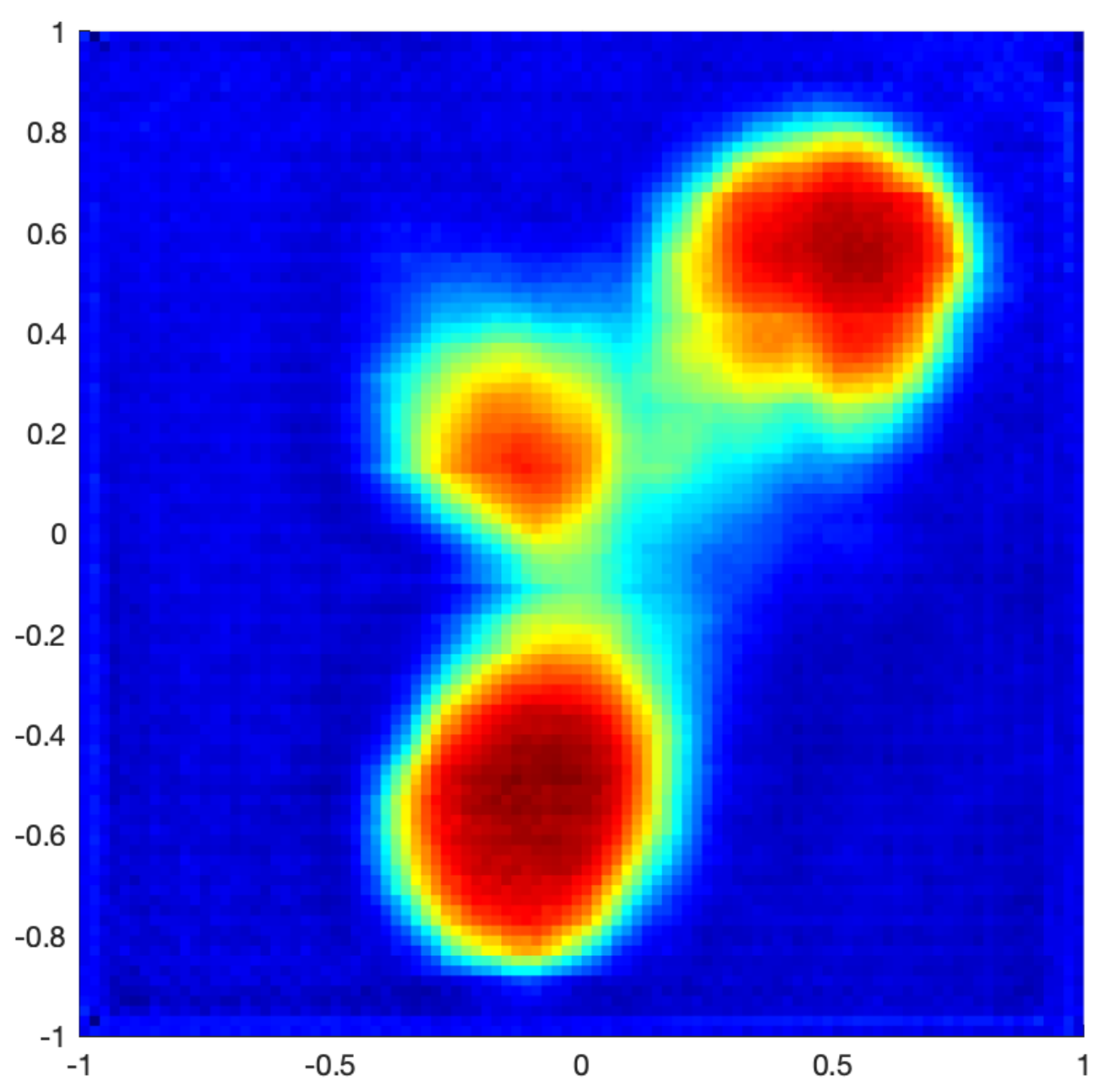}&
\includegraphics[width=1.2in]{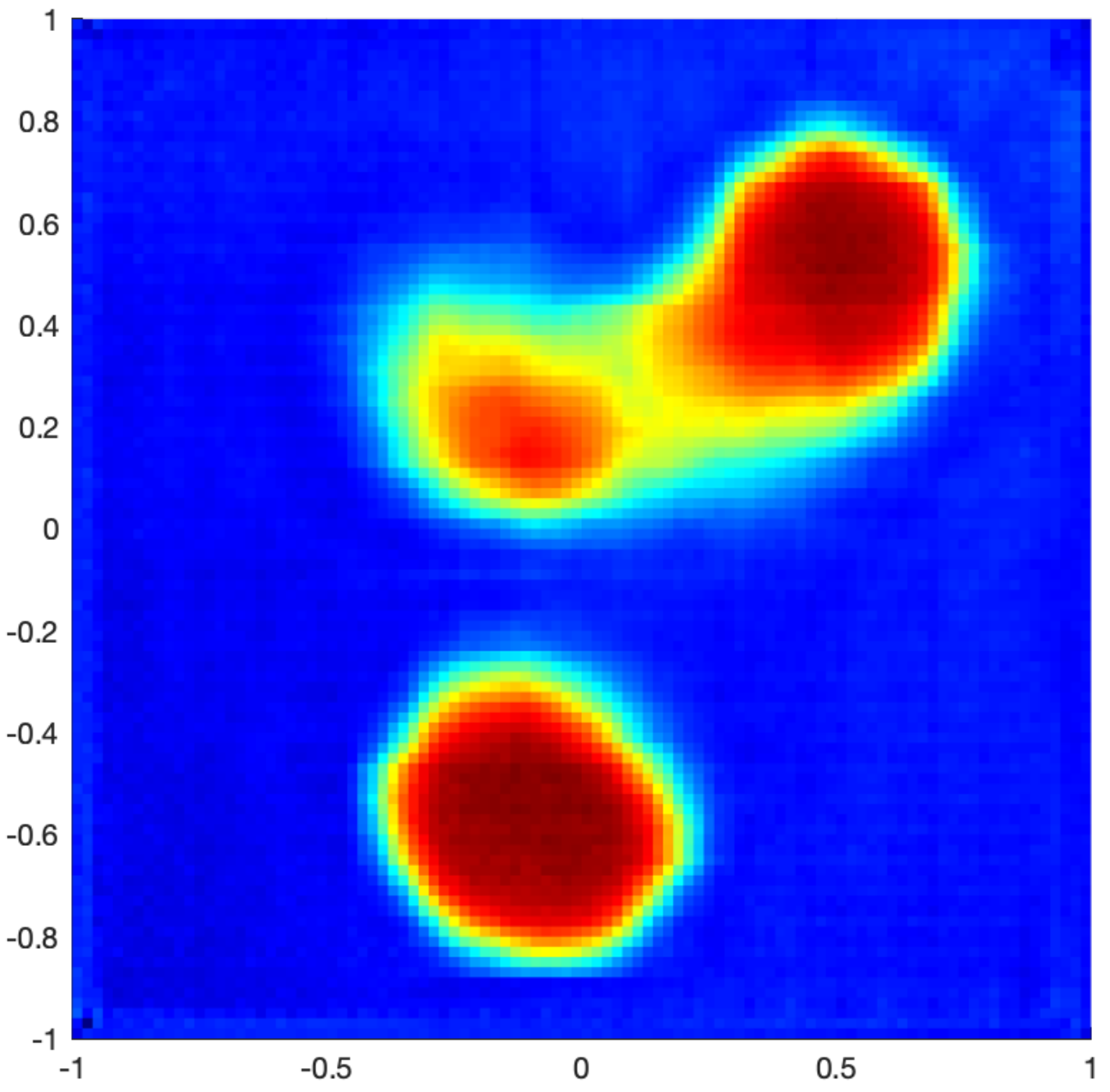}\\
\includegraphics[width=1.2in]{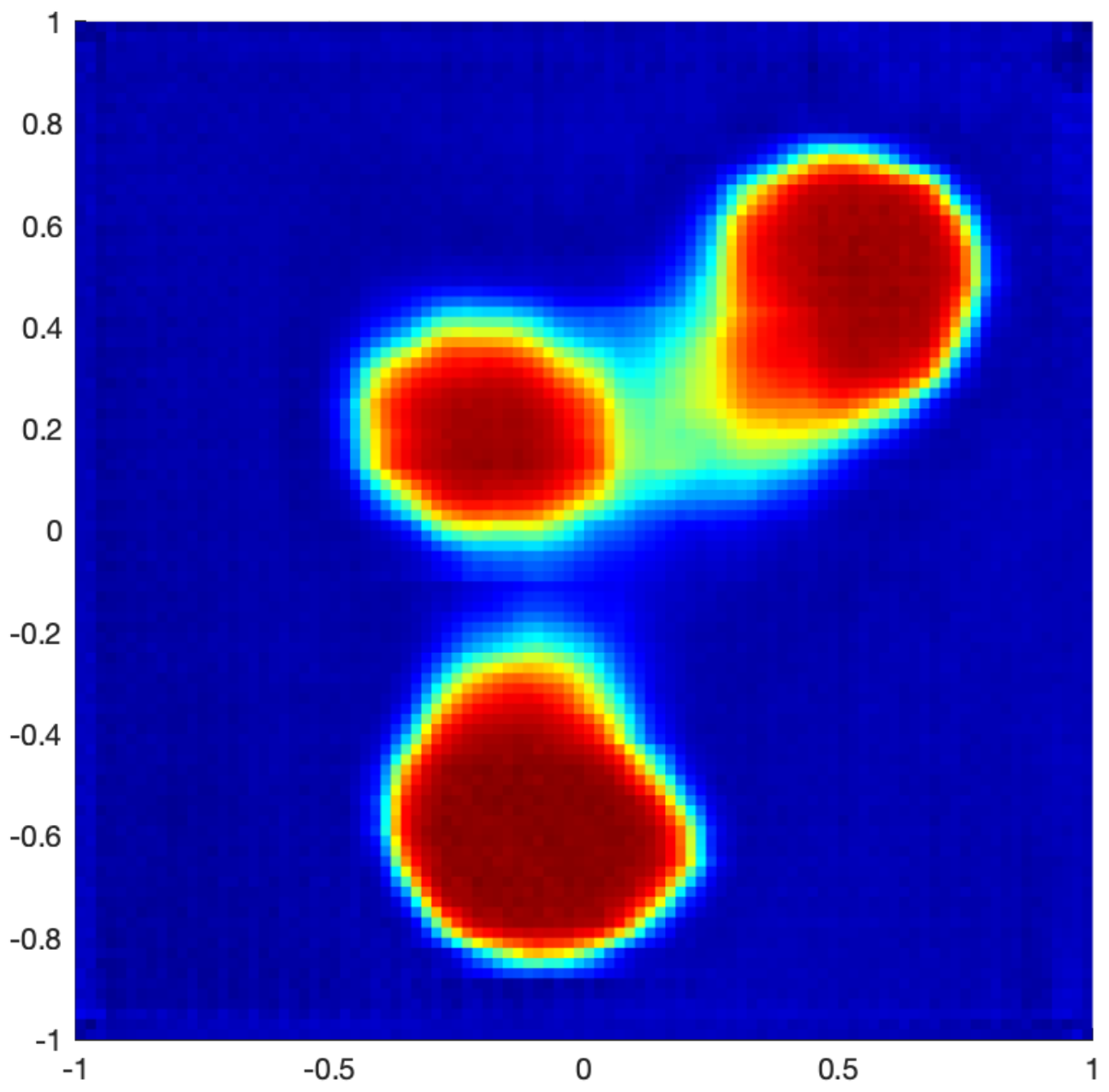}&
\includegraphics[width=1.2in]{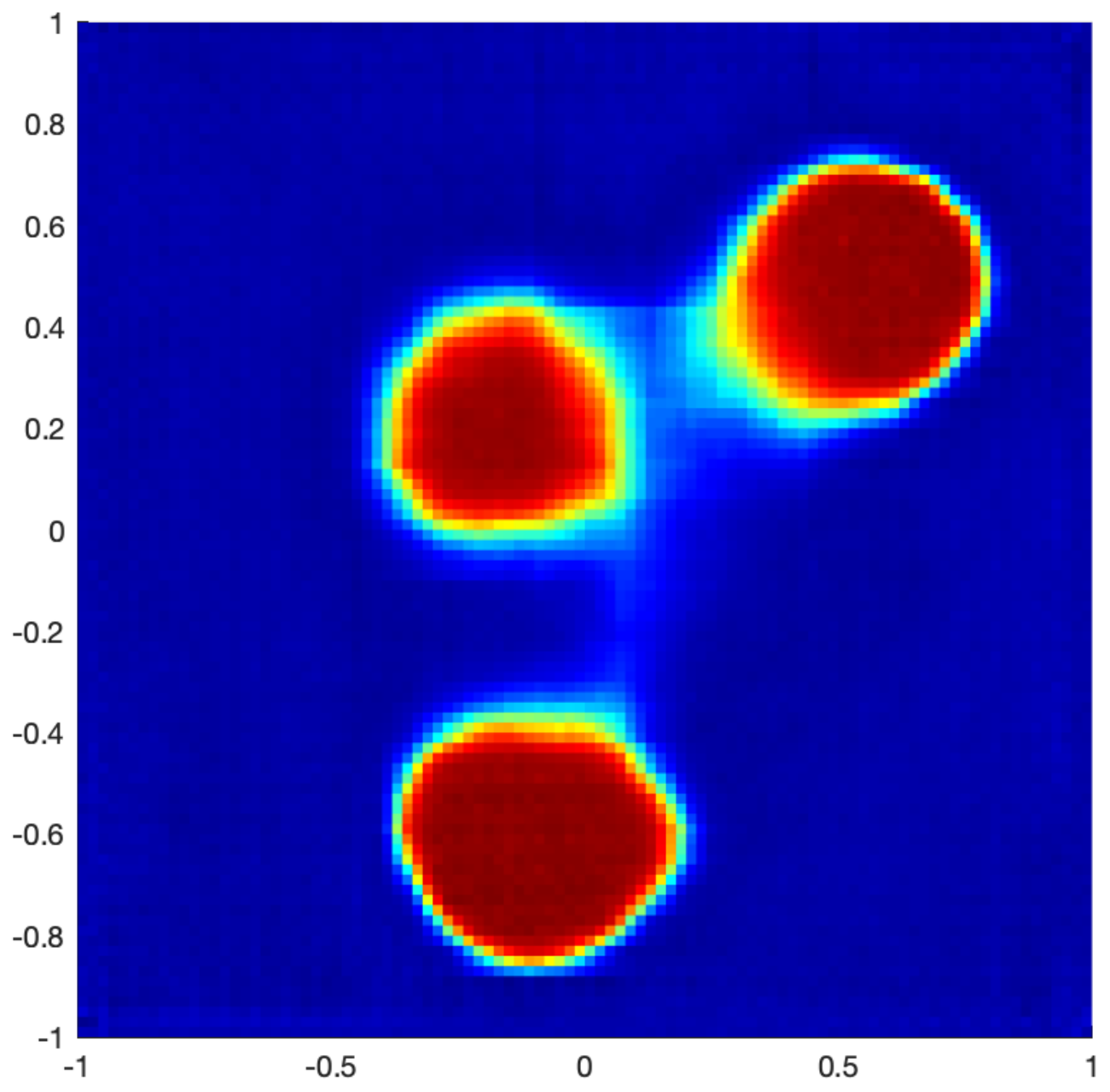}&
\includegraphics[width=1.2in]{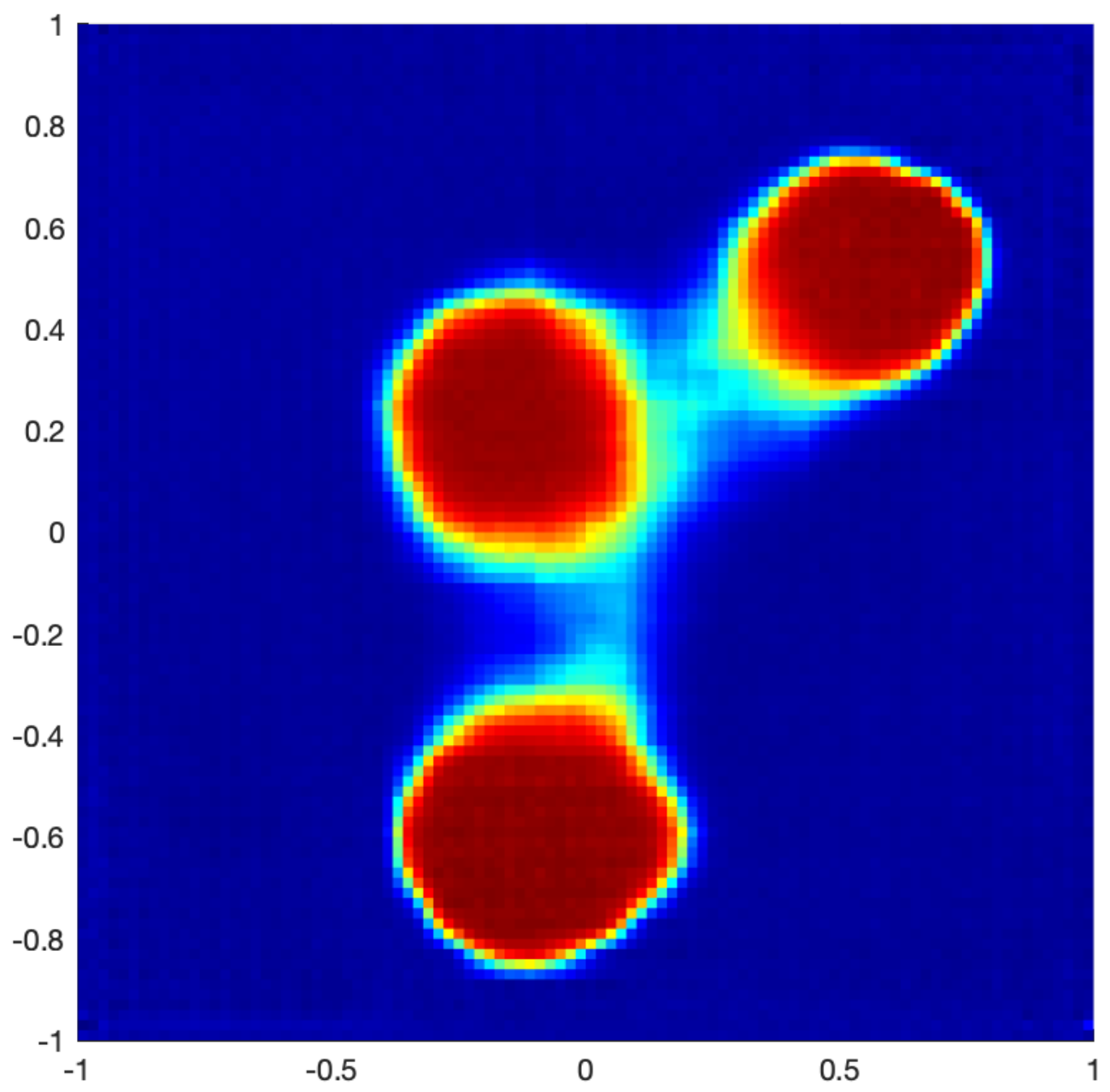}&
\includegraphics[width=1.2in]{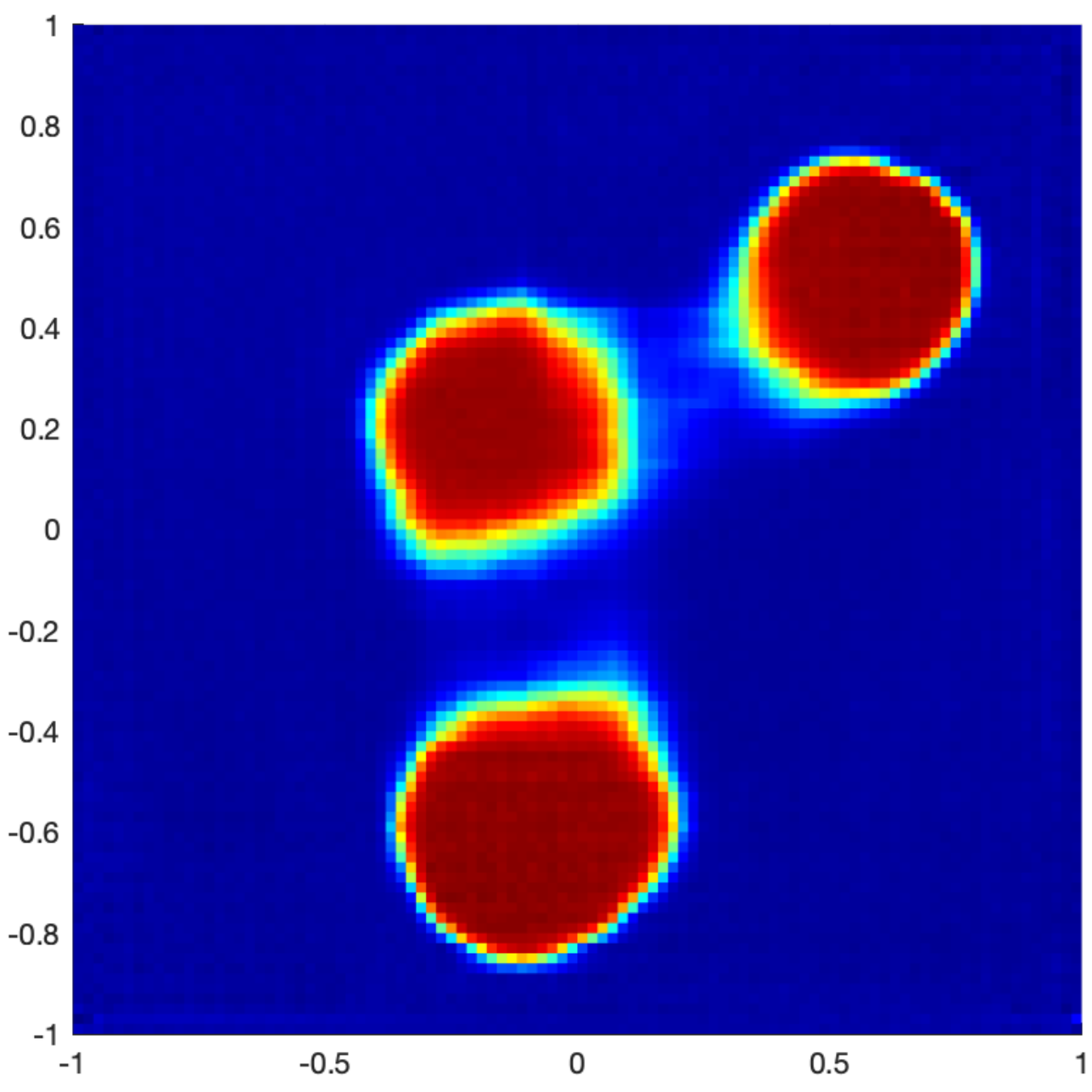}\\
\end{tabular}
  \caption{Some intermediate results during the training progress of  CNN-DDSM for  \textbf{Scenario 1} at iterations: $2000$, $5000$, $10000$, $25000$, $50000$, $75000$ and $100000$ (ordered from left to right and from top to bottom).} 
  \label{tab_cir_train}
\end{figure}

Figures \ref{tab_FN_cir_train} and \ref{tab_cir_train} show the reconstruction at certain iterations during the training progress in \textbf{Scenario 1}. It is very clear that both FNN-DDSM and CNN-DDSM will reconstruct the inclusions starting from the boundary to the center, which matches the nature of the EIT problem that only boundary data are available, and in general the closer to the boundary the easier the reconstruction. From the perspective of learning manner, the inclusions near the boundary can dominate the behavior of the boundary data which are relatively easier to be recognized and learned by DNNs during training, while the center inclusions have very minor effect on the boundary data which are more difficulty to be recognized. Similar behavior can be also found in prediction performance we will discuss later. 
Many conventional approaches for EIT such as \cite{2017JinXuZou} have shown that inclusions near the center and far away from the boundary is indeed very difficult to be detected since the boundary data are very unsensitive to its shape and location. But our DDSMs have some promising results for addressing this issue, which will be discussed in Section \ref{sec:further}. 

Now we present and discuss the reconstruction obtained by DDSMs for all three scenarios. For each one, we choose three inclusion samples in the test set and show the reconstruction in Figures \ref{tab_FN_3cir}-\ref{tab_ell}. We plot the reconstruction generated by $N=1,10,20$ pairs of Cauchy data, and for $N=20$ we include the noise $\delta=10\%$ and $20\%$ to test the numerical stability. Based on these figures, we can clearly see that both FNN-DDSM and CNN-DDSM with multiple but limited boundary Cauchy data can yield quite accurate reconstruction for all these inclusions having different geometry and topology. In particular, {\grc for some complicated geometry with concavity changes} such as the Case 3 of \textbf{Scenario 1}(Figures \ref{tab_FN_3cir}-\ref{tab_3cir}), the Cases 2,3 of \textbf{Scenario 2}(Figures \ref{tab_FN_5cir}-\ref{tab_5cir}) and Cases 2,3 of \textbf{Scenario 3}(Figures\ref{tab_FN_ell}-\ref{tab_ell}), both CNN-DDSM and FNN-DDSM can accurately capture the shape and the position of the inclusions. We highlight that our DDSMs show great potential to handle the case that the inclusions are not star-shape such as the Case 3 of \textbf{Scenario 1} (Figures \ref{tab_FN_3cir}-\ref{tab_3cir}) and the Cases 2,3 of the \textbf{Scenario 2} (Figures \ref{tab_FN_5cir}-\ref{tab_5cir}), which is very hard to achieve in general.  
Furthermore, we observe that the DDSMs are highly reliable with respect to the noise (up to 20\%) with 20 pairs. In some cases, even the 20\% noise has no obvious effect on the reconstruction. For other more challenging cases such as Case 2 of \textbf{Scenario 2} (Figures \ref{tab_FN_5cir}-\ref{tab_5cir}), the reconstruction with the 20\% noise can be still used as reasonable predictions to the true inclusions. 
We emphasize that such a large noise can totally destroy the reconstruction for many conventional approaches. So we believe the proposed DDSMs inherit and enhance the robustness feature of the DSM \cite{chow2014direct}, which is a considerable merit for solving the EIT problem that is extremely ill-posed and sensitive to noise.

Moreover, comparing the reconstruction from the FNN-DDSM and CNN-DDSM, we notice that the CNN-DDSM is a bit better than the FNN-DDSM which can be seen from the following two points. First, with the single Cauchy data, the reconstructions of the FNN-DDSM are too rough, which barely convey any information, while the reconstructions of the CNN-DDSM are able to contain the main features of the true coefficients. In particular, for the relatively simple case that the basic circular and elliptic components are disjoint, the CNN-DDSM can yield quite good reconstruction such as the Case 1 in each scenario, but the FNN-DDSM fails. Second, the comparison between the Cases 2,3 of \textbf{Scenario 2} (Figures \ref{tab_FN_5cir}-\ref{tab_5cir}) and the Case 3 of \textbf{Scenario 3} (Figures \ref{tab_FN_ell}-\ref{tab_ell}) shows that the CNN-DDSM can yield slightly better reconstruction at those portion near the domain center away from the boundary. {\grc All these gains of the CNN-DDSM are within our expectation since, as mentioned before, to predict the location of a point, it incorporates the information of the more neighbor points near this point} which may better reflect or approximate the format of the true high-dimensional index function. However we feel that the FNN-DDSM seems more stable with respect to noise which can be seen from the Case 3 of each Scenario, since, we think,
it only uses the information at the single point to predict its location and thus involves relatively less noise. Another attractive feature of the FNN-DDSM is that its output value has clear probabilistic interpretation. So we can directly read from the plots to conclude which portion is almost certainly inside (red) or outside (blue) the inclusion and which portion can be hardly determined due to lack of information.

\subsection{Further Investigation}
\label{sec:further}

\begin{figure}[htbp]
\centering
\begin{subfigure}{.32\textwidth}
     \includegraphics[width=2.2in]{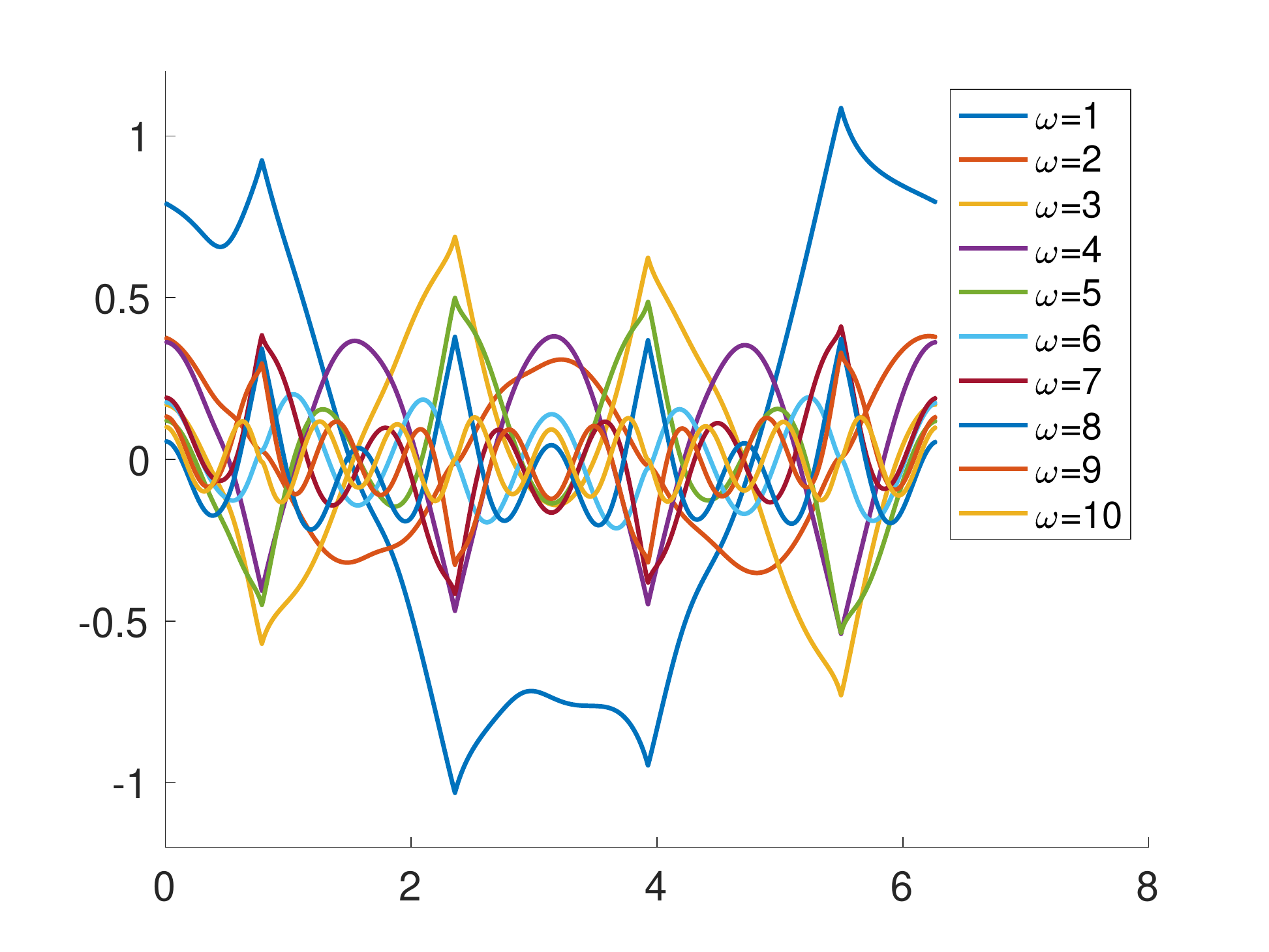}
\end{subfigure}
\begin{subfigure}{.32\textwidth}
     \includegraphics[width=2.2in]{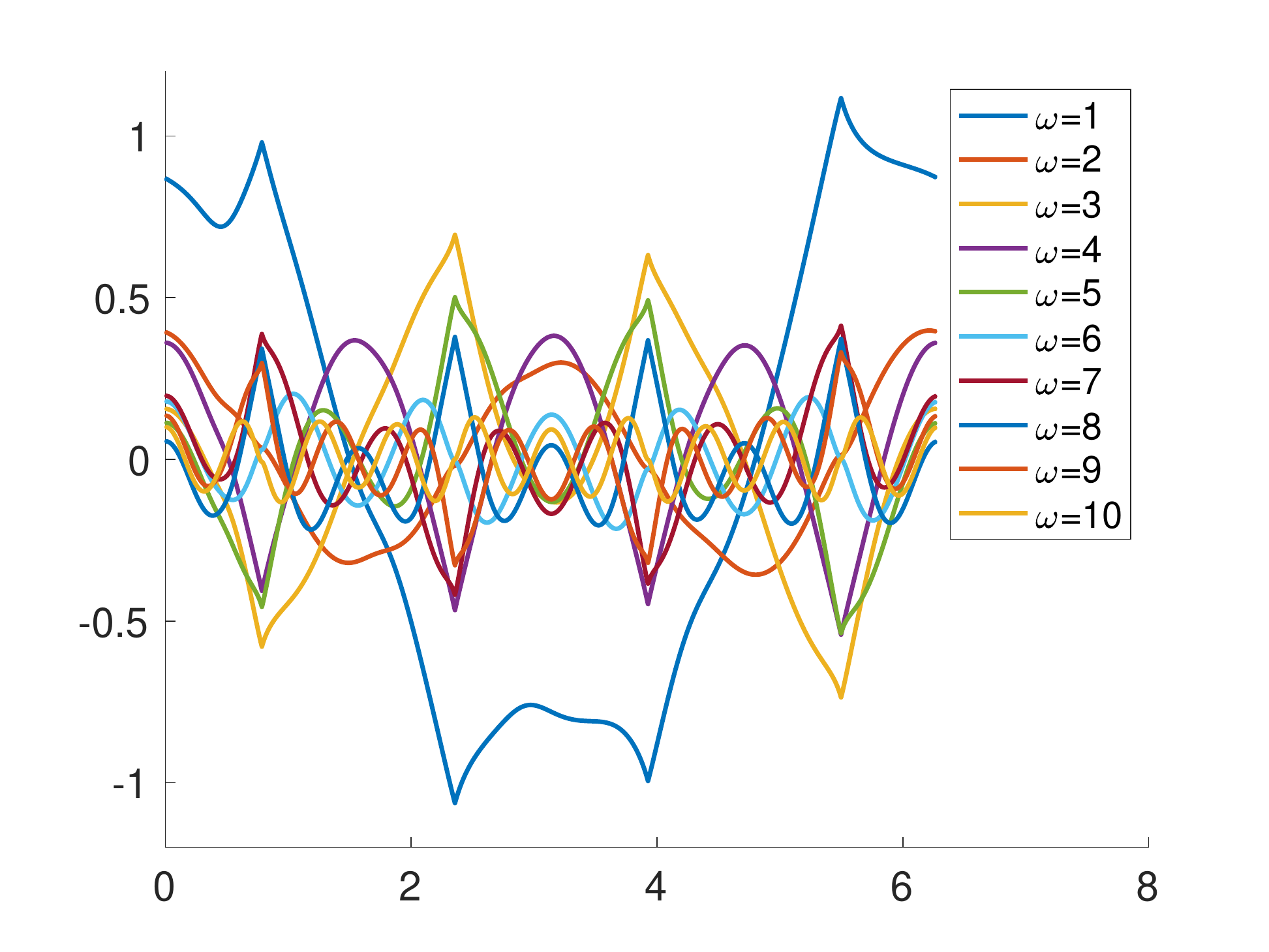}
\end{subfigure}
\begin{subfigure}{.32\textwidth}
     \includegraphics[width=2.2in]{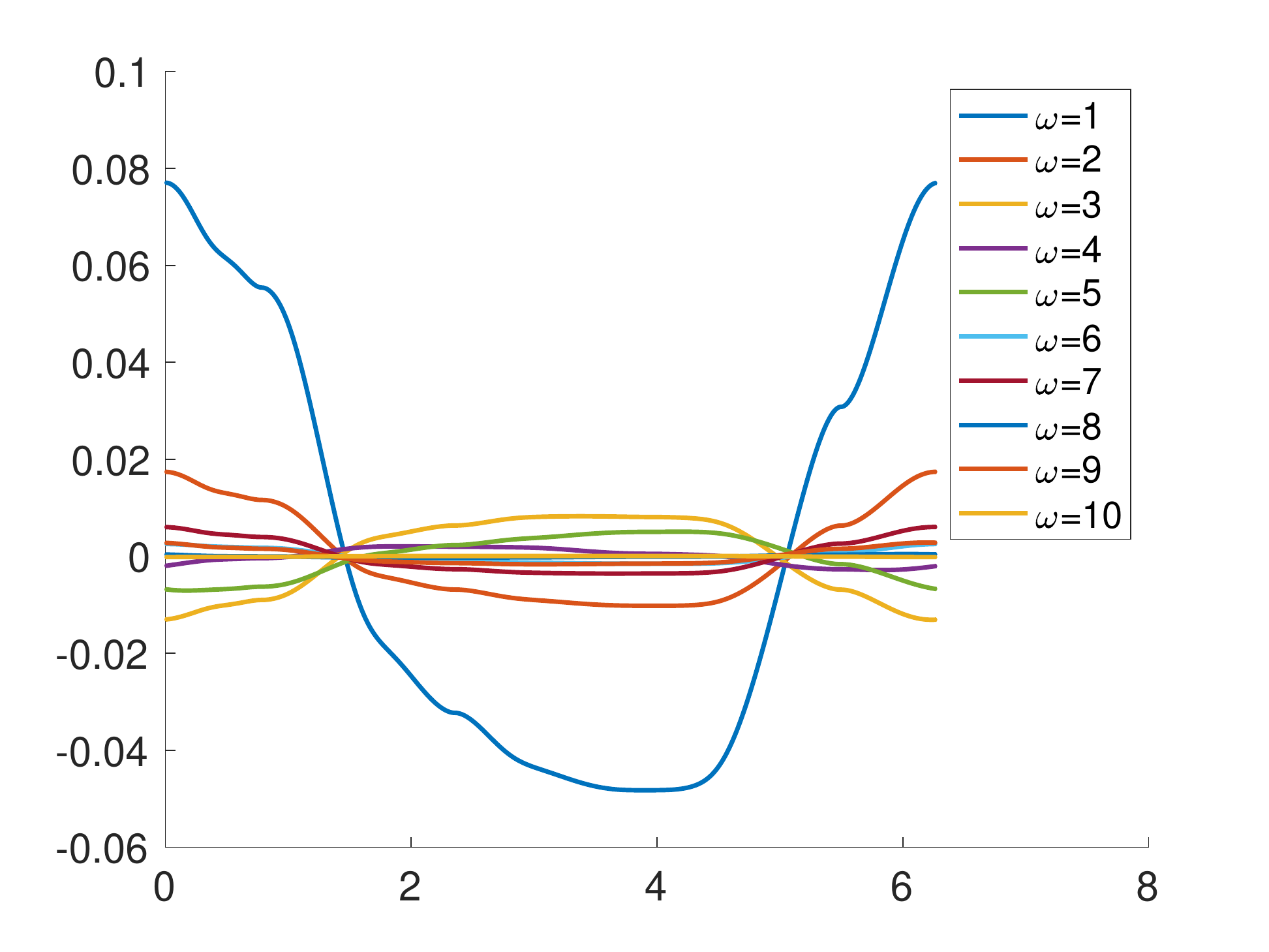}
\end{subfigure}
\caption{Plots of $u^c_{1,\omega}|_{\partial\Omega}$, $u^c_{2,\omega}|_{\partial\Omega}$ and $u^c_{1,\omega}|_{\partial\Omega}-u^c_{2,\omega}|_{\partial\Omega}$ versus the polar angle $\theta$ (of points on the boundary)}
  \label{fig:cir_compare} 
\end{figure}

\begin{figure}[htbp]
\centering
\begin{subfigure}{.32\textwidth}
     \includegraphics[width=2.2in]{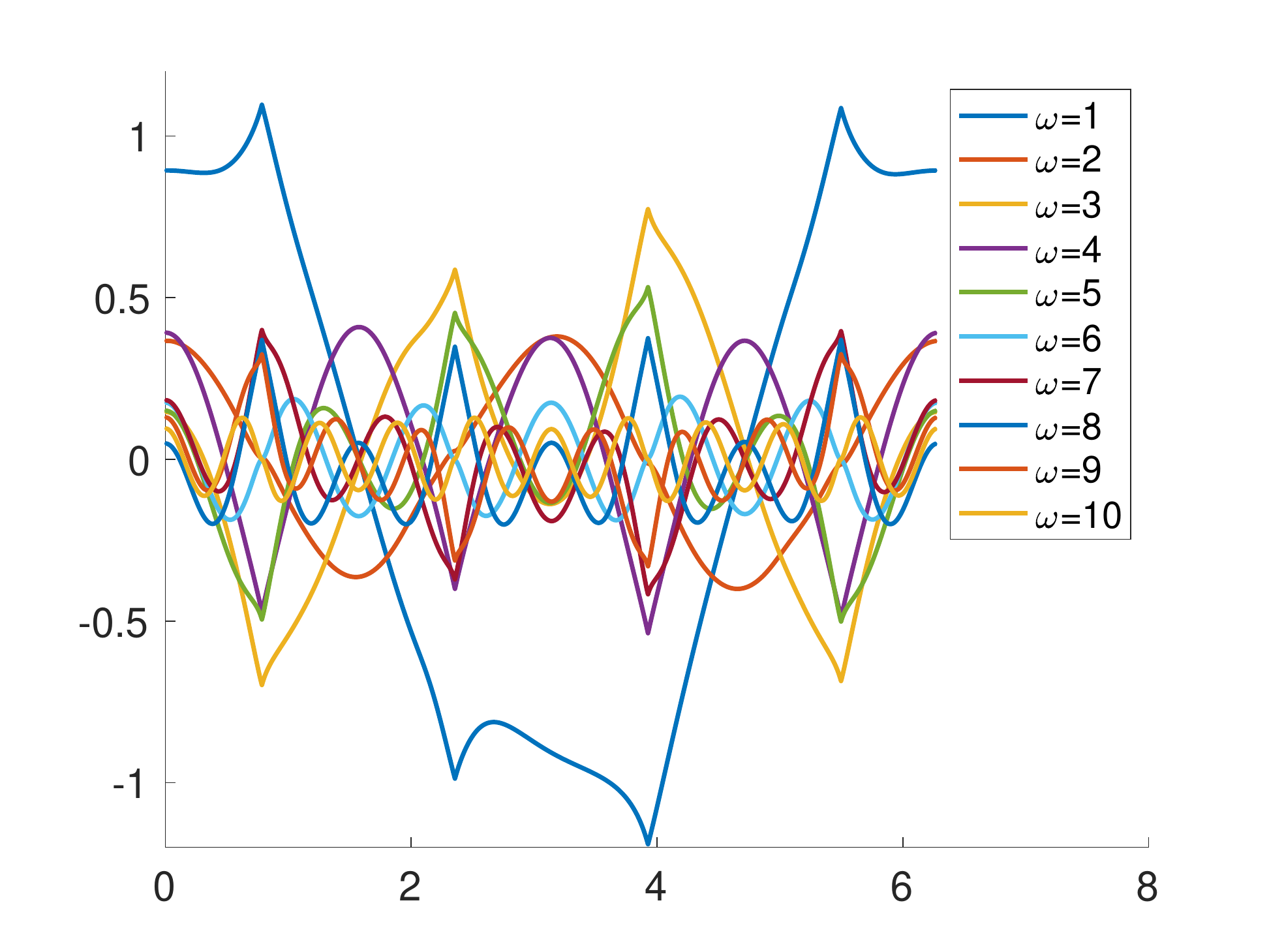}
\end{subfigure}
\begin{subfigure}{.32\textwidth}
     \includegraphics[width=2.2in]{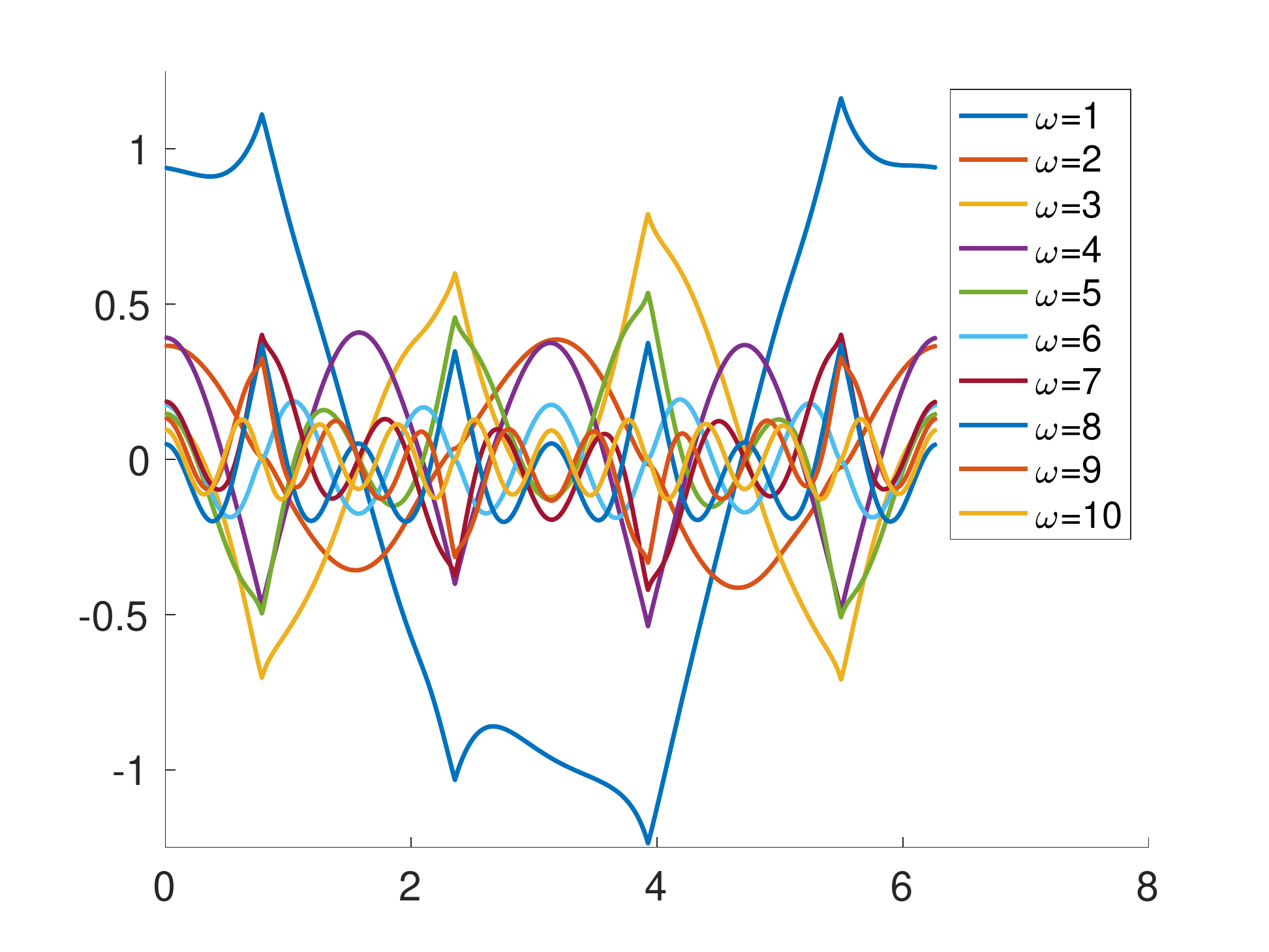}
\end{subfigure}
\begin{subfigure}{.32\textwidth}
     \includegraphics[width=2.2in]{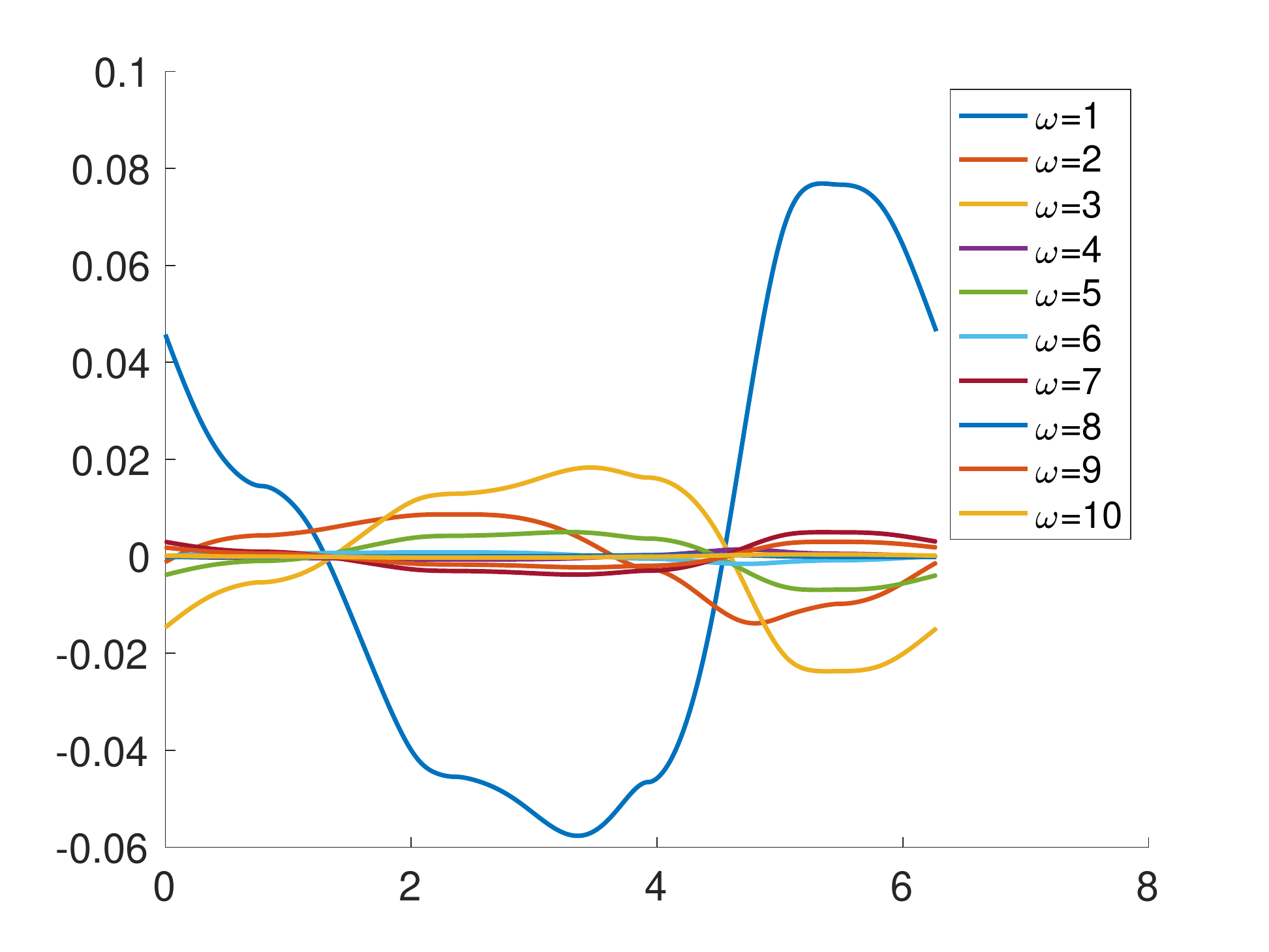}
\end{subfigure}
\caption{Plots of $u^e_{1,\omega}|_{\partial\Omega}$, $u^e_{2,\omega}|_{\partial\Omega}$ and $u^e_{1,\omega}|_{\partial\Omega}-u^e_{2,\omega}|_{\partial\Omega}$ versus the polar angle $\theta$ (of points on the boundary)}
  \label{fig:ell_compare} 
\end{figure}

\begin{figure}[htbp]
\begin{tabular}{ >{\centering\arraybackslash}m{0.9in} >{\centering\arraybackslash}m{0.9in} >{\centering\arraybackslash}m{0.9in}  >{\centering\arraybackslash}m{0.9in}  >{\centering\arraybackslash}m{0.9in}  >{\centering\arraybackslash}m{0.9in} }
\centering
True coefficients &
N=1, $\delta=0$&
N=10, $\delta=0$&
N=20, $\delta=0$&
N=20, $\delta=10\%$ &
N=20, $\delta=20\%$ \\
\includegraphics[width=1.1in]{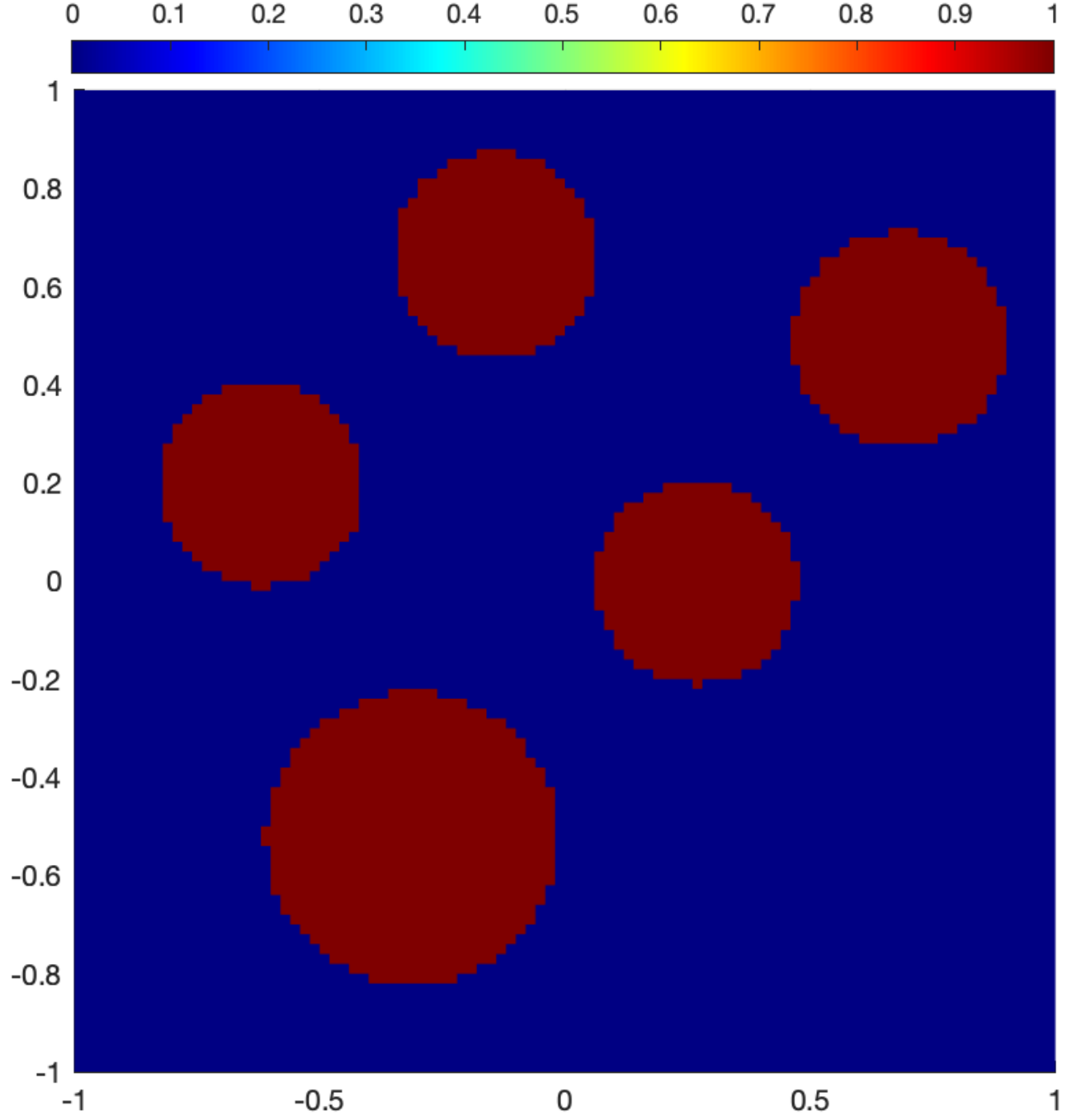}&
\includegraphics[width=1.1in]{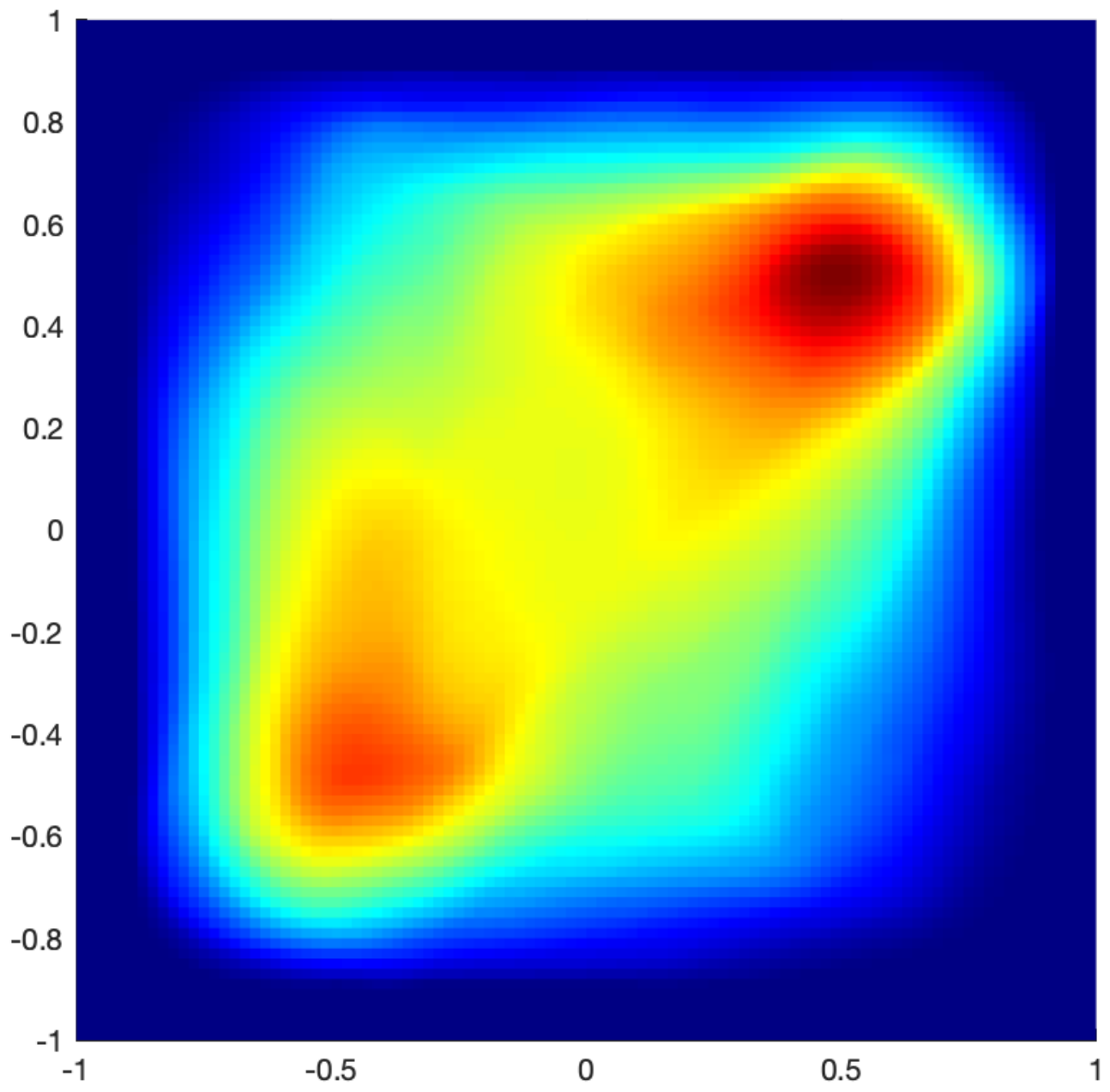}&
\includegraphics[width=1.1in]{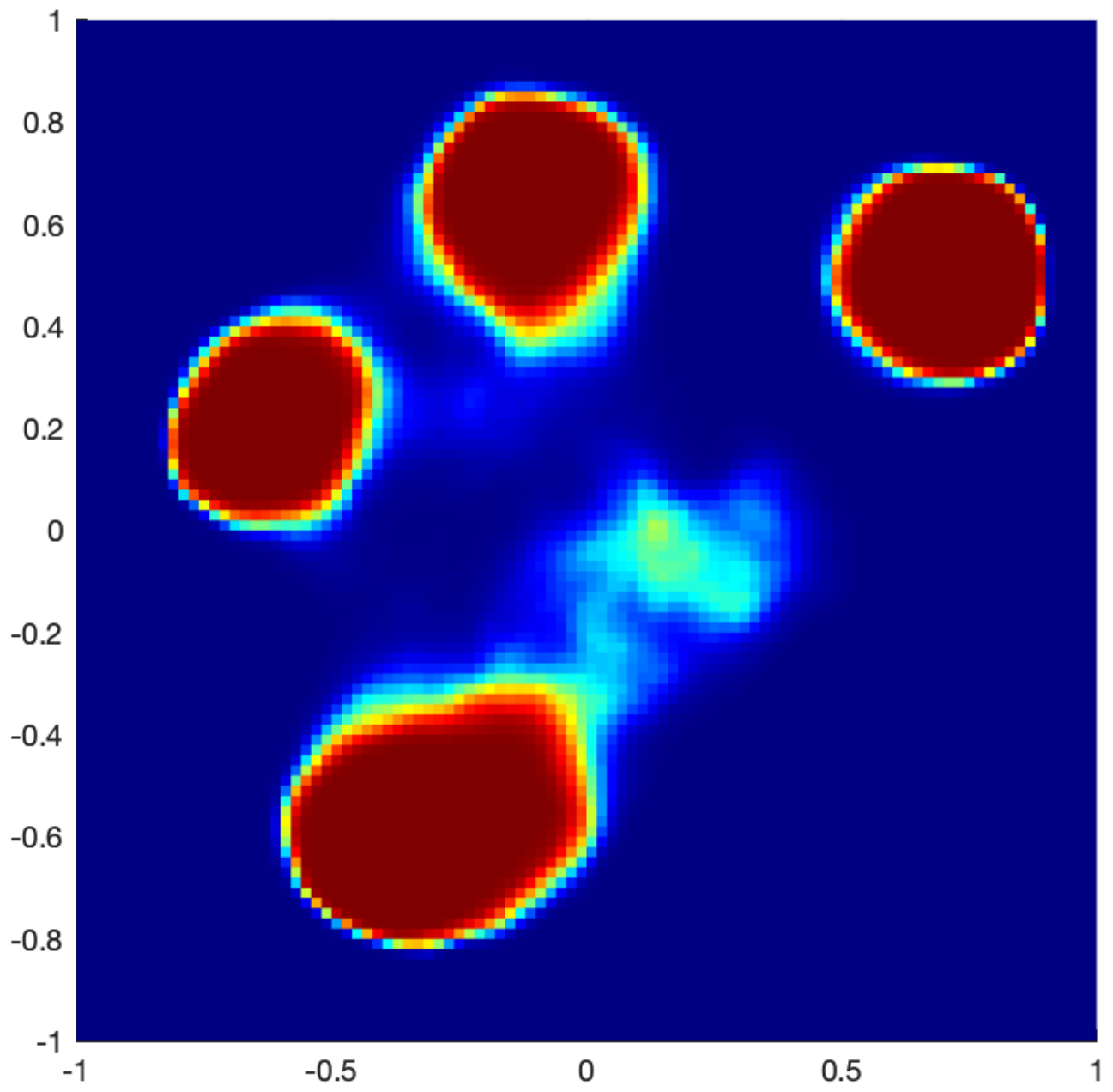}&
\includegraphics[width=1.1in]{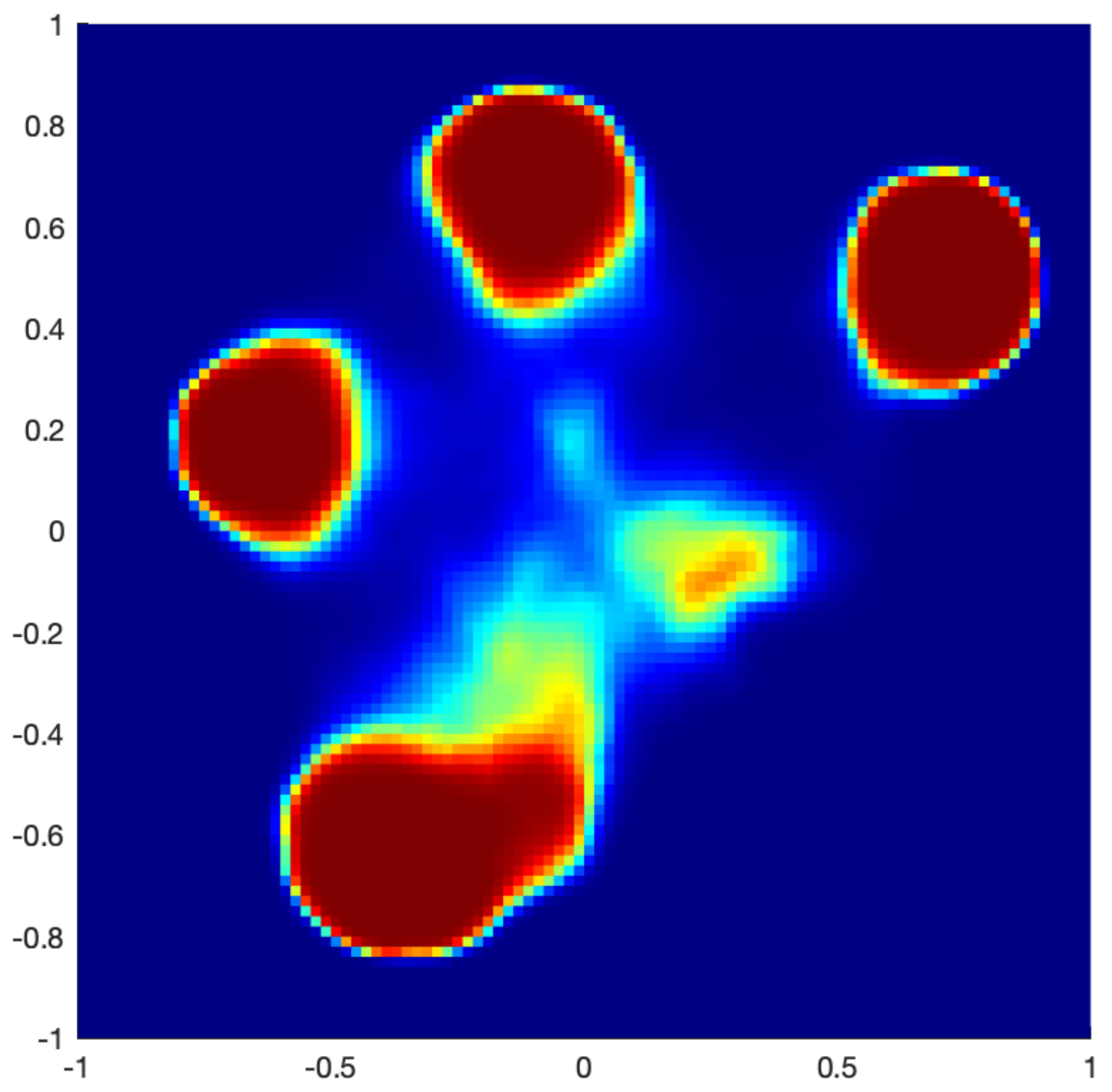}&
\includegraphics[width=1.1in]{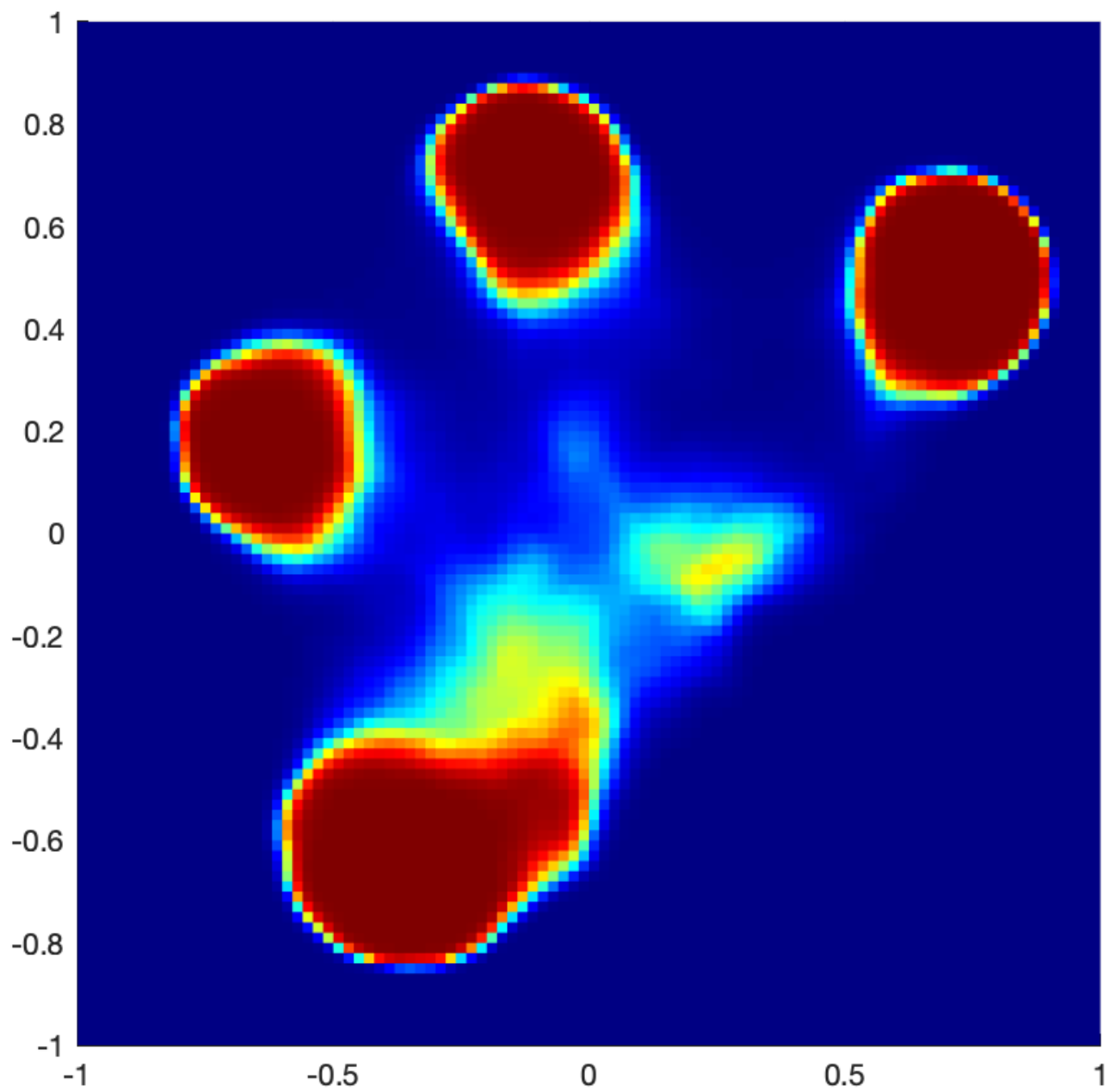}&
\includegraphics[width=1.1in]{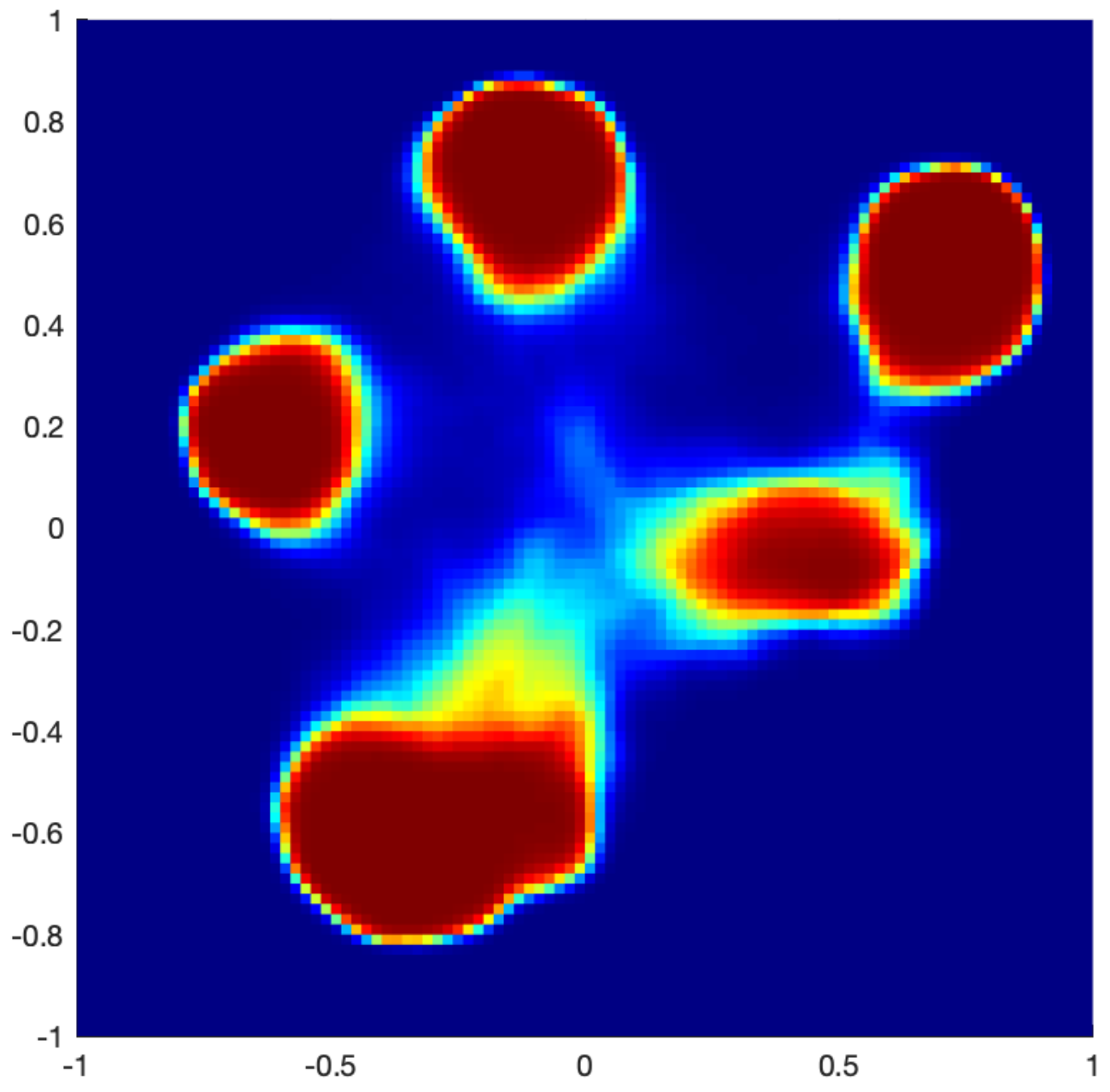}\\
\includegraphics[width=1.1in]{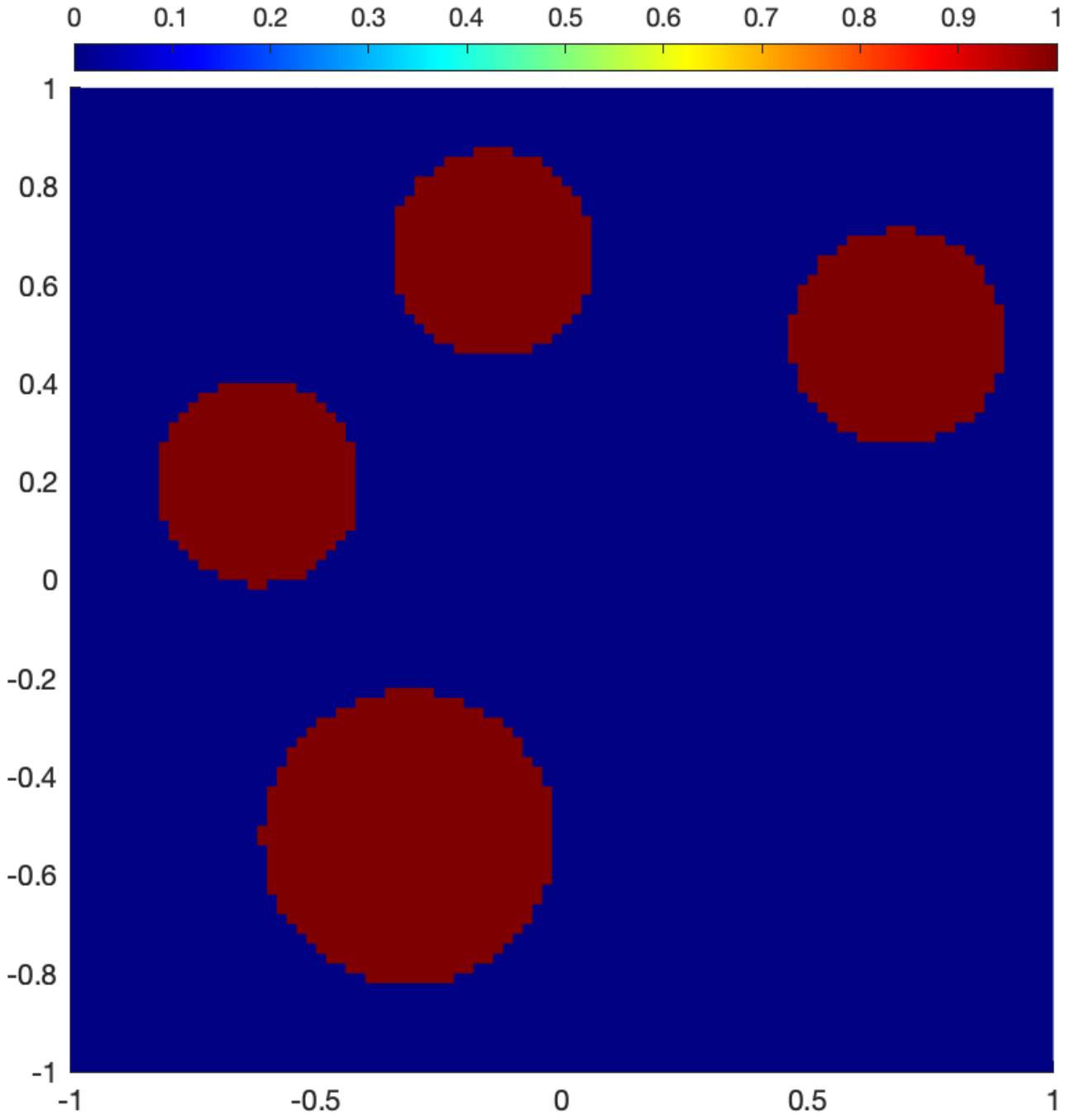}&
\includegraphics[width=1.1in]{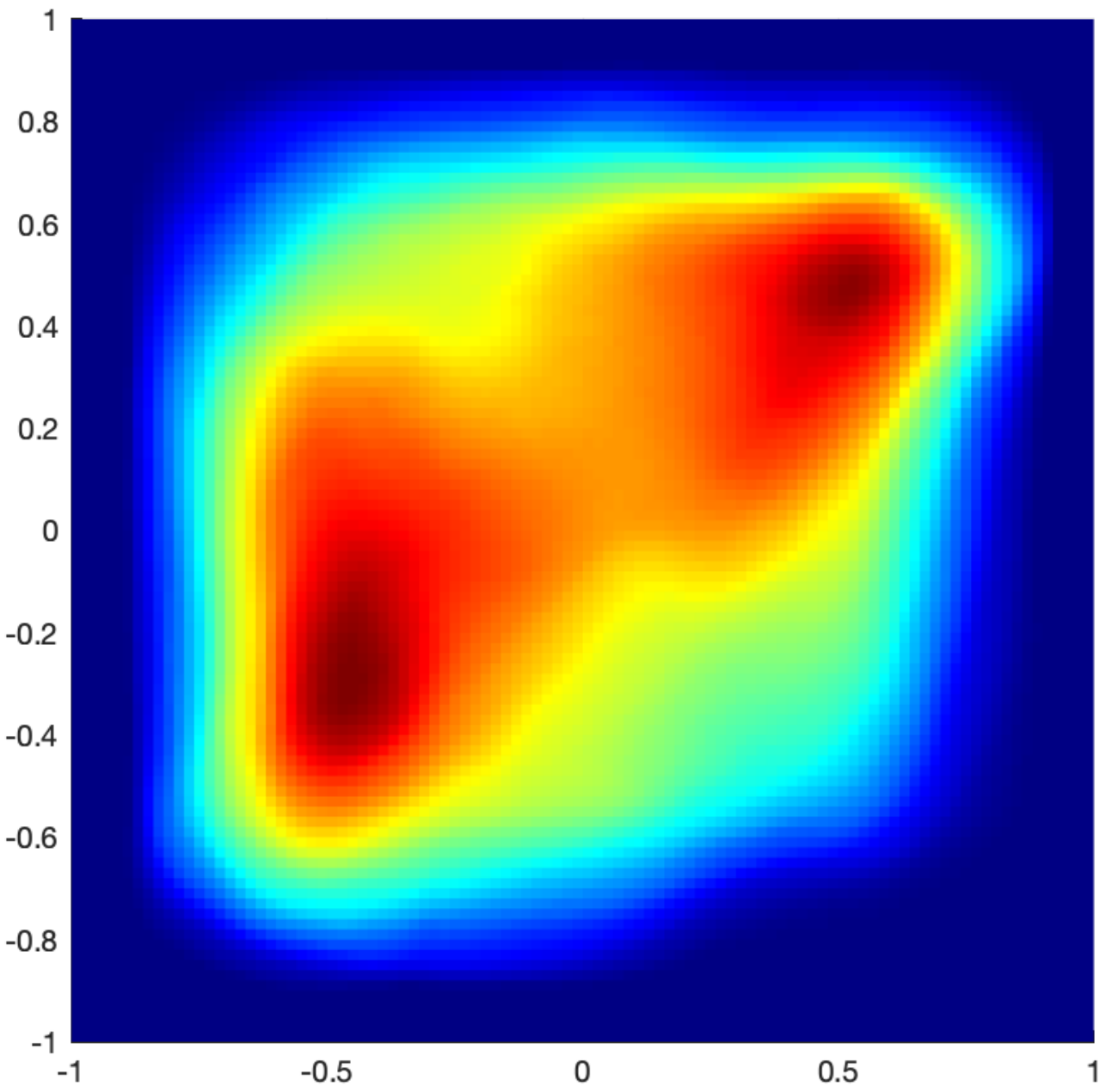}&
\includegraphics[width=1.1in]{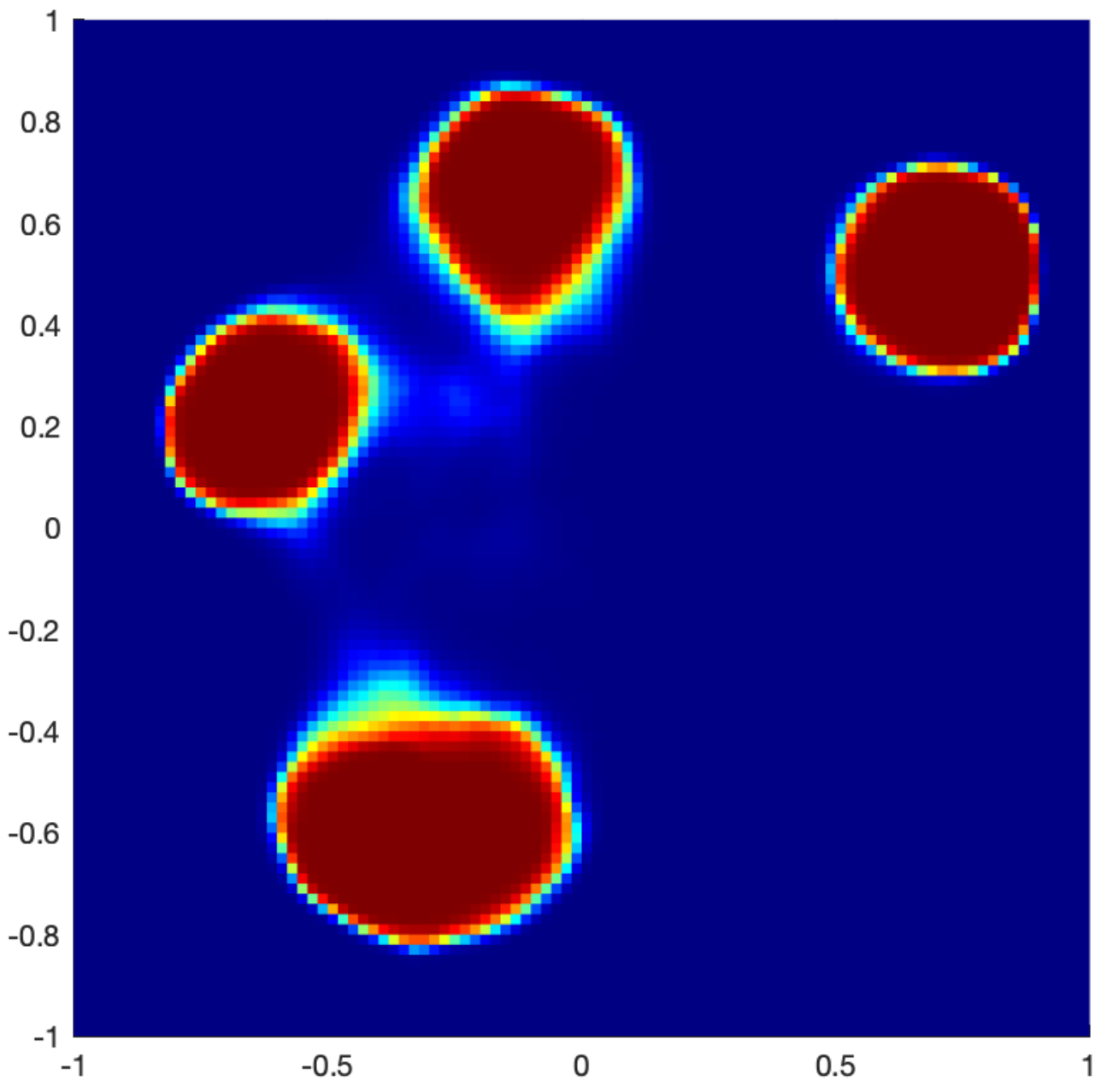}&
\includegraphics[width=1.1in]{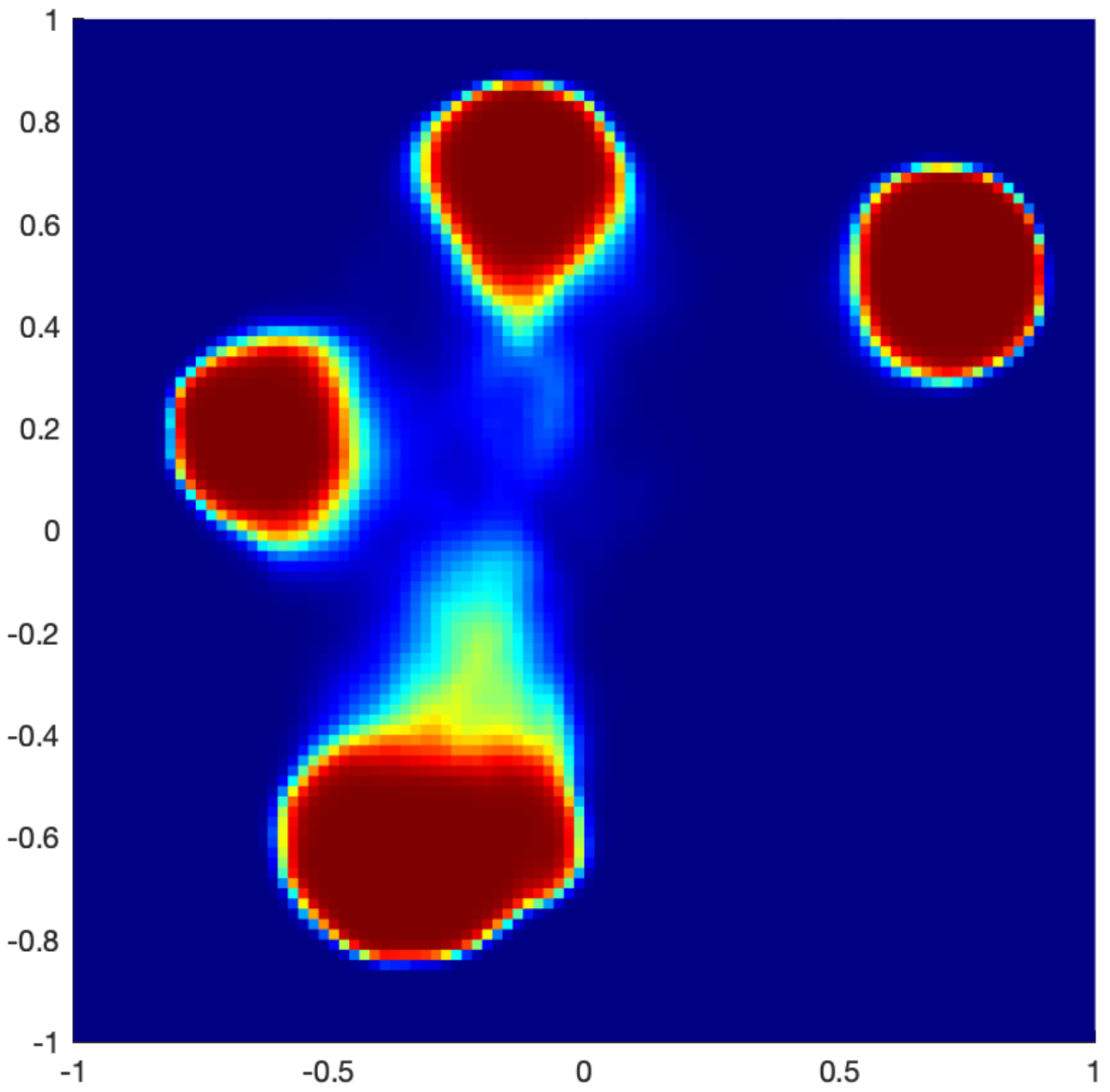}&
\includegraphics[width=1.1in]{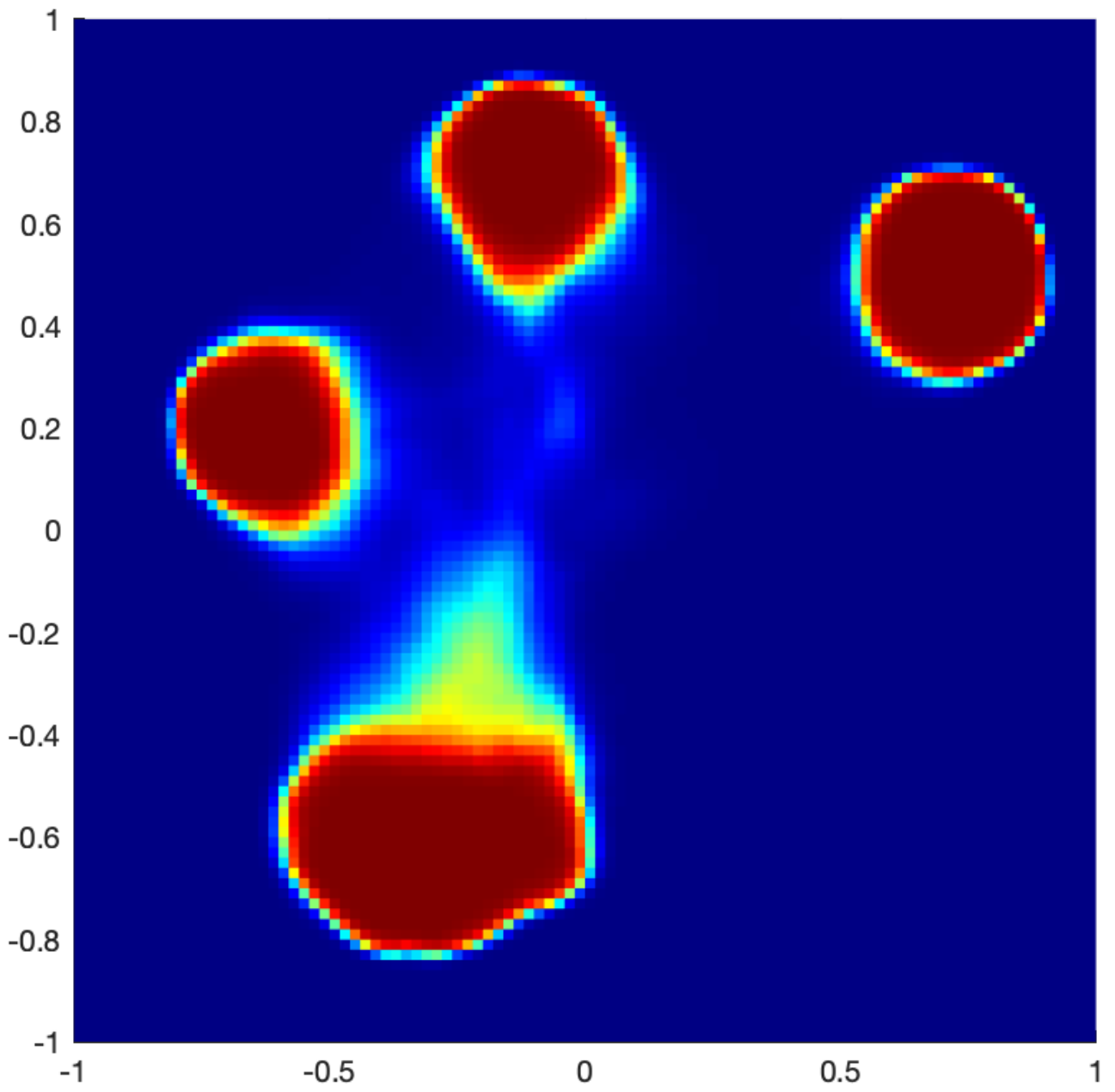}&
\includegraphics[width=1.1in]{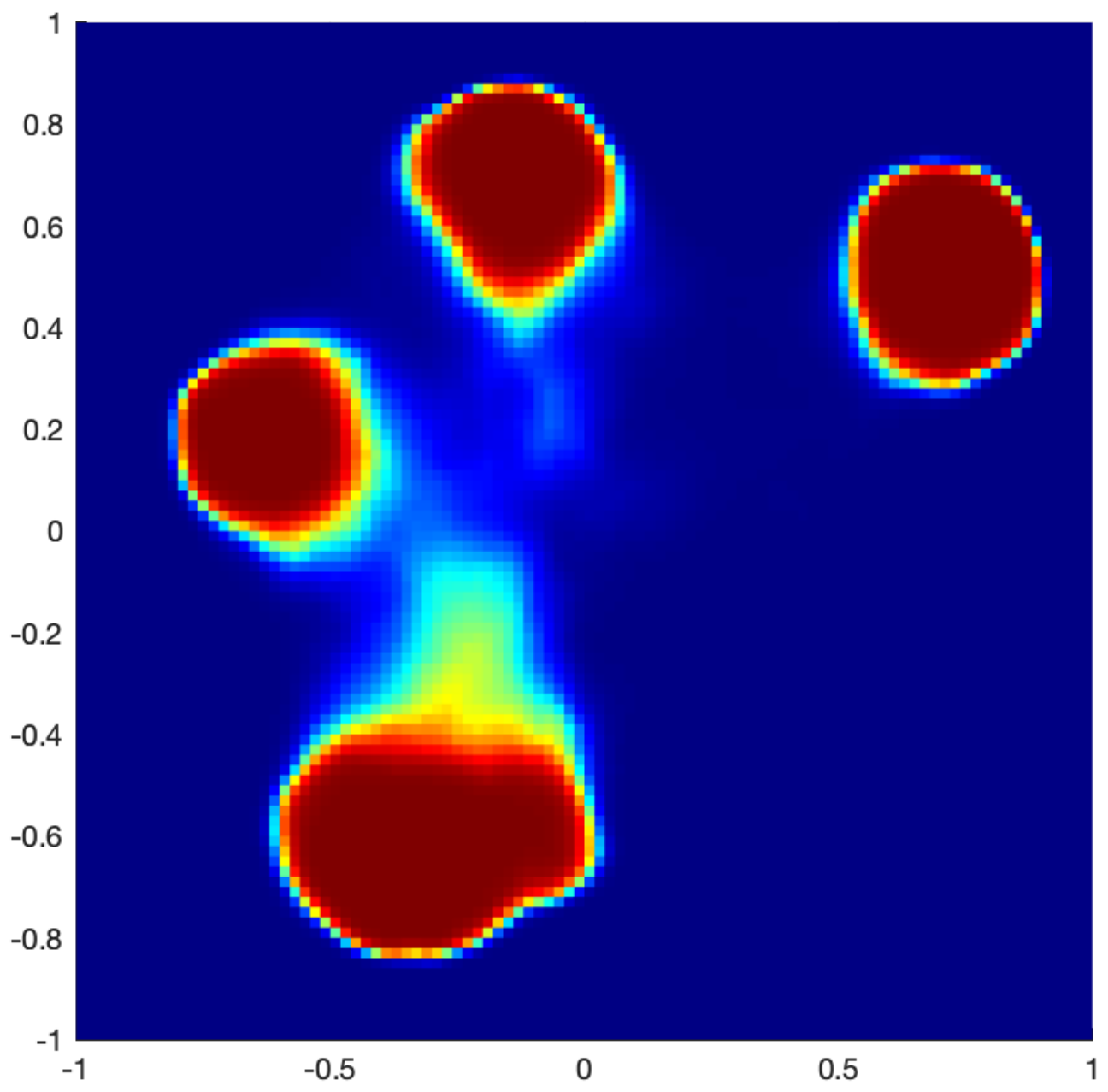}\\
\end{tabular}
  \caption{FNN-DDSM reconstruction for random circles: one circle is located closed to the center of domain and blocked from the boundary by other 4 circles (top) and this center circle is removed (bottom)} 
  \label{tab_FN_cir_comp}
\end{figure}

\begin{figure}[htbp]
\begin{tabular}{ >{\centering\arraybackslash}m{0.9in} >{\centering\arraybackslash}m{0.9in} >{\centering\arraybackslash}m{0.9in}  >{\centering\arraybackslash}m{0.9in}  >{\centering\arraybackslash}m{0.9in}  >{\centering\arraybackslash}m{0.9in} }
\centering
True coefficients &
N=1, $\delta=0$&
N=10, $\delta=0$&
N=20, $\delta=0$&
N=20, $\delta=10\%$ &
N=20, $\delta=20\%$ \\
\includegraphics[width=1.1in]{cir_comp1_0-eps-converted-to.pdf}&
\includegraphics[width=1.1in]{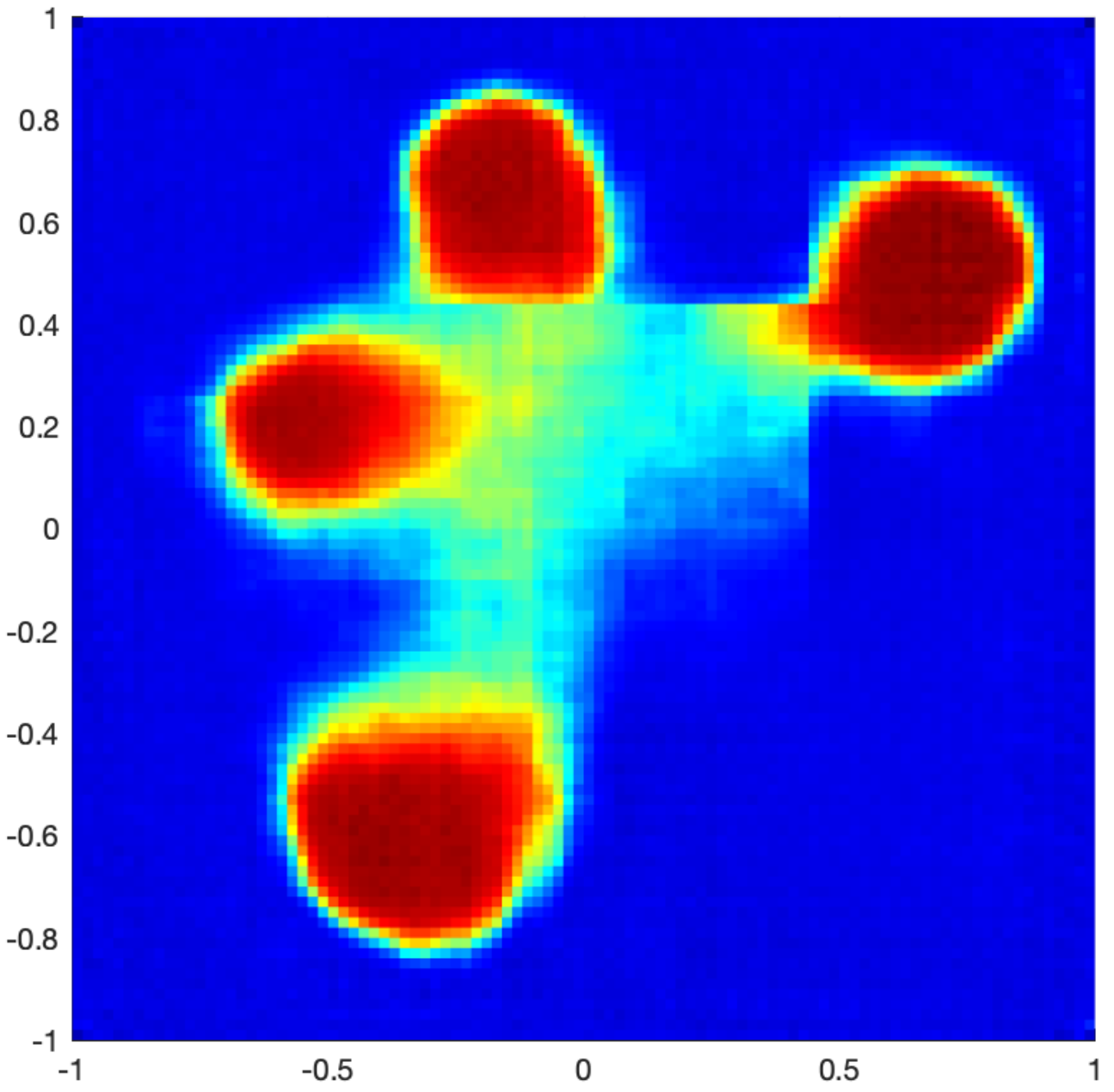}&
\includegraphics[width=1.1in]{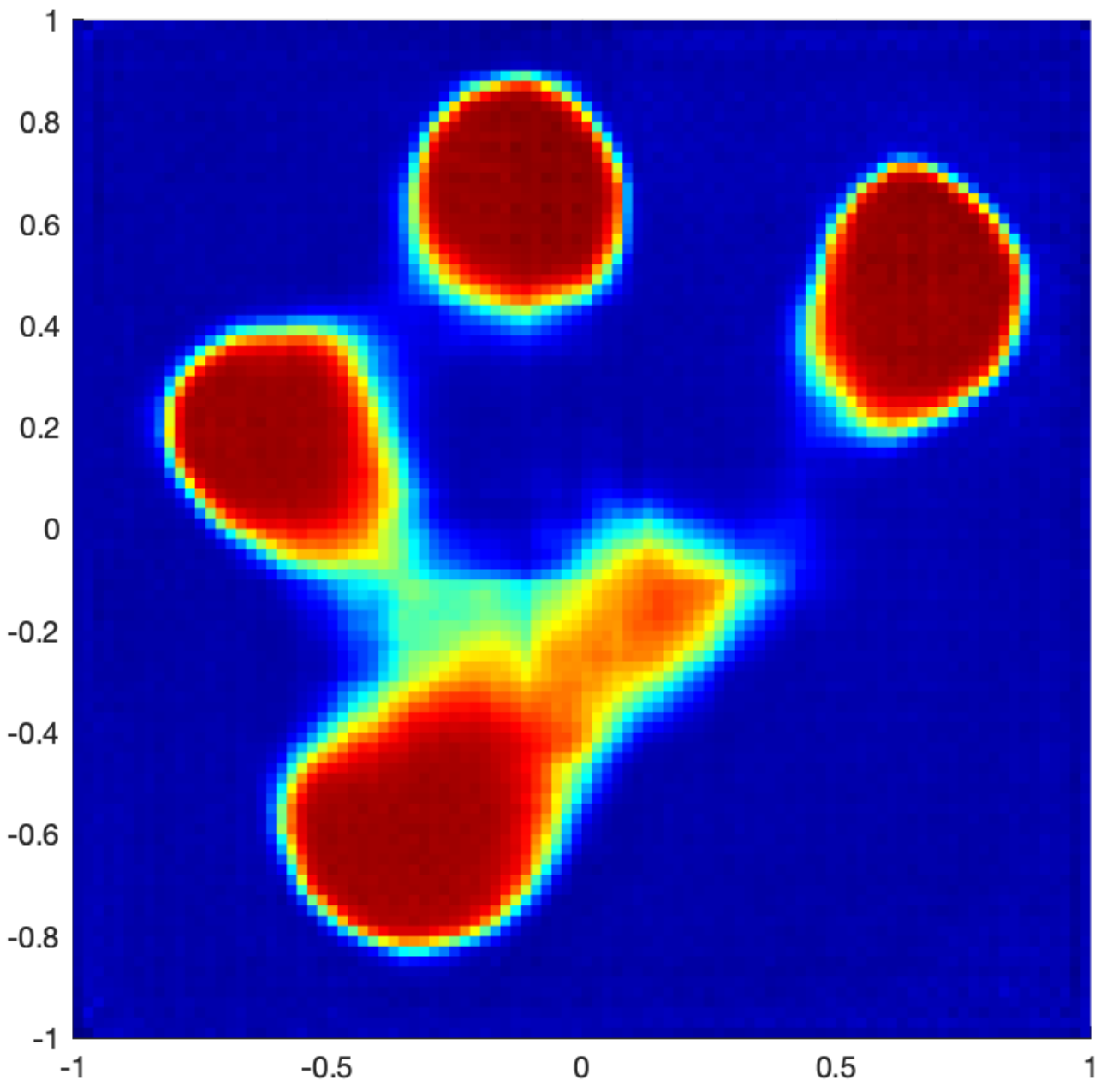}&
\includegraphics[width=1.1in]{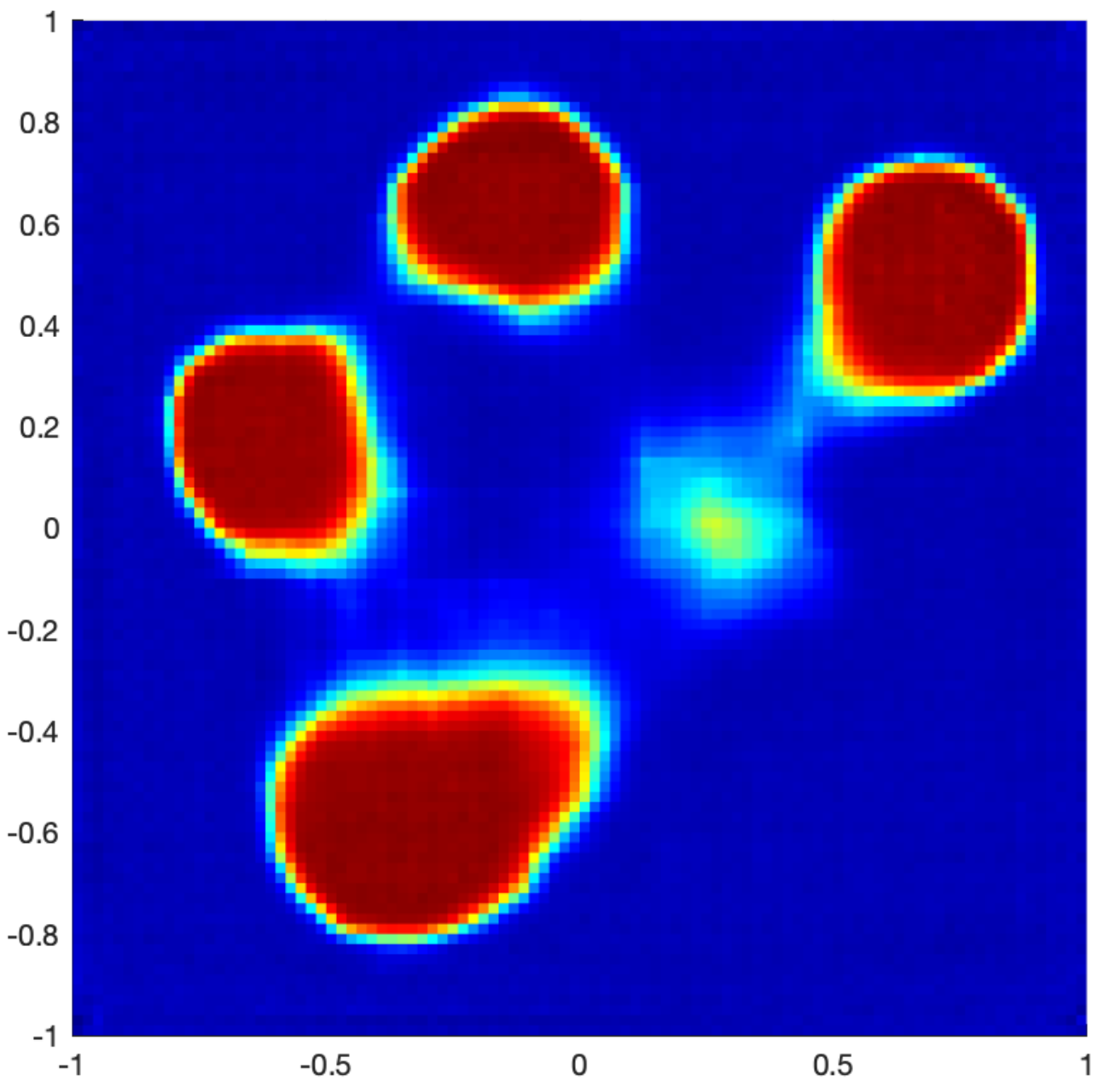}&
\includegraphics[width=1.1in]{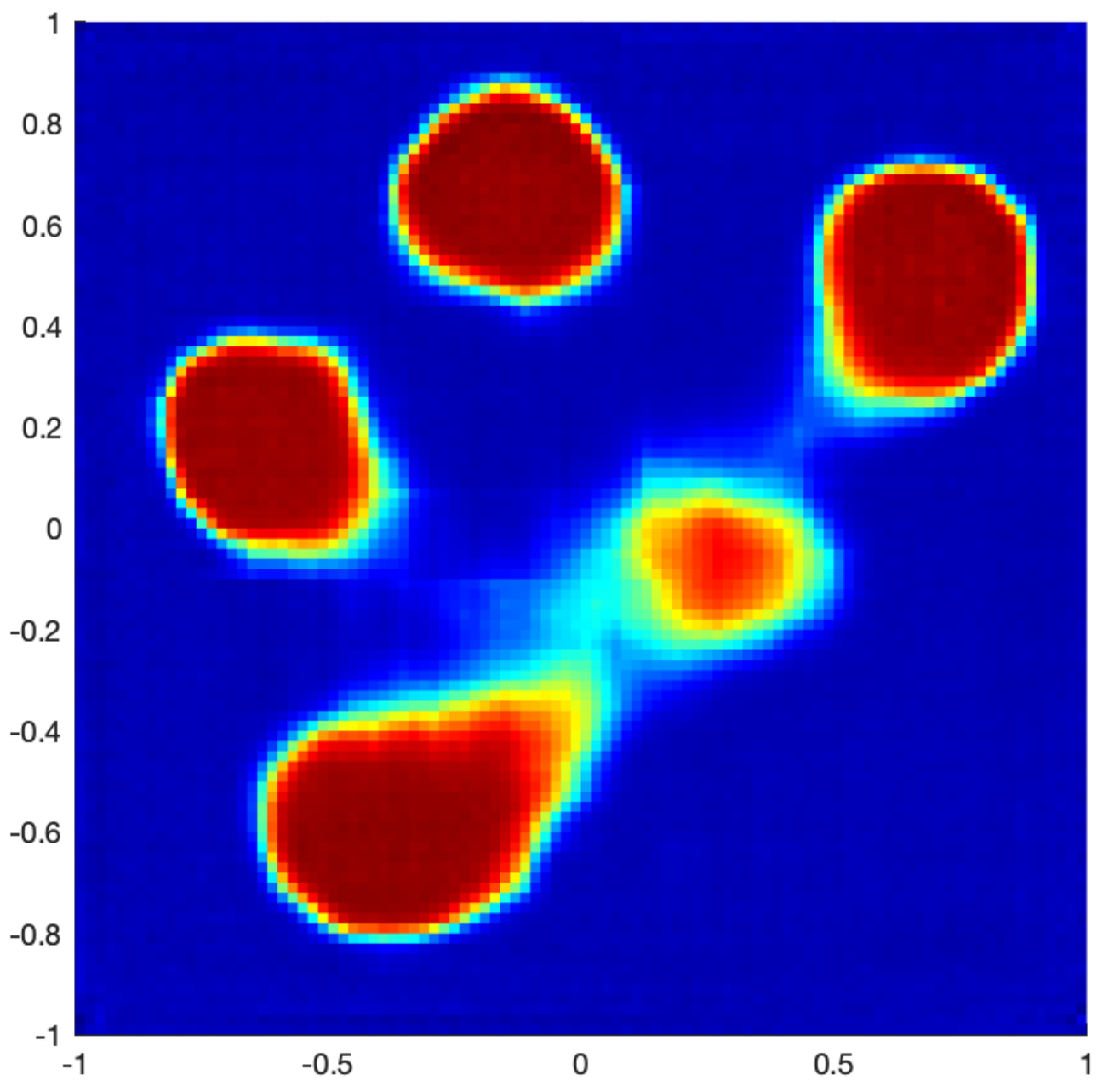}&
\includegraphics[width=1.1in]{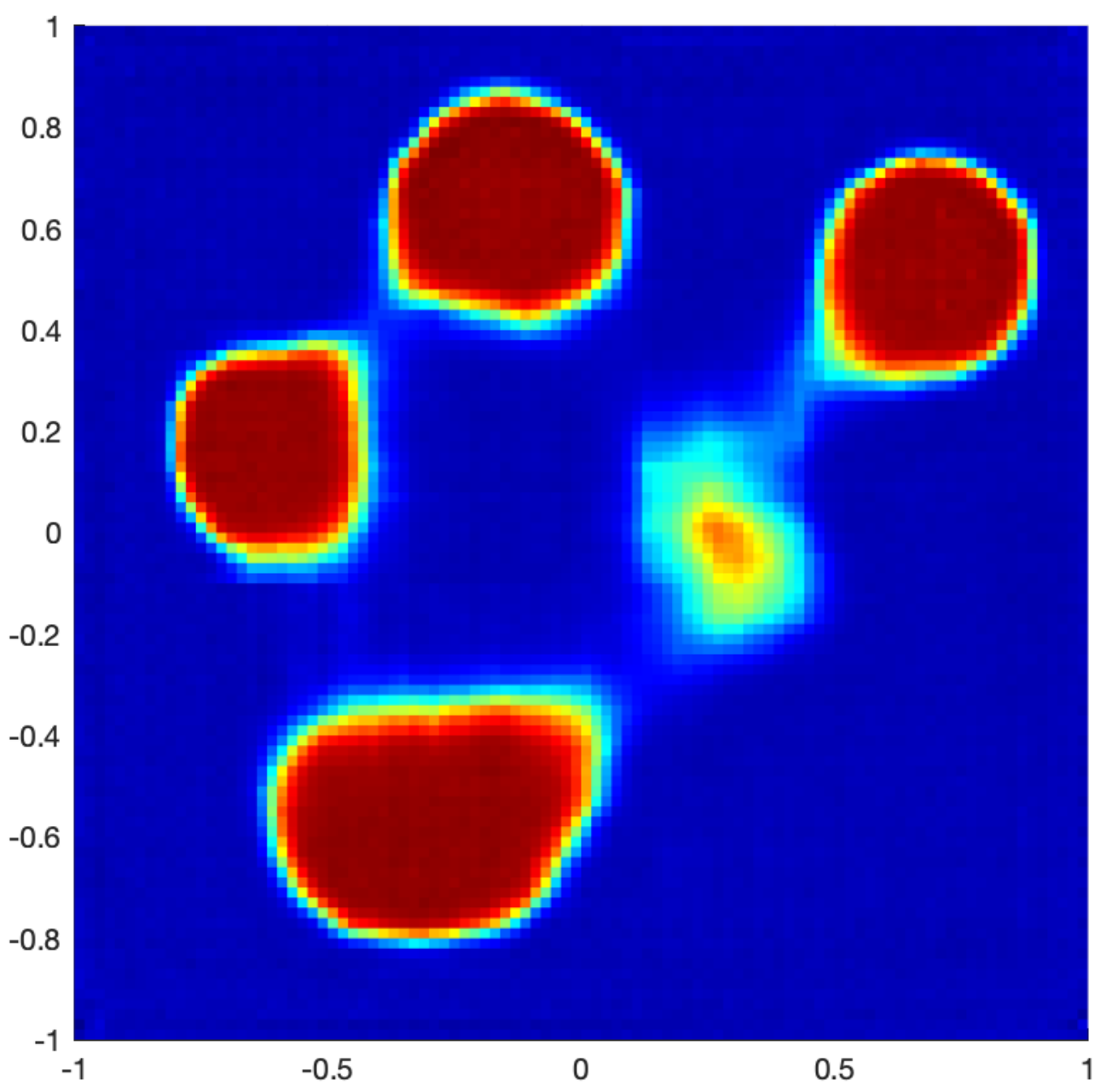}\\
\includegraphics[width=1.1in]{cir_comp2_0-eps-converted-to.pdf}&
\includegraphics[width=1.1in]{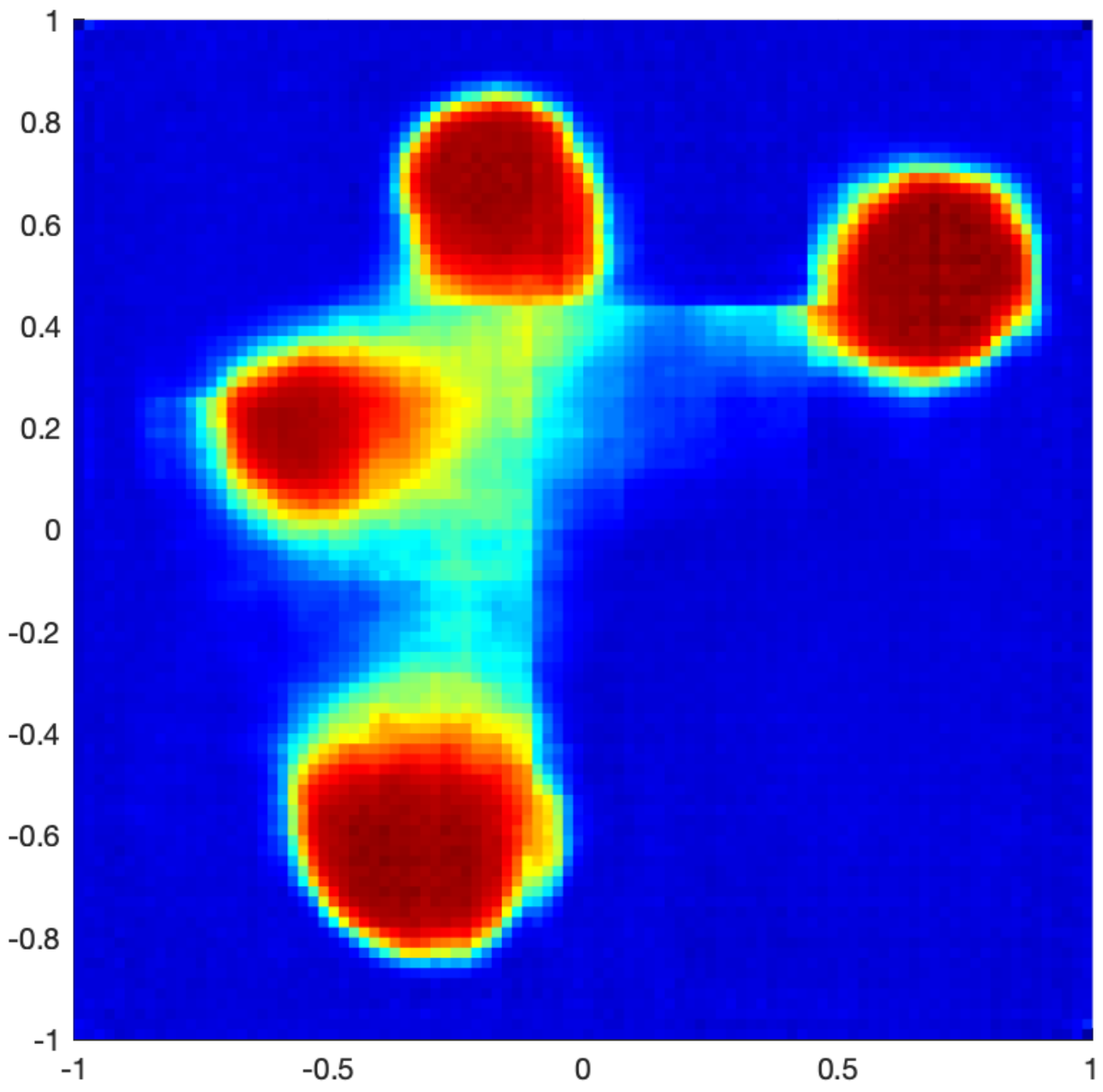}&
\includegraphics[width=1.1in]{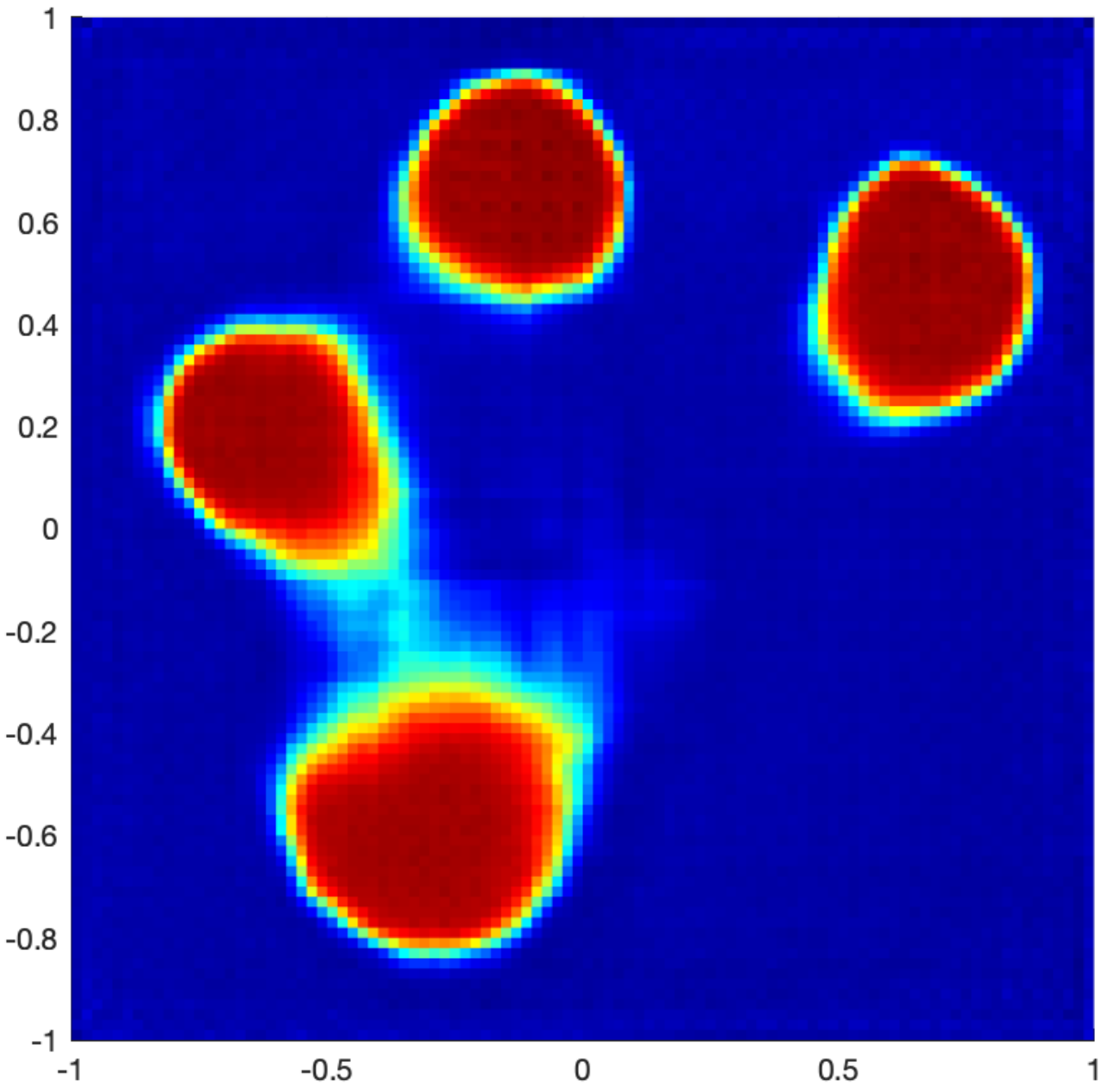}&
\includegraphics[width=1.1in]{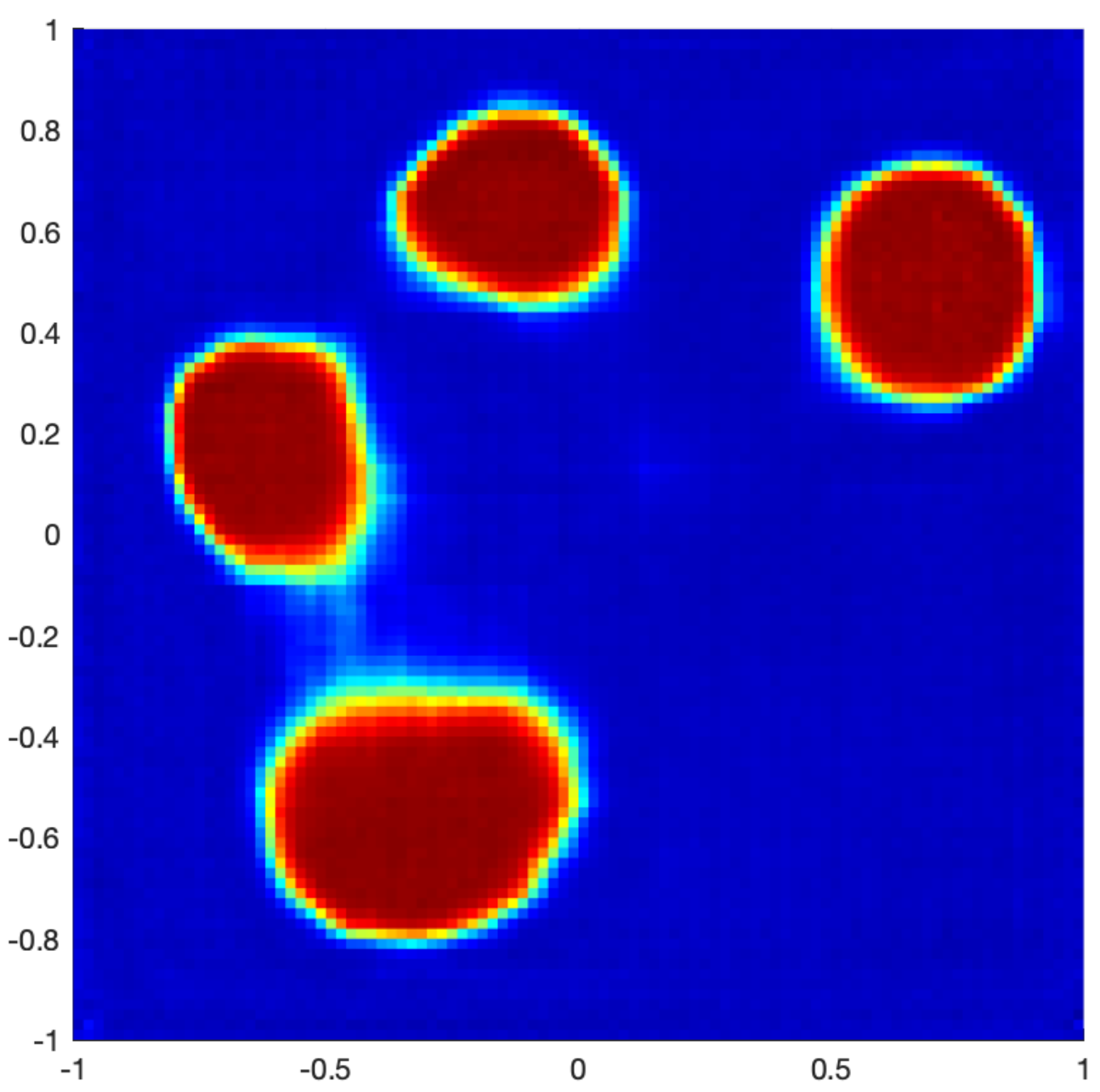}&
\includegraphics[width=1.1in]{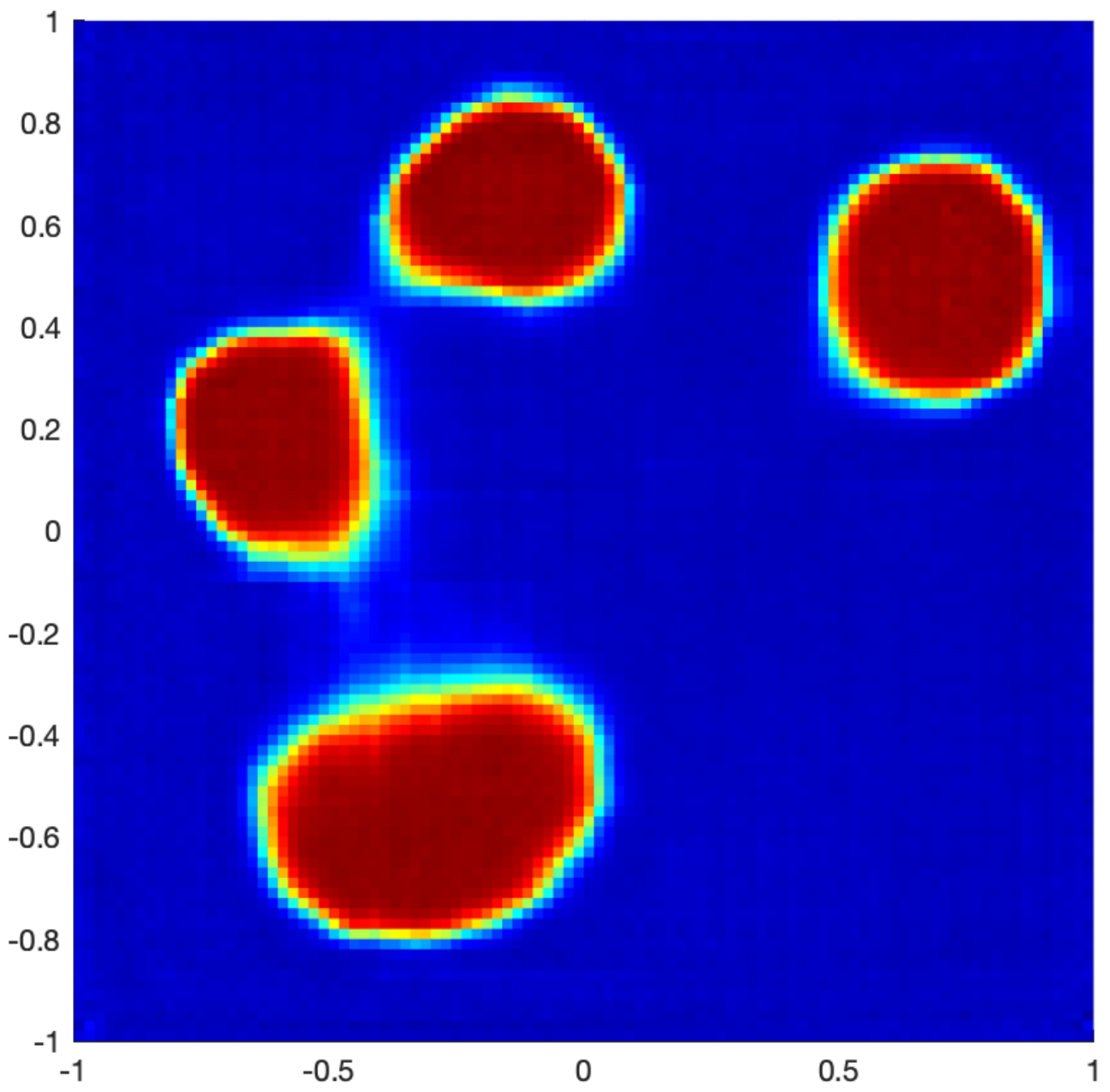}&
\includegraphics[width=1.1in]{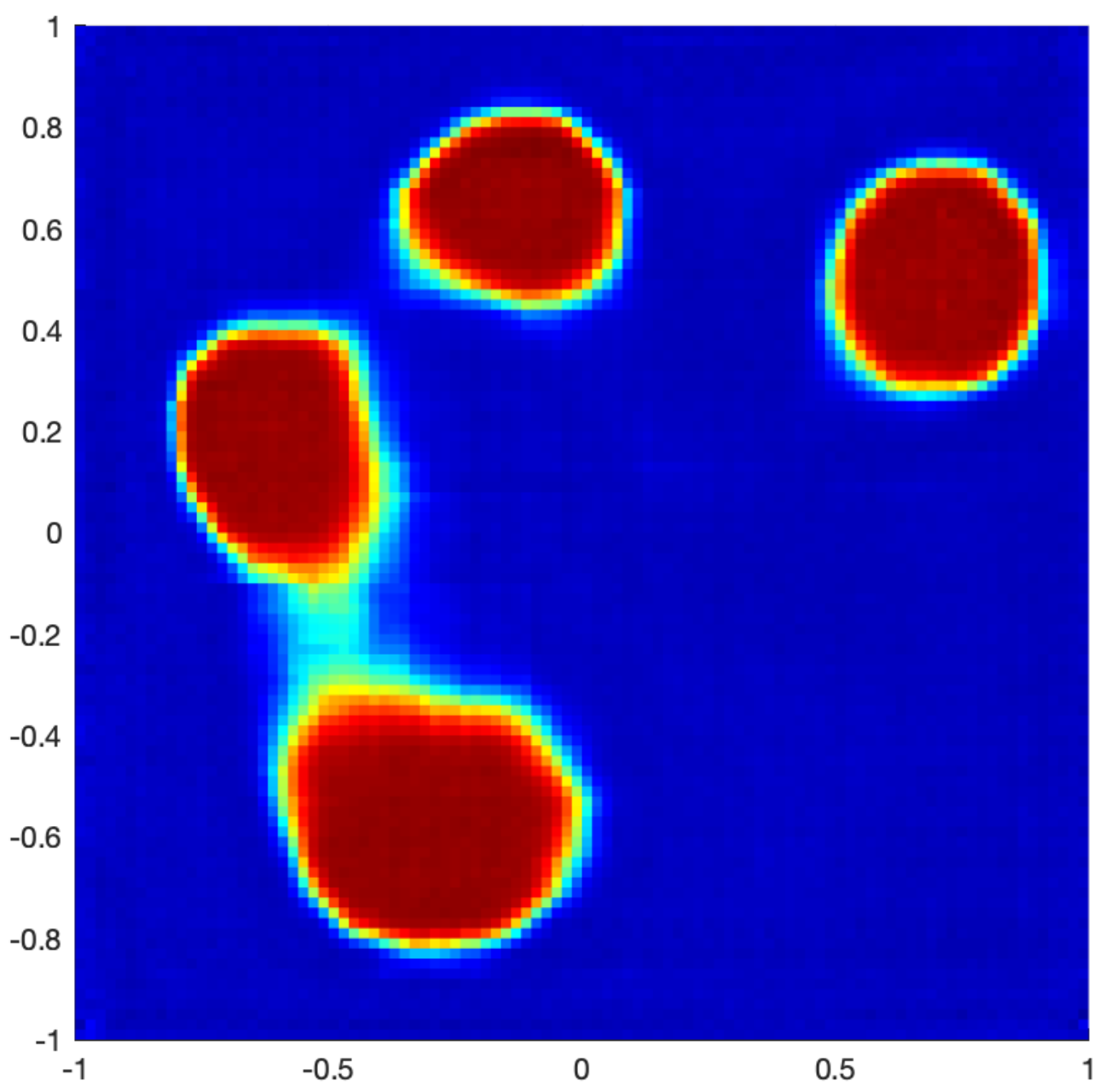}\\
\end{tabular}
  \caption{CNN-DDSM reconstruction for random circles: one circle is located closed to the center of domain and blocked from the boundary by other 4 circles (top) and this center circle is removed (bottom)} 
  \label{tab_cir_comp}
\end{figure}

\begin{figure}[htbp]
\begin{tabular}{ >{\centering\arraybackslash}m{0.9in} >{\centering\arraybackslash}m{0.9in} >{\centering\arraybackslash}m{0.9in}  >{\centering\arraybackslash}m{0.9in}  >{\centering\arraybackslash}m{0.9in}  >{\centering\arraybackslash}m{0.9in} }
\centering
True coefficients &
N=1, $\delta=0$&
N=10, $\delta=0$&
N=20, $\delta=0$&
N=20, $\delta=10\%$ &
N=20, $\delta=20\%$ \\
\includegraphics[width=1.1in]{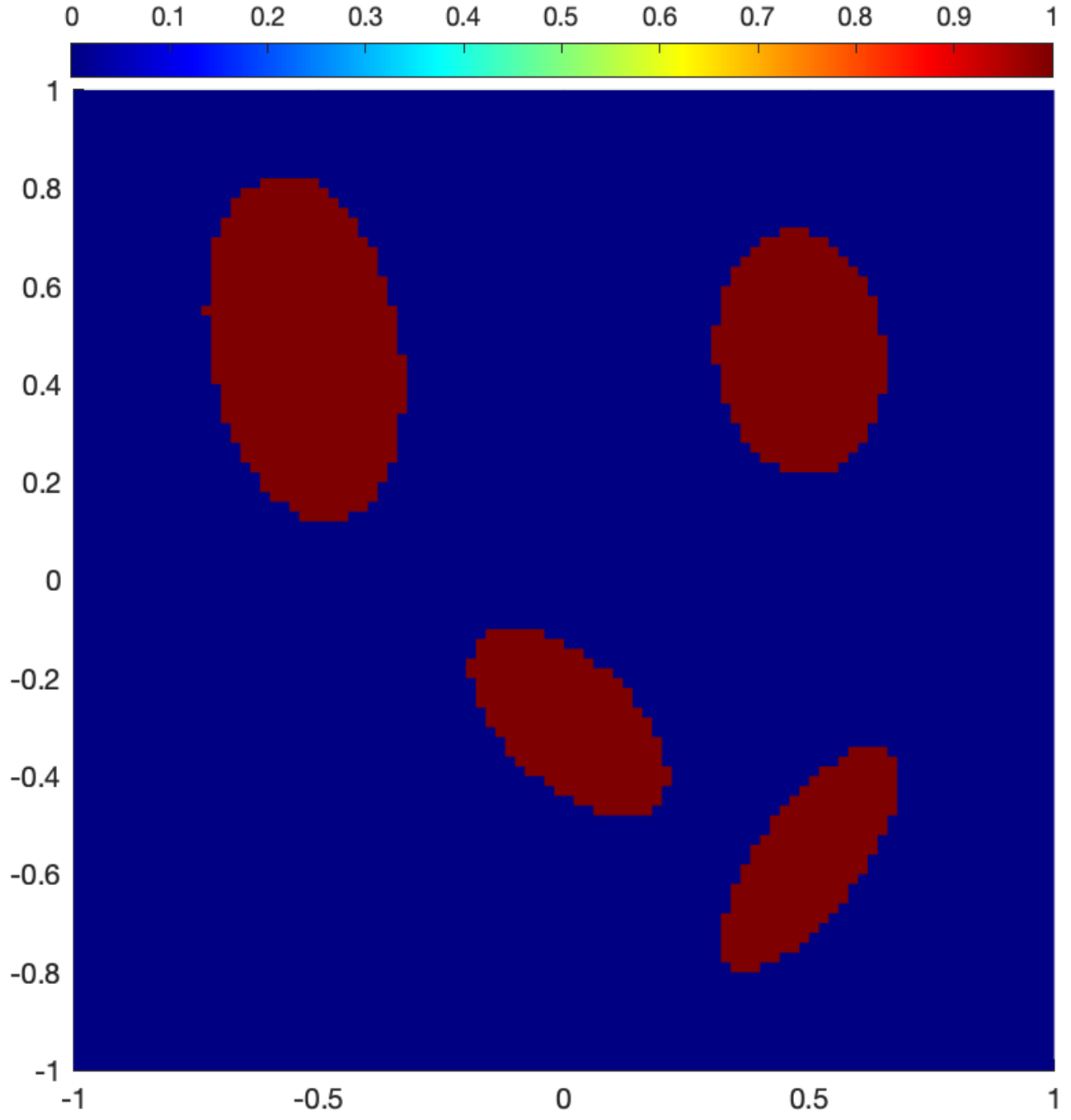}&
\includegraphics[width=1.1in]{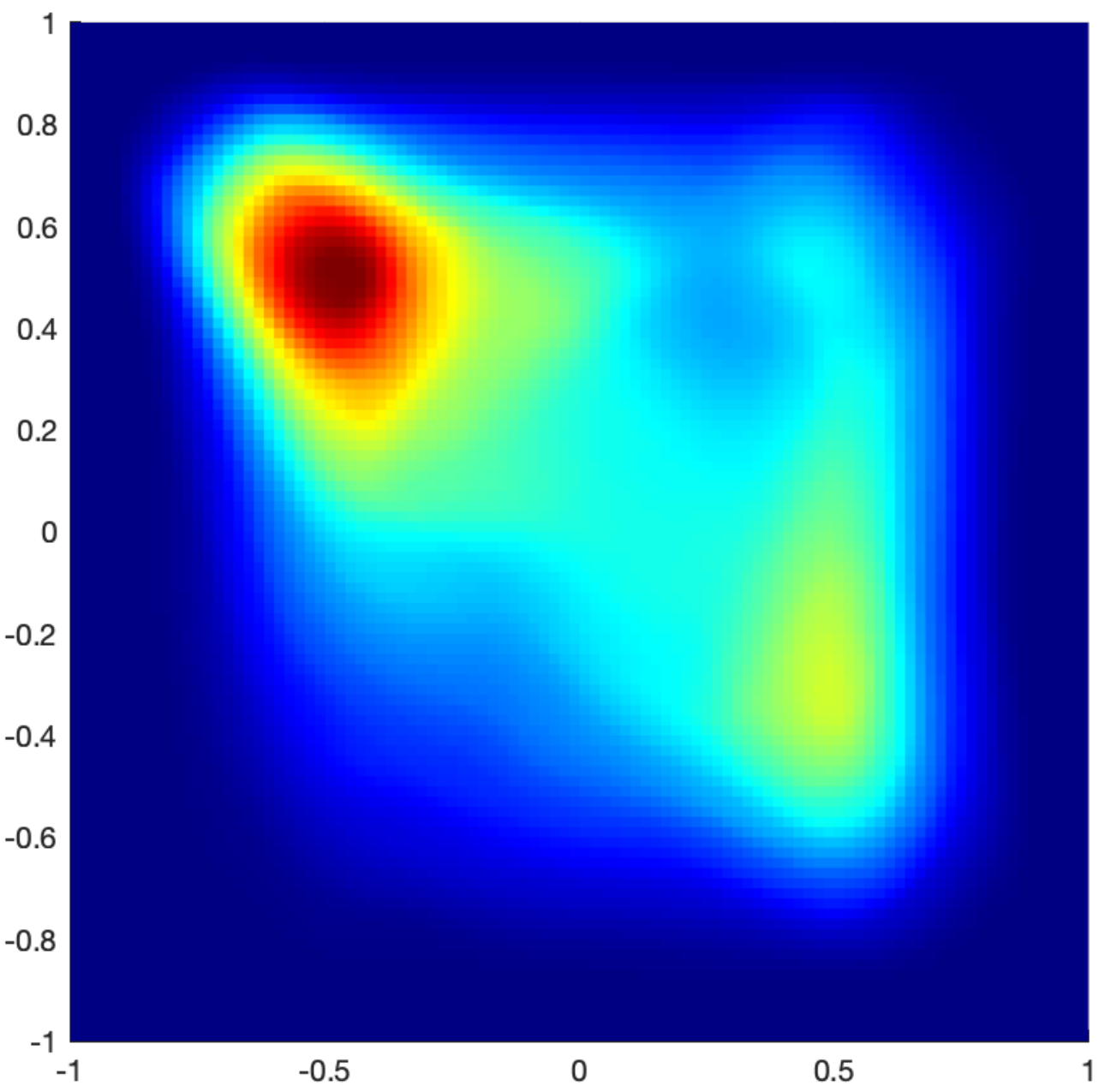}&
\includegraphics[width=1.1in]{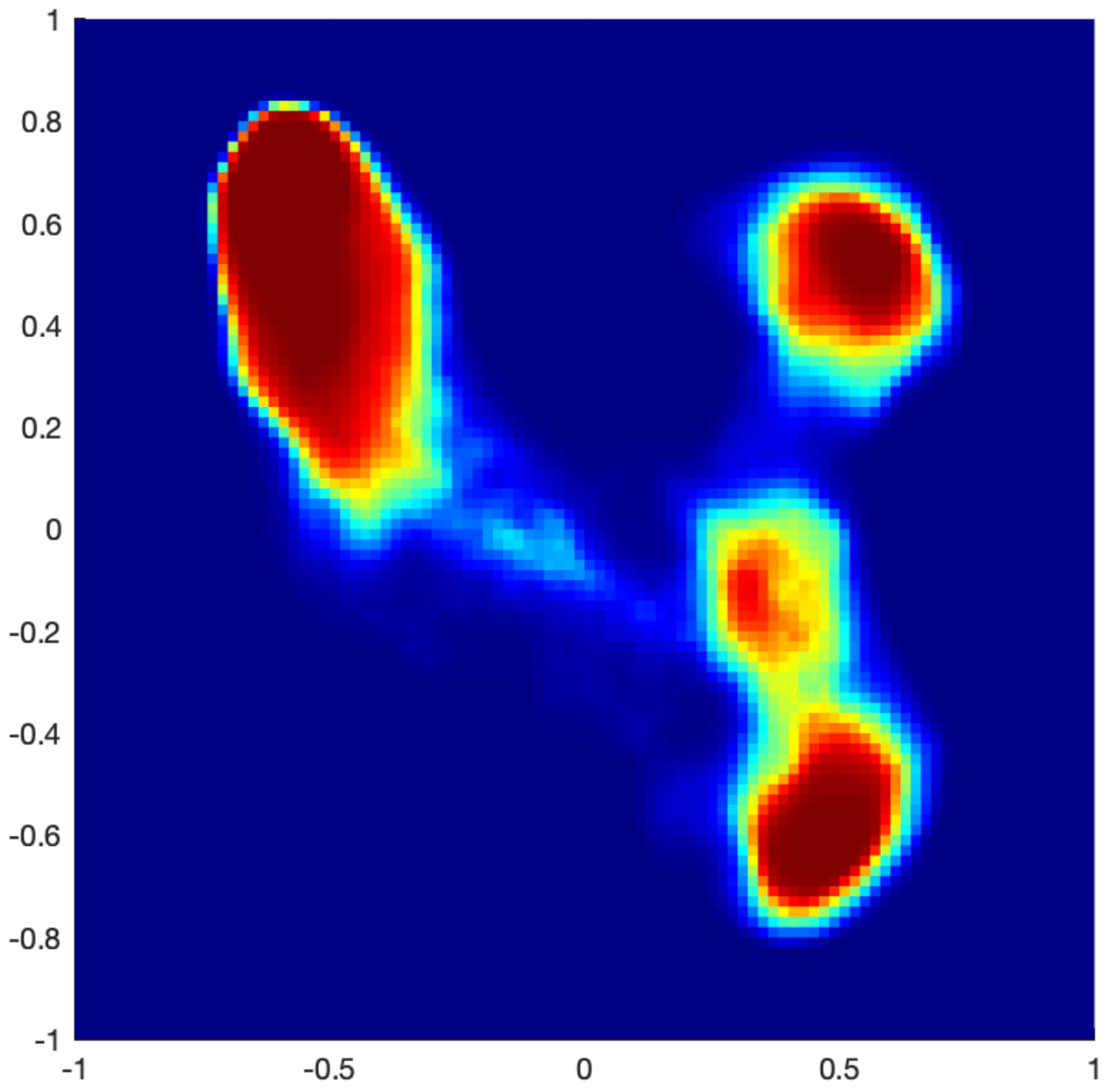}&
\includegraphics[width=1.1in]{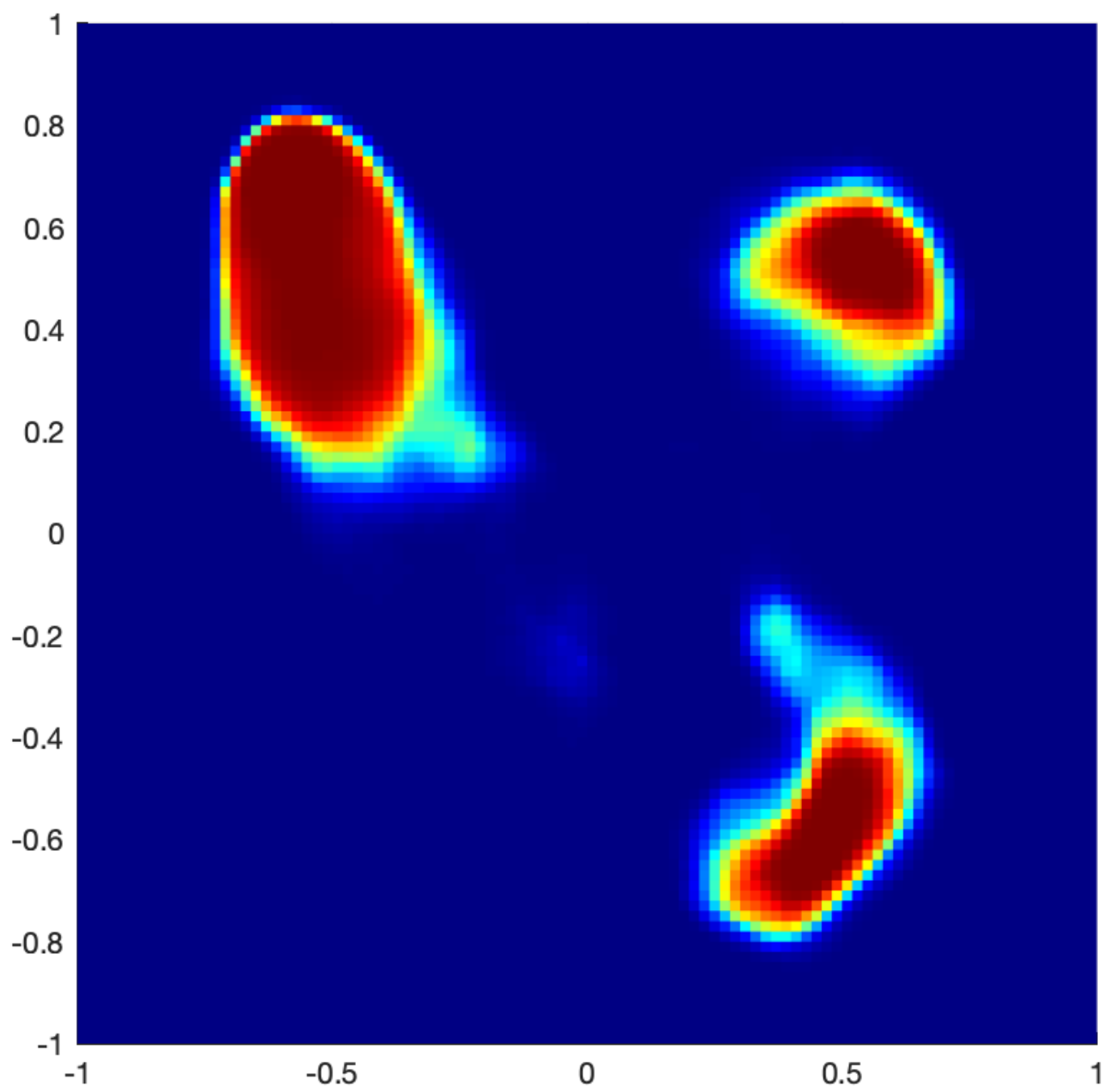}&
\includegraphics[width=1.1in]{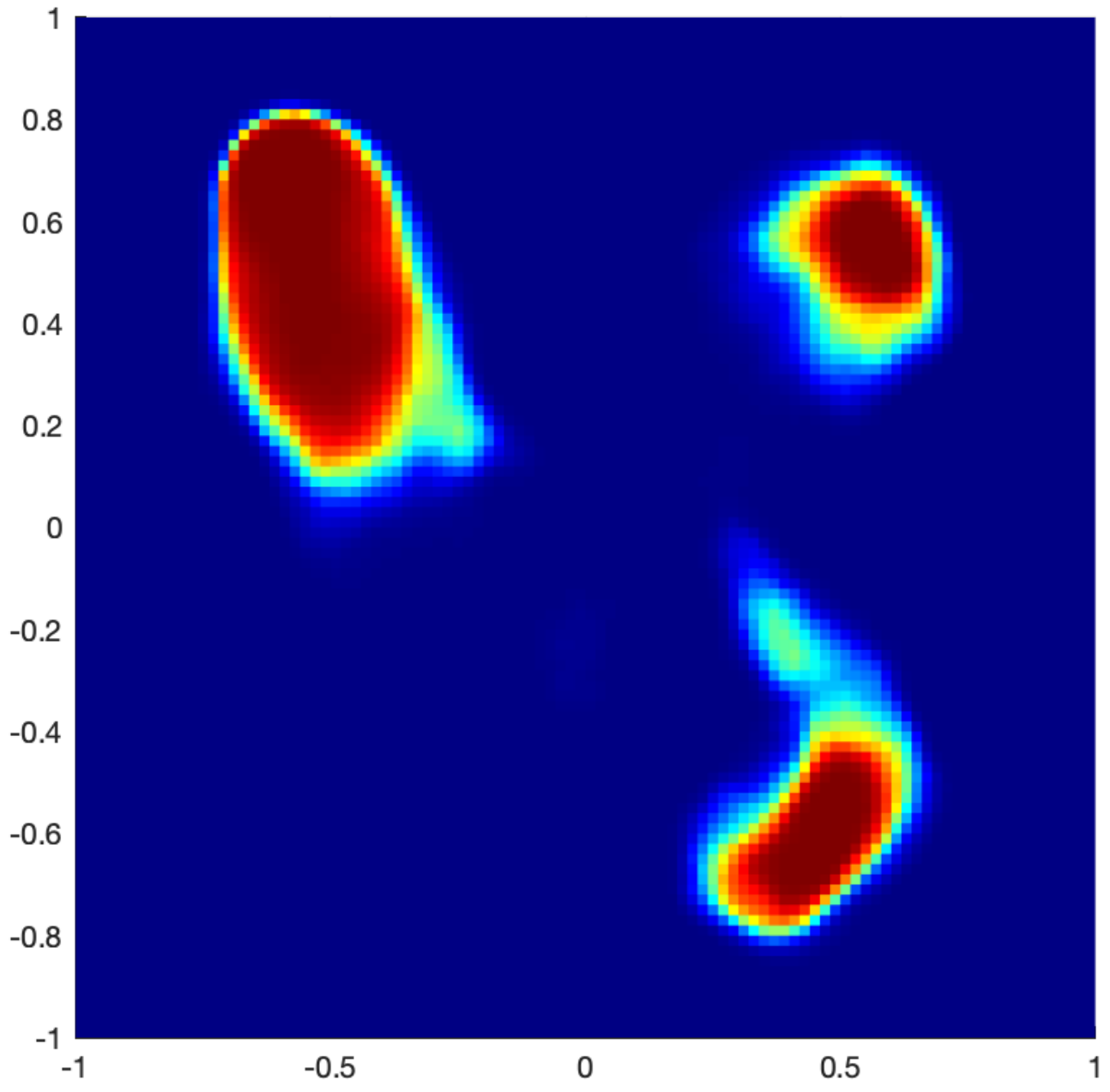}&
\includegraphics[width=1.1in]{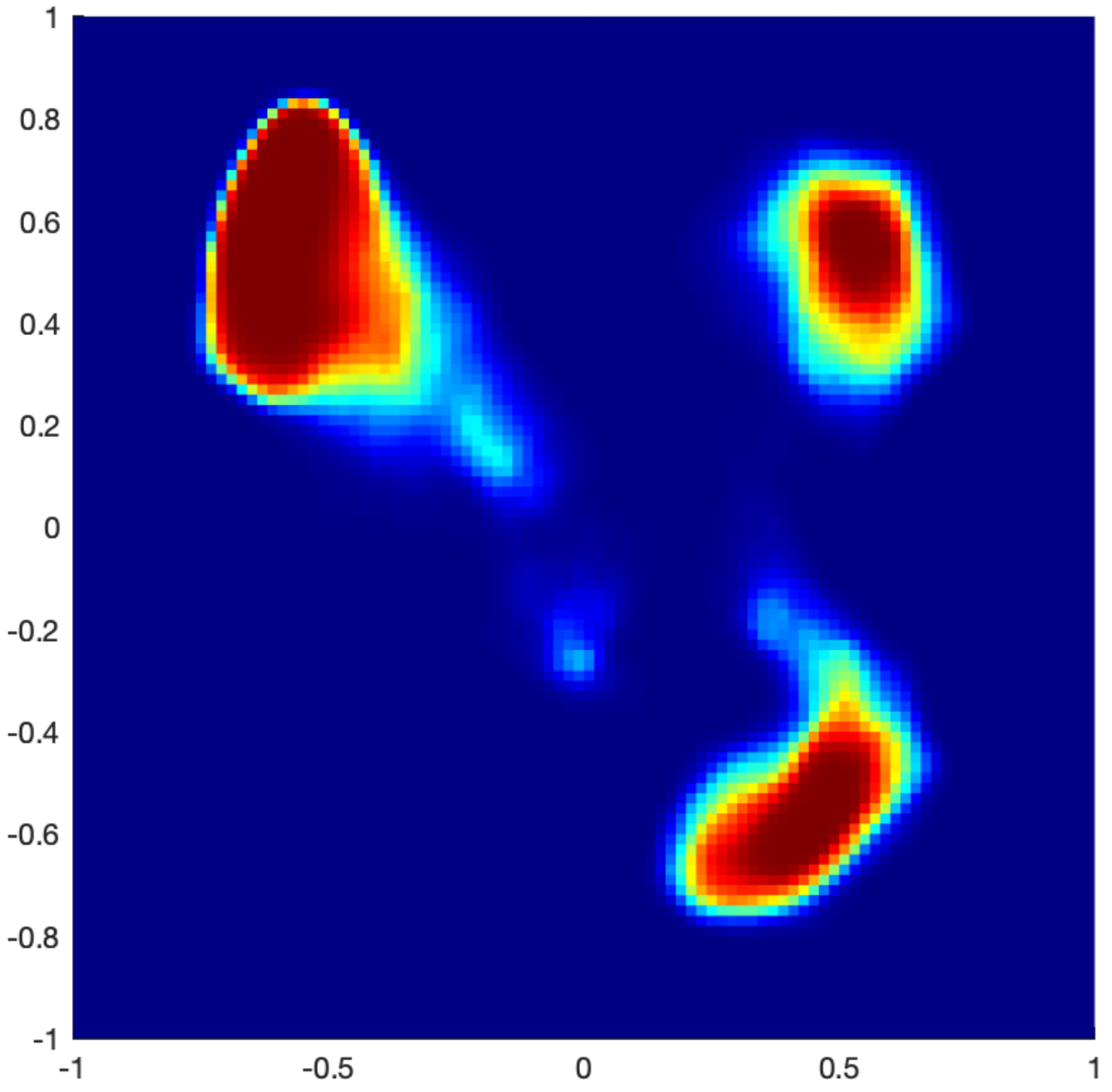}\\
\includegraphics[width=1.1in]{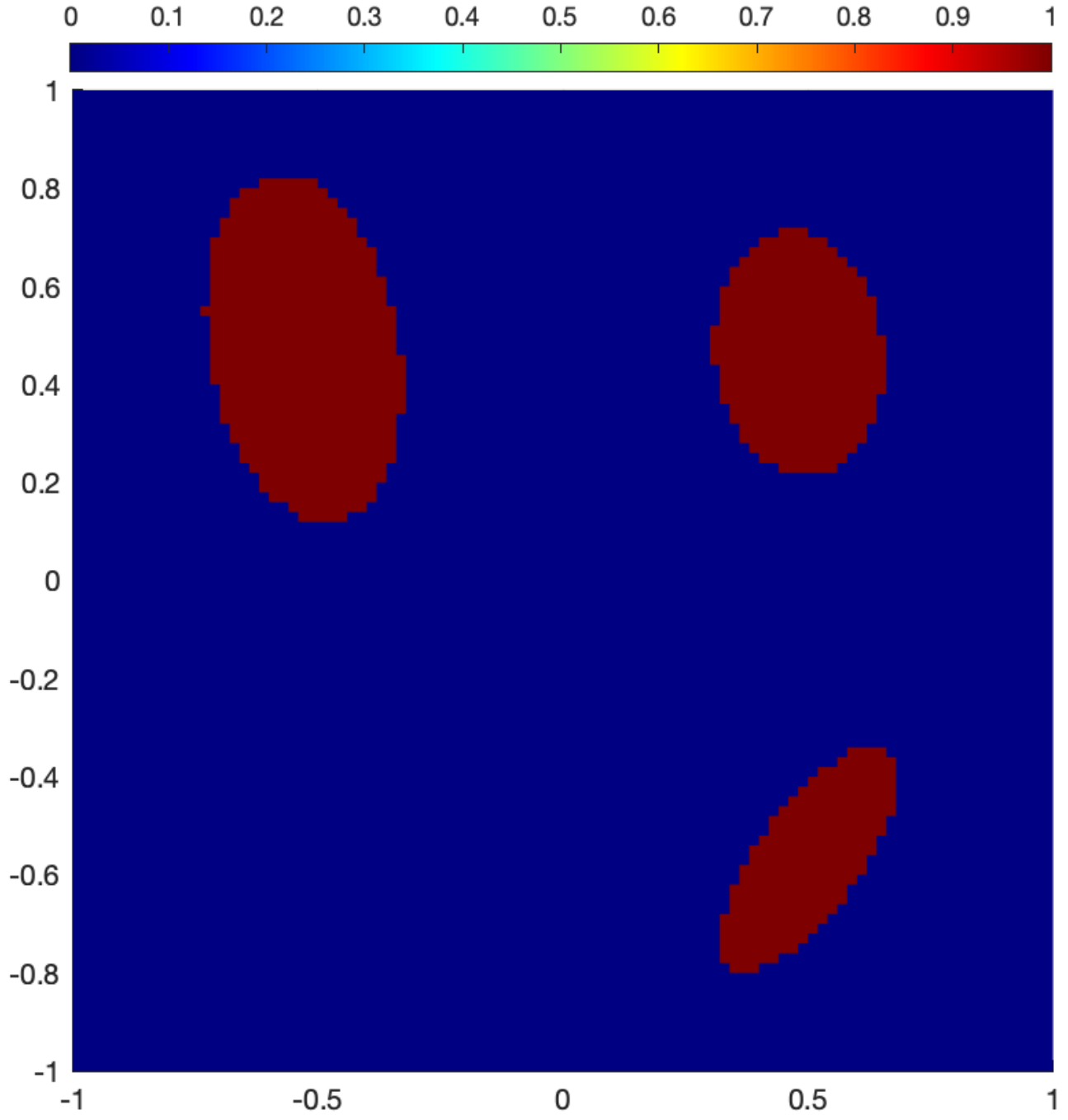}&
\includegraphics[width=1.1in]{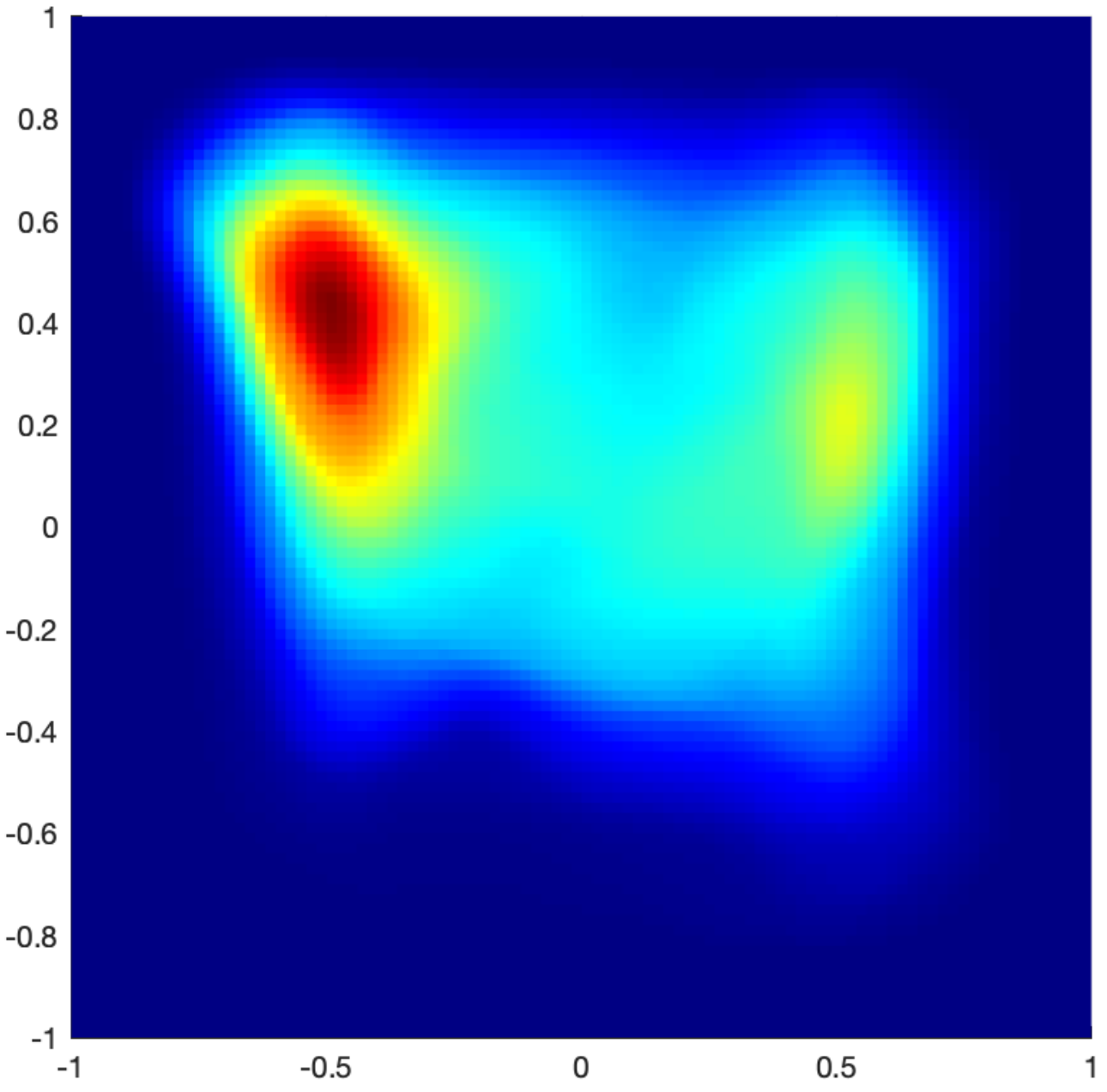}&
\includegraphics[width=1.1in]{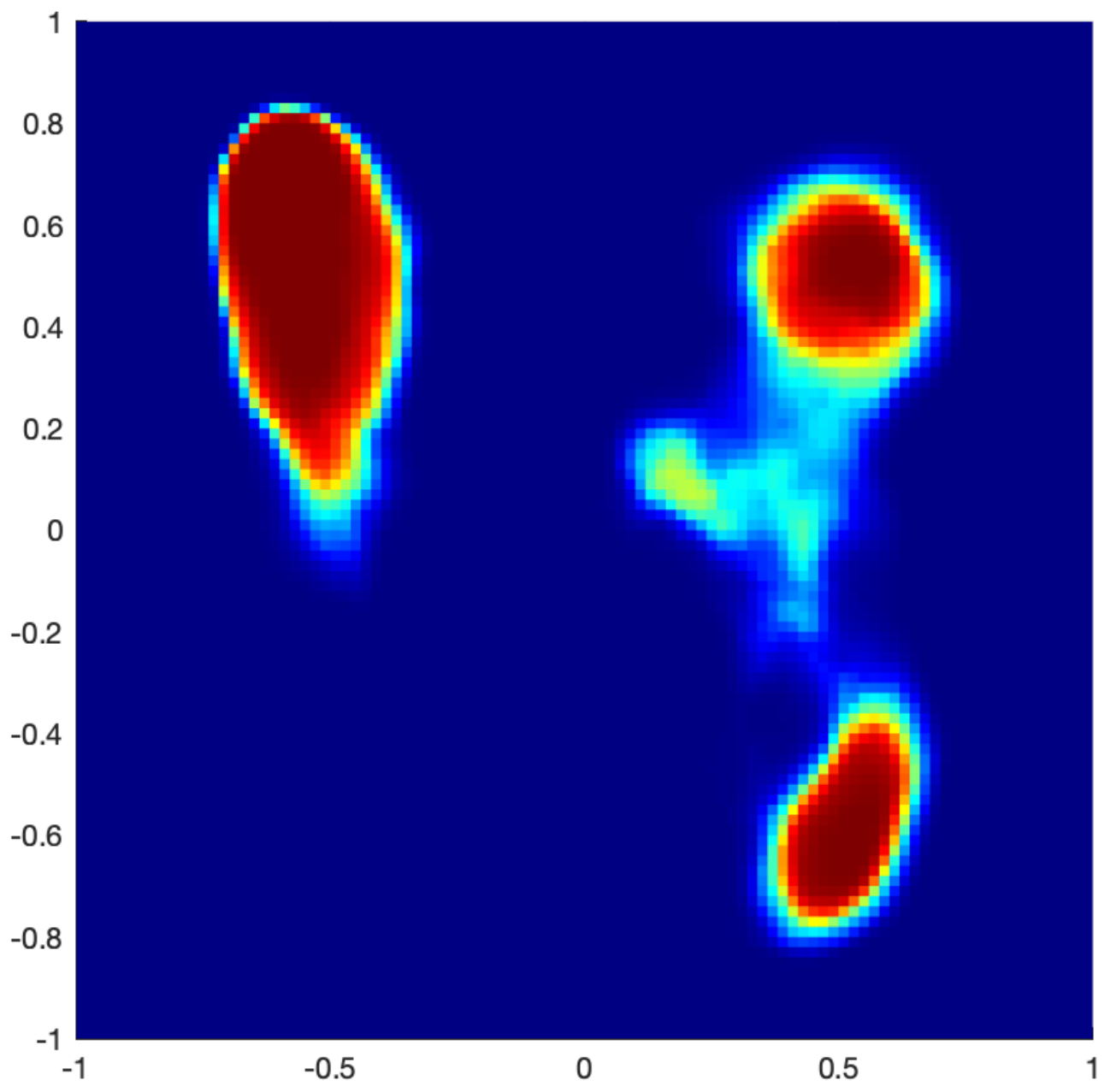}&
\includegraphics[width=1.1in]{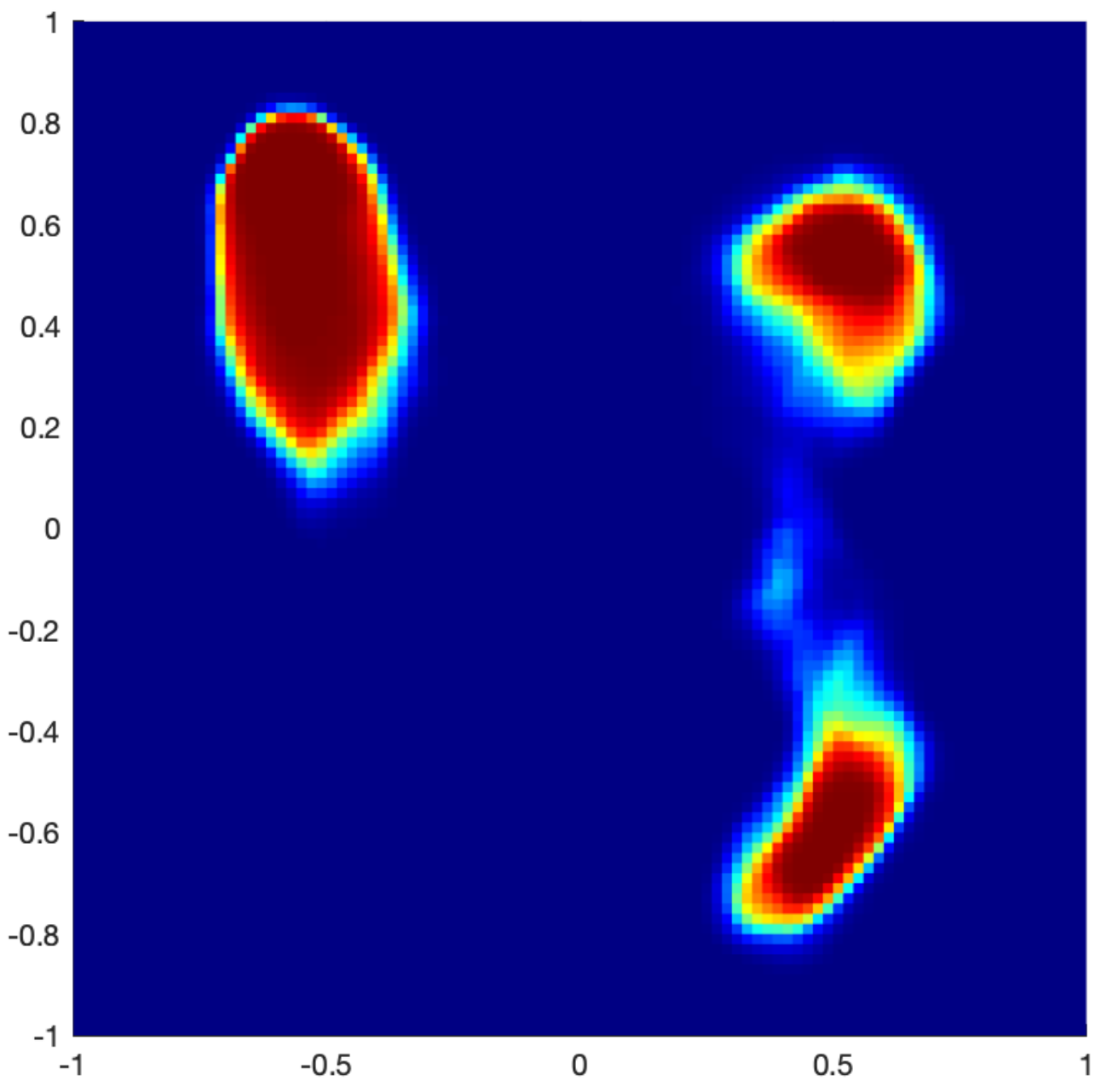}&
\includegraphics[width=1.1in]{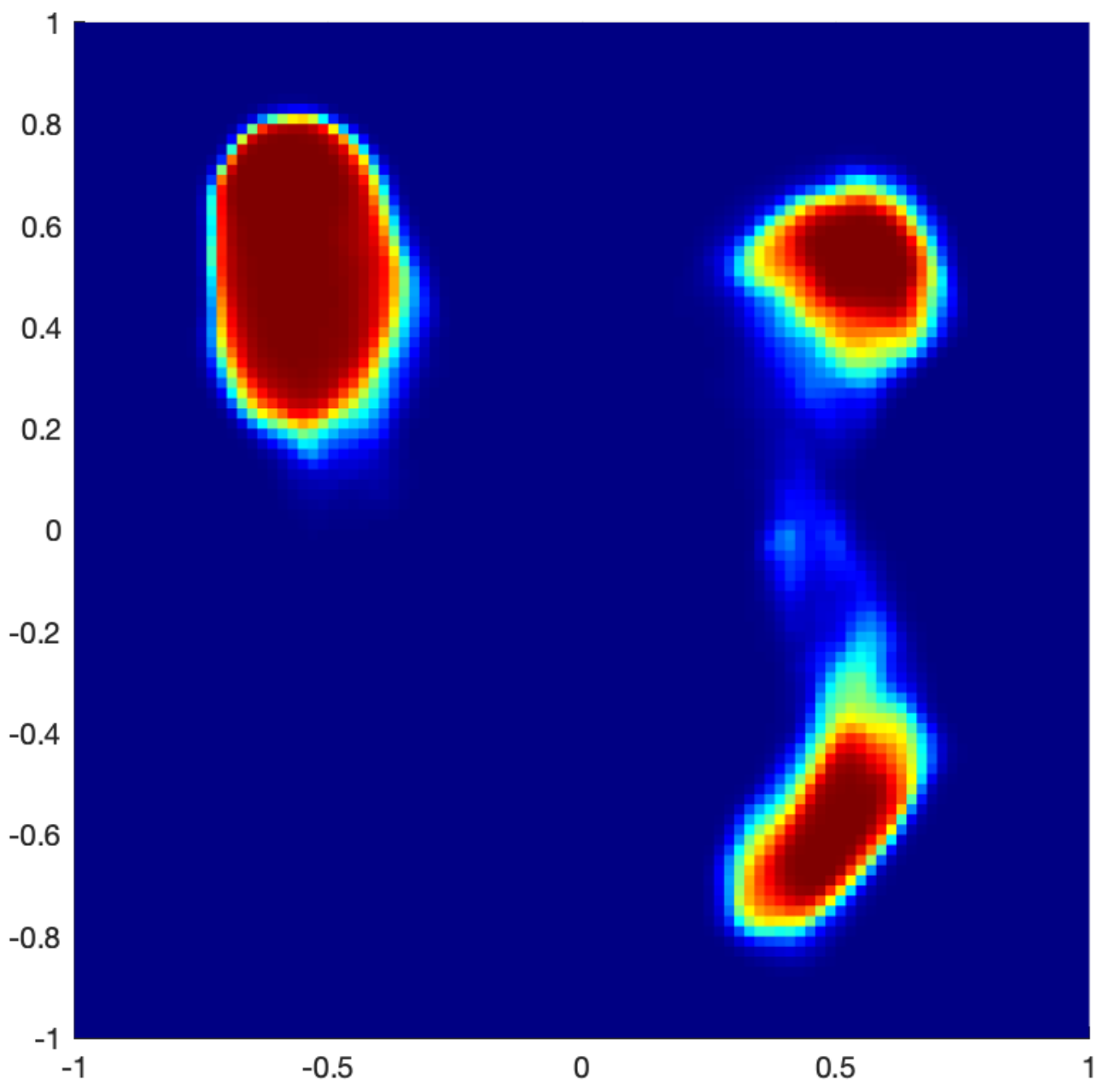}&
\includegraphics[width=1.1in]{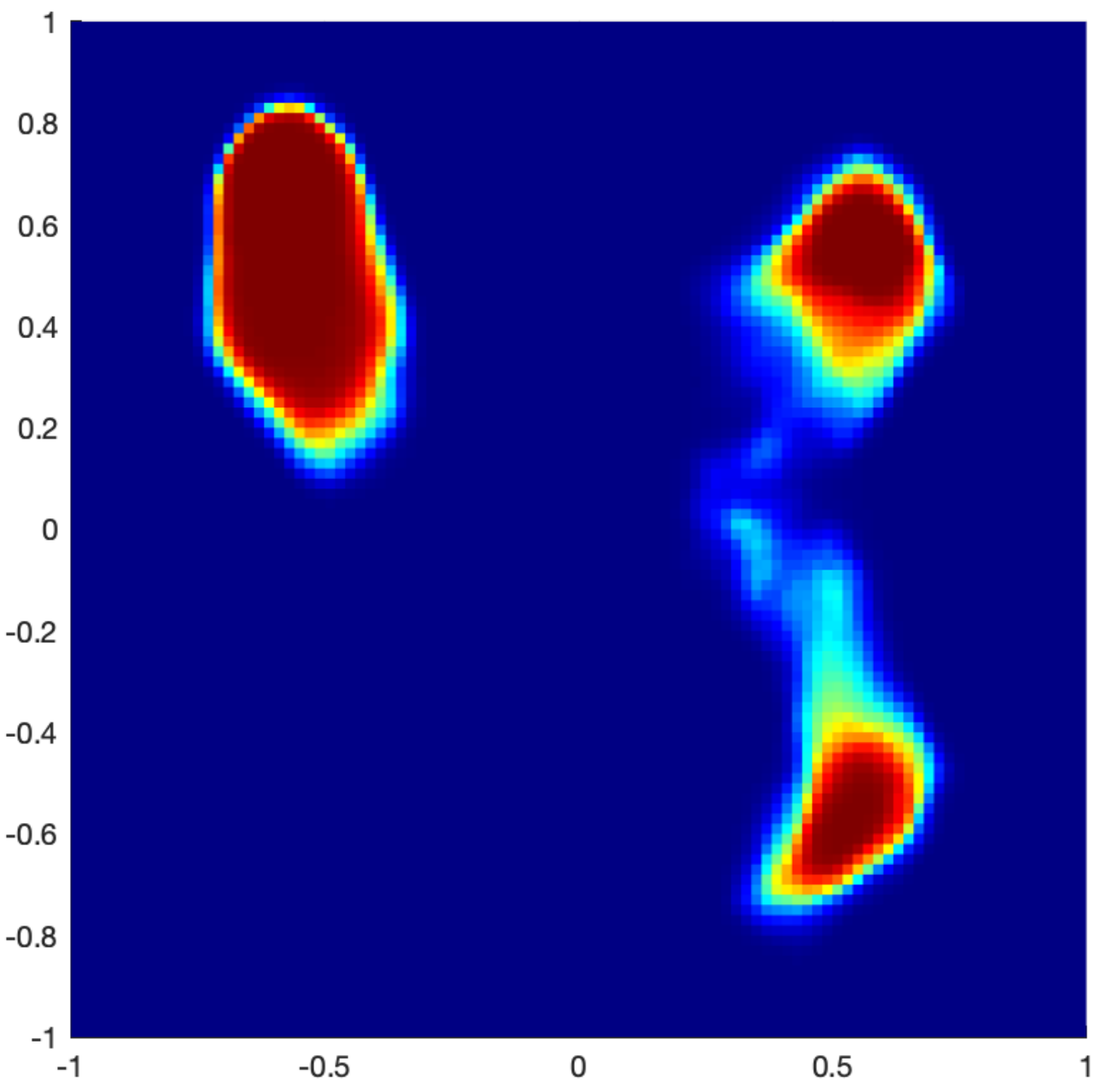}\\
\end{tabular}
  \caption{FNN-DDSM reconstruction for random ellipses: one ellipse is located closed to the center of domain and blocked from the boundary by other 3 ellipses (top) and this center ellipse is removed (bottom)} 
  \label{tab_FN_ell_comp}
\end{figure}

\begin{figure}[htbp]
\begin{tabular}{ >{\centering\arraybackslash}m{0.9in} >{\centering\arraybackslash}m{0.9in} >{\centering\arraybackslash}m{0.9in}  >{\centering\arraybackslash}m{0.9in}  >{\centering\arraybackslash}m{0.9in}  >{\centering\arraybackslash}m{0.9in} }
\centering
True coefficients &
N=1, $\delta=0$&
N=10, $\delta=0$&
N=20, $\delta=0$&
N=20, $\delta=10\%$ &
N=20, $\delta=20\%$ \\
\includegraphics[width=1.1in]{ell_comp1_0-eps-converted-to.pdf}&
\includegraphics[width=1.1in]{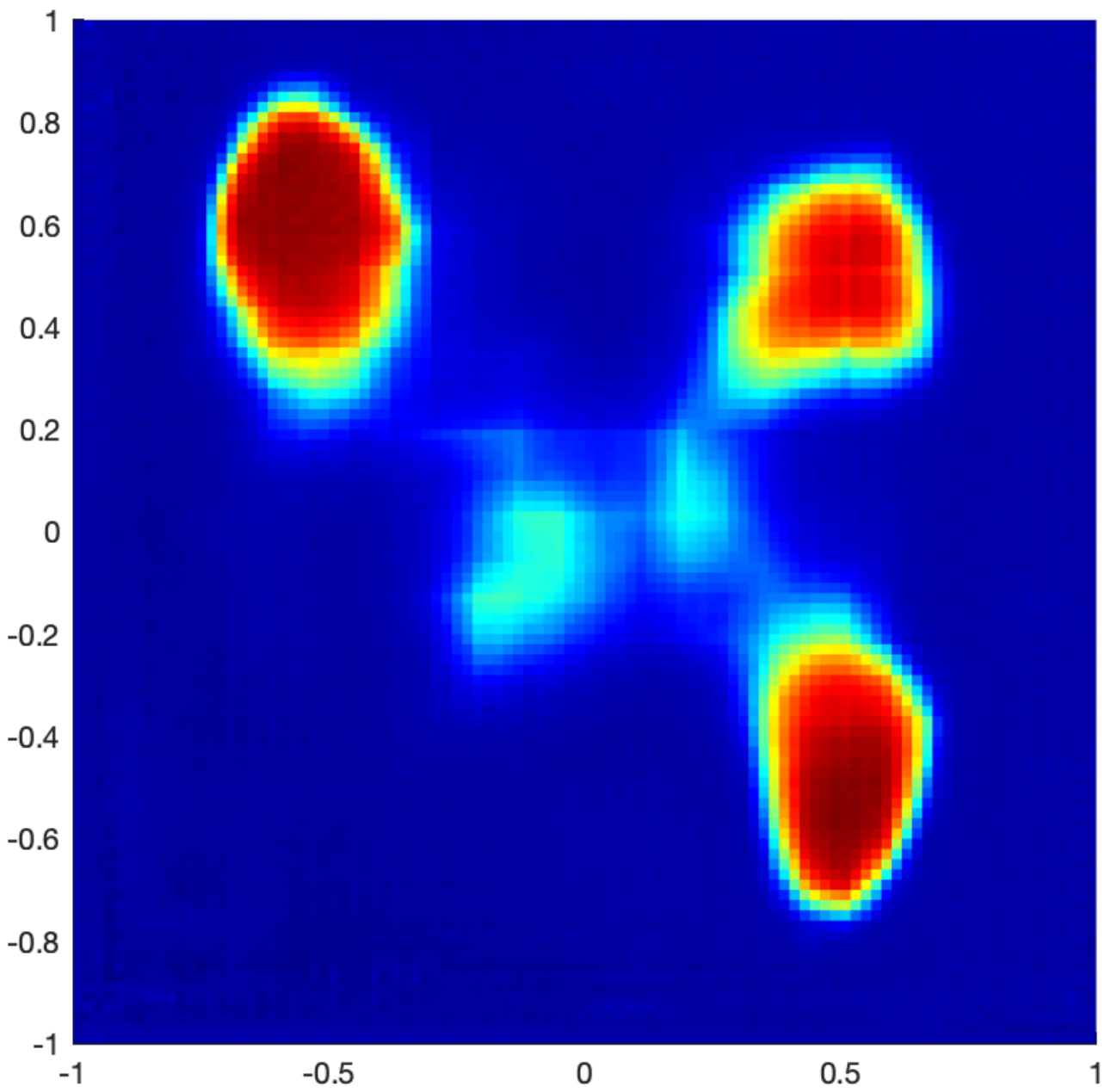}&
\includegraphics[width=1.1in]{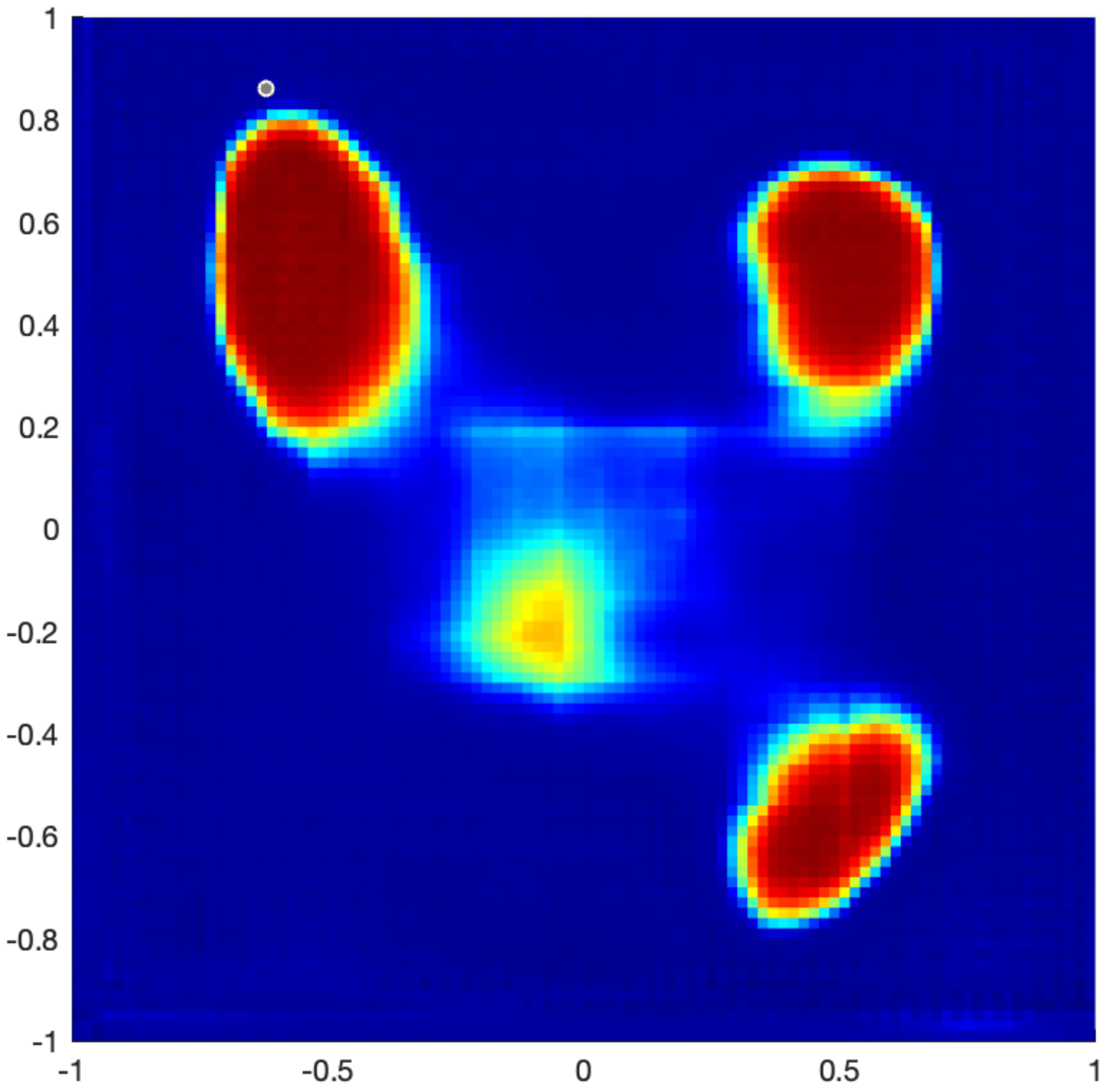}&
\includegraphics[width=1.1in]{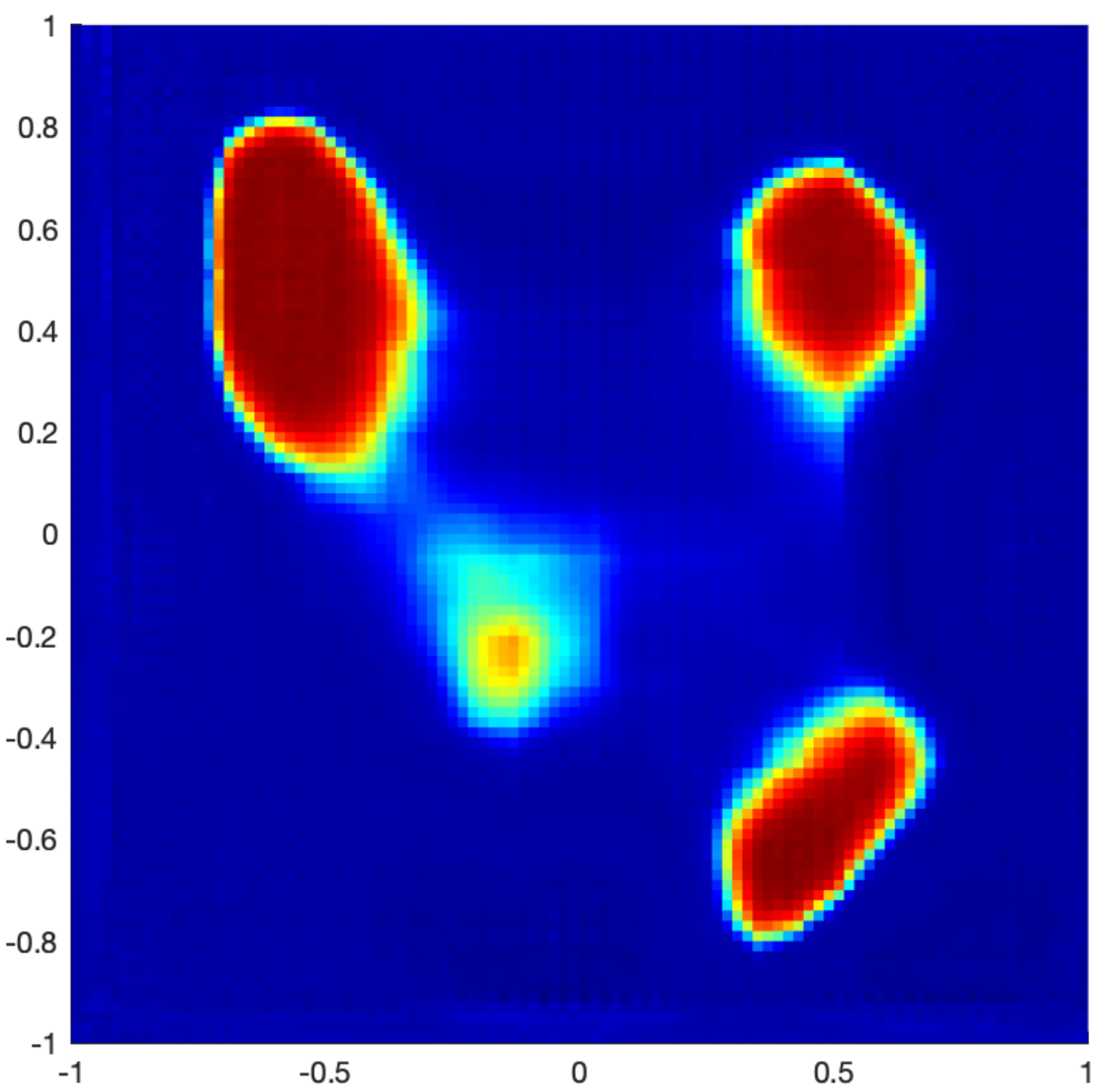}&
\includegraphics[width=1.1in]{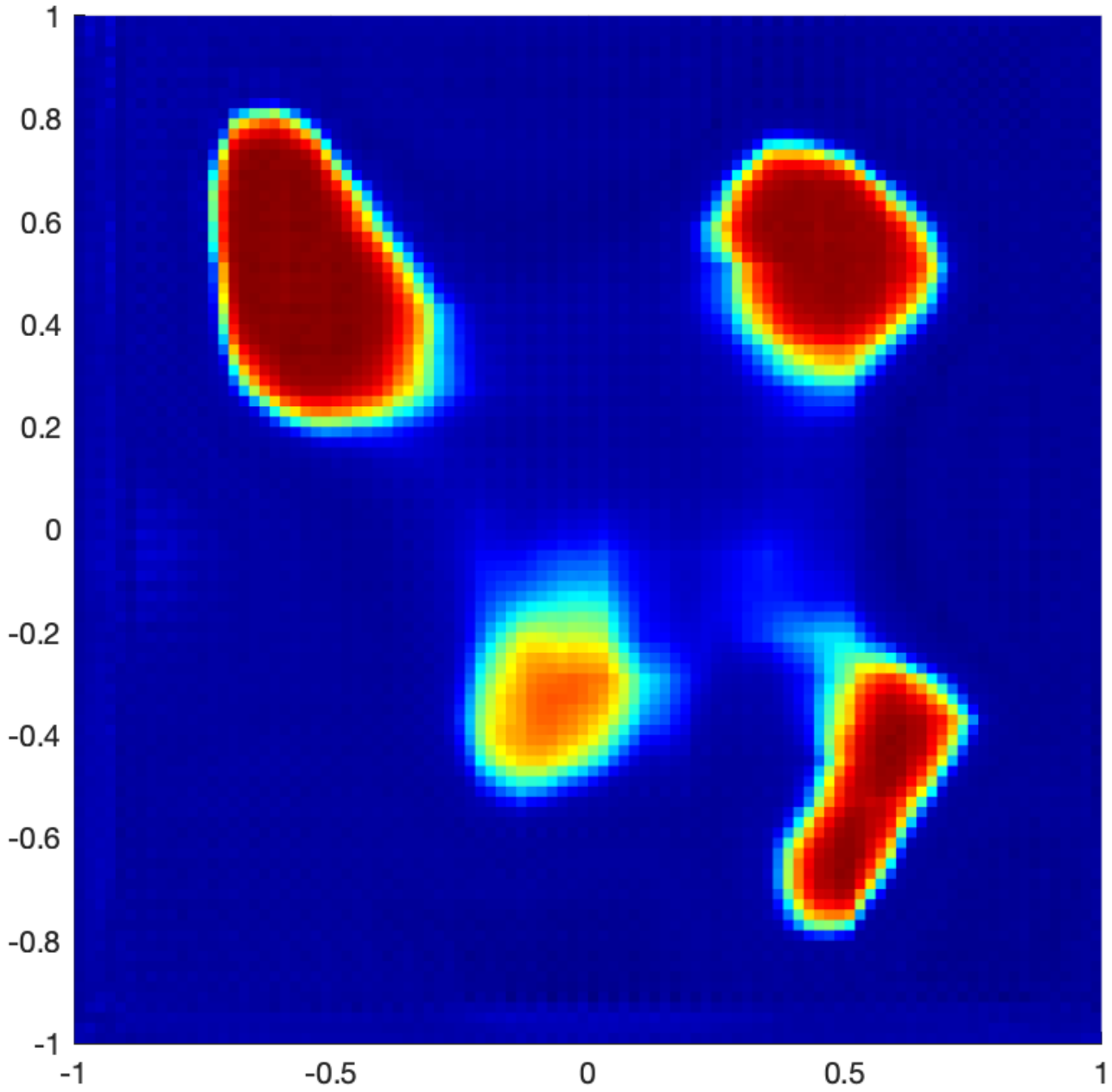}&
\includegraphics[width=1.1in]{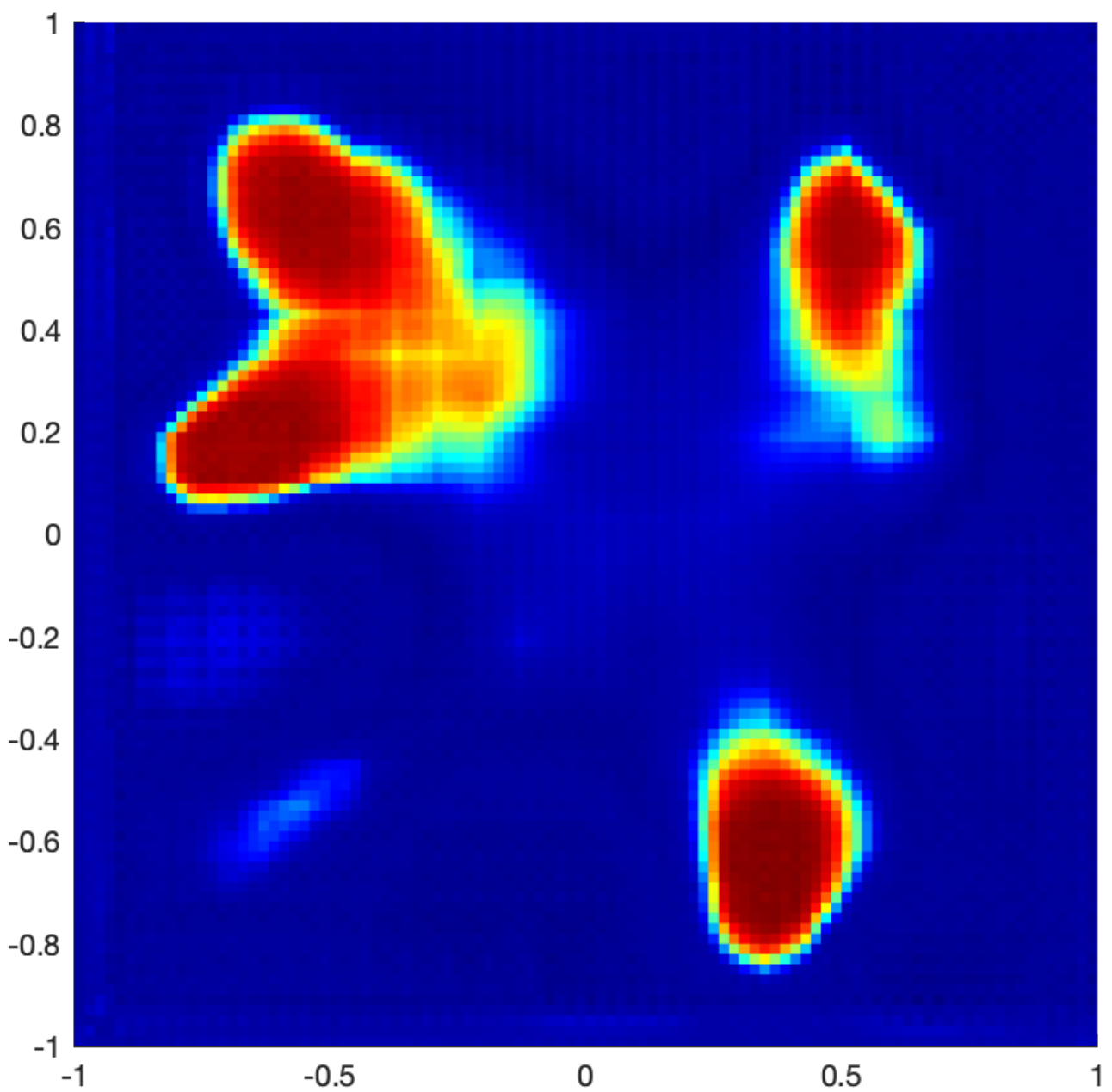}\\
\includegraphics[width=1.1in]{ell_comp2_0-eps-converted-to.pdf}&
\includegraphics[width=1.1in]{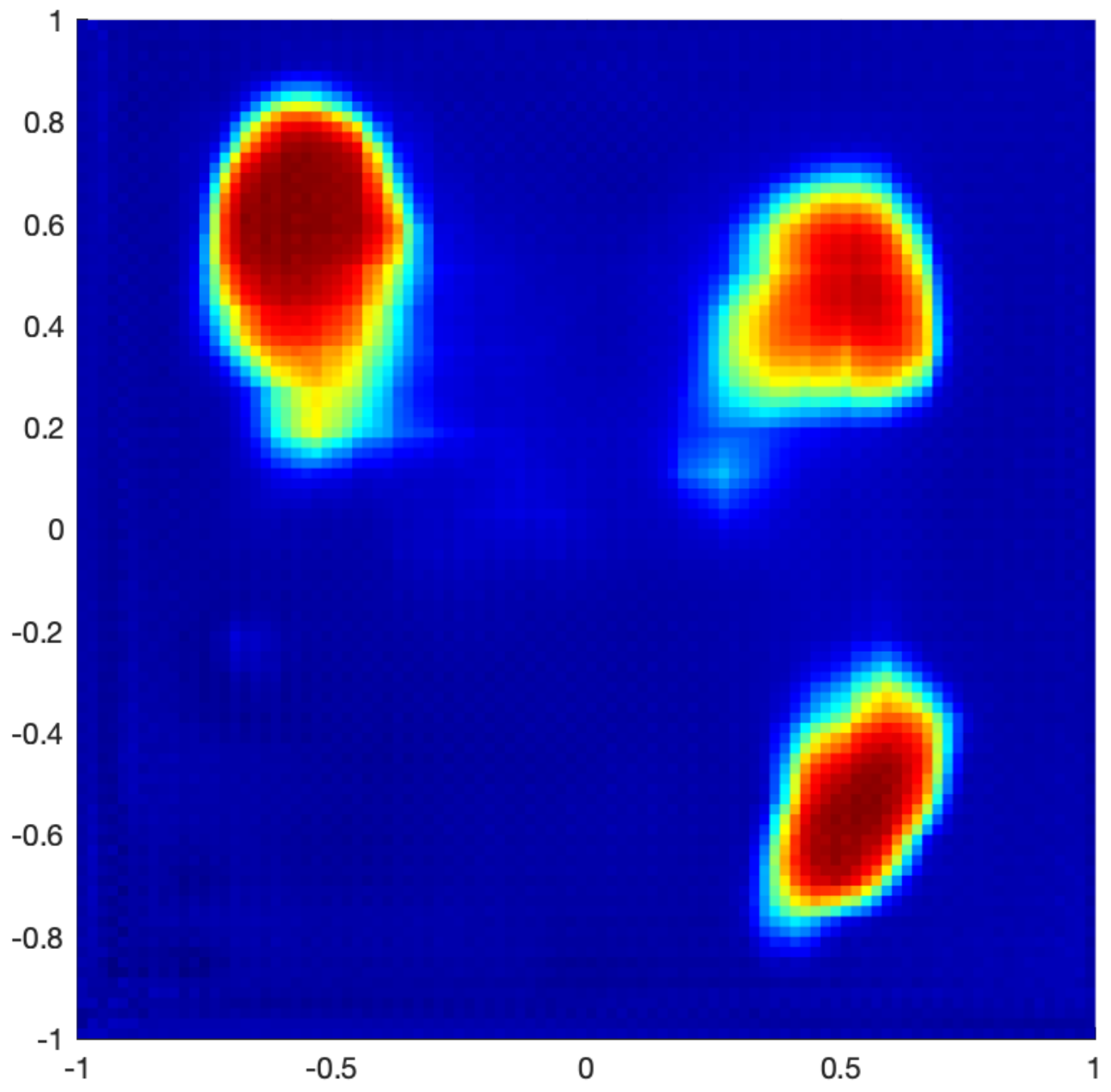}&
\includegraphics[width=1.1in]{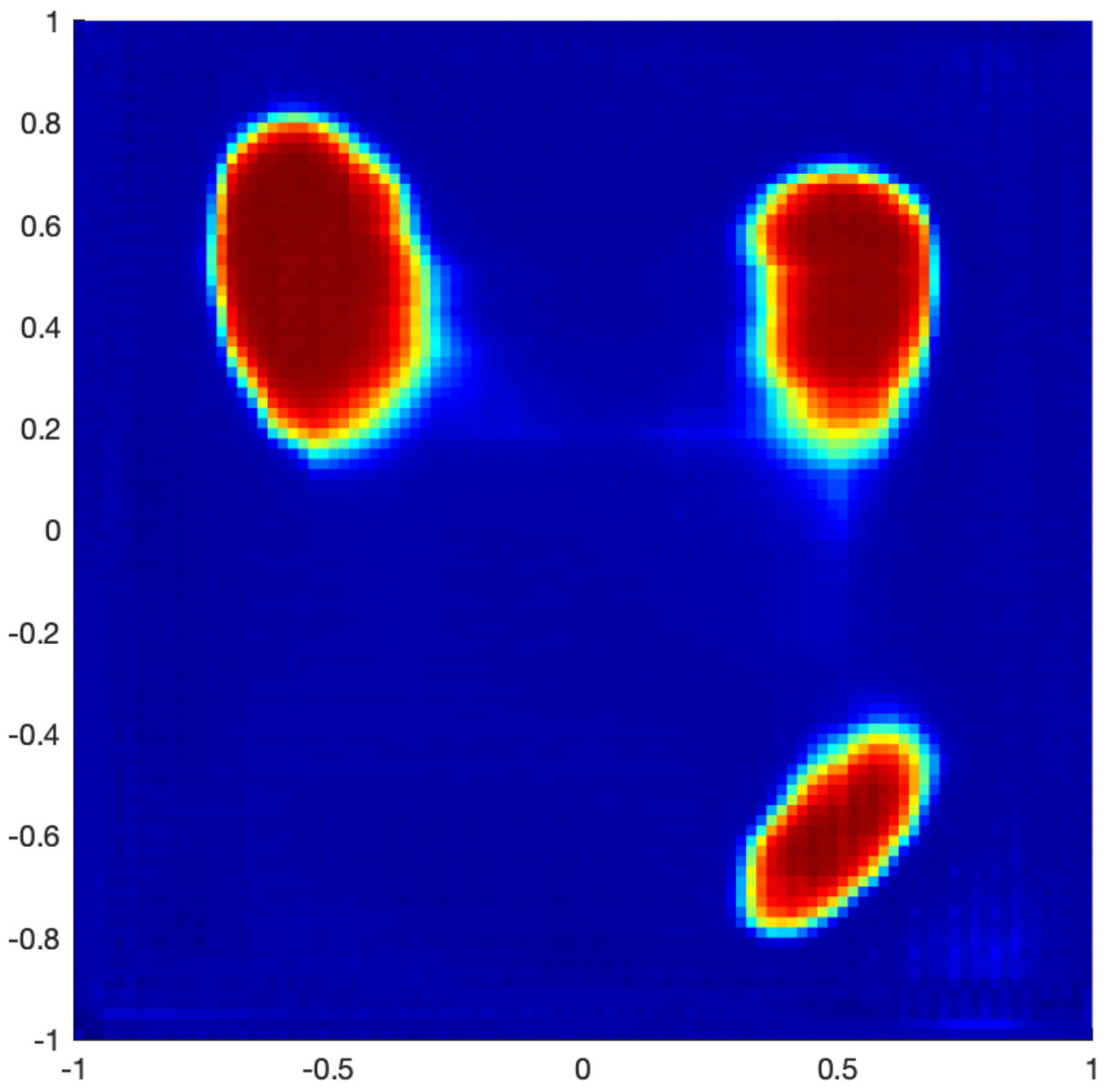}&
\includegraphics[width=1.1in]{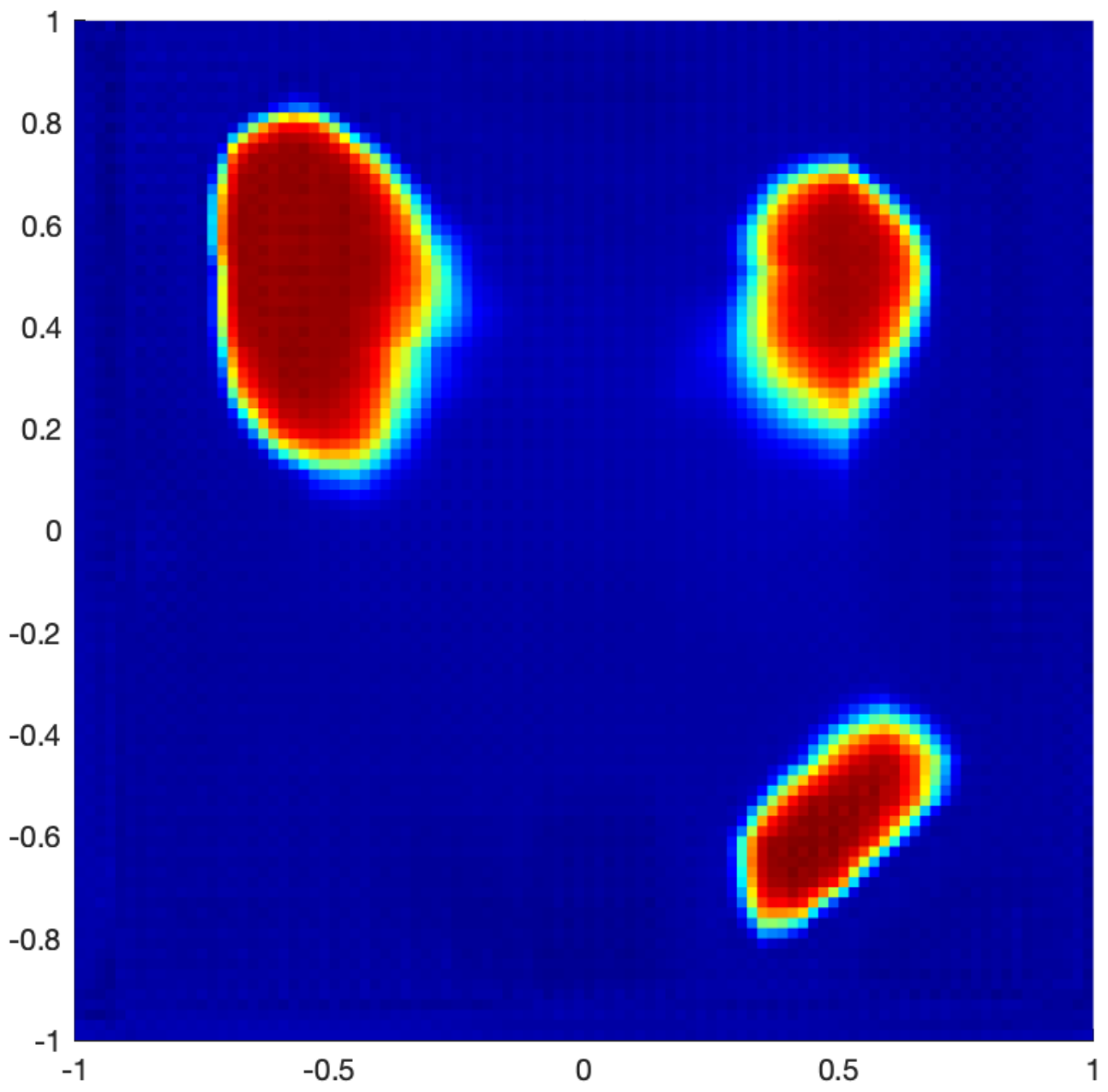}&
\includegraphics[width=1.1in]{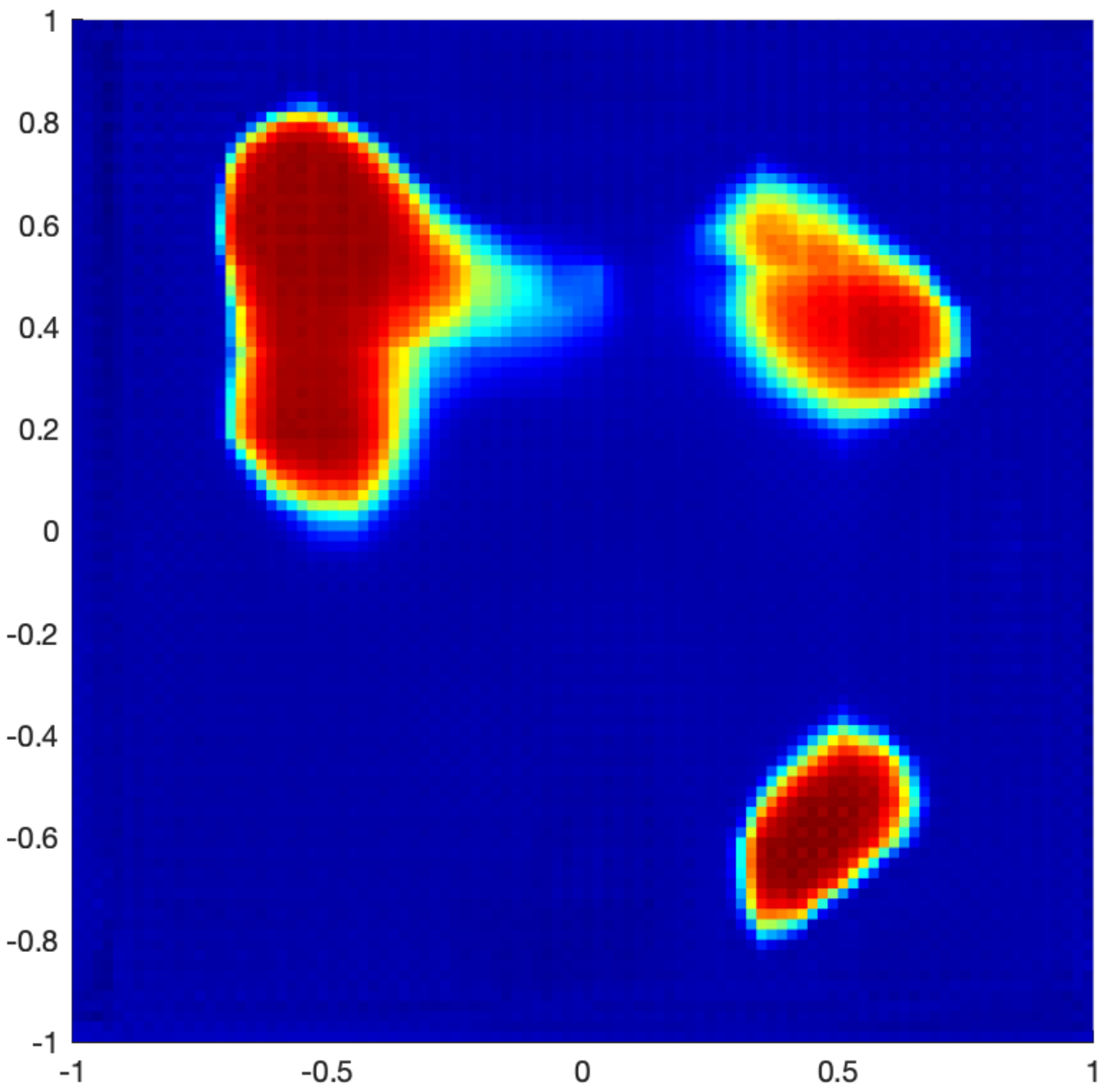}&
\includegraphics[width=1.1in]{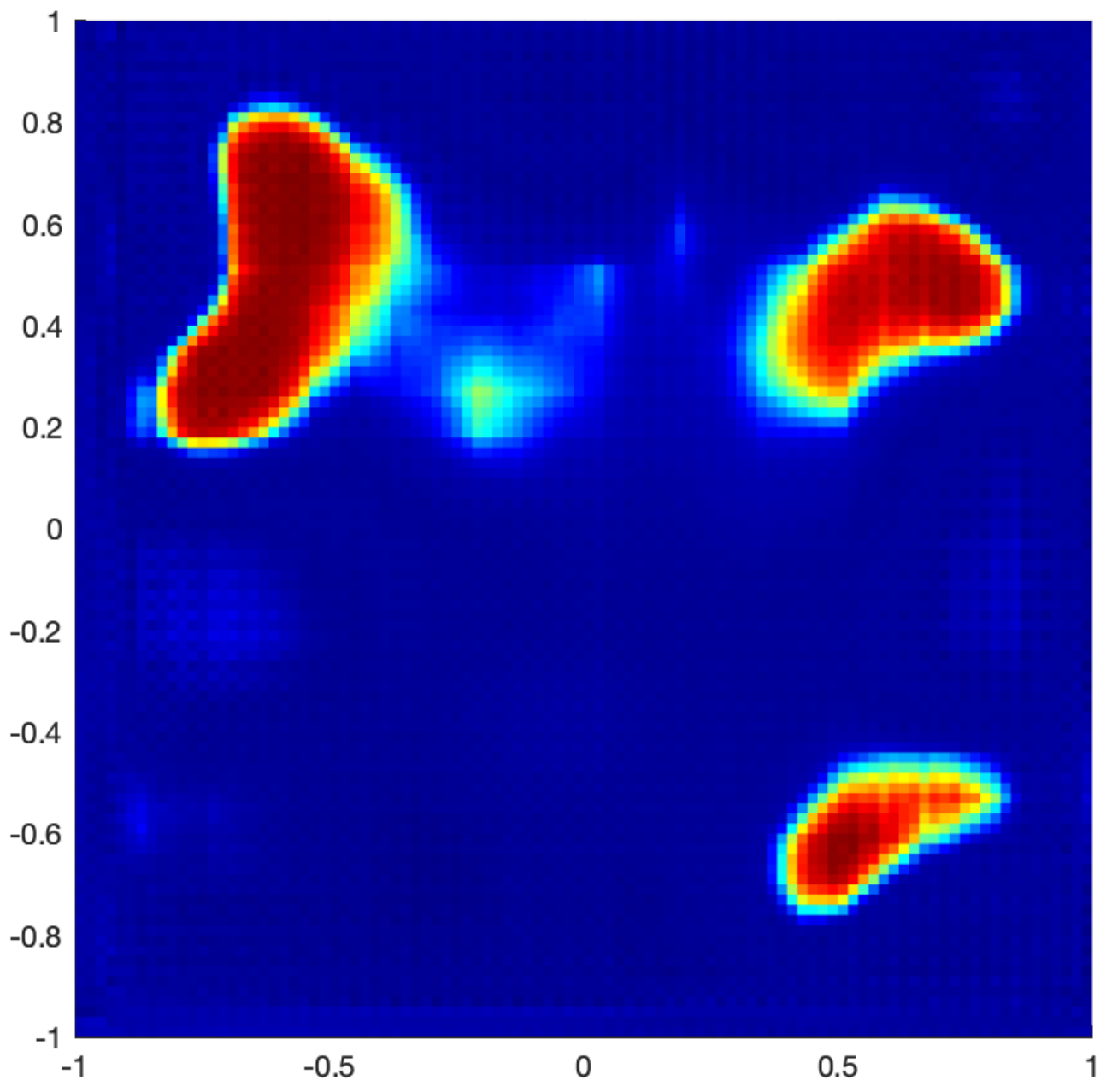}\\
\end{tabular}
  \caption{CNN-DDSM reconstruction for random ellipses: one ellipse is located closed to the center of domain and blocked from the boundary by other 3 ellipses (top) and this center ellipse is removed (bottom)} 
  \label{tab_ell_comp}
\end{figure}

The results in the previous subsection have shown the effectiveness and robustness of the proposed DDSMs. In this subsection, we investigate the sensitivity of the DDSMs with respect to the Cauchy data on the boundary and study their prediction performance on some inclusions that their geometry and topology are out of the scope of training sample set-up, i.e., those inclusions can not be the union of circles or ellipses,

We modify the true coefficients in Case 1 of \textbf{Scenario 2} and \textbf{Scenario 3} (first subplot in the first row in Figure \ref{tab_FN_5cir}-\ref{tab_ell})
by artificially moving one of the circles and ellipses to somewhere near the domain center away from the boundary and blocked by the surrounding circles or ellipses to get the first subplot in the first row in Figures \ref{tab_FN_cir_comp}-\ref{tab_ell_comp}. For the purpose of comparison, we also consider another situation that this center circle/ellipse is completely removed as shown in the first subplot in the second rows in Figures \ref{tab_FN_cir_comp}-\ref{tab_ell_comp}. We denote these two different coefficient distribution by $\sigma^s_{1}$ and $\sigma^s_2$ and compute the corresponding Dirichlet data $u^s_{1,\omega}|_{\partial\Omega}$ and $u^s_{2,\omega}|_{\partial\Omega}$ according to the Neumann data in \eqref{g_cos}. Let $s=c$ indicate the case of circles and $s=e$ indicate the case of ellipses. We plot $u^s_{1,\omega}|_{\partial\Omega}$, $u^s_{2,\omega}|_{\partial\Omega}$ and $u^s_{1,\omega}|_{\partial\Omega}-u^s_{2,\omega}|_{\partial\Omega}$ versus the polar angle $\theta\in[0.2\pi]$ in Figures \ref{fig:cir_compare} and \ref{fig:ell_compare}. 
We observe that with the same Neumann data $g_{\omega}$, the two Dirichlet data $u^s_{1,\omega}|_{\partial\Omega}$, $u^s_{2,\omega}|_{\partial\Omega}$ are very close with each other, although their coefficient distribution $\sigma^s_{1}, \sigma^s_2$ are very different. 
We can see that the relative difference between $u^s_{1,\omega}|_{\partial\Omega}$ and $u^s_{2,\omega}|_{\partial\Omega}$
\begin{equation}
\label{rela_diff}
\frac{\| u^s_{1,\omega} - u^s_{2,\omega} \|_{L^2(\partial\Omega)} }{ \| u^s_{1,\omega}\|_{L^2(\partial\Omega)}} \leqslant 6\%, ~~ \text{for} ~ \omega=1,2,...,5, ~~~ \text{and} ~~~\frac{\| u^s_{1,\omega} - u^s_{2,\omega} \|_{L^2(\partial\Omega)} }{ \| u^s_{1,\omega}\|_{L^2(\partial\Omega)}} \leqslant 2\%, ~~ \text{for} ~ \omega \geqslant 5.
\end{equation}
is in the same magnitude of regular noise level. It means that the small inclusions near the domain center have very subtle effect on the boundary Cauchy data, which makes them easily hidden in the domain and very difficult to be detected. In fact, such difference are small enough to be considered as just noise in the data for many conventional approaches. Fortunately, Figures \ref{tab_FN_cir_comp} and \ref{tab_cir_comp} show that both the FNN-DDSM and CNN-DDSM are able to sense such small change in the data and reflect it in a correct manner in the prediction. In Figures \ref{tab_FN_ell_comp} and \ref{tab_ell_comp}, the CNN-DDSM still captures the center inclusion quite well but FNN-DDSM barely gives {\jjhn this} information. These {\jjhn phenomena} agree with our observation from the previous results that the CNN-DDSM is more sensitive to the center inclusions. Although the reconstruction of the center inclusions is not as accurate as the surrounding inclusions near the boundary, we think it is still satisfactory to certain extent given that the difference between the Cauchy data of the two different coefficient distributions is very small. Furthermore, we note that this kind of sensitivity to data is also stable with respect to the noise, namely we can observe the reconstructed center inclusions even with 20\% noise for circles and 10\% noise for ellipses. More interestingly, it seems sometimes that the noise actually enhance the reconstruction of the center inclusions instead of undermining it, for example, in Figures \ref{tab_FN_cir_comp}, \ref{tab_cir_comp}, \ref{tab_ell_comp}, the center inclusion reconstructed with 10\% or 20\% noise is certainly more clear than no noise.

{\grc Due to the severe instability of the EIT problem, a small perturbation of the data may yield completely wrong reconstruction.}
In our opinion, 
a good algorithm should be, on one hand sensitive to the true(correct) data perturbation, i.e., recognize the inclusion information hidden in the Cauchy data as much as possibly, on the other hand, insensitive to the noise interruption, i.e., the reconstruction is not affected too much by noise. According to the numerical experiments, we believe the proposed DDSMs have this kind of feature. 

We remark that the inclusion shapes are never known exactly in practice and the a-priori knowledge of the shapes may not be always correct. It is important for a DNN based algorithm to have stable performance on a large variety of inclusion shapes no matter whether they are {\grc within the training set-up(library)}. Therefore, we study the performance of the DDSMs on some typical inclusion shapes {\grc which are out of the scope of the training sample set-up.} 
It is more challenging than predicting the inclusions just in the test set since it requires that the DNNs truly learn and fit the non-linear mapping from the Cauchy data set to the coefficient distribution instead of just a certain projection of the mapping on some low dimensional subspaces of data. 

\begin{figure}[htbp]
\begin{tabular}{ >{\centering\arraybackslash}m{0.9in} >{\centering\arraybackslash}m{0.9in} >{\centering\arraybackslash}m{0.9in}  >{\centering\arraybackslash}m{0.9in}  >{\centering\arraybackslash}m{0.9in}  >{\centering\arraybackslash}m{0.9in} }
\centering
True coefficients &
N=1, $\delta=0$&
N=10, $\delta=0$&
N=20, $\delta=0$&
N=20, $\delta=10\%$ &
N=20, $\delta=20\%$ \\
\includegraphics[width=1.1in]{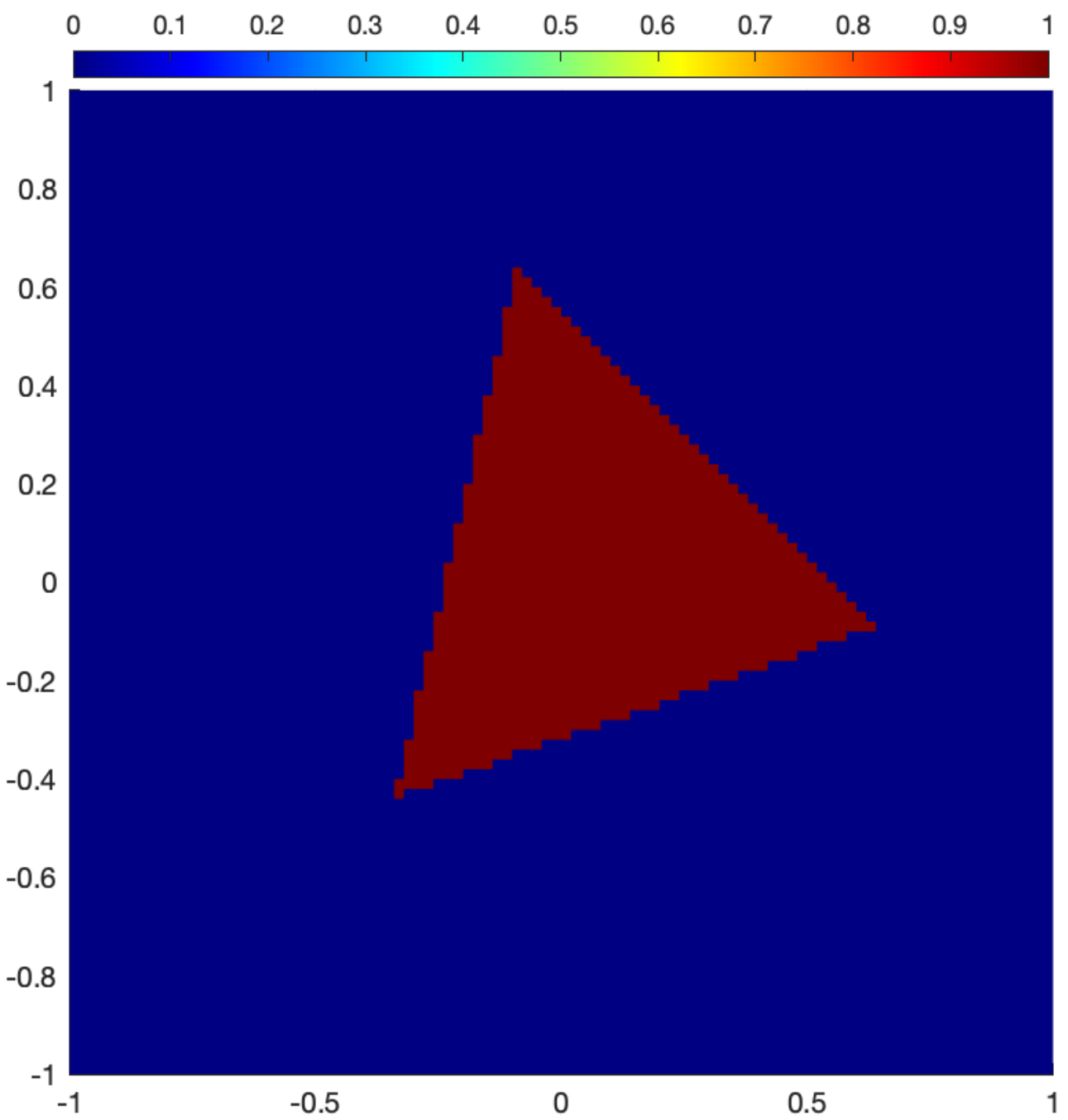}&
\includegraphics[width=1.1in]{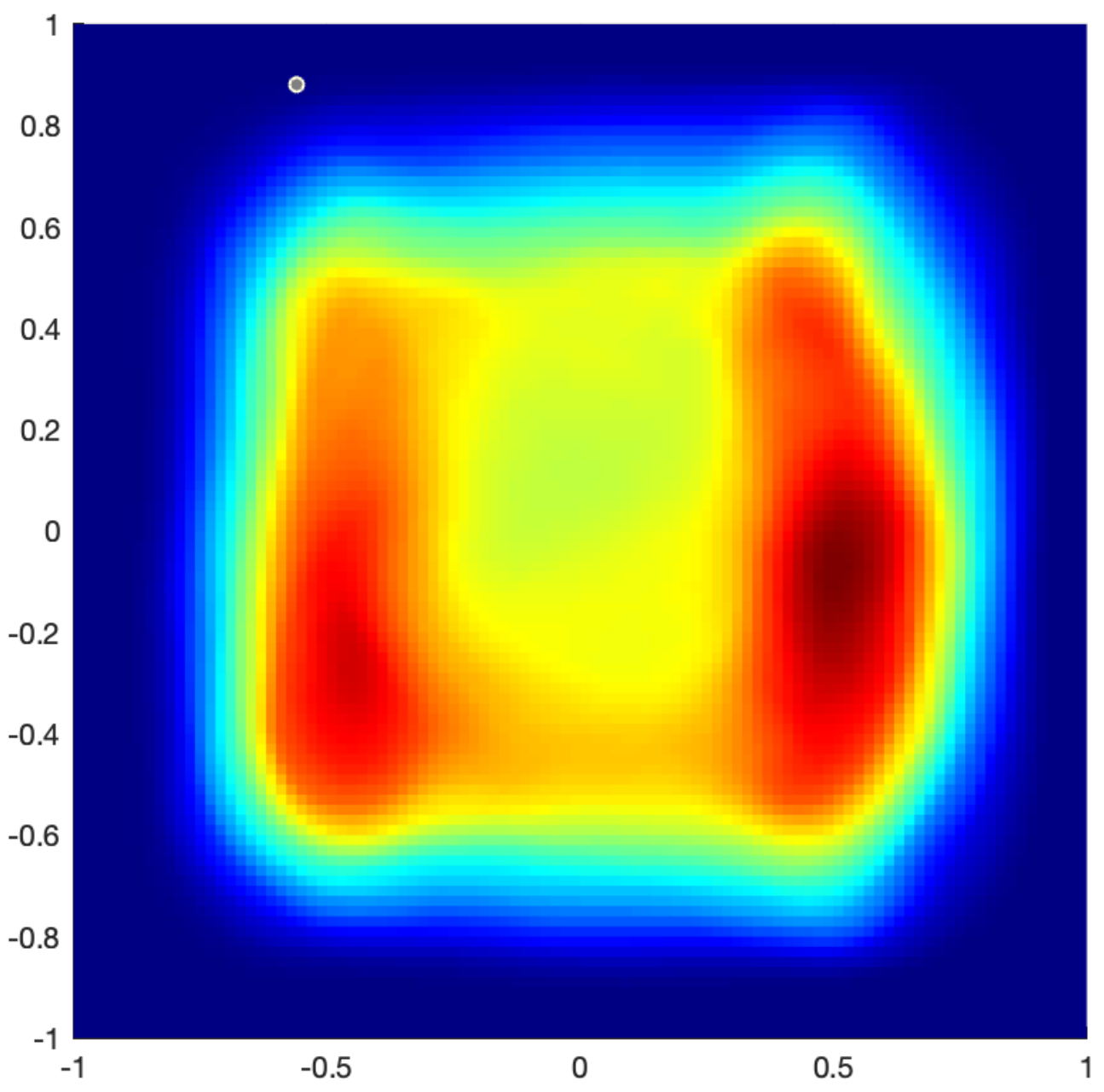}&
\includegraphics[width=1.1in]{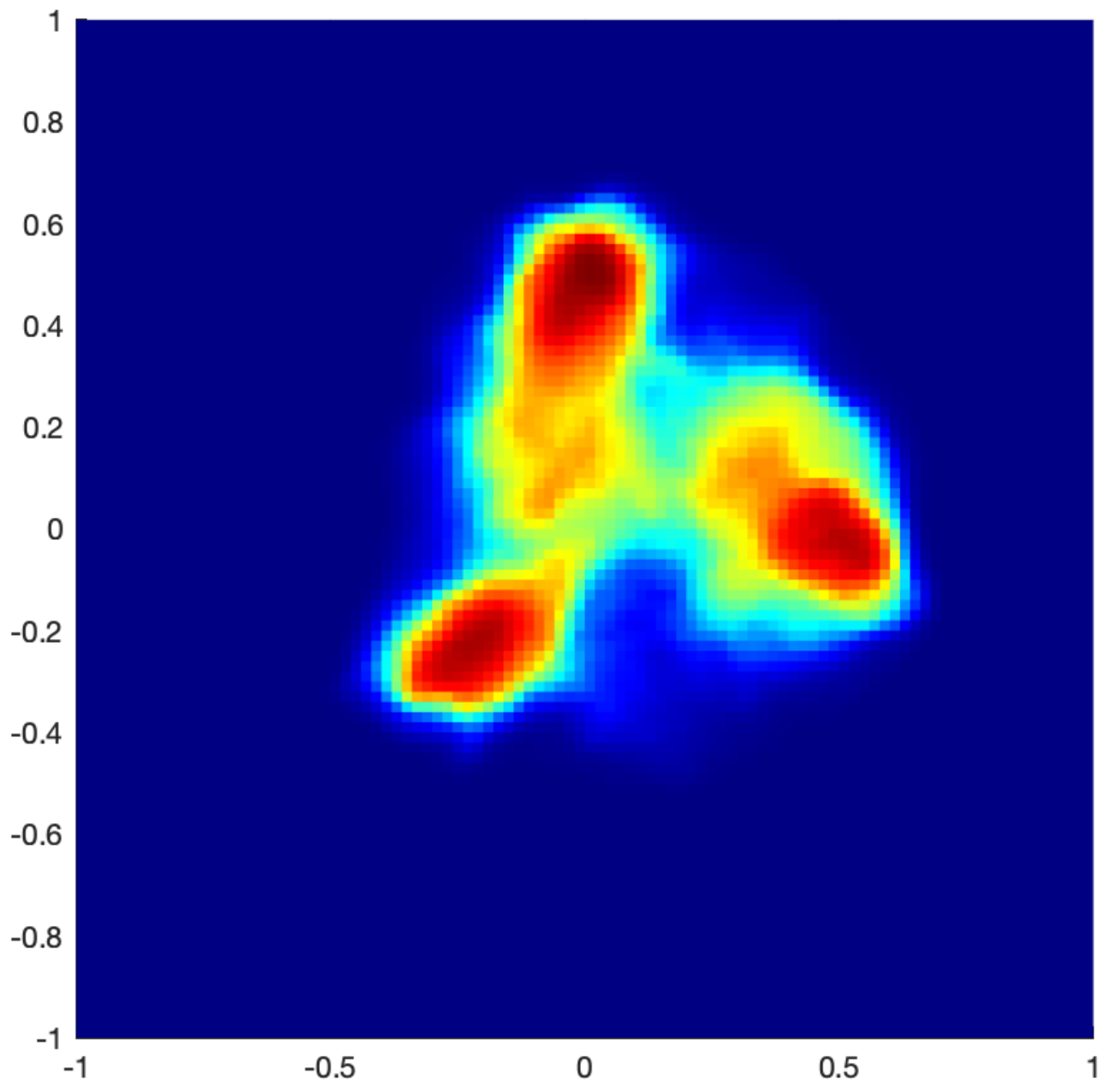}&
\includegraphics[width=1.1in]{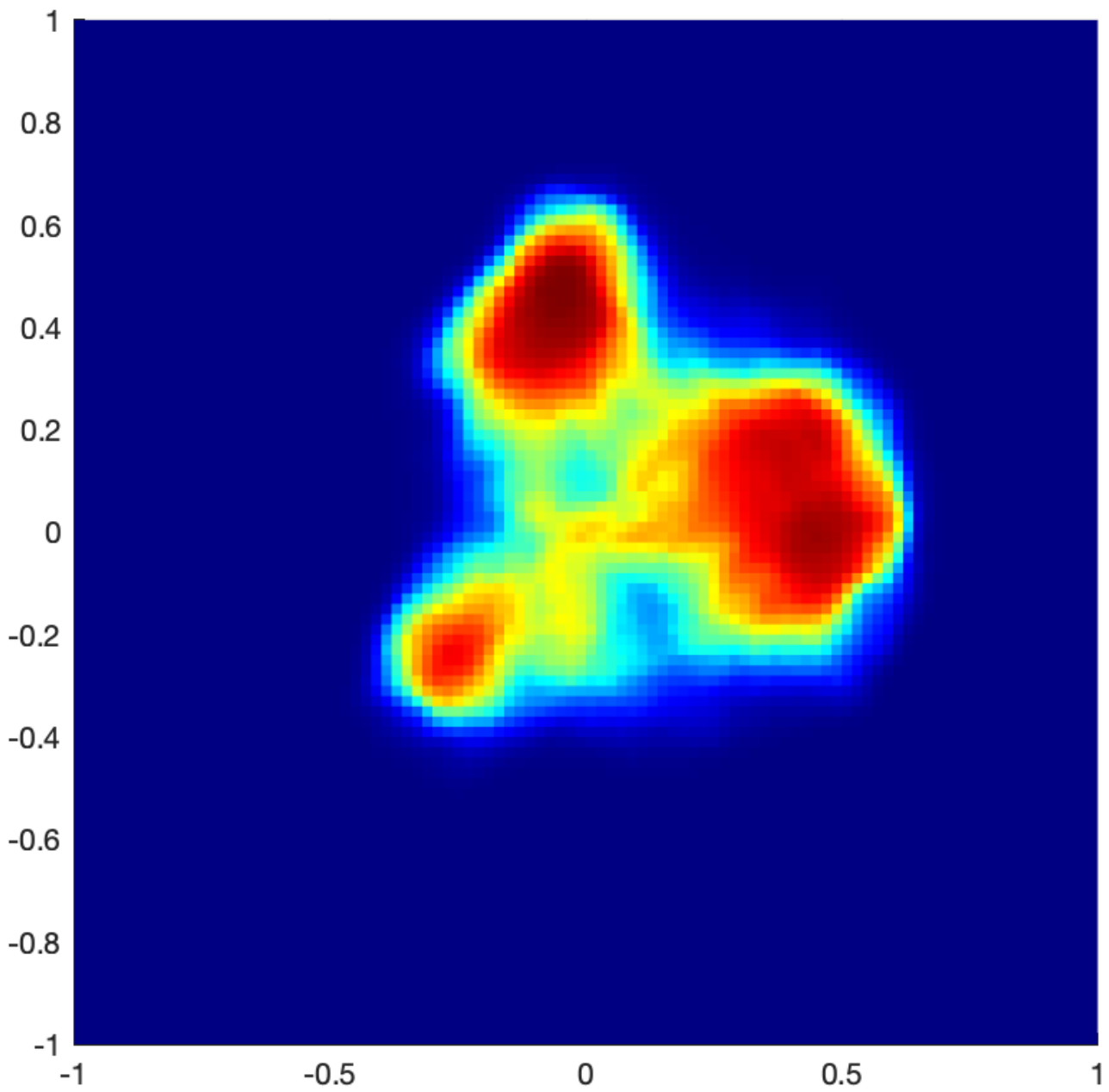}&
\includegraphics[width=1.1in]{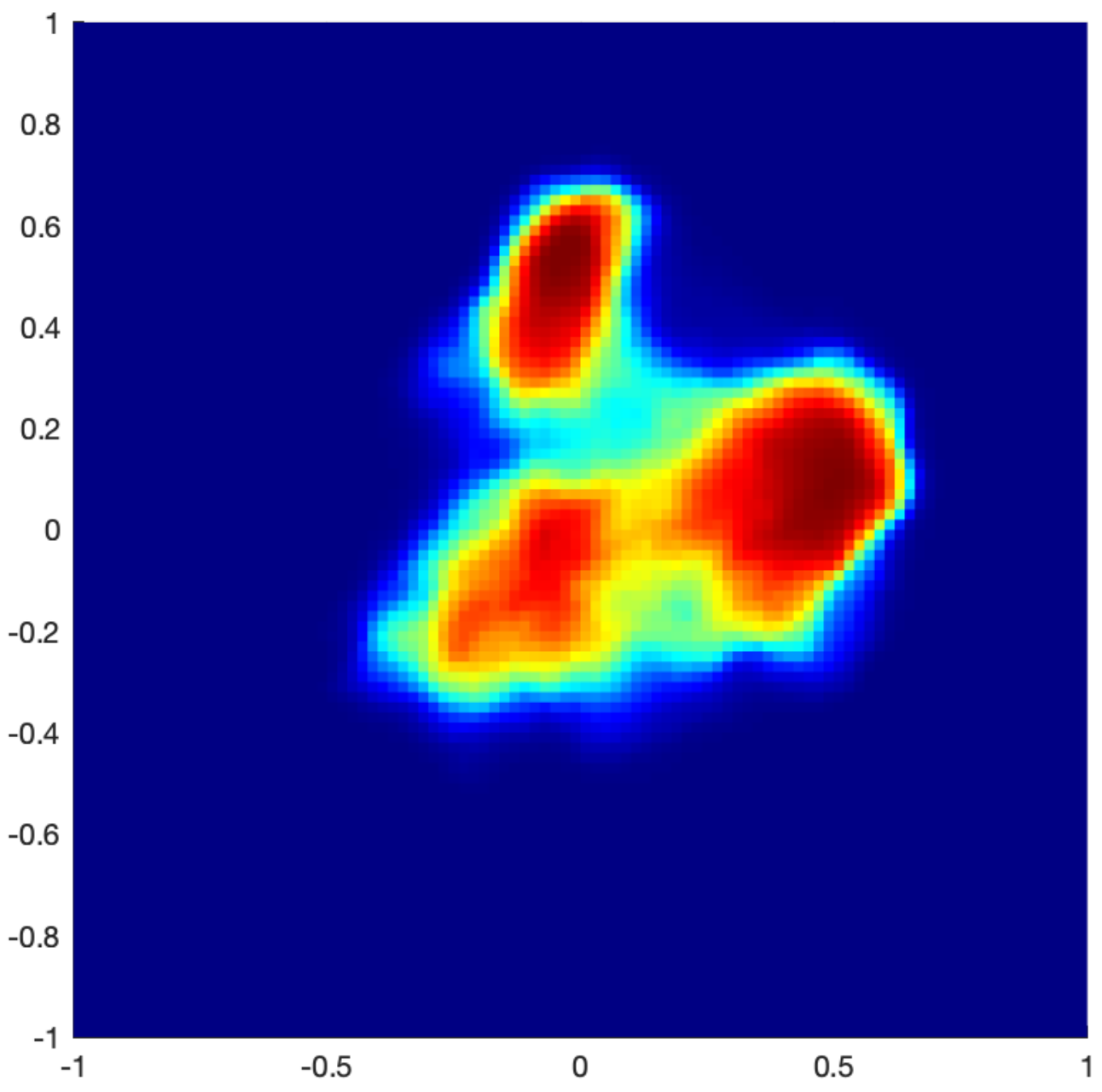}&
\includegraphics[width=1.1in]{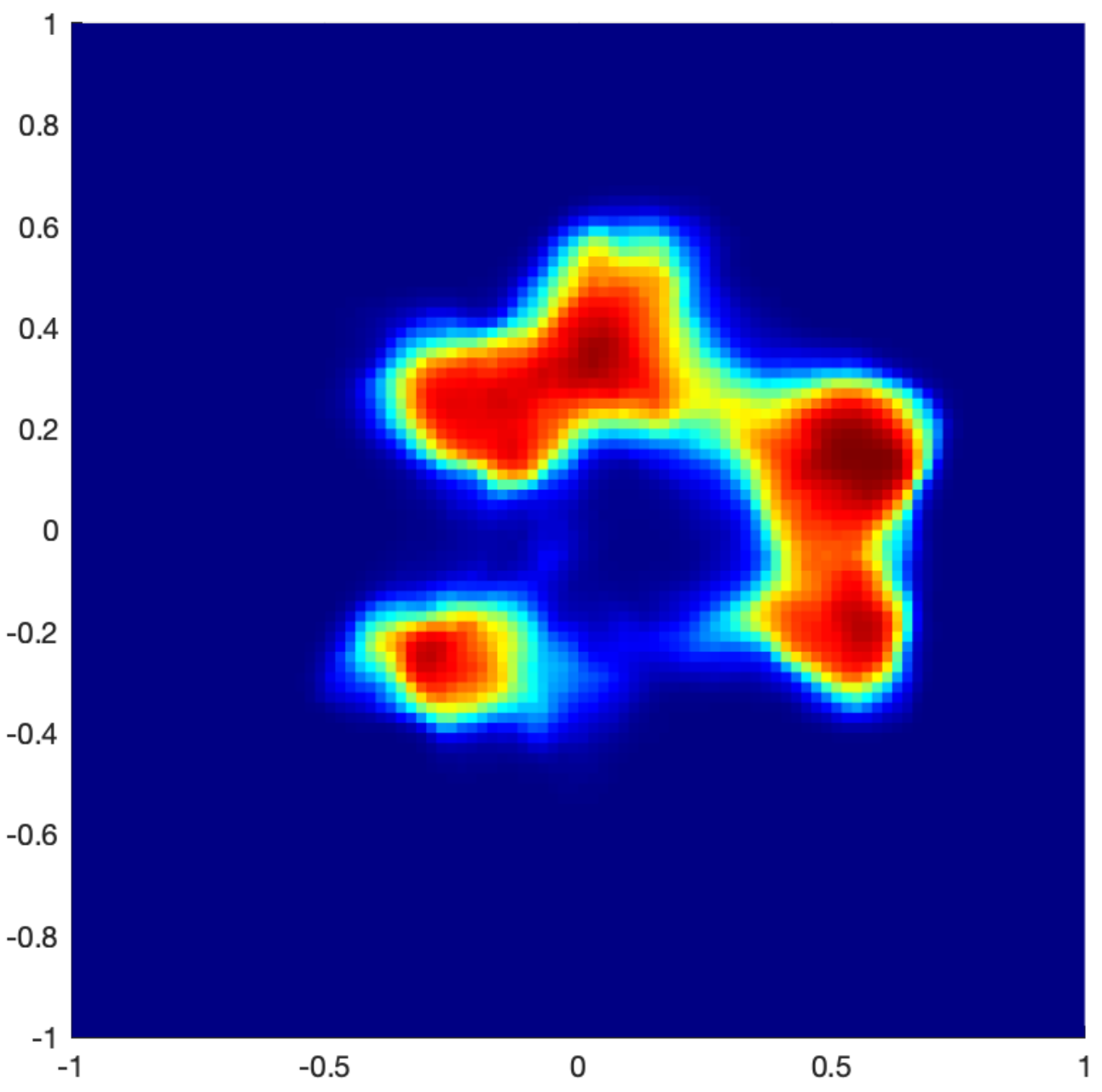}\\
\includegraphics[width=1.1in]{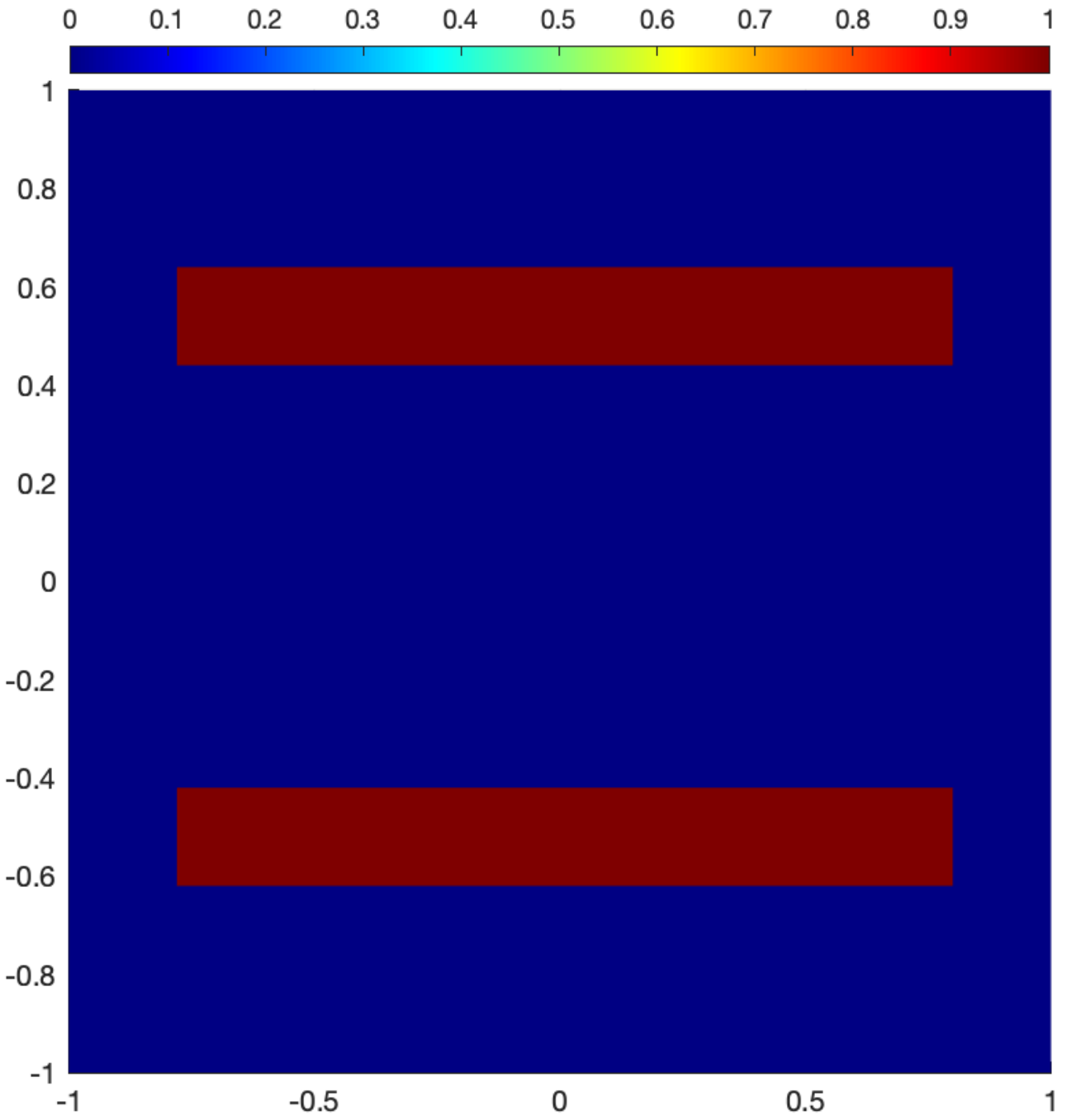}&
\includegraphics[width=1.1in]{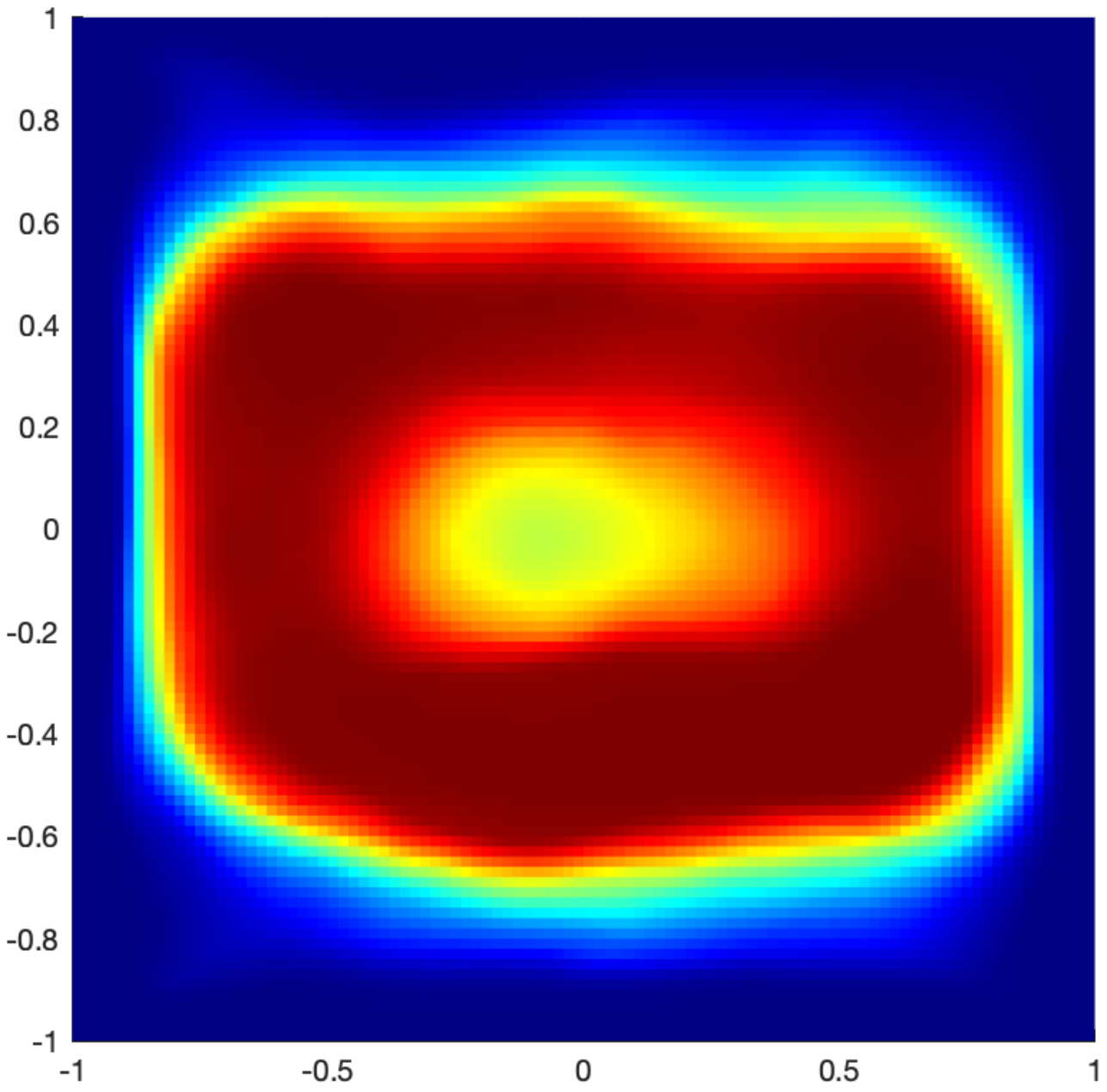}&
\includegraphics[width=1.1in]{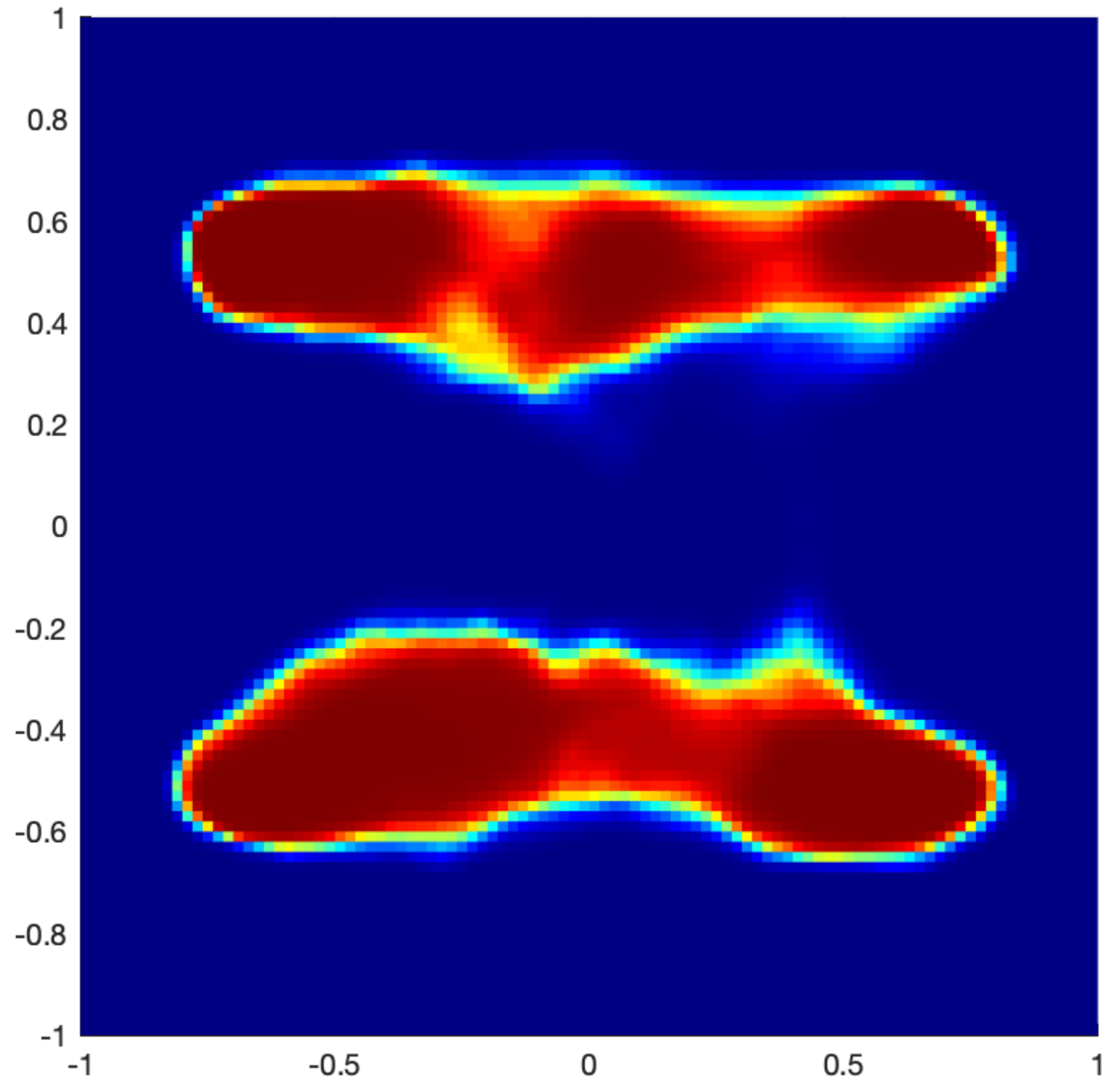}&
\includegraphics[width=1.1in]{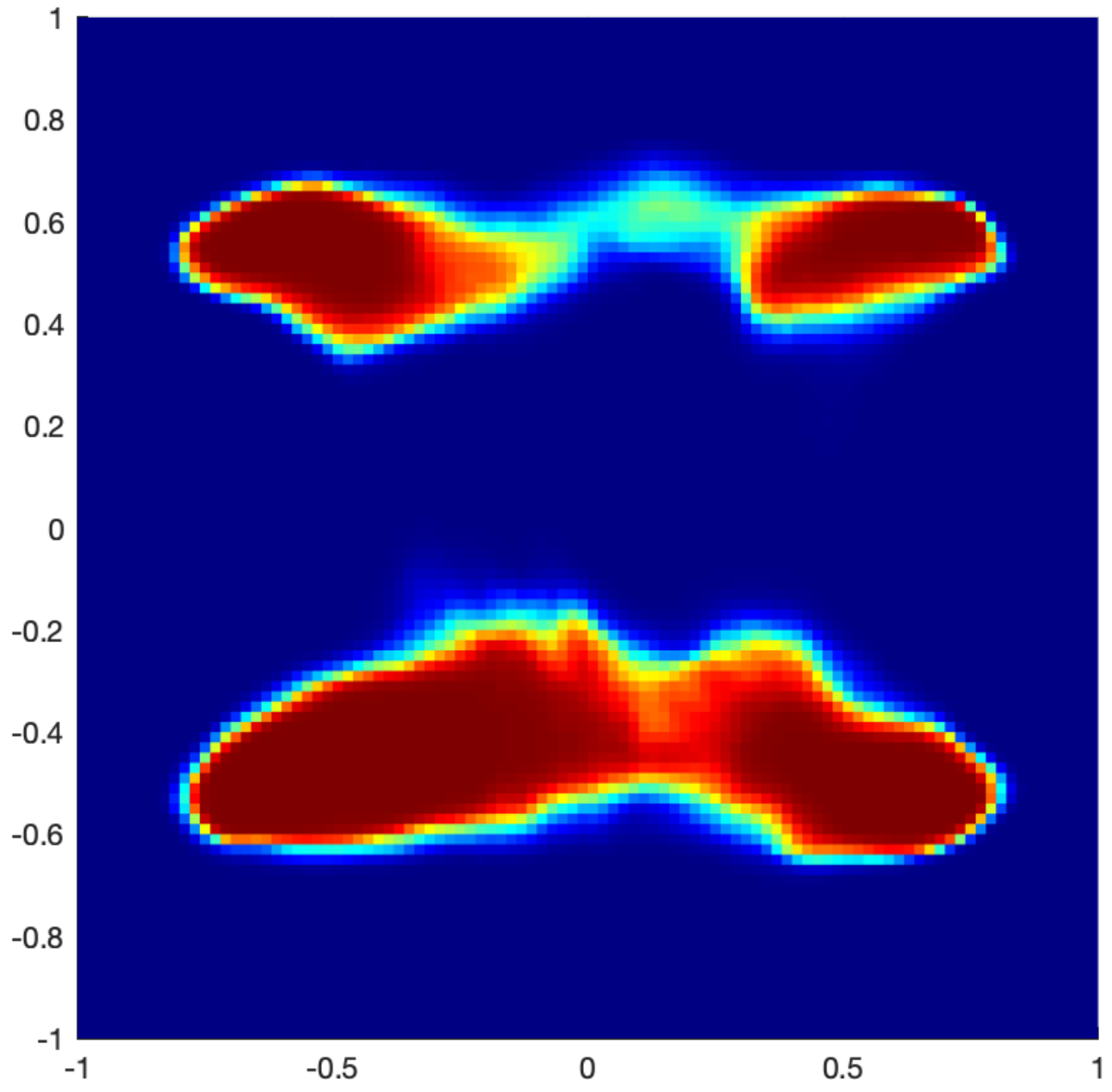}&
\includegraphics[width=1.1in]{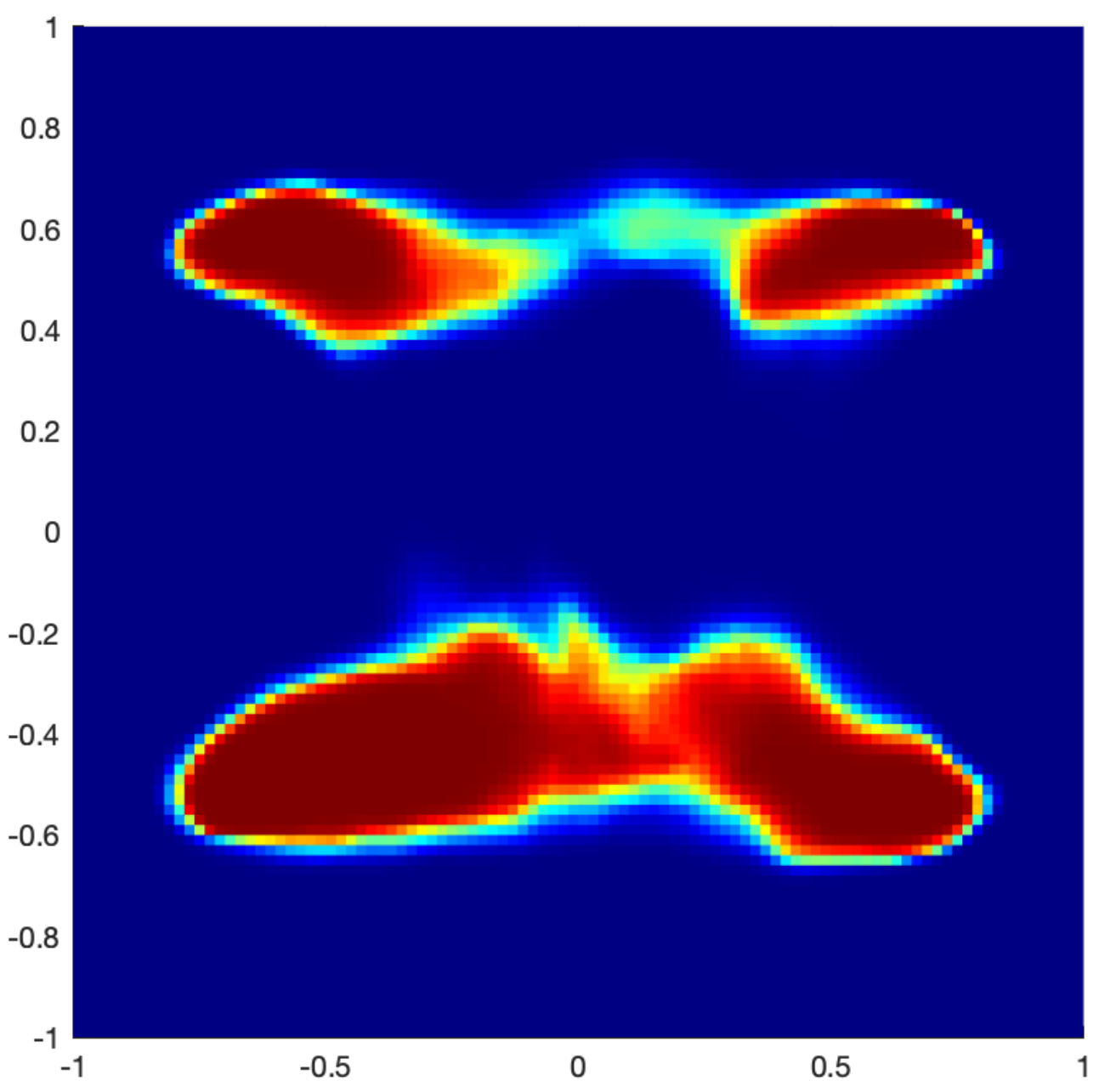}&
\includegraphics[width=1.1in]{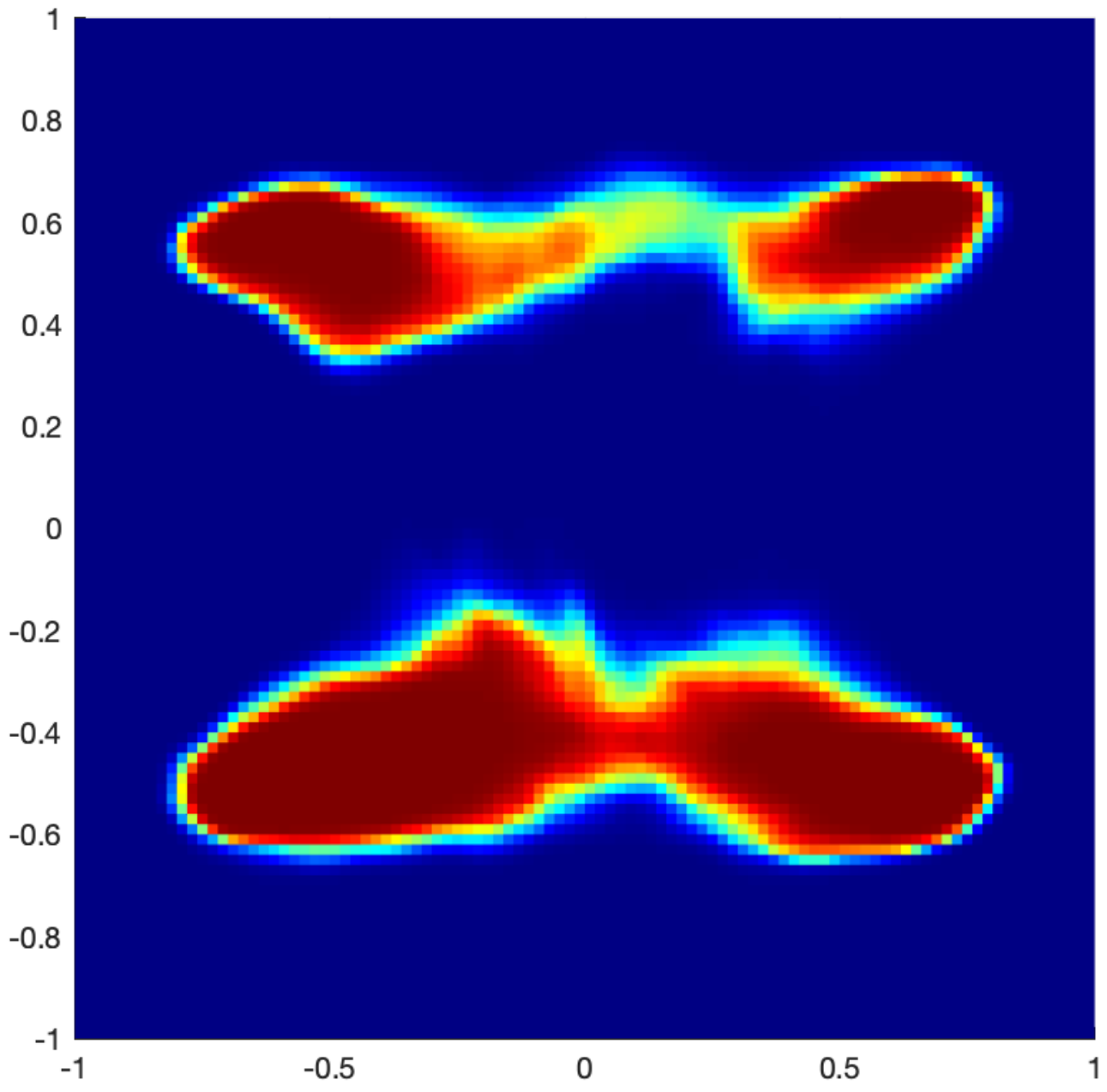}\\
\includegraphics[width=1.1in]{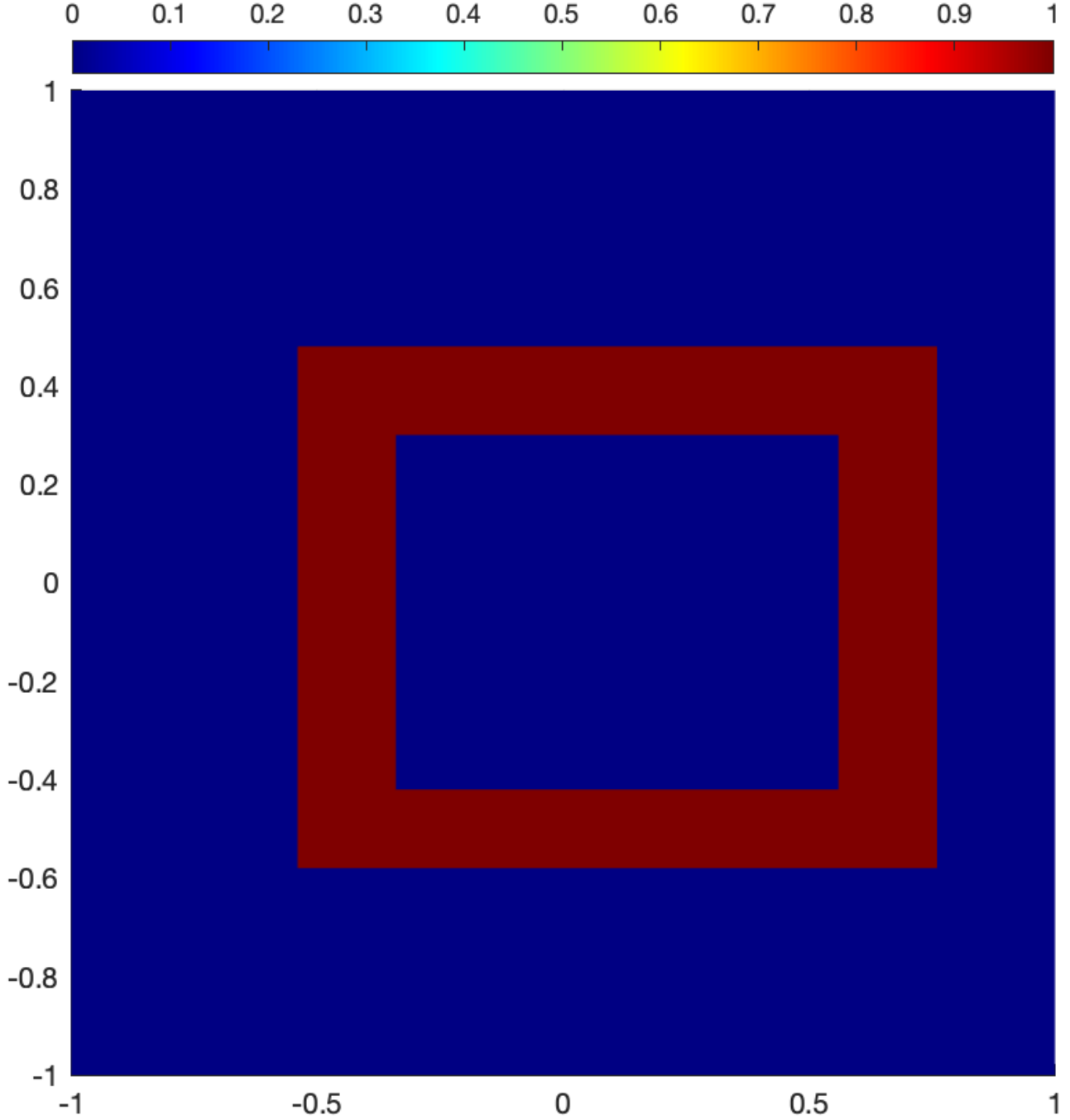}&
\includegraphics[width=1.1in]{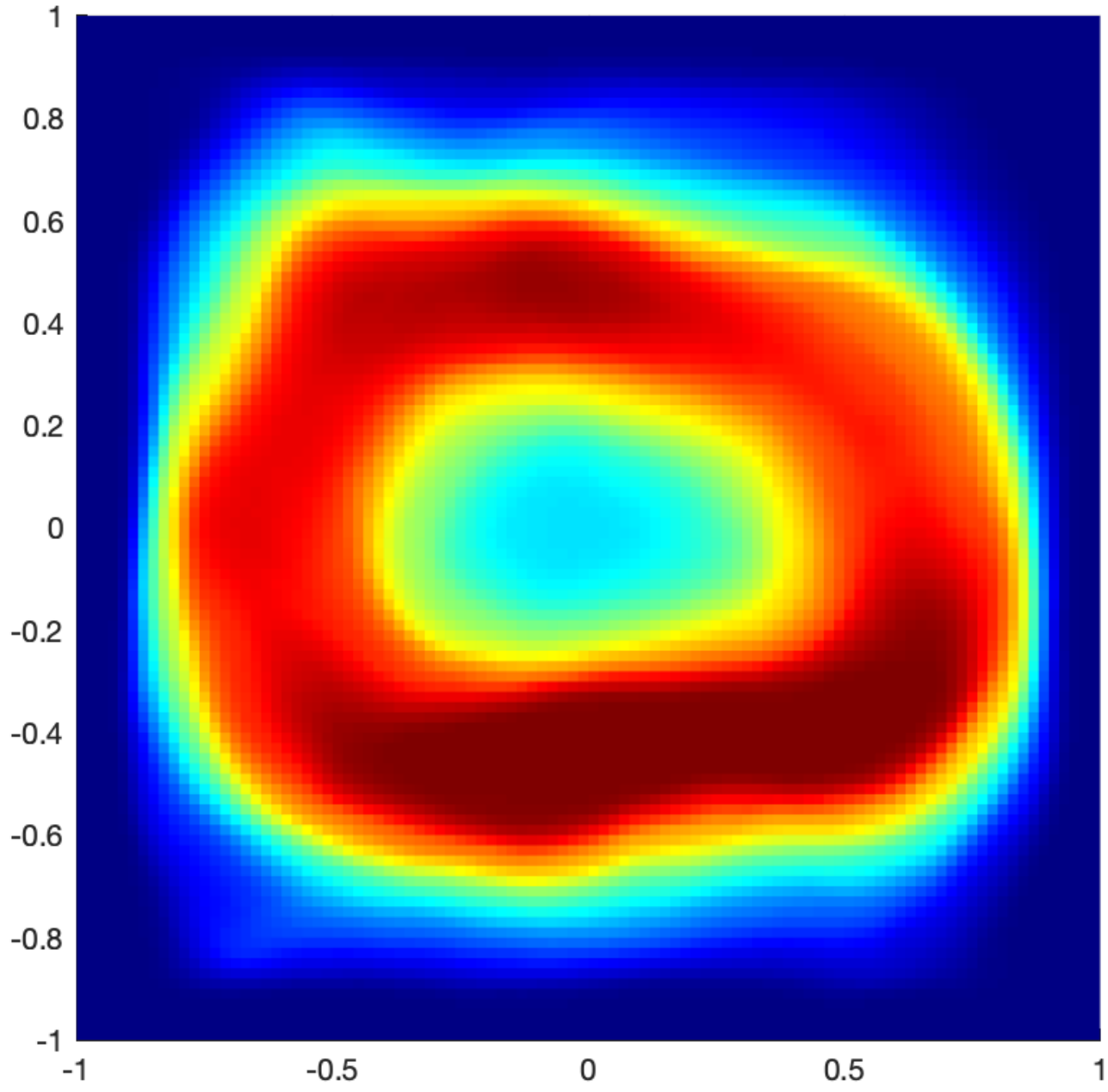}&
\includegraphics[width=1.1in]{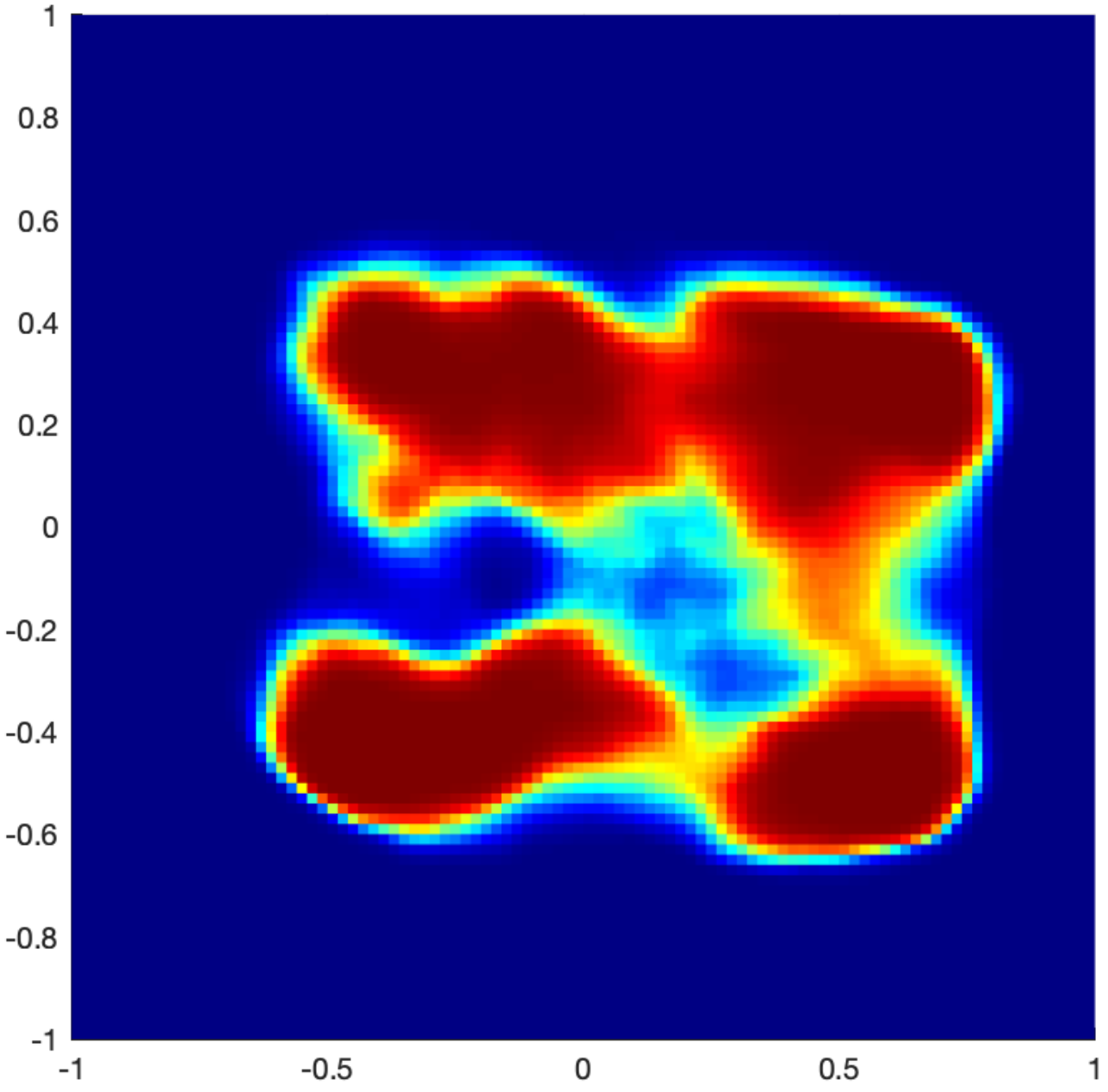}&
\includegraphics[width=1.1in]{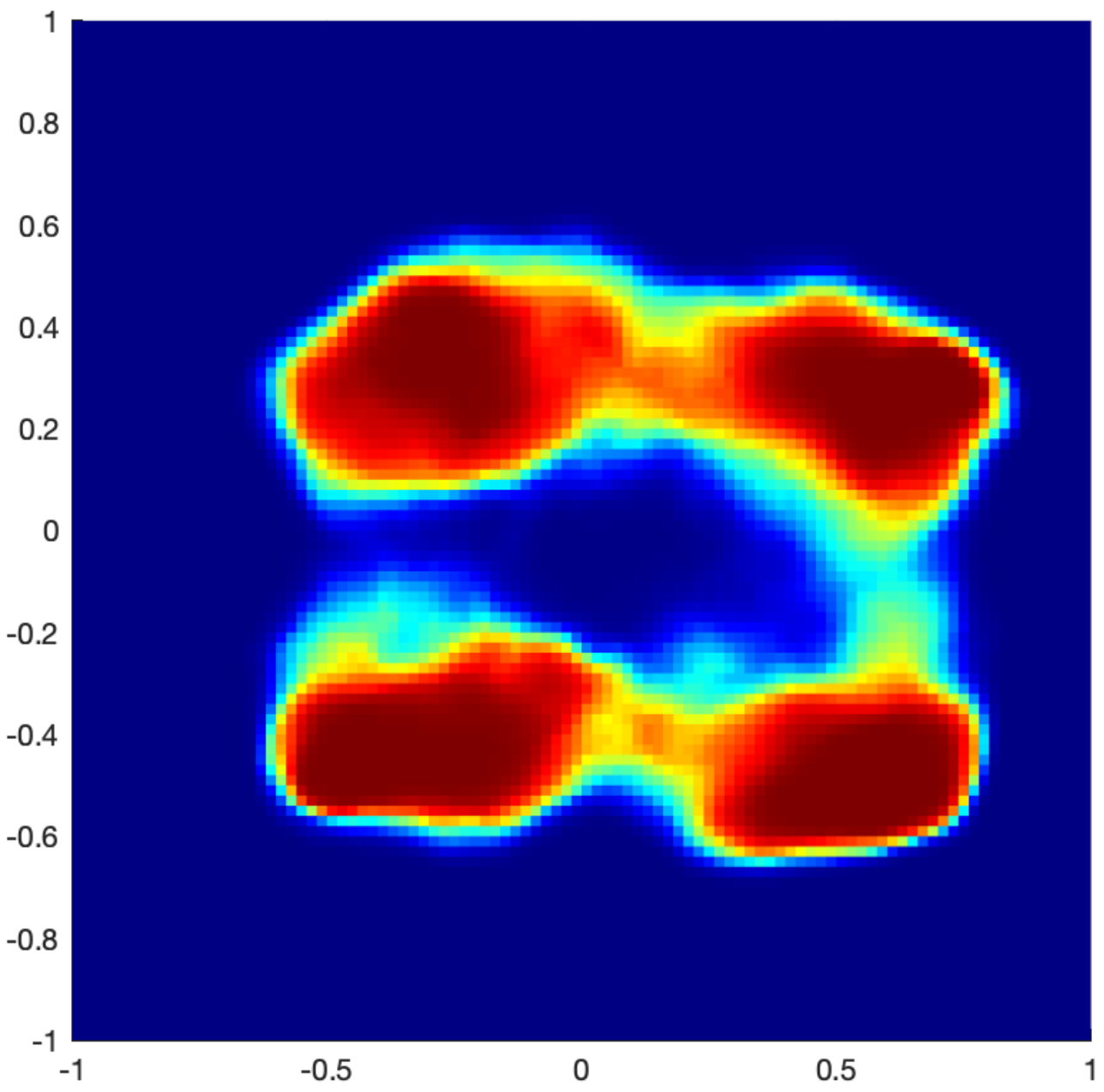}&
\includegraphics[width=1.1in]{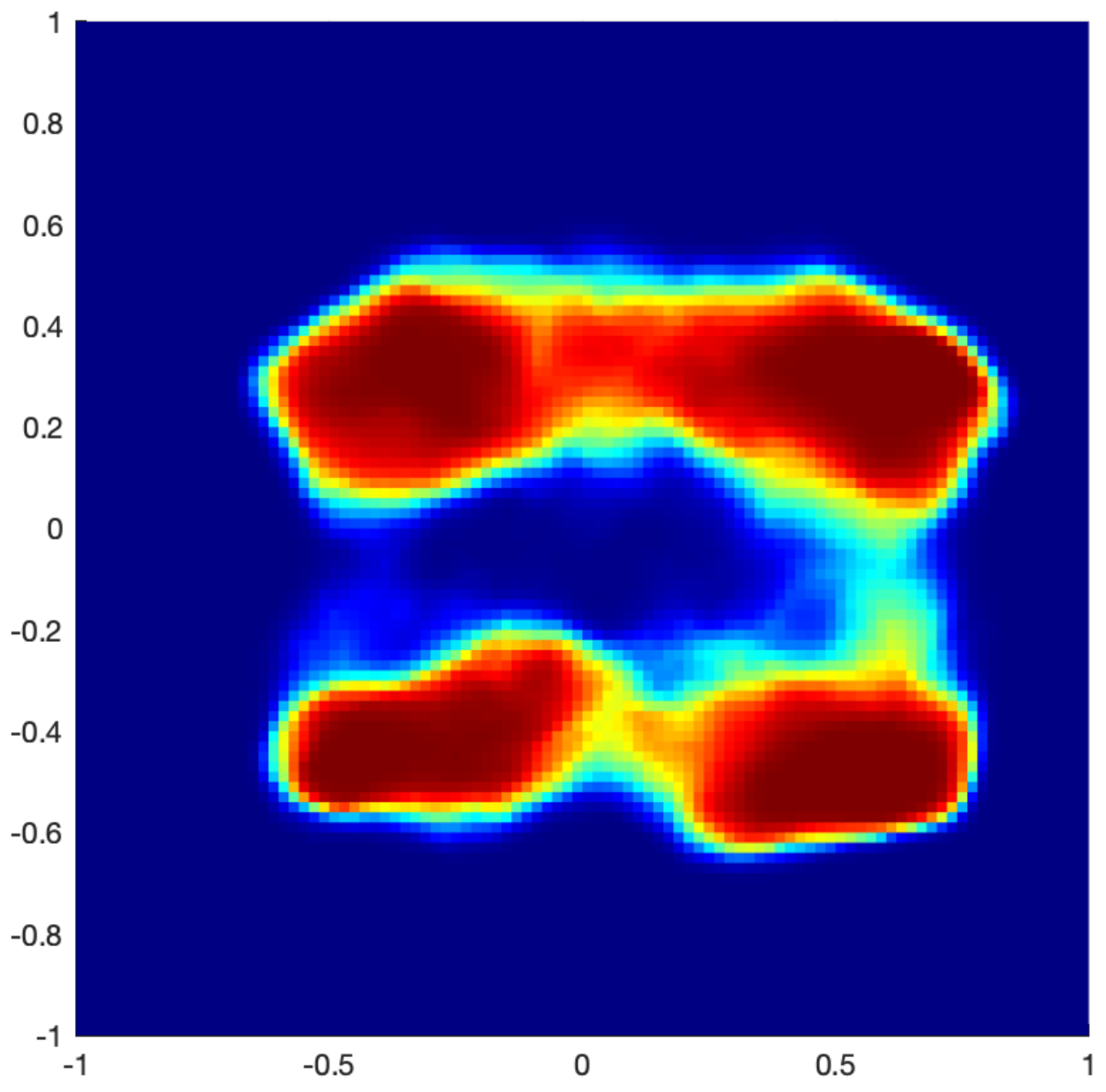}&
\includegraphics[width=1.1in]{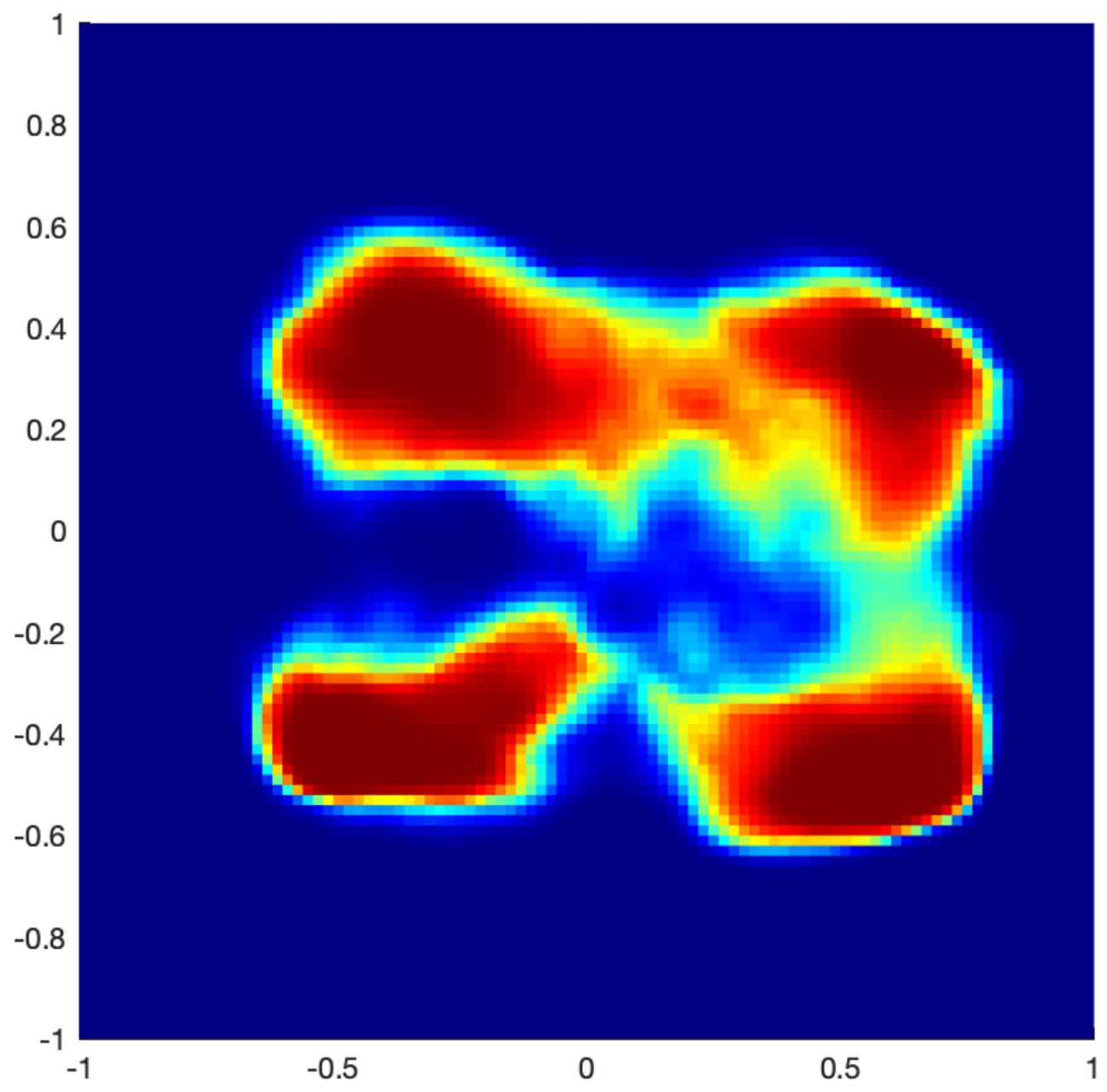}\\
\end{tabular}
  \caption{FNN-DDSM reconstruction for 3 special inclusion shapes: one triangle (top), two long rectangular bars (middle) and a rectangular annulus (top)} 
  \label{tab_FN_nonell}
\end{figure}

\begin{figure}[htbp]
\begin{tabular}{ >{\centering\arraybackslash}m{0.9in} >{\centering\arraybackslash}m{0.9in} >{\centering\arraybackslash}m{0.9in}  >{\centering\arraybackslash}m{0.9in}  >{\centering\arraybackslash}m{0.9in}  >{\centering\arraybackslash}m{0.9in} }
\centering
True coefficients &
N=1, $\delta=0$&
N=10, $\delta=0$&
N=20, $\delta=0$&
N=20, $\delta=10\%$ &
N=20, $\delta=20\%$ \\
\includegraphics[width=1.1in]{tri_0-eps-converted-to.pdf}&
\includegraphics[width=1.1in]{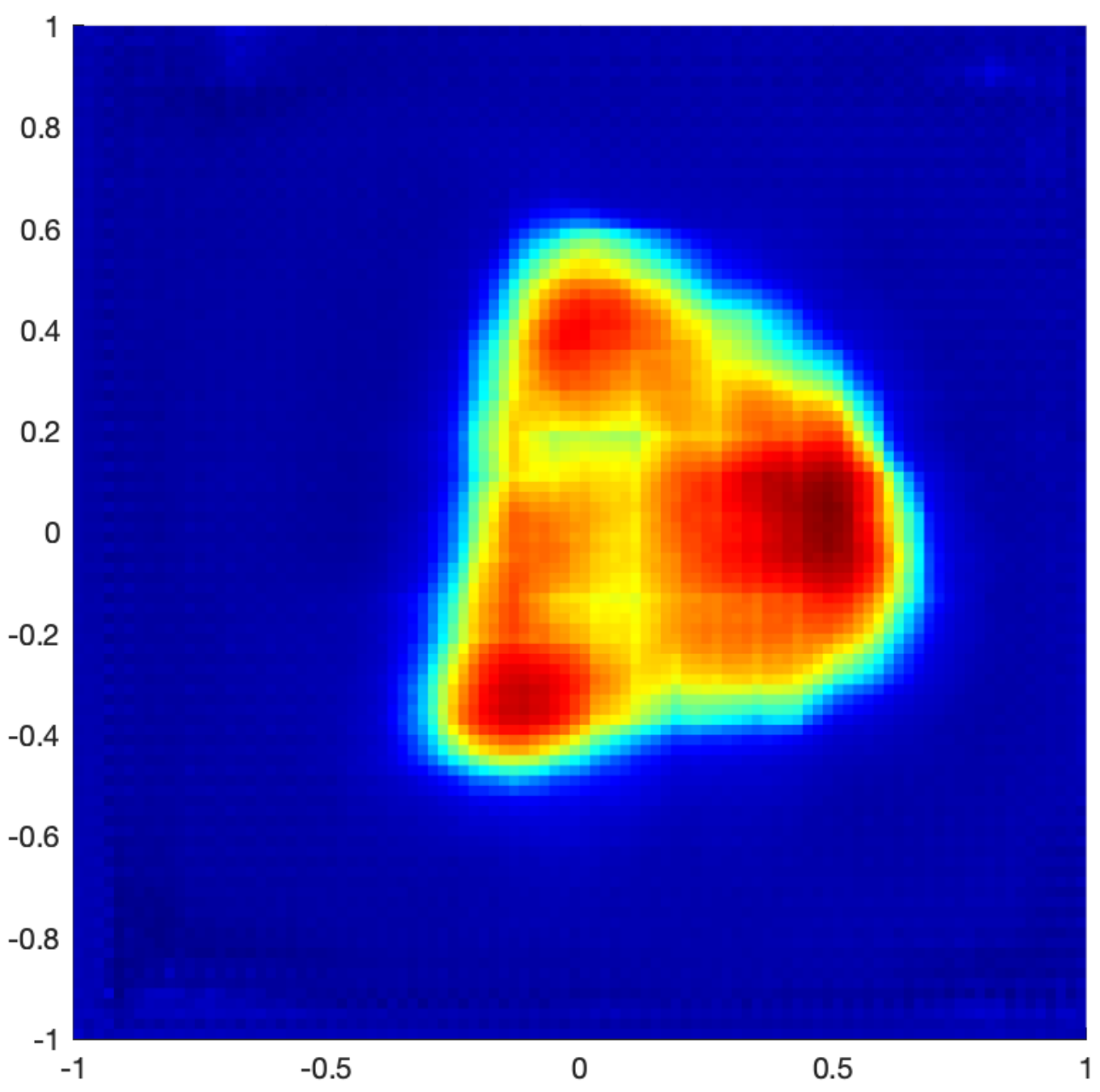}&
\includegraphics[width=1.1in]{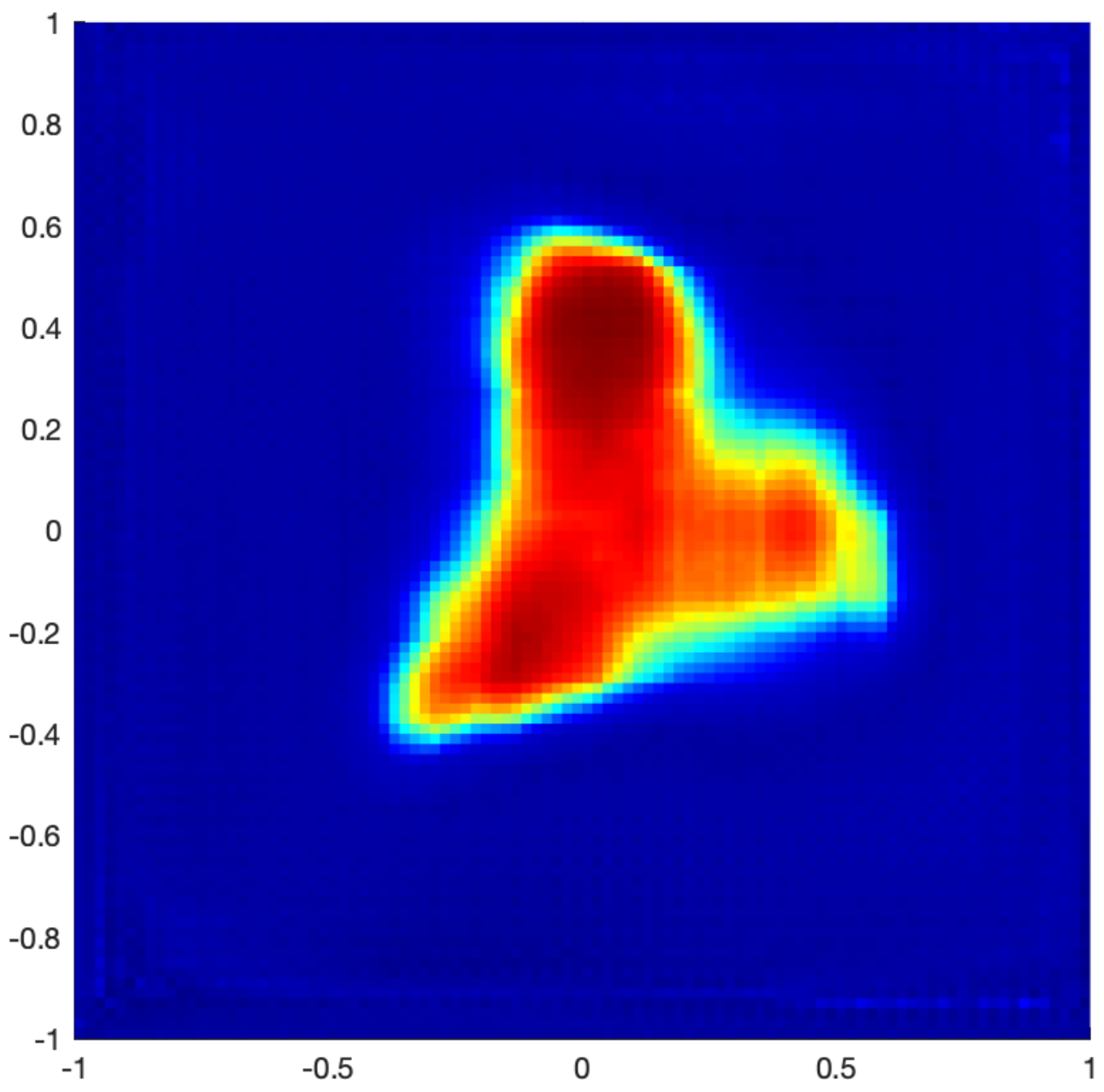}&
\includegraphics[width=1.1in]{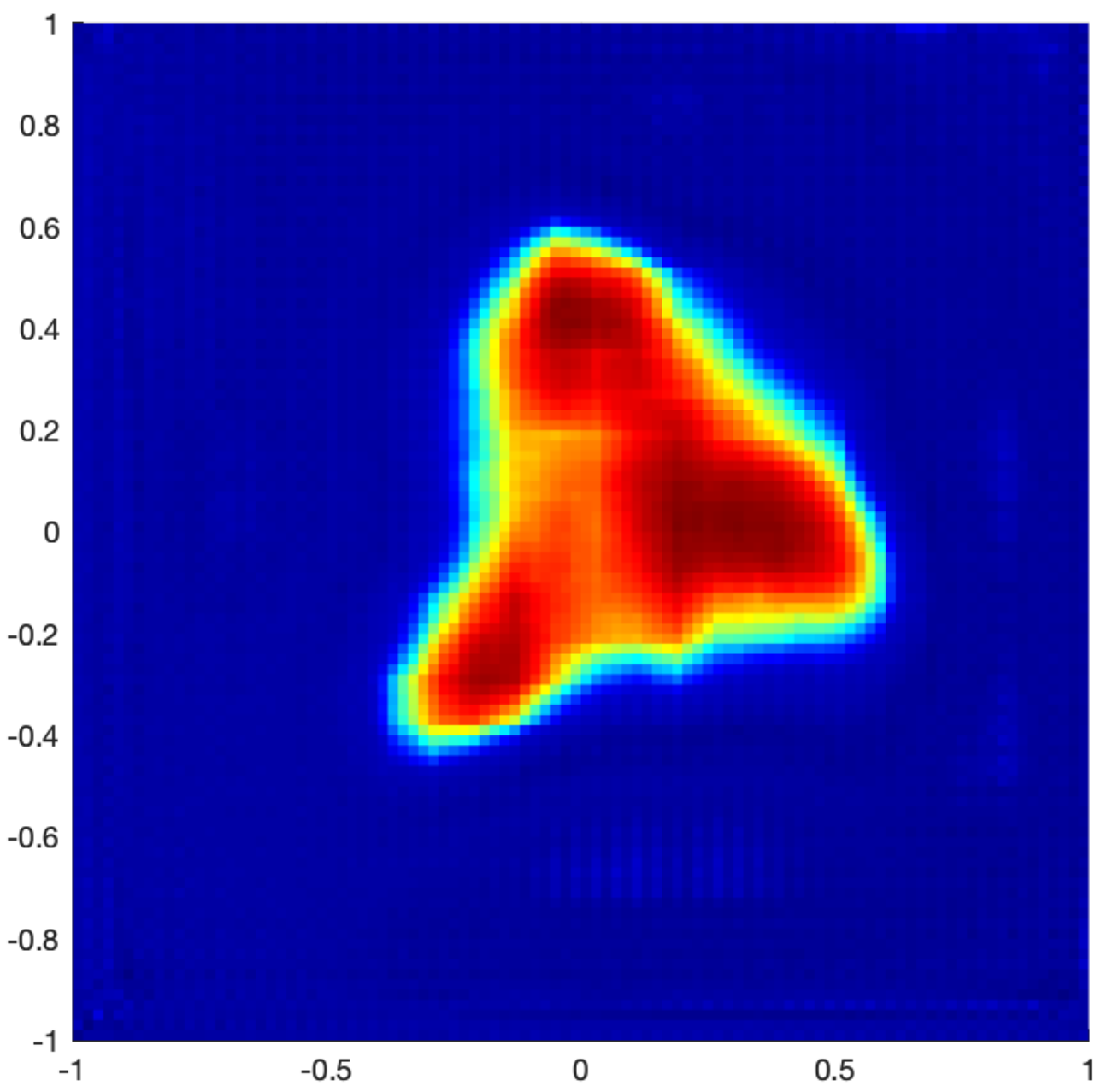}&
\includegraphics[width=1.1in]{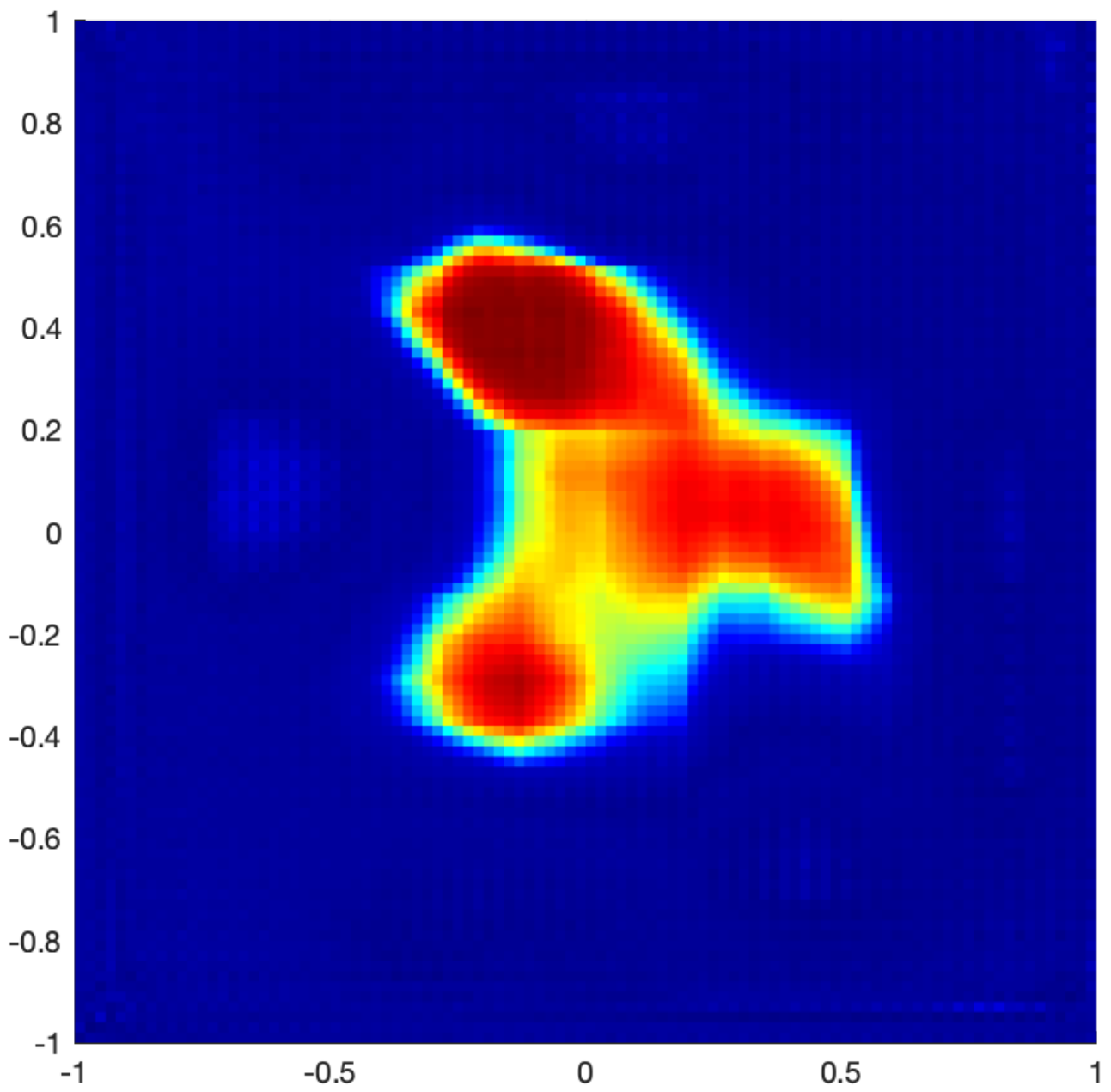}&
\includegraphics[width=1.1in]{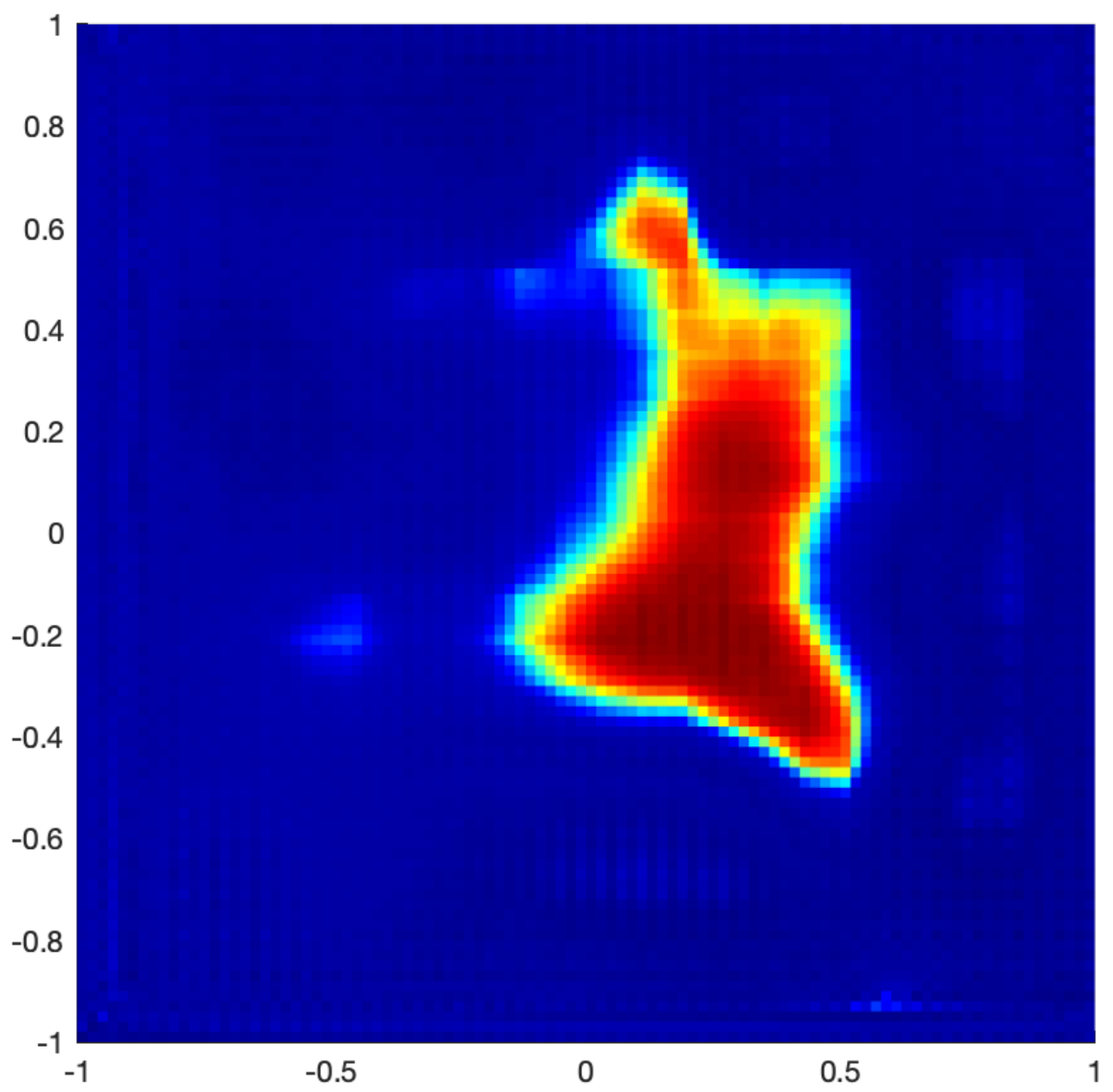}\\
\includegraphics[width=1.1in]{rec_case1_0-eps-converted-to.pdf}&
\includegraphics[width=1.1in]{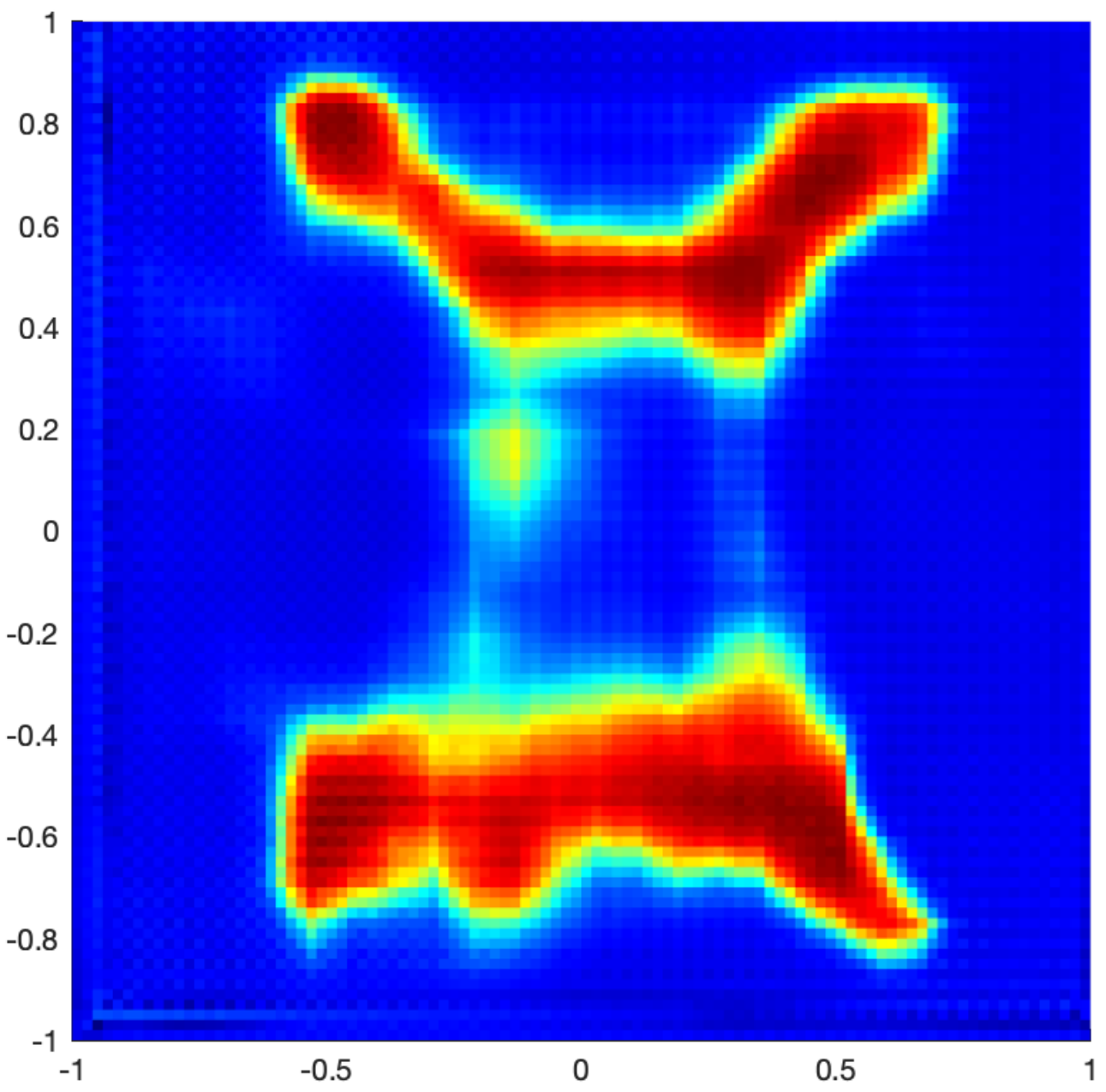}&
\includegraphics[width=1.1in]{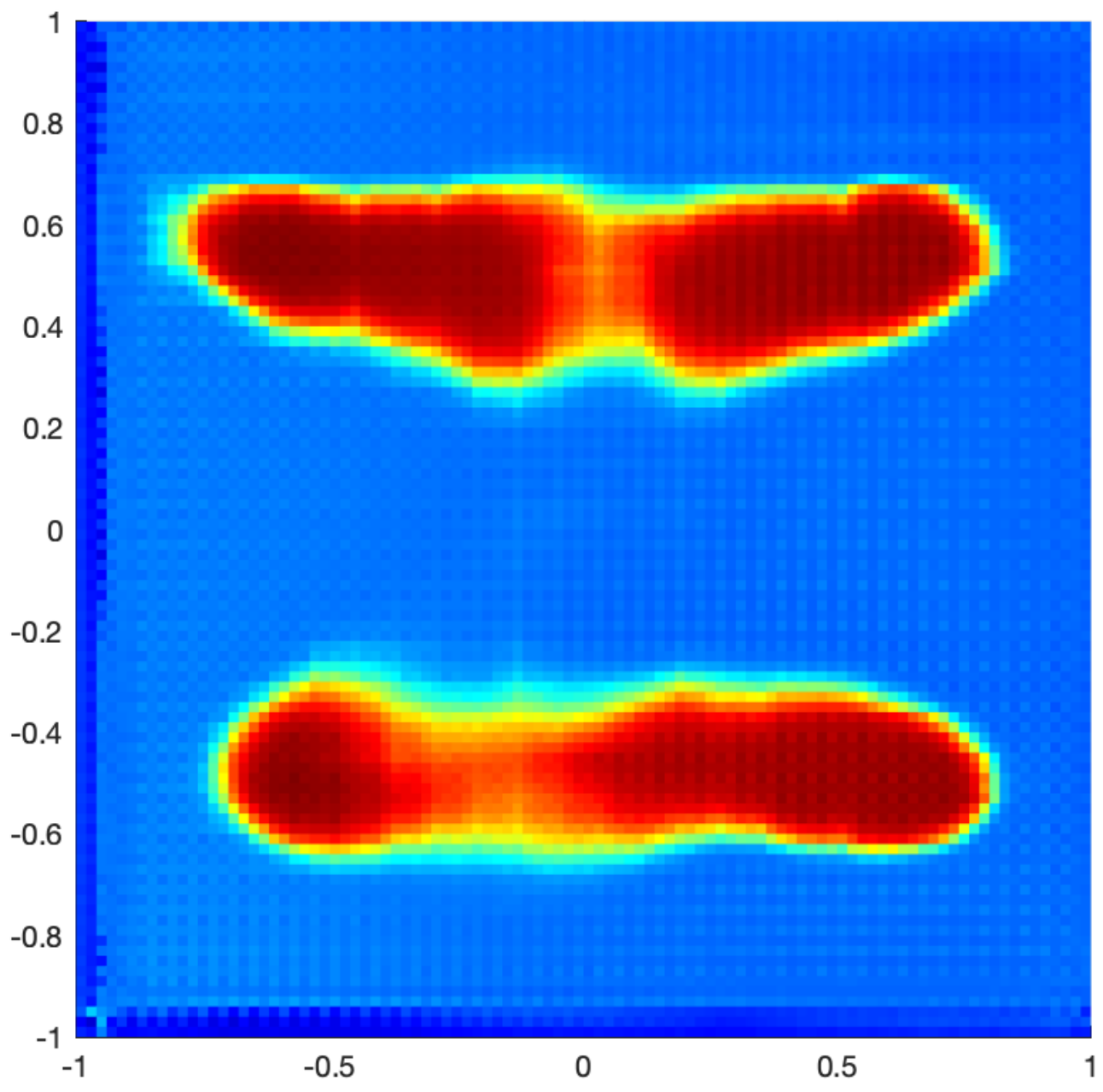}&
\includegraphics[width=1.1in]{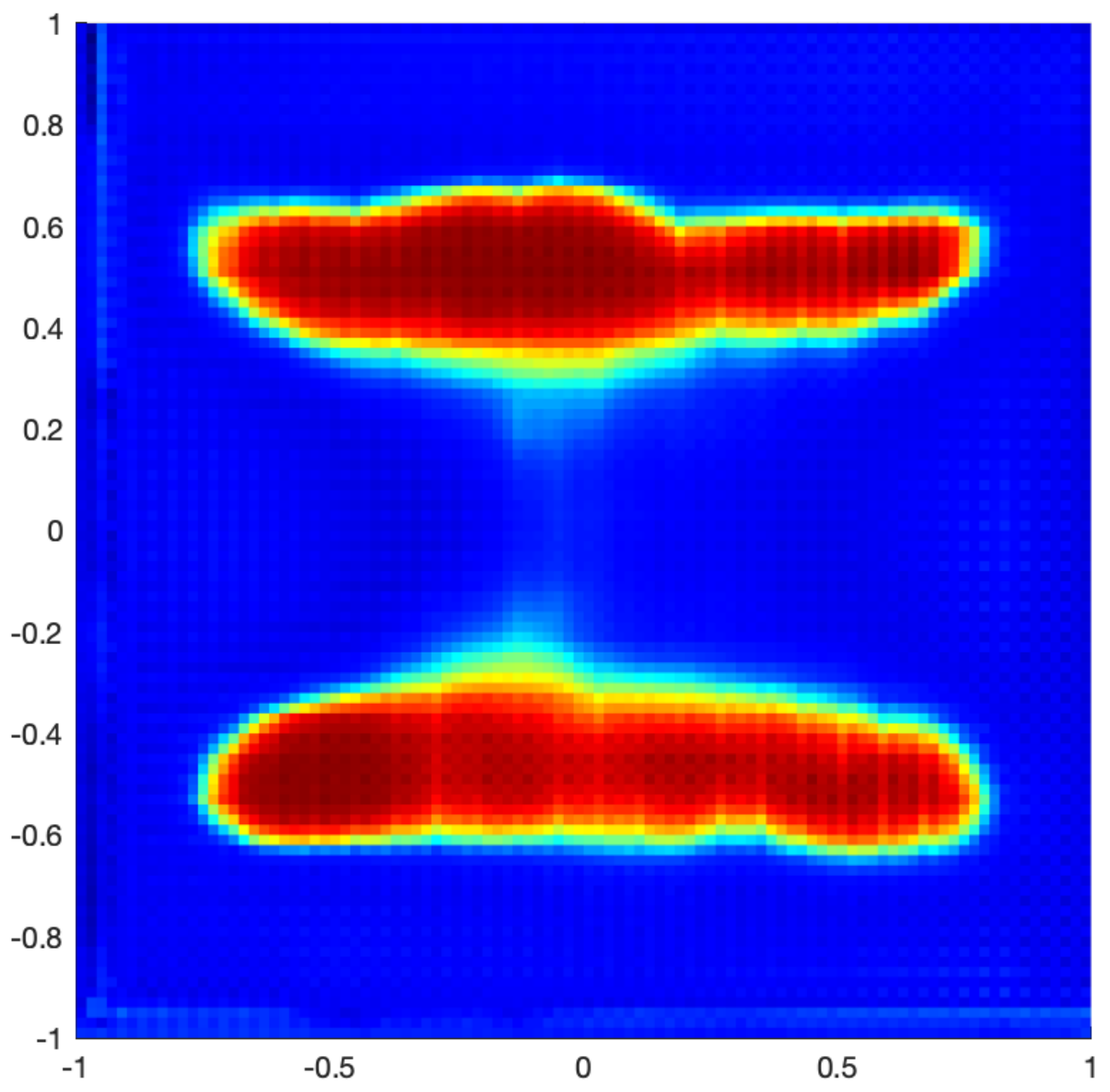}&
\includegraphics[width=1.1in]{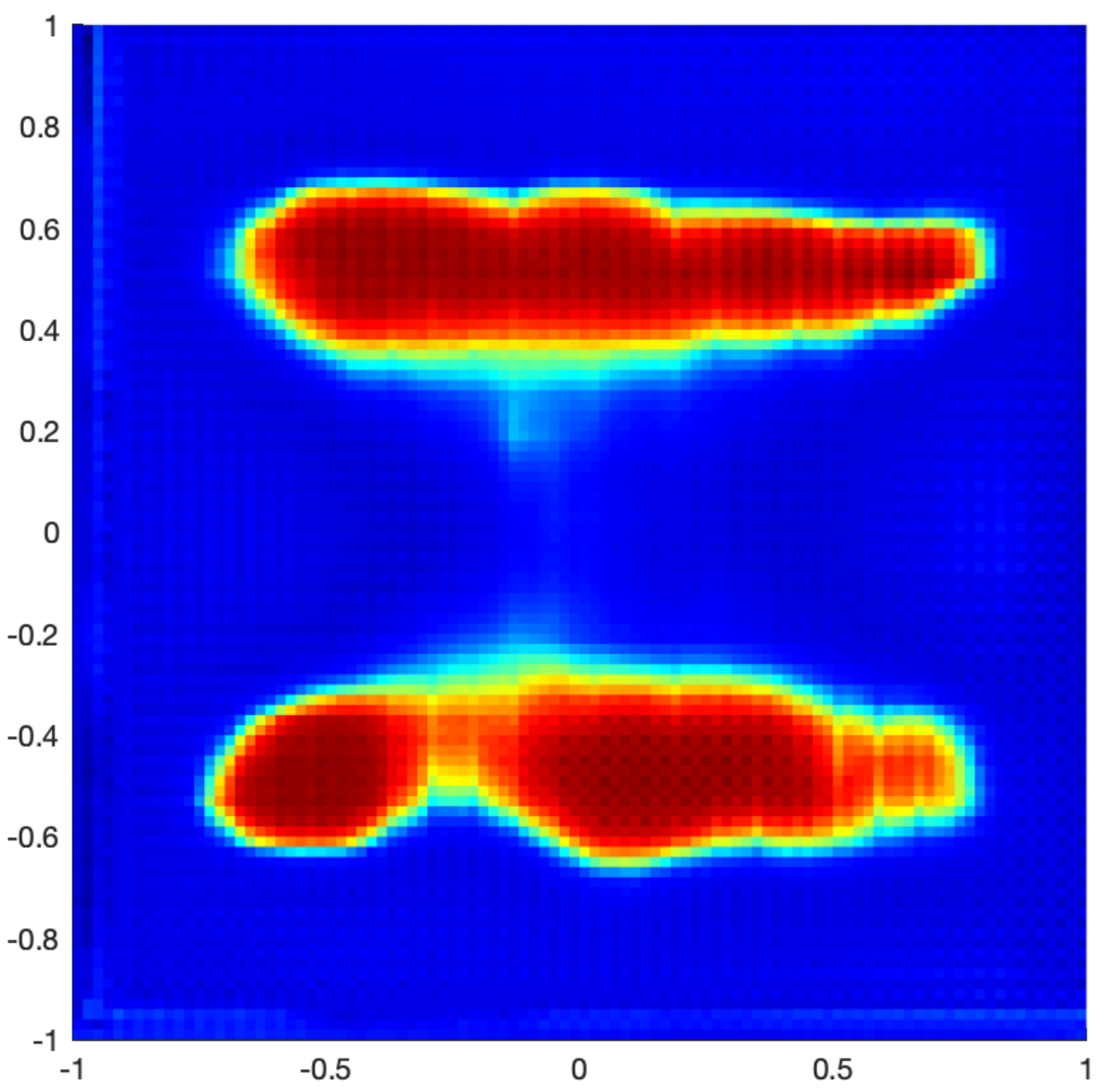}&
\includegraphics[width=1.1in]{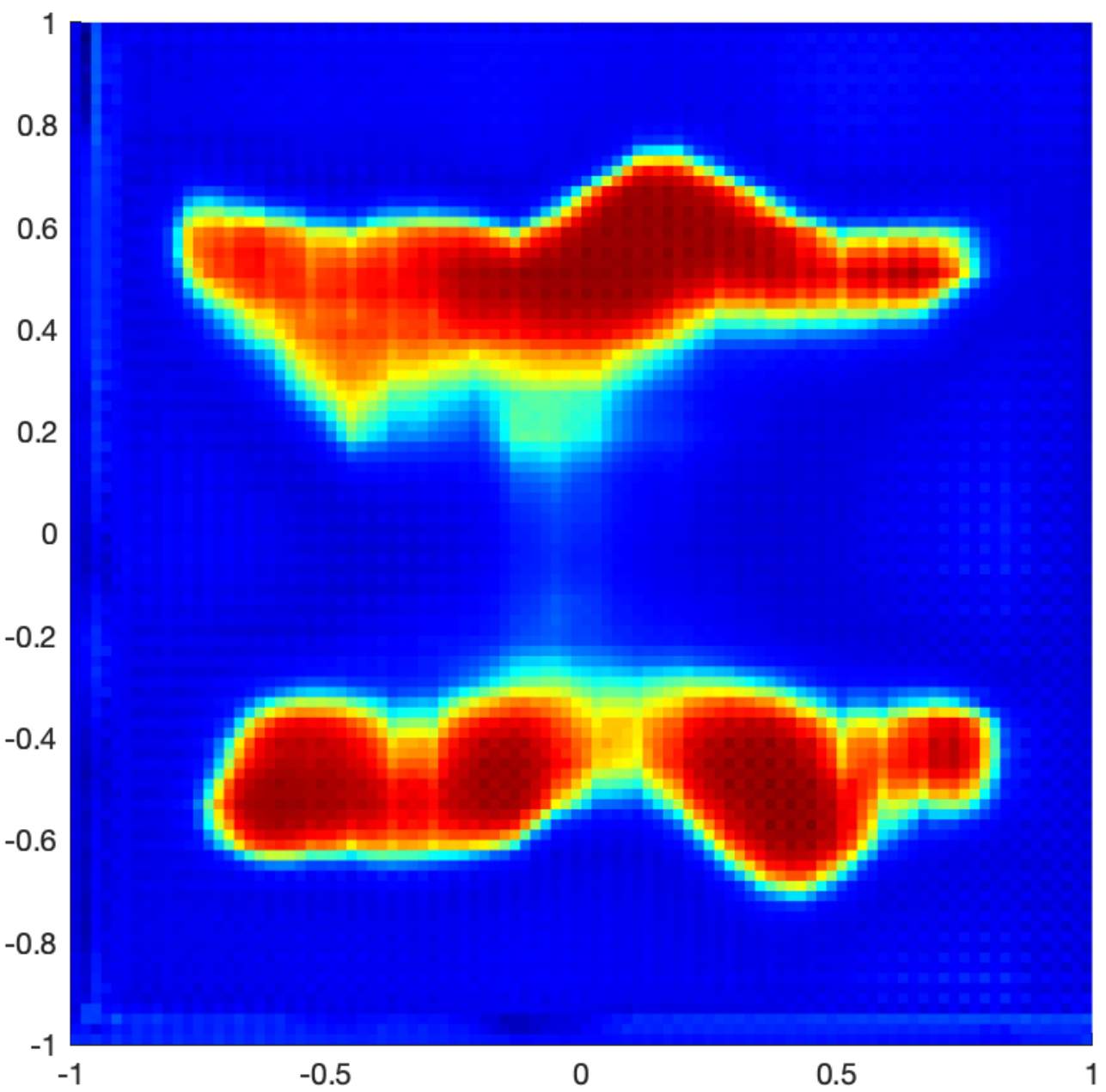}\\
\includegraphics[width=1.1in]{rec_case2_0-eps-converted-to.pdf}&
\includegraphics[width=1.1in]{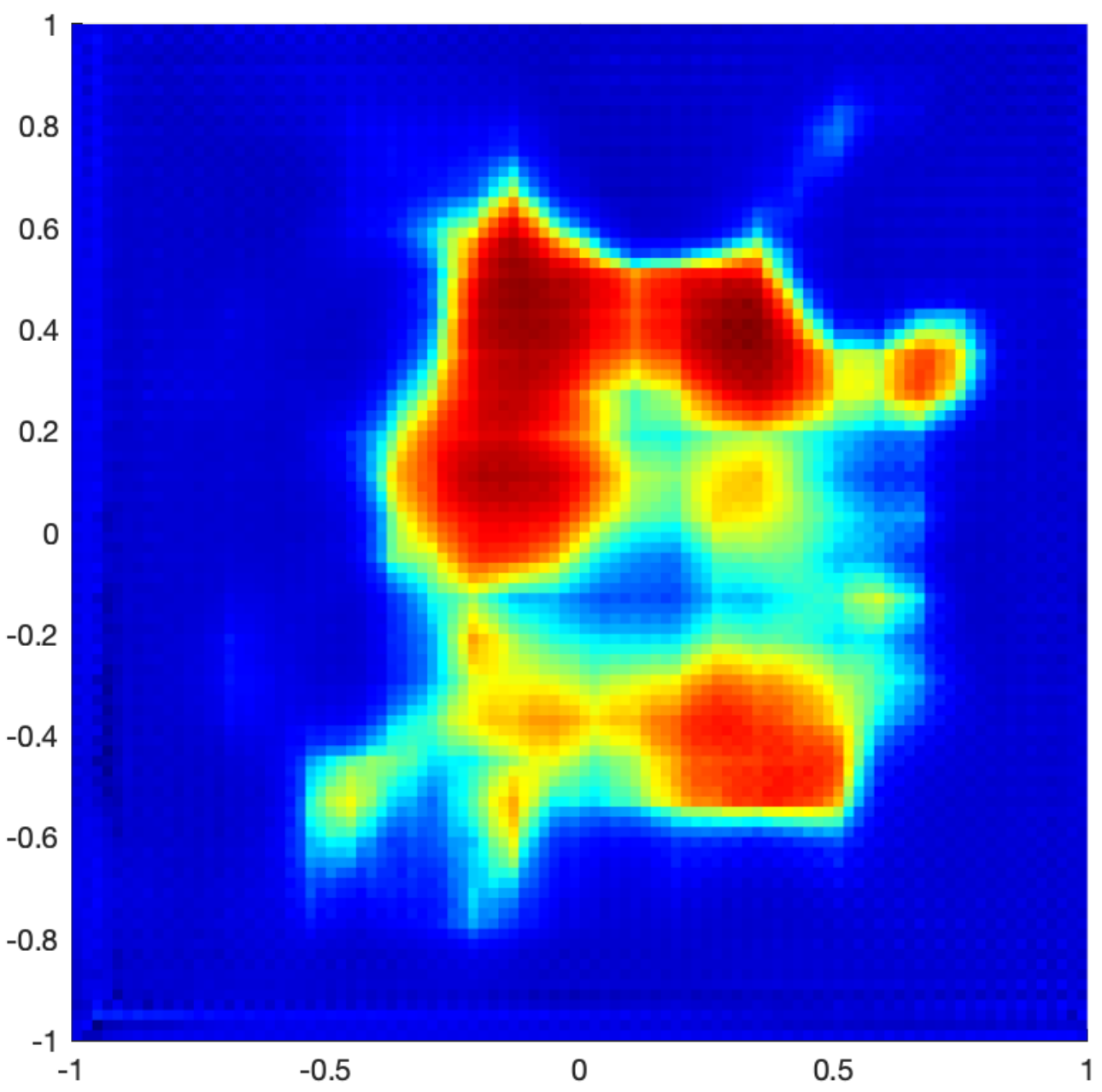}&
\includegraphics[width=1.1in]{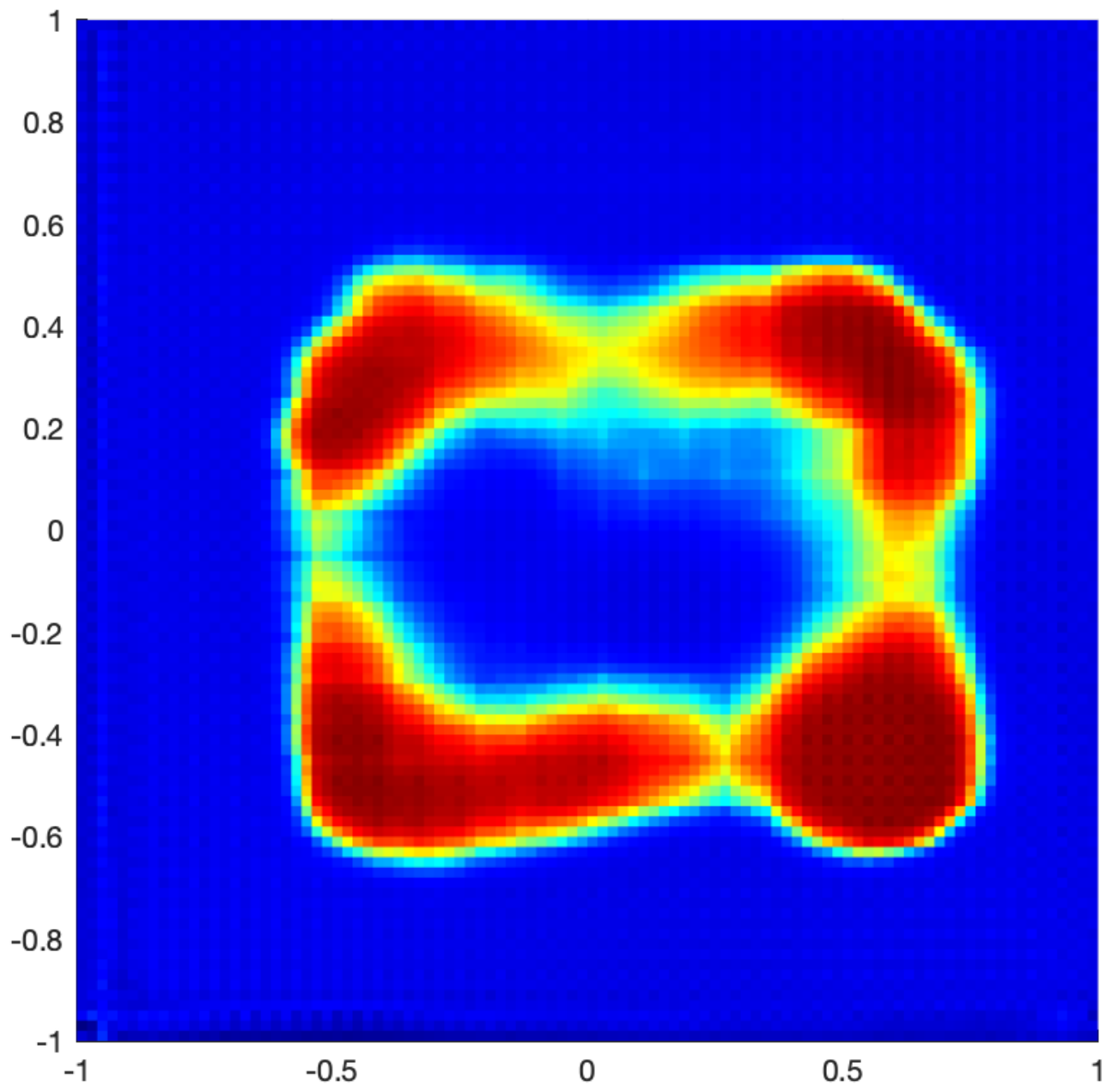}&
\includegraphics[width=1.1in]{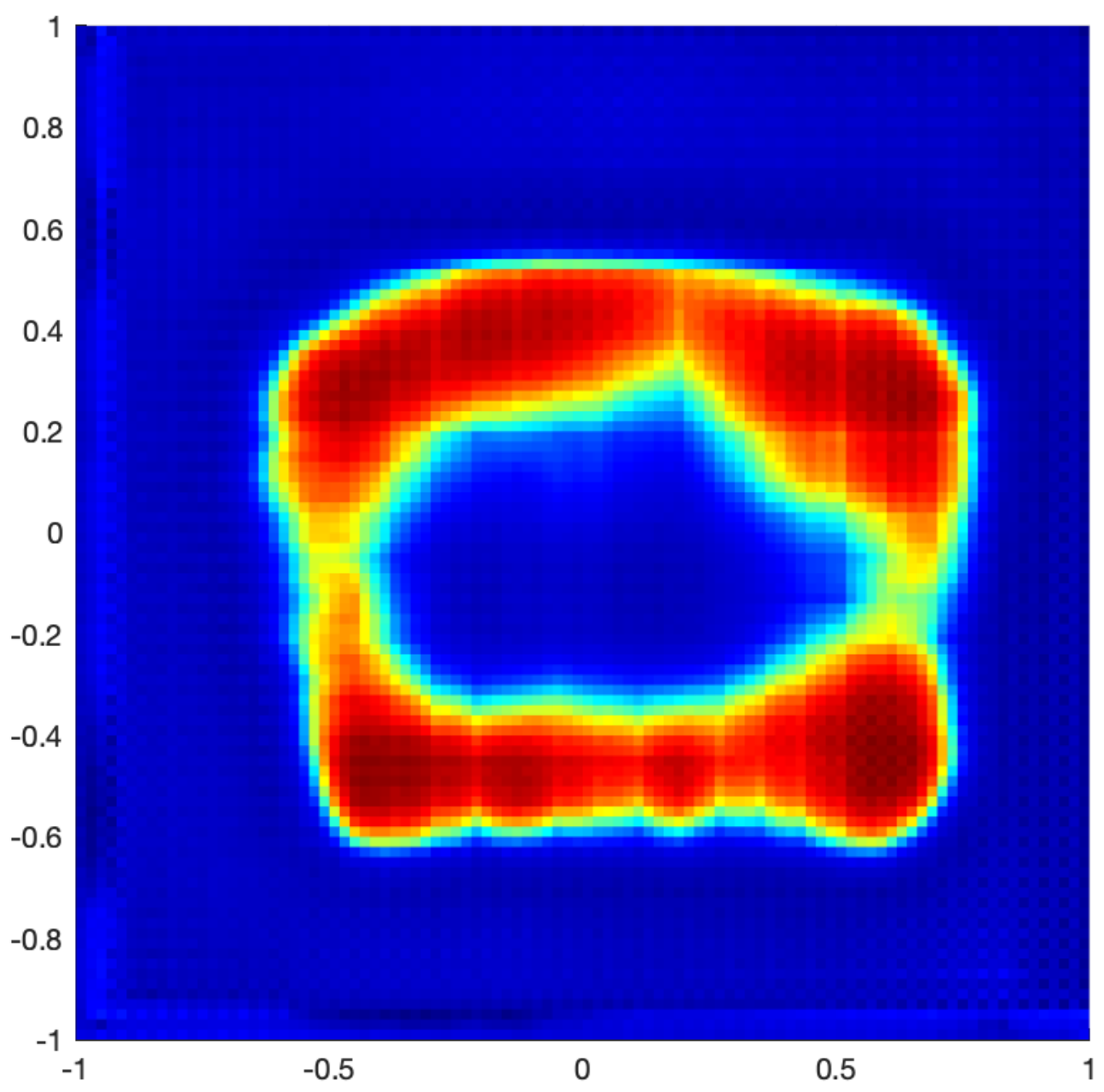}&
\includegraphics[width=1.1in]{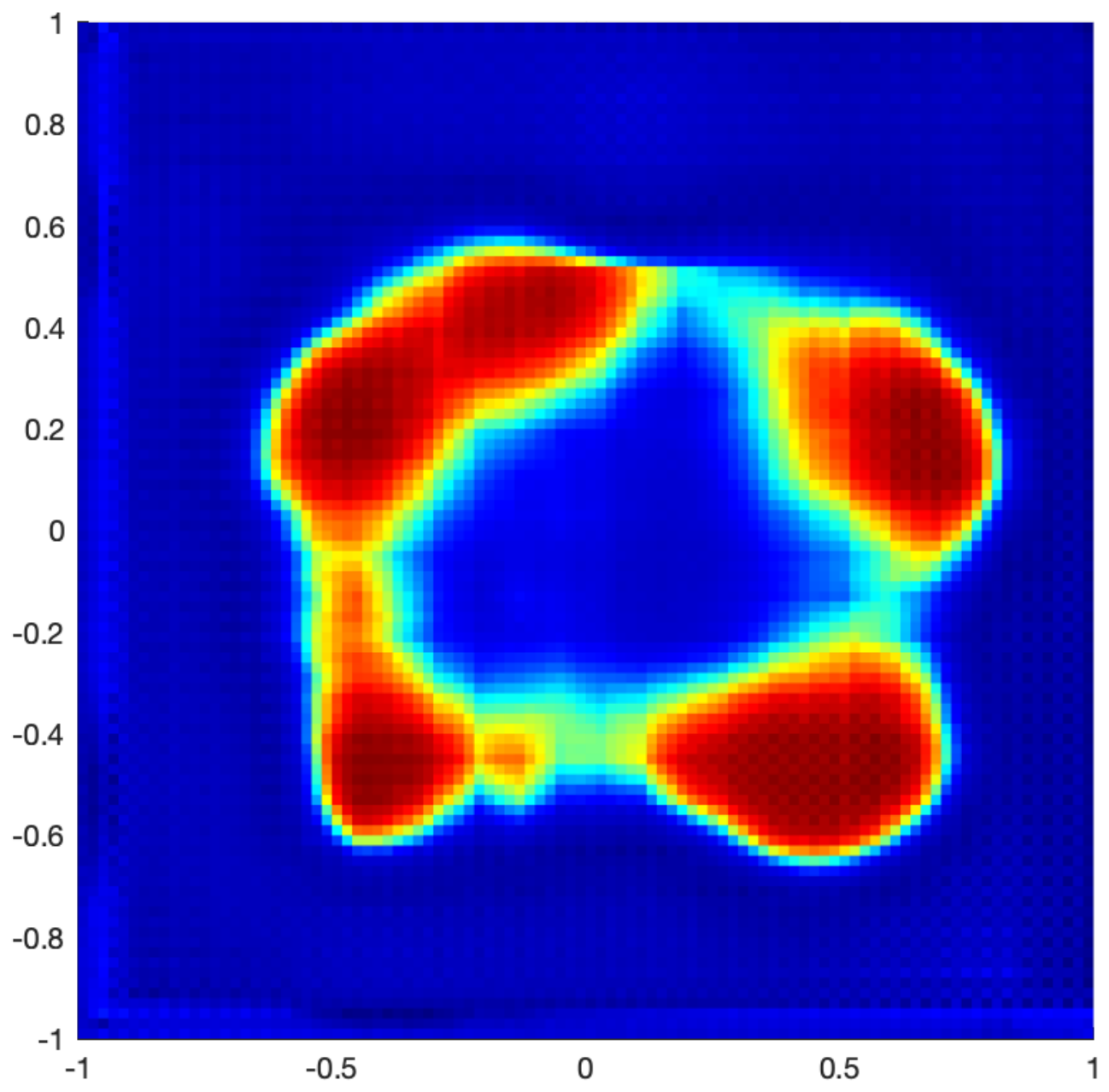}&
\includegraphics[width=1.1in]{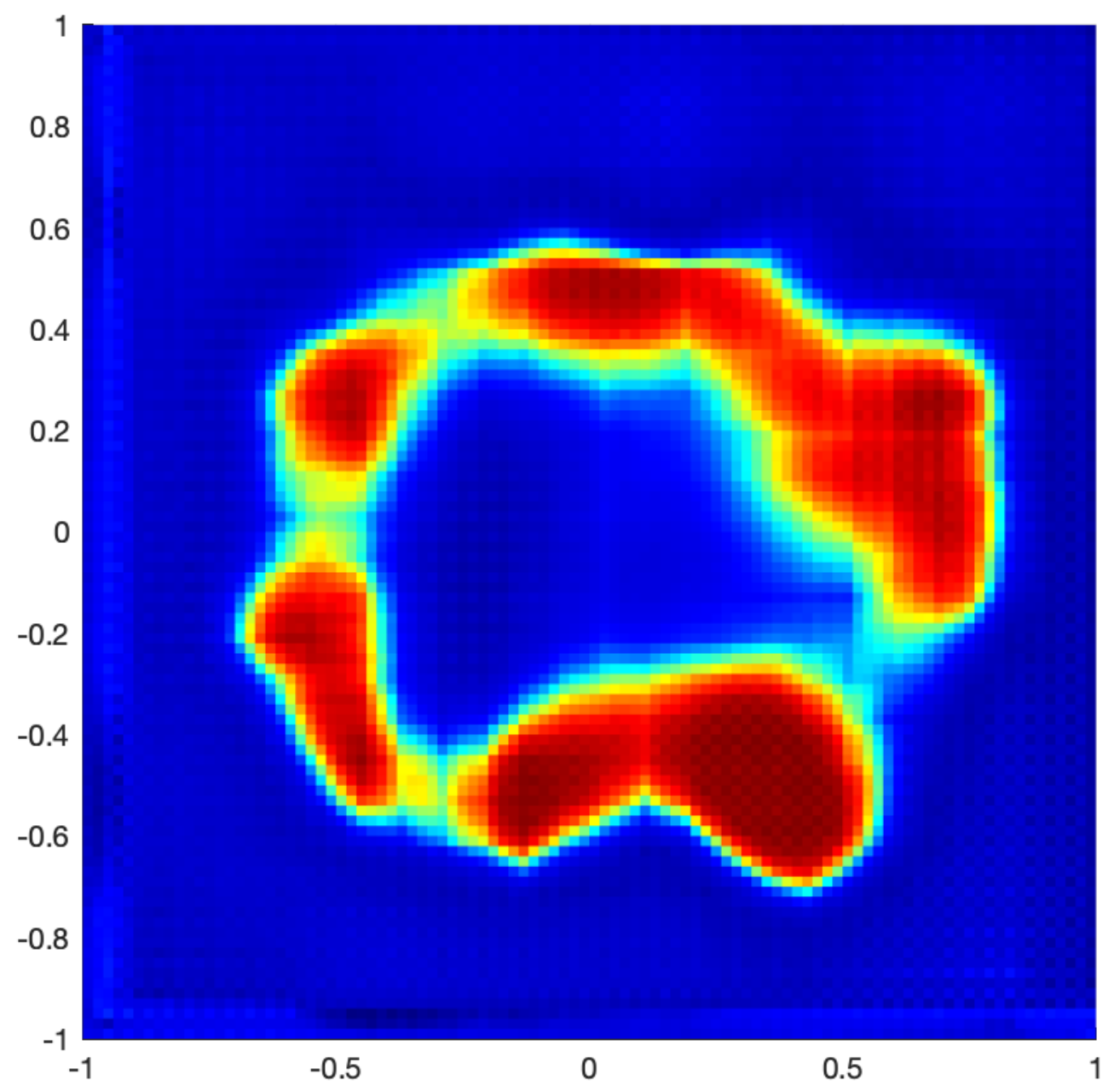}\\
\end{tabular}
  \caption{CNN-DDSM reconstruction for 3 special inclusion shapes: one triangle (top), two long rectangular bars (middle) and a rectangular annulus (top)} 
  \label{tab_nonell}
\end{figure}

For this purpose, we focus on the FNN-DDSM and CNN-DDSM trained by the data of 4 ellipses and show the reconstruction on some typical inclusion shapes: a triangle, two rectangular bar and a rectangular annulus in Figures \ref{tab_FN_nonell} and \ref{tab_nonell} which are indeed far away from union of ellipses and circles. 
Some similar ones were also used in \cite{chow2014direct,chow2015direct} for the conventional DSM. 
For the triangle shown in the first row in Figures \ref{tab_FN_nonell} and \ref{tab_nonell}, we can clearly see that the CNN-DDSM reconstructs the triangular shape while the FNN-DDSM {\grc only captures the three angles}. The reconstructions by CNN-DDSM with zero noise are all quite satisfactory. Even with 10\% noise the reconstruction still gives us some rough information about the triangular shape. The case of two rectangular bars is provided in the second row in Figures \ref{tab_FN_nonell} and \ref{tab_nonell}. For the CNN-DDSM,  except the reconstruction by the only single Cauchy data pair, all the other reconstructions recover the shape and position very well. The FNN-DDSN with $10$ pairs of Cauchy data and no noise also results in a good performance, but the reconstructions with other settings {\grc are relatively worse}. 
As mentioned in \cite{chow2014direct}, the most challenging case is the rectangular annulus since there is a hole in the inclusion which is hardly detected by the boundary data shown in the last row in Figures \ref{tab_FN_nonell} and \ref{tab_nonell}. For the CNN-DDSM, the reconstruction becomes much more accurate as the number of Cauchy data pairs increases. Even with 20\% noise, we can still obtain an annulus at the right location although its shape deteriorates.  Contrary to the CNN-DDSM, the FNN-DDSM captures the feature of the annulus accurately with single pair of Cauchy data {\grc but losses two edges for other settings.}

To summarize our findings for predicting inclusions different from the training samples, the CNN-DDSM performs quite satisfactorily and much better than the FNN-DDSM, which might be due to the more general structure of {\jjhn its index function trained by the CNN}. We believe that the behavior of DDSMs can be improved if the inclusions with more kinds of shapes such as the triangles and rectangles are added to the training set (users' library) to train DNNs. 




\section{Concluding remarks}
\label{sec:con}

In this work, based on the DSM invented in \cite{chow2014direct} we propose two approaches to construct DNNs for solving EIT problems and the DNN based DSM is called DDSMs. {\grc Our basic idea is to use a large amount of inclusion samples together with boundary measurements to learn the index function in DSM of which the construction is unknown by classical mathematical derivation for multiple Cauchy data pairs and general-shaped domain}. The first DNN we have proposed is a FNN directly approximating the index function pointwisely by taking the spatial variable $x$ and the gradient of Cauchy difference functions $\nabla\phi^{\omega}(x)$ as the input, which is called FNN-DDSM. A remarkable feature is that its output has very clear probabilistic meaning, i.e., the chance of $x$ inside or outside the inclusions, and we hope it can motivate the research between probability and DSM. The second one is a CNN approximating a so-called index functional from a stack of Cauchy difference functions to the inclusion distribution, which is called CNN-DDSM. It is worth mentioning that the CNN-DDSM is a further generalization of the conventional DSM \cite{chow2014direct} since it relaxes the assumption that the location of each point $x$ only lies on the data at this point. Additionally, it has strong connection to the image segmentation problems.  The proposed two DDSMs inherit the features of both classical optimization-type methods (accuracy) and the DSM \cite{chow2014direct} (efficiency) based on their offline-online decomposition structure, that is, the costly optimization procedure only needs to be done once and in offline {\jjh phase}, and the prediction is done by fast direct evaluation and in online {\jjh phase}. 

We have also carried out extensive numerical experiments to show that the proposed DDSMs are effective, accurate and robust even for very complicated geometry. In particular, the DDSMs are highly stable with respect to large noise which, we believe, is due to noise smoothing procedure in the generation of the Cauchy difference functions \eqref{eq:defindexfun}. Comparing these two methods on the testing set within the scope of geometry set-up of the training set, we conclude that the prediction of these two DDSMs are almost comparable. But the CNN-DDSM is slightly more sensitive to the center inclusion than FNN-DDSM while the FNN-DDSM seems more stable with respect to noise. When applying them to the inclusions of which the geometry is out of the scope of the training set, the performance of the CNN-DDSM is much better than the FNN-DDSM, which may be due to the fact that the CNN-DDSM uses more neighborhood information to predict the location of each point. The observations from the numerical experiments suggest that our DDSMs have great potential to improve if more types of inclusions are included in the training set.


\end{document}